\documentclass{article}

\usepackage{a4wide,cite,latexsym,amsfonts,amssymb,exscale,enumerate}
\usepackage[centertags,sumlimits,intlimits,namelimits,reqno]{amsmath}
\usepackage{amsthm}
\usepackage{hyperref,picins}

\usepackage{LPsymb}

\usepackage{pstricks}
\usepackage{pst-grad}
\usepackage{pst-plot}
\usepackage[tiling]{pst-fill}
\usepackage{pstcol}

\usepackage{graphicx}
\usepackage[all, knot]{xy}

\usepackage{fancyheadings}
\newcommand{\cat}[1]{\ensuremath{\mbox{\bfseries {\upshape {#1}}}}}

\newcommand{\BOX}{\hbox {$\sqcap$ \kern -1em $\sqcup$}}

\renewcommand{\to}{\rightarrow}
\newcommand{\maps}{\colon}

\newcommand{\id}{{\rm id}}
\newcommand{\del}{\partial}
\newcommand{\End}{{\rm End}}

\newcommand{\im}{{\rm im\ }}
\newcommand{\coim}{{\rm coim\ }}
\newcommand{\chr}{{\rm char\ }}

\newcommand{\spann}{{\rm span}}
\newcommand{\rk}{{\rm rk\ }}
\def\bigboxtimes{\mathop{\boxtimes}\limits}

\newcommand{\F}{{\mathbb F}}

\newcommand{\mT}{{\rm T} }

\newcommand{\scs}{\scriptstyle}

\theoremstyle{definition}
\newtheorem{thm}{Theorem}[section]
\newtheorem{cor}[thm]{Corollary}

\newtheorem{lem}[thm]{Lemma}
\newtheorem{rem}[thm]{Remark}
\newtheorem{prop}[thm]{Proposition}
\newtheorem{defn}[thm]{Definition}
\newtheorem{example}[thm]{Example}

\newcommand{\be}{\begin{equation}}
\newcommand{\ee}{\end{equation}}
\newcommand{\ba}{\begin{eqnarray}}
\newcommand{\ea}{\end{eqnarray}}
\newcommand{\ban}{\begin{eqnarray*}}
\newcommand{\ean}{\end{eqnarray*}}
\newcommand{\barr}{\begin{array}}
\newcommand{\earr}{\end{array}}

\SelectTips{cm}{}

\psset{linewidth=0.3pt,dimen=middle}
\psset{xunit=.70cm,yunit=0.70cm}
\newgray{whitegray}{.80}
\psset{arrowsize=1pt 5,arrowlength=.6,arrowinset=.7}

\numberwithin{equation}{section}

\def\emph#1{{\sl #1\/}}
\def\ie{{\sl i.e.\/}}

\def\etc{{\sl etc.\/}}
\def\cf{{\sl c.f.\/}}

\let\hat=\widehat
\let\tilde=\widetilde

\let\phi=\varphi
\let\theta=\vartheta
\let\epsilon=\varepsilon

\usepackage{bbm}

\def\N{{\mathbbm N}}
\def\R{{\mathbbm R}}
\def\Z{{\mathbbm Z}}

\def\H{{\mathbbm H}}

\def\cal#1{\mathcal{#1}}%
\def\1{\mathbbm{1}}%
\def\tr{\mathrm{tr}}%
\def\nn{\notag}

\def\llbracket{\left[\hbox to-0.3em{\hss}\left[}
\def\rrbracket{\right]\hbox to-0.3em{\hss}\right]}
\def\tangle#1{\llbracket #1\rrbracket}
\def\la{\langle}
\def\ra{\rangle}

\def\theauthor{}
\def\empty{}
\def\theaffiliation{}
\def\preprint#1{
  \thispagestyle{plain}
  \def\theauthor{#1}
  \ifx\theauthor\empty
  \else
    \begin{flushright}{\small #1\par}\end{flushright}
  \fi
  \begin{center}}
\def\title#1{
  {\LARGE #1\par}\vskip 1em}
\def\author#1{
  \ifx\theaffiliation\empty
  \else
    \par
  \fi
  \def\theauthor{#1}\def\theaffiliation{}}
\def\email#1{
  \vskip 1em{\large\theauthor\footnote{\small email: {\tt #1}}\par}\vskip .5em}
\def\affiliation#1{
  \ifx\theaffiliation\empty
    \def\theaffiliation{second}
  \else
    \par and\par
  \fi
  {\small\sl #1}}
\def\date#1{
  \vskip 1em{(#1)\par}\end{center}\vskip 2em}

\newcommand{\hopfresolutions}{
 \xy
  (-5,0)*+{
  \psset{xunit=.3cm,yunit=.3cm}
 \begin{pspicture}[.4](2,2)
 \psbezier[linewidth=.8pt](.5,0)(.5,.5)(1.5,.5)(1.5,1)
 \psbezier[linewidth=.8pt](1.5,2)(1.5,1.5)(.5,1.5)(.5,1)
 \pspolygon[linecolor=white,fillstyle=solid,fillcolor=white](.9,0)(.9,2)(1.1,2)(1.1,0)(.9,0)
 \psbezier[linewidth=.8pt](1.5,0)(1.5,.5)(.5,.5)(.5,1)
 \psbezier[linewidth=.8pt](.5,2)(.5,1.5)(1.5,1.5)(1.5,1)
 \psbezier[linewidth=.8pt](1.5,0)(1.5,-.25)(2,-.25)(2,0)
 \psbezier[linewidth=.8pt](1.5,2)(1.5,2.25)(2,2.25)(2,2)
 \psline[linewidth=.8pt](2,0)(2,2)
 \psbezier[linewidth=.8pt](.5,0)(.5,-.25)(0,-.25)(0,0)
 \psbezier[linewidth=.8pt](.5,2)(.5,2.25)(0,2.25)(0,2)
 \psline[linewidth=.8pt](0,0)(0,2)
\end{pspicture}
  }="t";
  (5,10)*+{\psset{xunit=.3cm,yunit=.3cm}
  \begin{pspicture}[.4](2,2)
 \psbezier[linewidth=.8pt](.5,0)(.5,.5)(1.5,.5)(1.5,1)
 \pspolygon[linecolor=white,fillstyle=solid,fillcolor=white](.9,0)(.9,2)(1.1,2)(1.1,0)(.9,0)
 \psbezier[linewidth=.8pt](1.5,0)(1.5,.5)(.5,.5)(.5,1)
 \psbezier[linewidth=.8pt](.5,2)(.5,1.5)(1.5,1.5)(1.5,2)
 \psbezier[linewidth=.8pt](.5,1)(.5,1.5)(1.5,1.5)(1.5,1)
 \psbezier[linewidth=.8pt](1.5,0)(1.5,-.25)(2,-.25)(2,0)
 \psbezier[linewidth=.8pt](1.5,2)(1.5,2.25)(2,2.25)(2,2)
 \psline[linewidth=.8pt](2,0)(2,2)
 \psbezier[linewidth=.8pt](.5,0)(.5,-.25)(0,-.25)(0,0)
 \psbezier[linewidth=.8pt](.5,2)(.5,2.25)(0,2.25)(0,2)
 \psline[linewidth=.8pt](0,0)(0,2)
  \end{pspicture}
  }="0";
  (5,-10)*+{\psset{xunit=.3cm,yunit=.3cm}
  \begin{pspicture}[.4](2,2)
 \psbezier[linewidth=.8pt](.5,0)(.5,.5)(1.5,.5)(1.5,1)
 \pspolygon[linecolor=white,fillstyle=solid,fillcolor=white](.9,0)(.9,2)(1.1,2)(1.1,0)(.9,0)
 \psbezier[linewidth=.8pt](1.5,0)(1.5,.5)(.5,.5)(.5,1)
 \pscurve[linewidth=.8pt](.5,2)(.75,1.5)(.5,1)
 \pscurve[linewidth=.8pt](1.5,2)(1.25,1.5)(1.5,1)
 \psbezier[linewidth=.8pt](1.5,0)(1.5,-.25)(2,-.25)(2,0)
 \psbezier[linewidth=.8pt](1.5,2)(1.5,2.25)(2,2.25)(2,2)
 \psline[linewidth=.8pt](2,0)(2,2)
 \psbezier[linewidth=.8pt](.5,0)(.5,-.25)(0,-.25)(0,0)
 \psbezier[linewidth=.8pt](.5,2)(.5,2.25)(0,2.25)(0,2)
 \psline[linewidth=.8pt](0,0)(0,2)
\end{pspicture}}="1";
  (20,15)*+{\psset{xunit=.3cm,yunit=.3cm}
    \begin{pspicture}[.4](2,2)
 \psbezier[linewidth=.8pt](.5,2)(.5,1.5)(1.5,1.5)(1.5,2)
 \psbezier[linewidth=.8pt](.5,1)(.5,1.5)(1.5,1.5)(1.5,1)
 \psbezier[linewidth=.8pt](.5,0)(.5,.5)(1.5,.5)(1.5,0)
 \psbezier[linewidth=.8pt](.5,1)(.5,.5)(1.5,.5)(1.5,1)
 \psbezier[linewidth=.8pt](1.5,0)(1.5,-.25)(2,-.25)(2,0)
 \psbezier[linewidth=.8pt](1.5,2)(1.5,2.25)(2,2.25)(2,2)
 \psline[linewidth=.8pt](2,0)(2,2)
 \psbezier[linewidth=.8pt](.5,0)(.5,-.25)(0,-.25)(0,0)
 \psbezier[linewidth=.8pt](.5,2)(.5,2.25)(0,2.25)(0,2)
 \psline[linewidth=.8pt](0,0)(0,2)
  \end{pspicture}
  }="00";
  (20,5)*+{
  \psset{xunit=.3cm,yunit=.3cm}
    \begin{pspicture}[.4](2,2)
 \psbezier[linewidth=.8pt](.5,2)(.5,1.5)(1.5,1.5)(1.5,2)
 \psbezier[linewidth=.8pt](.5,1)(.5,1.5)(1.5,1.5)(1.5,1)
 \pscurve[linewidth=.8pt](.5,1)(.75,.5)(.5,0)
 \pscurve[linewidth=.8pt](1.5,1)(1.25,.5)(1.5,0)
 \psbezier[linewidth=.8pt](1.5,0)(1.5,-.25)(2,-.25)(2,0)
 \psbezier[linewidth=.8pt](1.5,2)(1.5,2.25)(2,2.25)(2,2)
 \psline[linewidth=.8pt](2,0)(2,2)
 \psbezier[linewidth=.8pt](.5,0)(.5,-.25)(0,-.25)(0,0)
 \psbezier[linewidth=.8pt](.5,2)(.5,2.25)(0,2.25)(0,2)
 \psline[linewidth=.8pt](0,0)(0,2)
  \end{pspicture}}="01";
  (20,-5)*+{
  \psset{xunit=.3cm,yunit=.3cm}
  \begin{pspicture}[.4](2,2)
 \pscurve[linewidth=.8pt](.5,2)(.75,1.5)(.5,1)
 \pscurve[linewidth=.8pt](1.5,2)(1.25,1.5)(1.5,1)
  \psbezier[linewidth=.8pt](.5,0)(.5,.5)(1.5,.5)(1.5,00)
 \psbezier[linewidth=.8pt](.5,1)(.5,.5)(1.5,.5)(1.5,1)
 \psbezier[linewidth=.8pt](1.5,0)(1.5,-.25)(2,-.25)(2,0)
 \psbezier[linewidth=.8pt](1.5,2)(1.5,2.25)(2,2.25)(2,2)
 \psline[linewidth=.8pt](2,0)(2,2)
 \psbezier[linewidth=.8pt](.5,0)(.5,-.25)(0,-.25)(0,0)
 \psbezier[linewidth=.8pt](.5,2)(.5,2.25)(0,2.25)(0,2)
 \psline[linewidth=.8pt](0,0)(0,2)
\end{pspicture}}="10";
  (20,-15)*+{
    \psset{xunit=.3cm,yunit=.3cm}
  \begin{pspicture}[.4](2,2)
 \pscurve[linewidth=.8pt](.5,2)(.65,1.5)(.5,1)
 \pscurve[linewidth=.8pt](1.5,2)(1.35,1.5)(1.5,1)
  \pscurve[linewidth=.8pt](.5,1)(.65,.5)(.5,0)
 \pscurve[linewidth=.8pt](1.5,1)(1.35,.5)(1.5,0)
 \psbezier[linewidth=.8pt](1.5,0)(1.5,-.25)(2,-.25)(2,0)
 \psbezier[linewidth=.8pt](1.5,2)(1.5,2.25)(2,2.25)(2,2)
 \psline[linewidth=.8pt](2,0)(2,2)
 \psbezier[linewidth=.8pt](.5,0)(.5,-.25)(0,-.25)(0,0)
 \psbezier[linewidth=.8pt](.5,2)(.5,2.25)(0,2.25)(0,2)
 \psline[linewidth=.8pt](0,0)(0,2)
\end{pspicture}}="11";
    {\ar^0 "t";"0"};
    {\ar_1 "t";"1"};
    {\ar^0 "0";"00"};
    {\ar_1 "0";"01"};
    {\ar^0 "1";"10"};
    {\ar_1 "1";"11"};
 \endxy
 }

\newcommand{\tangleresolutions}{
 \xy
  (-5,0)*+{
  \psset{xunit=.3cm,yunit=.3cm}
 \begin{pspicture}[.4](2,2)
 \psbezier[linewidth=.8pt](.5,0)(.5,.5)(1.5,.5)(1.5,1)
 \psbezier[linewidth=.8pt](1.5,2)(1.5,1.5)(.5,1.5)(.5,1)
 \pspolygon[linecolor=white,fillstyle=solid,fillcolor=white](.9,0)(.9,2)(1.1,2)(1.1,0)(.9,0)
 \psbezier[linewidth=.8pt](1.5,0)(1.5,.5)(.5,.5)(.5,1)
 \psbezier[linewidth=.8pt](.5,2)(.5,1.5)(1.5,1.5)(1.5,1)
\end{pspicture}
  }="t";
  (5,10)*+{\psset{xunit=.3cm,yunit=.3cm}
  \begin{pspicture}[.4](2,2)
 \psbezier[linewidth=.8pt](.5,0)(.5,.5)(1.5,.5)(1.5,1)
 \pspolygon[linecolor=white,fillstyle=solid,fillcolor=white](.9,0)(.9,2)(1.1,2)(1.1,0)(.9,0)
 \psbezier[linewidth=.8pt](1.5,0)(1.5,.5)(.5,.5)(.5,1)
 \psbezier[linewidth=.8pt](.5,2)(.5,1.5)(1.5,1.5)(1.5,2)
 \psbezier[linewidth=.8pt](.5,1)(.5,1.5)(1.5,1.5)(1.5,1)
  \end{pspicture}
  }="0";
  (5,-10)*+{\psset{xunit=.3cm,yunit=.3cm}
  \begin{pspicture}[.4](2,2)
 \psbezier[linewidth=.8pt](.5,0)(.5,.5)(1.5,.5)(1.5,1)
 \pspolygon[linecolor=white,fillstyle=solid,fillcolor=white](.9,0)(.9,2)(1.1,2)(1.1,0)(.9,0)
 \psbezier[linewidth=.8pt](1.5,0)(1.5,.5)(.5,.5)(.5,1)
 \pscurve[linewidth=.8pt](.5,2)(.75,1.5)(.5,1)
 \pscurve[linewidth=.8pt](1.5,2)(1.25,1.5)(1.5,1)
\end{pspicture}}="1";
  (20,15)*+{\psset{xunit=.3cm,yunit=.3cm}
    \begin{pspicture}[.4](2,2)
 \psbezier[linewidth=.8pt](.5,2)(.5,1.5)(1.5,1.5)(1.5,2)
 \psbezier[linewidth=.8pt](.5,1)(.5,1.5)(1.5,1.5)(1.5,1)
 \psbezier[linewidth=.8pt](.5,0)(.5,.5)(1.5,.5)(1.5,0)
 \psbezier[linewidth=.8pt](.5,1)(.5,.5)(1.5,.5)(1.5,1)
  \end{pspicture}
  }="00";
  (20,5)*+{
  \psset{xunit=.3cm,yunit=.3cm}
    \begin{pspicture}[.4](2,2)
 \psbezier[linewidth=.8pt](.5,2)(.5,1.5)(1.5,1.5)(1.5,2)
 \psbezier[linewidth=.8pt](.5,1)(.5,1.5)(1.5,1.5)(1.5,1)
 \pscurve[linewidth=.8pt](.5,1)(.75,.5)(.5,0)
 \pscurve[linewidth=.8pt](1.5,1)(1.25,.5)(1.5,0)
  \end{pspicture}}="01";
  (20,-5)*+{
  \psset{xunit=.3cm,yunit=.3cm}
  \begin{pspicture}[.4](2,2)
 \pscurve[linewidth=.8pt](.5,2)(.75,1.5)(.5,1)
 \pscurve[linewidth=.8pt](1.5,2)(1.25,1.5)(1.5,1)
  \psbezier[linewidth=.8pt](.5,0)(.5,.5)(1.5,.5)(1.5,00)
 \psbezier[linewidth=.8pt](.5,1)(.5,.5)(1.5,.5)(1.5,1)
\end{pspicture}}="10";
  (20,-15)*+{
    \psset{xunit=.3cm,yunit=.3cm}
  \begin{pspicture}[.4](2,2)
 \pscurve[linewidth=.8pt](.5,2)(.65,1.5)(.5,1)
 \pscurve[linewidth=.8pt](1.5,2)(1.35,1.5)(1.5,1)
  \pscurve[linewidth=.8pt](.5,1)(.65,.5)(.5,0)
 \pscurve[linewidth=.8pt](1.5,1)(1.35,.5)(1.5,0)
\end{pspicture}}="11";
{\ar^0 "t";"0"};{\ar_1 "t";"1"};{\ar^0 "0";"00"};{\ar_1
"0";"01"};{\ar^0 "1";"10"};{\ar_1 "1";"11"};
 \endxy
 }

\newcommand{\saddlel}{
 \pscustom[fillstyle=solid,fillcolor=lightgray]{
      \psline(.5,2.2)(.5,5)
      \psline(.5,5)(-.5,5)
      \psbezier(-.5,5)(-.5,4)(-4.5,3)(-4.5,0)
      \psline(-4.5,0)(-3.5,0)
      \psbezier(-3.5,0)(-3.5,1)(-2,3)(-1,3)
      \psbezier(-1,3)(-.5,3)(-.5,2.2)(-1,2.2)
 }
 \pscustom[fillstyle=solid,fillcolor=lightgray]{
      \psbezier(-1,2.2)(-1.5,2.2)(-3.5,3)(-3.5,5)
      \psline(-4.5,5)
      \psbezier(-4.5,5)(-4.5,2)(-.5,1)(-.5,0)
      \psline(.5,0)
      \psline(.5,2.2)
  }
 }

 \newcommand{\saddler}{
 \pscustom[fillstyle=solid,fillcolor=lightgray]{
      \psline(-.5,2.2)(-.5,5)
      \psline(-.5,5)(.5,5)
      \psbezier(.5,5)(.5,4)(4.5,3)(4.5,0)
      \psline(4.5,0)(3.5,0)
      \psbezier(3.5,0)(3.5,1)(2,3)(1,3)
      \psbezier(1,3)(.5,3)(.5,2.2)(1,2.2)
 }
 \pscustom[fillstyle=solid,fillcolor=lightgray]{
      \psbezier(1,2.2)(1.5,2.2)(3.5,3)(3.5,5)
      \psline(4.5,5)
      \psbezier(4.5,5)(4.5,2)(.5,1)(.5,0)
      \psline(-.5,0)
      \psline(-.5,2.2)
  }
 }

\newcommand{\mediummultl}{
  \pscustom[fillstyle=solid,fillcolor=lightgray]{
    \psline(-0.5,2.5)(0.5,2.5)
    \psbezier(0.5,2.5)(1.5,1.2)(2.5,1.2)(3.5,2.5)
    \psline(3.5,2.5)(4.5,2.5)
    \psbezier(4.5,2.5)(4.0,1.2)(3.0,1.6)(2.5,0)
    \psline(2.5,0)(1.5,0)
    \psbezier(1.5,0)(1.0,1.6)(0.0,1.2)(-0.5,2.5)
  }
}
\newcommand{\mediumcomultl}{
  \pscustom[fillstyle=solid,fillcolor=lightgray]{
    \psline(-0.5,0)(0.5,0)
    \psbezier(0.5,0)(1.5,1.3)(2.5,1.3)(3.5,0)
    \psline(3.5,0)(4.5,0)
    \psbezier(4.5,0)(4.0,1.3)(3.0,0.9)(2.5,2.5)
    \psline(2.5,2.5)(1.5,2.5)
    \psbezier(1.5,2.5)(1.0,0.9)(0.0,1.3)(-0.5,0)
  }
}

\newcommand{\widecrossl}{
  \pscustom[fillstyle=solid,fillcolor=lightgray]{
  \psline(-2.5,0)(-1.5,0)
  \psbezier(-1.5,0)(-1.5,.5)(2.5,1.7)(2.5,2.5)
  \psline(1.5,2.5)
  \psbezier(1.5,2.5)(1.5,2)(-2.5,.8)(-2.5,0)
  }
  \pscustom[fillstyle=solid,fillcolor=lightgray]{
  \psline(2.5,0)(1.5,0)
  \psbezier(1.5,0)(1.5,.5)(-2.5,1.7)(-2.5,2.5)
  \psline(-1.5,2.5)
  \psbezier(-1.5,2.5)(-1.5,2)(2.5,.8)(2.5,0)
  }
}

\newcommand{\Ronei}{
\psset{xunit=.25cm,yunit=.25cm}\begin{pspicture}[.25](2,1.5)
  \psbezier[linewidth=.6pt](.5,0)(.5,.25)(.85,.25)(.85,.5)
  \psbezier[linewidth=.6pt](.85,.5)(.85,.75)(.5,.75)(.5,1)
  \psbezier[linewidth=.6pt](1.5,0)(1.5,.25)(1.15,.25)(1.15,.5)
  \psbezier[linewidth=.6pt](1.5,1)(1.5,.75)(1.15,.75)(1.15,.5)
  \psbezier[linewidth=.6pt](.5,1)(.5,1.45)(1.5,1.45)(1.5,1)
\end{pspicture}
}

\newcommand{\Roneii}{
\psset{xunit=.25cm,yunit=.25cm}
\begin{pspicture}[.25](2,1.5)
 \psbezier[linewidth=.6pt](1.5,0)(1.5,.5)(.5,.5)(.5,1)
 \pspolygon[linecolor=white,fillstyle=solid,fillcolor=white](.85,0)(.85,2)(1.15,2)(1.15,0)(.9,0)
  \psbezier[linewidth=.6pt](.5,1)(.5,1.45)(1.5,1.45)(1.5,1)
  \psbezier[linewidth=.6pt](.5,0)(.5,.5)(1.5,.5)(1.5,1)
\end{pspicture}
}

\newcommand{\Reidemeister}{
\psset{xunit=.5cm,yunit=.5cm}
\begin{pspicture}[.4](2,1.5)
  \psbezier[linewidth=.6pt](.5,0)(.5,.5)(1.5,.5)(1.5,1)
 \pspolygon[linecolor=white,fillstyle=solid,fillcolor=white](.9,0)(.9,2)(1.1,2)(1.1,0)(.9,0)
  \psbezier[linewidth=.6pt](.5,1)(.5,1.45)(1.5,1.45)(1.5,1)
  \psbezier[linewidth=.6pt](1.5,0)(1.5,.5)(.5,.5)(.5,1)
\end{pspicture}
  \; \leftrightsquigarrow \;
  \begin{pspicture}[.4](2,1.5)
  \psbezier[linewidth=.6pt](.5,0)(.5,.25)(.85,.25)(.85,.5)
  \psbezier[linewidth=.6pt](.85,.5)(.85,.75)(.5,.75)(.5,1)
  \psbezier[linewidth=.6pt](1.5,0)(1.5,.25)(1.15,.25)(1.15,.5)
  \psbezier[linewidth=.6pt](1.5,1)(1.5,.75)(1.15,.75)(1.15,.5)
  \psbezier[linewidth=.6pt](.5,1)(.5,1.45)(1.5,1.45)(1.5,1)
\end{pspicture}
 \; \leftrightsquigarrow \;
\begin{pspicture}[.4](2,1.5)
 \psbezier[linewidth=.6pt](1.5,0)(1.5,.5)(.5,.5)(.5,1)
 \pspolygon[linecolor=white,fillstyle=solid,fillcolor=white](.9,0)(.9,2)(1.1,2)(1.1,0)(.9,0)
  \psbezier[linewidth=.6pt](.5,1)(.5,1.45)(1.5,1.45)(1.5,1)
  \psbezier[linewidth=.6pt](.5,0)(.5,.5)(1.5,.5)(1.5,1)
\end{pspicture}
\qquad
 \begin{pspicture}[.4](2,2)
  \psline[linewidth=.6pt](.5,0)(.5,2)
  \psline[linewidth=.6pt](1.5,0)(1.5,2)
\end{pspicture}
 \; \leftrightsquigarrow \;
 \begin{pspicture}[.4](2,2)
 \psbezier[linewidth=.6pt](1.5,0)(1.5,.5)(.5,.5)(.5,1)
 \psbezier[linewidth=.6pt](1.5,2)(1.5,1.5)(.5,1.5)(.5,1)
 \pspolygon[linecolor=white,fillstyle=solid,fillcolor=white](.9,0)(.9,2)(1.1,2)(1.1,0)(.9,0)
  \psbezier[linewidth=.6pt](.5,0)(.5,.5)(1.5,.5)(1.5,1)
 \psbezier[linewidth=.6pt](.5,2)(.5,1.5)(1.5,1.5)(1.5,1)
\end{pspicture}
 \qquad
 \begin{pspicture}[.4](3,2)
 \psbezier[linewidth=.8pt](1.5,2)(1.5,1.5)(.5,1.5)(.5,1)
 \psbezier[linewidth=.8pt](1.5,0)(1.5,.5)(.5,.5)(.5,1)
   \pspolygon[linecolor=white,fillstyle=solid,fillcolor=white](.8,0)(.8,2)(1.1,2)(1.1,0)(.8,0)
   \psbezier[linewidth=.8pt](2.5,2)(2.5,1.25)(.5,.5)(.5,0)
   \pspolygon[linecolor=white,fillstyle=solid,fillcolor=white](1.5,0)(1.5,2)(1.8,2)(1.8,0)(1.5,0)
 \psbezier[linewidth=.8pt](.5,2)(.5,1.5)(2.5,.75)(2.5,0)
\end{pspicture}
\; \leftrightsquigarrow \;
\begin{pspicture}[.4](3,2)
   \psbezier[linewidth=.8pt](1.5,0)(1.5,.5)(2.5,.5)(2.5,1)
   \psbezier[linewidth=.8pt](1.5,2)(1.5,1.5)(2.5,1.5)(2.5,1)
   \pspolygon[linecolor=white,fillstyle=solid,fillcolor=white](2.2,0)(2.2,2)(1.9,2)(1.9,0)(2.2,0)
   \psbezier[linewidth=.8pt](2.5,2)(2.5,1.5)(.5,.75)(.5,0)
   \pspolygon[linecolor=white,fillstyle=solid,fillcolor=white](1.5,0)(1.5,2)(1.2,2)(1.2,0)(1.5,0)
 \psbezier[linewidth=.8pt](.5,2)(.5,1.25)(2.5,.5)(2.5,0)
\end{pspicture}
}

 \newcommand{\holec}{
      \pscustom[fillstyle=gradient,
    gradbegin=white, gradend=gray,gradmidpoint=0,gradangle=70]{
        \psbezier(1.5,2)(1.5,.6)(.4,1.1)(.5,0)
        \psbezier(.5,0)(.4,-.25)(-.4,-.25)(-.5,0)
        \psbezier(-0.5,0)(-.4,1.1)(-1.5,.6)(-1.5,2)
        \psline(-.5,2)
        \psbezier(-.5,2)(-.55,1.2)(0.55,1.2)(.5,2)
        \psline(1.5,2)
    }
     \begin{psclip}{
 \pspolygon[linestyle=none](.5,0)(.5,.3)(-.5,.3)(-.5,0)(.5,0)
 }
 \psellipse[linestyle=dotted](0,0)(.5,0.2)
 \end{psclip}
 \pscustom[fillstyle=gradient,
    gradbegin=white, gradend=gray,gradmidpoint=0,gradangle=110,linestyle=none]{
        \psbezier(1.5,2)(1.5,3.4)(.4,2.9)(.5,4)
        \psline(-0.5,4)
        \psbezier(-0.5,4)(-.4,2.9)(-1.5,3.4)(-1.5,2)
        \psbezier(-1.5,2)(-1.4,1.75)(-.6,1.75)(-.5,2)
        \psbezier(-.5,2)(-.55,2.8)(0.55,2.8)(.5,2)
        \psbezier(.5,2)(.6,1.75)(1.4,1.75)(1.5,2)
    }
 \psbezier(1.5,2)(1.5,3.4)(.4,2.9)(.5,4)
 \psbezier(-0.5,4)(-.4,2.9)(-1.5,3.4)(-1.5,2)
        \psbezier(-.5,2)(-.55,2.8)(0.55,2.8)(.5,2)
  \psellipse[fillcolor=lightgray,fillstyle=gradient,
        gradbegin=lightgray, gradend=gray,gradmidpoint=1,gradangle=110](0,4)(.5,.2)
  }

\newcommand{\multl}{
  \pscustom[fillcolor=lightgray, fillstyle=solid]{
        \psbezier(1.5,2.5)(1.5,1.1)(.4,1.6)(.5,0)
        \psline(-0.5,0)
        \psbezier(-0.5,0)(-.4,1.6)(-1.5,1.1)(-1.5,2.5)
        \psline(-.5,2.5)
        \psbezier(-.5,2.5)(-.6,1.5)(0.6,1.5)(.5,2.5)
        \psline(1.5,2.5)
    }
}

\newcommand{\comultl}{
  \pscustom[fillcolor=lightgray, fillstyle=solid]{
        \psbezier(1.5,0)(1.5,1.4)(.4,.9)(.5,2.5)
        \psline(-0.5,2.5)
        \psbezier(-0.5,2.5)(-.4,.9)(-1.5,1.4)(-1.5,0)
        \psline(-.5,0)
        \psbezier(-.5,0)(-.6,1)(0.6,1)(.5,0)
        \psline(1.5,0)
    }
}

\newcommand{\ctl}{
  \begin{psclip}{
    \pscustom{
        \psline(-.58,2)(-.58,0)
        \psline(-.58,0)(.42,0)
        \psline(.42,0)(.42,2)
        \psellipse(-.08,2)(.5,.2)
    }
  }
    \pspolygon[fillcolor=lightgray,fillstyle=gradient,
    gradbegin=lightgray,gradend=gray,gradmidpoint=1,gradangle=110](-.58,0)(-.58,2.4)(.42,2.4)(.42,0)(-.58,0)
 \end{psclip}
 \pscustom[fillcolor=lightgray,fillstyle=gradient,
        gradbegin=white, gradend=gray,gradmidpoint=0,gradangle=88]{
    \psline(-.58,2)(-.58,0)
    \psbezier(-.58,0)(-.48,.5)(-.48,.7)(-.08,1)
    \psbezier(-.08,1)(.32,.7)(.32,.5)(.42,0)
    \psellipse(-.08,2)(.5,.2)
 }
 \psellipse[fillcolor=lightgray,fillstyle=gradient,
        gradbegin=lightgray, gradend=gray,gradmidpoint=1,gradangle=110](-.08,2)(.5,.2)
}

\newcommand{\ltc}{
    \pspolygon[fillcolor=lightgray,fillstyle=gradient,
    gradbegin=lightgray,gradend=gray,gradmidpoint=1,gradangle=60](.58,2)(.58,.4)(-.42,.4)(-.42,2)(.58,2)
 \pscustom[fillcolor=lightgray,fillstyle=gradient,
        gradbegin=white, gradend=gray,gradmidpoint=0,gradangle=88]{
    \psline(.58,0)(.58,2)
    \psbezier(.58,2)(.48,1.5)(.48,1.3)(.08,1)
    \psbezier(.08,1)(-.32,1.3)(-.32,1.5)(-.42,2)
    \psline(-.42,0)
    \psbezier(-.42,0)(-.32,-.25)(.48,-.25)(.58,0)
 }
 \begin{psclip}{
 \pspolygon[linestyle=none](.58,0)(.58,.3)(-.42,.3)(-.42,0)(.58,0)
 }
 \psellipse[linestyle=dotted](.08,0)(.5,0.2)
 \end{psclip}
}

 \newcommand{\birthl}{
 \pscustom[fillcolor=lightgray, fillstyle=solid]{
        \psbezier(-.5,0)(-.5,.9)(0.5,.9)(.5,0)
        \psline(-.5,0)
    }
 }

  \newcommand{\deathl}{
 \pscustom[fillcolor=lightgray, fillstyle=solid]{
        \psbezier(-.5,0)(-.5,-.9)(0.5,-.9)(.5,0)
        \psline(-.5,0)    }
 }

\newcommand{\identl}{
    \pspolygon[fillcolor=lightgray,fillstyle=solid](-.5,0)(.5,0)(.5,2.5)(-.5,2.5)(-.5,0)
}

\newcommand{\medidentl}{
    \pspolygon[fillcolor=lightgray,fillstyle=solid](-.5,0)(.5,0)(.5,2)(-.5,2)(-.5,0)
}

\newcommand{\crossl}{
    \pspolygon[fillcolor=lightgray,fillstyle=solid](-1.5,0)(-.5,0)(1.5,2.5)(.5,2.5)(-1.5,0)
    \pspolygon[fillcolor=lightgray,fillstyle=solid](1.5,0)(.5,0)(-1.5,2.5)(-.5,2.5)(1.5,0)
}

\newcommand{\curverightl}{
  \pscustom[fillcolor=lightgray, fillstyle=solid]{
        \psbezier(1.5,2.5)(1.5,1.5)(.4,1.3)(.5,0)
        \psline(-0.5,0)
        \psbezier(-.5,0)(-.6,1.3)(.5,1.5)(.5,2.5)
        \psline(1.5,2.5)
    }
}

\newcommand{\curveleftl}{
  \pscustom[fillcolor=lightgray, fillstyle=solid]{
        \psbezier(-1.5,2.5)(-1.5,1.5)(-.4,1.3)(-.5,0)
        \psline(0.5,0)
        \psbezier(.5,0)(.6,1.3)(-.5,1.5)(-.5,2.5)
        \psline(-1.5,2.5)
    }
}

\newgray{mygray}{0.65}
\newgray{mylightgray}{0.85}

\newcommand{\bigsaddle}{
  \pscustom[fillcolor=mygray,fillstyle=solid]{
    \psbezier(0.7,1.8)(0.7,2.0)(0.4,2.1)(0.2,2.1)
    \psline(0.2,2.1)(0.2,1.0)(0.7,1.0)(0.7,1.8)
  }
  \pscustom[fillcolor=mygray,fillstyle=solid]{
    \psbezier(2.0,1.8)(2.0,2.0)(2.5,2.1)(2.7,2.1)
    \psline(2.7,2.1)(2.7,0.6)(2.5,0.55)(2.0,0.6)(2.0,1.8)
  }
  \pscustom[fillcolor=mylightgray,fillstyle=solid]{
    \psline(0.0,1.5)(0.0,0.0)
    \psbezier(0.0,0.0)(0.5,0.2)(2.0,0.2)(2.5,0.0)
    \psline(2.5,0.0)(2.5,1.5)
    \psbezier(2.5,1.5)(2.2,1.5)(2.0,1.6)(2.0,1.8)
    \psbezier(2.0,1.8)(1.8,1.0)(0.9,1.0)(0.7,1.8)
    \psbezier(0.7,1.8)(0.7,1.6)(0.3,1.5)(0.0,1.5)
  }
  \psline[linestyle=dotted,dotsep=2pt](0.2,1.52)(0.2,0.6)
  \psbezier[linestyle=dotted,dotsep=2pt](0.2,0.6)(0.7,0.4)(2.2,0.4)(2.5,0.55)
  \psline[linestyle=dotted,dotsep=2pt](2.5,0.55)(2.7,0.6)
}

\newcommand{\bigsaddleorient}{
  \pscustom[fillcolor=mygray,fillstyle=solid]{
    \psbezier(0.7,1.8)(0.7,2.0)(0.4,2.1)(0.2,2.1)
    \psline(0.2,2.1)(0.2,1.0)(0.7,1.0)(0.7,1.8)
  }
  \pscustom[fillcolor=mygray,fillstyle=solid]{
    \psbezier(2.0,1.8)(2.0,2.0)(2.5,2.1)(2.7,2.1)
    \psline(2.7,2.1)(2.7,0.6)(2.5,0.55)(2.0,0.6)(2.0,1.8)
  }
  \pscustom[fillcolor=mylightgray,fillstyle=solid]{
    \psline(0.0,1.5)(0.0,0.0)
    \psbezier(0.0,0.0)(0.5,0.2)(1.1,0.2)(1.25,0.2)
    \psbezier(1.25,0.2)(1.4,0.2)(2.0,0.2)(2.5,0.0)
    \psline(2.5,0.0)(2.5,1.5)
    \psbezier(2.5,1.5)(2.2,1.5)(2.0,1.6)(2.0,1.8)
    \psbezier(2.0,1.8)(1.8,1.0)(0.9,1.0)(0.7,1.8)
    \psbezier(0.7,1.8)(0.7,1.6)(0.3,1.5)(0.0,1.5)
  }
  \psline[linestyle=dotted,dotsep=2pt](0.2,1.52)(0.2,0.6)
  \psbezier[linestyle=dotted,dotsep=2pt](0.2,0.6)(0.7,0.4)(1.3,0.4)(1.45,0.4)
  \psbezier[linestyle=dotted,dotsep=2pt,arrowsize=4pt 6,arrowinset=0.7]{<-}(1.45,0.4)(1.6,0.4)(2.2,0.4)(2.5,0.55)
  \psline[linestyle=dotted,dotsep=2pt](2.5,0.55)(2.7,0.6)
  \psbezier[arrowsize=4pt 6,arrowinset=0.7]{->}(0.0,0.0)(0.5,0.2)(1.1,0.2)(1.25,0.2)
  \psbezier(1.25,0.2)(1.4,0.2)(2.0,0.2)(2.5,0.0)
  \psbezier[arrowsize=4pt 6,arrowinset=0.7]{<-}(2.0,1.8)(2.0,2.0)(2.5,2.1)(2.7,2.1)
  \psbezier[arrowsize=4pt 6,arrowinset=0.7]{<-}(0.7,1.8)(0.7,1.6)(0.3,1.5)(0.0,1.5)
  }

\newcommand{\bigsaddleb}{
  \pscustom[fillcolor=mygray,fillstyle=solid]{
    \psbezier(0.7,1.8)(0.7,2.0)(0.4,2.1)(0.2,2.1)
    \psline(0.2,2.1)(0.2,1.0)(0.7,1.0)(0.7,1.8)
  }
  \pscustom[fillcolor=mygray,fillstyle=solid]{
    \psbezier(2.0,1.8)(2.0,2.0)(2.5,2.1)(2.7,2.1)
    \psline(2.7,2.1)(2.7,0.6)(2.5,0.55)(2.0,0.6)(2.0,1.8)
  }
  \pscustom[fillcolor=mylightgray,fillstyle=solid]{
    \psline(0.0,1.5)(0.0,0.0)
    \psbezier(0.0,0.0)(0.5,0.2)(2.0,0.2)(2.5,0.0)
    \psline(2.5,0.0)(2.5,1.5)
    \psbezier(2.5,1.5)(2.2,1.5)(2.0,1.6)(2.0,1.8)
    \psbezier(2.0,1.8)(1.8,1.0)(0.9,1.0)(0.7,1.8)
    \psbezier(0.7,1.8)(0.7,1.6)(0.3,1.5)(0.0,1.5)
  }
  \psline[linestyle=dotted,dotsep=2pt](0.2,1.52)(0.2,0.0)
  \psbezier[linestyle=dotted,dotsep=2pt](0.2,0.6)(0.7,0.4)(2.2,0.4)(2.5,0.55)
  \psline[linestyle=dotted,dotsep=2pt](2.5,0.55)(2.7,0.6)
}

\newcommand{\bigsaddlet}{
  \pscustom[fillcolor=mygray,fillstyle=solid]{
    \psbezier(0.7,1.8)(0.7,2.0)(0.4,2.1)(0.2,2.1)
    \psline(0.2,2.1)(0.2,1.0)(0.7,1.0)(0.7,1.8)
  }
  \pscustom[fillcolor=mygray,fillstyle=solid]{
    \psbezier(2.0,1.8)(2.0,2.0)(2.5,2.1)(2.7,2.1)
    \psline(2.7,2.1)(2.7,0.6)(2.5,0.55)(2.0,0.6)(2.0,1.8)
  }
  \pscustom[fillcolor=mylightgray,fillstyle=solid]{
    \psline(0.0,1.5)(0.0,0.0)
    \psbezier(0.0,0.0)(0.5,0.2)(2.0,0.2)(2.5,0.0)
    \psline(2.5,0.0)(2.5,1.5)
    \psbezier(2.5,1.5)(2.2,1.5)(2.0,1.6)(2.0,1.8)
    \psbezier(2.0,1.8)(1.8,1.0)(0.9,1.0)(0.7,1.8)
    \psbezier(0.7,1.8)(0.7,1.6)(0.3,1.5)(0.0,1.5)
  }
  \psline[linestyle=dotted,dotsep=2pt](0.2,1.52)(0.2,0.6)
  \psbezier[linestyle=dotted,dotsep=2pt](0.2,0.6)(0.7,0.4)(2.2,0.4)(2.5,0.55)
  \psline[linestyle=dotted,dotsep=2pt](2.5,0.55)(2.7,0.6)

  \psbezier(0.0,1.5)(0.3,1.2)(0.4,1.2)(0.9,1.4)
  \psbezier[linestyle=dotted,dotsep=2pt](0.9,1.4)(0.8,1.2)(0.5,0.8)(0.2,0.6)
  \psbezier(0.0,0.0)(0.3,0.2)(0.9,1.1)(1.0,1.3)
  \psbezier[linestyle=dotted,dotsep=2pt](1.0,1.3)(0.9,1.1)(0.5,0.8)(0.2,0.6)
  \psbezier(0.0,0.0)(0.4,0.2)(1.25,0.8)(1.35,1.2)
  \psbezier[linestyle=dotted,dotsep=2pt](1.35,1.2)(1.45,0.8)(2.4,0.7)(2.7,0.6)
  \psbezier(2.5,0.0)(2.2,0.2)(1.7,1.1)(1.6,1.27)
  \psbezier[linestyle=dotted,dotsep=2pt](1.6,1.27)(1.7,1.0)(2.4,0.8)(2.7,0.6)
  \psbezier(2.5,1.5)(2.2,1.2)(2.1,1.2)(1.72,1.33)
  \psbezier[linestyle=dotted,dotsep=2pt](1.72,1.33)(1.8,1.2)(2.4,0.8)(2.7,0.6)
}

\newcommand{\bigarc}{
  \pscustom[fillcolor=mygray,fillstyle=solid]{
    \psbezier(0.7,1.8)(0.7,2.0)(0.4,2.1)(0.2,2.1)
    \psline(0.2,2.1)(0.2,1.0)(0.7,1.0)(0.7,1.8)
  }
  \pscustom[fillcolor=mylightgray,fillstyle=solid]{
    \psline(0.0,1.5)(0.0,0.0)
    \psbezier(0.0,0.0)(0.2,0.0)(0.7,0.1)(0.7,0.3)
    \psline(0.7,0.3)(0.7,1.8)
    \psbezier(0.7,1.8)(0.7,1.6)(0.3,1.5)(0.0,1.5)
  }
  \psline[linestyle=dotted,dotsep=2pt](0.2,1.52)(0.2,0.6)
  \psbezier[linestyle=dotted,dotsep=2pt](0.7,0.3)(0.7,0.5)(0.4,0.6)(0.2,0.6)
}

\newcommand{\bigarcb}{
  \pscustom[fillcolor=mygray,fillstyle=solid]{
    \psbezier(0.7,1.8)(0.7,2.0)(0.4,2.1)(0.2,2.1)
    \psline(0.2,2.1)(0.2,1.0)(0.7,1.0)(0.7,1.8)
  }
  \pscustom[fillcolor=mylightgray,fillstyle=solid]{
    \psline(0.0,1.5)(0.0,0.0)
    \psbezier(0.0,0.0)(0.2,0.0)(0.7,0.1)(0.7,0.3)
    \psline(0.7,0.3)(0.7,1.8)
    \psbezier(0.7,1.8)(0.7,1.6)(0.3,1.5)(0.0,1.5)
  }
  \psline[linestyle=dotted,dotsep=2pt](0.2,1.52)(0.2,0.0)
  \psbezier[linestyle=dotted,dotsep=2pt](0.7,0.3)(0.7,0.5)(0.4,0.6)(0.2,0.6)
}

\newcommand{\bigarcc}{
  \pscustom[fillcolor=mygray,fillstyle=solid]{
    \psbezier(0.7,1.8)(0.7,2.0)(0.4,2.1)(0.2,2.1)
    \psline(0.2,2.1)(0.2,1.0)(0.7,1.0)(0.7,1.8)
  }
  \pscustom[fillcolor=mylightgray,fillstyle=solid]{
    \psline(0.0,1.5)(0.0,0.0)
    \psbezier(0.0,0.0)(0.2,0.0)(0.7,0.1)(0.7,0.3)
    \psline(0.7,0.3)(0.7,1.8)
    \psbezier(0.7,1.8)(0.7,1.6)(0.3,1.5)(0.0,1.5)
  }
  \psline[linestyle=dotted,dotsep=2pt](0.2,1.52)(0.2,0.6)
  \psline[linestyle=dotted,dotsep=2pt](0.2,0.6)(0.0,0.55)
  \psbezier[linestyle=dotted,dotsep=2pt](0.7,0.3)(0.7,0.5)(0.4,0.6)(0.2,0.6)
}

\newcommand{\bigarct}{
  \pscustom[fillcolor=mygray,fillstyle=solid]{
    \psbezier(0.7,1.8)(0.7,2.0)(0.4,2.1)(0.2,2.1)
    \psline(0.2,2.1)(0.2,1.0)(0.7,1.0)(0.7,1.8)
  }
  \pscustom[fillcolor=mylightgray,fillstyle=solid]{
    \psline(0.0,1.5)(0.0,0.0)
    \psbezier(0.0,0.0)(0.2,0.0)(0.7,0.1)(0.7,0.3)
    \psline(0.7,0.3)(0.7,1.8)
    \psbezier(0.7,1.8)(0.7,1.6)(0.3,1.5)(0.0,1.5)
  }
  \psline[linestyle=dotted,dotsep=2pt](0.2,1.52)(0.2,0.6)
  \psline[linestyle=dotted,dotsep=2pt](0.2,0.6)(0.0,0.55)
  \psbezier[linestyle=dotted,dotsep=2pt](0.7,0.3)(0.7,0.5)(0.4,0.6)(0.2,0.6)
  \psbezier(0.0,1.5)(0.3,1.2)(0.4,1.05)(0.7,1.05)
  \psbezier[linestyle=dotted,dotsep=2pt](0.7,1.05)(0.4,0.95)(0.5,0.9)(0.2,0.6)
  \psline[linestyle=dotted,dotsep=2pt](0.2,0.6)(0.0,0.75)
  \psline[linestyle=dotted,dotsep=2pt](0.2,0.6)(0.0,0.7)
  \psline[linestyle=dotted,dotsep=2pt](0.2,0.6)(0.0,0.65)
}

\newcommand{\bigtopleft}{
  \pscustom[fillcolor=mygray,fillstyle=solid]{
    \psbezier(0.7,2.4)(0.7,2.6)(0.4,2.7)(0.2,2.7)
    \psline(0.2,2.7)(0.2,1.0)(0.7,1.0)(0.7,2.4)
  }
  \pscustom[fillcolor=mylightgray,fillstyle=solid]{
    \psline(0.0,2.1)(0.0,0.0)
    \psbezier(0.0,0.0)(0.2,0.0)(0.7,0.1)(0.7,0.3)
    \psline(0.7,0.3)(0.7,2.4)
    \psbezier(0.7,2.4)(0.7,2.2)(0.3,2.1)(0.0,2.1)
  }
  \psline[linestyle=dotted,dotsep=2pt](0.2,2.12)(0.2,0.6)
  \psbezier[linestyle=dotted,dotsep=2pt](0.7,0.3)(0.7,0.5)(0.4,0.6)(0.2,0.6)
}

\newcommand{\bigtopleftb}{
  \pscustom[fillcolor=mygray,fillstyle=solid]{
    \psbezier(0.7,2.4)(0.7,2.6)(0.4,2.7)(0.2,2.7)
    \psline(0.2,2.7)(0.2,1.0)(0.7,1.0)(0.7,2.4)
  }
  \pscustom[fillcolor=mylightgray,fillstyle=solid]{
    \psline(0.0,2.1)(0.0,0.0)
    \psbezier(0.0,0.0)(0.2,0.0)(0.7,0.1)(0.7,0.3)
    \psline(0.7,0.3)(0.7,2.4)
    \psbezier(0.7,2.4)(0.7,2.2)(0.3,2.1)(0.0,2.1)
  }
  \psline[linestyle=dotted,dotsep=2pt](0.2,2.12)(0.2,0.0)
  \psbezier[linestyle=dotted,dotsep=2pt](0.7,0.3)(0.7,0.5)(0.4,0.6)(0.2,0.6)
}

\newcommand{\bigtopright}{
  \pscustom[fillcolor=mygray,fillstyle=solid]{
    \psbezier(1.2,2.5)(1.2,2.7)(0.95,2.7)(0.85,2.7)
    \psbezier(0.85,2.7)(0.75,2.7)(0.5,2.7)(0.5,2.5)
    \psline(0.5,2.5)(0.5,2.0)(1.2,2.0)(1.2,2.5)
  }
  \pscustom[fillcolor=mygray,fillstyle=solid]{
    \psline(0.0,0.3)(0.3,1.6)(0.95,1.6)
    \psbezier(0.95,1.4)(0.85,1.4)(0.7,1.3)(0.7,0.6)
    \psbezier(0.7,0.6)(0.5,0.6)(0.0,0.5)(0.0,0.3)
  }
  \pscustom[fillcolor=mylightgray,fillstyle=solid]{
    \psbezier(0.5,2.5)(0.5,2.0)(0.0,1.3)(0.0,0.3)
    \psbezier(0.0,0.3)(0.0,0.1)(0.3,0.0)(0.5,0.0)
    \psbezier(0.5,0.0)(0.5,0.7)(0.55,1.2)(0.75,1.2)
    \psbezier(0.75,1.2)(0.95,1.2)(1.0,0.7)(1.0,0.0)
    \psbezier(1.0,0.0)(1.2,0.0)(1.7,0.1)(1.7,0.3)
    \psbezier(1.7,0.3)(1.7,1.3)(1.2,2.0)(1.2,2.5)
    \psbezier(1.2,2.5)(1.2,2.3)(0.95,2.3)(0.85,2.3)
    \psbezier(0.85,2.3)(0.75,2.3)(0.5,2.3)(0.5,2.5)
  }
  \psbezier[linestyle=dotted,dotsep=2pt](1.7,0.3)(1.7,0.5)(1.4,0.6)(1.2,0.6)
  \psbezier[linestyle=dotted,dotsep=2pt](1.2,0.6)(1.2,1.3)(1.05,1.6)(0.95,1.6)
  \psbezier[linestyle=dotted,dotsep=2pt](0.95,1.6)(0.85,1.6)(0.7,1.3)(0.7,0.6)
  \psbezier[linestyle=dotted,dotsep=2pt](0.7,0.6)(0.5,0.6)(0.0,0.5)(0.0,0.3)
}

\newcommand{\bigtoprightb}{
  \pscustom[fillcolor=mygray,fillstyle=solid]{
    \psbezier(1.2,2.5)(1.2,2.7)(0.95,2.7)(0.85,2.7)
    \psbezier(0.85,2.7)(0.75,2.7)(0.5,2.7)(0.5,2.5)
    \psline(0.5,2.5)(0.5,2.0)(1.2,2.0)(1.2,2.5)
  }
  \pscustom[fillcolor=mygray,fillstyle=solid]{
    \psline(0.0,0.3)(0.3,1.6)(0.95,1.6)
    \psbezier(0.95,1.4)(0.85,1.4)(0.7,1.3)(0.7,0.6)
    \psbezier(0.7,0.6)(0.5,0.6)(0.0,0.5)(0.0,0.3)
  }
  \pscustom[fillcolor=mylightgray,fillstyle=solid]{
    \psbezier(0.5,2.5)(0.5,2.0)(0.0,1.3)(0.0,0.3)
    \psbezier(0.0,0.3)(0.0,0.1)(0.3,0.0)(0.5,0.0)
    \psbezier(0.5,0.0)(0.5,0.7)(0.55,1.2)(0.75,1.2)
    \psbezier(0.75,1.2)(0.95,1.2)(1.0,0.7)(1.0,0.0)
    \psbezier(1.0,0.0)(1.2,0.0)(1.7,0.1)(1.7,0.3)
    \psbezier(1.7,0.3)(1.7,1.3)(1.2,2.0)(1.2,2.5)
    \psbezier(1.2,2.5)(1.2,2.3)(0.95,2.3)(0.85,2.3)
    \psbezier(0.85,2.3)(0.75,2.3)(0.5,2.3)(0.5,2.5)
  }
  \psbezier[linestyle=dotted,dotsep=2pt](1.7,0.3)(1.7,0.5)(1.4,0.6)(1.2,0.6)
  \psline[linestyle=dotted,dotsep=2pt](1.2,0.6)(1.2,0.0)
  \psbezier[linestyle=dotted,dotsep=2pt](1.2,0.6)(1.2,1.3)(1.05,1.6)(0.95,1.6)
  \psbezier[linestyle=dotted,dotsep=2pt](0.95,1.6)(0.85,1.6)(0.7,1.3)(0.7,0.6)
  \psbezier[linestyle=dotted,dotsep=2pt](0.7,0.6)(0.5,0.6)(0.0,0.5)(0.0,0.3)
}

\newcommand{\bigbottom}{
  \pscustom[fillcolor=mygray,fillstyle=solid]{
    \psbezier(0.2,2.1)(0.7,1.9)(2.2,1.9)(2.7,2.1)
    \psbezier(2.7,2.1)(2.7,1.4)(2.75,0.9)(2.95,0.9)
    \psbezier(2.95,0.9)(3.15,0.9)(3.2,1.4)(3.2,2.1)
    \psbezier(3.2,2.1)(3.4,2.1)(3.7,2.0)(3.7,1.8)
    \psline(3.7,1.8)(3.7,0.3)
    \psbezier(3.7,0.3)(3.7,0.5)(2.9,0.6)(2.7,0.6)
    \psbezier(2.7,0.6)(2.2,0.6)(0.7,0.4)(0.2,0.6)
    \psline(0.2,0.6)(0.2,2.1)
  }
  \pscustom[fillcolor=mylightgray,fillstyle=solid]{
    \psline(0.0,1.5)(0.0,0.0)
    \psbezier(0.0,0.0)(0.5,0.2)(2.0,0.0)(2.5,0.0)
    \psbezier(2.5,0.0)(2.7,0.0)(3.7,0.1)(3.7,0.3)
    \psline(3.7,0.3)(3.7,1.8)
    \psbezier(3.7,1.8)(3.7,1.6)(3.2,1.5)(3.0,1.5)
    \psbezier(3.0,1.5)(3.0,0.8)(2.95,0.3)(2.75,0.3)
    \psbezier(2.75,0.3)(2.55,0.3)(2.5,0.8)(2.5,1.5)
    \psbezier(2.5,1.5)(2.0,1.7)(0.5,1.7)(0.0,1.5)
  }
  \psline[linestyle=dotted,dotsep=2pt](0.2,2.1)(0.2,0.6)
  \psbezier[linestyle=dotted,dotsep=2pt](0.2,0.6)(0.7,0.4)(2.2,0.6)(2.7,0.6)
  \psbezier[linestyle=dotted,dotsep=2pt](2.7,0.6)(2.9,0.6)(3.7,0.5)(3.7,0.3)
  \psbezier[linestyle=dotted,dotsep=2pt](2.95,0.9)(3.15,0.9)(3.2,1.4)(3.2,2.1)
}

\newcommand{\bigcomposed}{
  \pscustom[fillcolor=mygray,fillstyle=solid]{
    \psbezier(0.7,1.8)(0.7,2.0)(0.4,2.1)(0.2,2.1)
    \psline(0.2,2.1)(0.2,1.0)(0.7,1.0)(0.7,1.8)
  }
  \pscustom[fillcolor=mygray,fillstyle=solid]{
    \psbezier(2.0,1.8)(2.0,2.0)(2.5,2.1)(2.7,2.1)
    \psbezier(2.7,2.1)(2.9,2.1)(3.3,2.0)(3.3,1.8)
    \psline(3.3,1.8)(3.3,1.0)(2.0,1.0)(2.0,1.8)
  }
  \pscustom[fillcolor=mylightgray,fillstyle=solid]{
    \psline(0.0,1.5)(0.0,0.0)
    \psbezier(0.0,0.0)(0.5,0.2)(2.0,0.0)(2.5,0.0)
    \psbezier(2.5,0.0)(2.7,0.0)(3.3,0.1)(3.3,0.3)
    \psline(3.3,0.3)(3.3,1.8)
    \psbezier(3.3,1.8)(3.3,1.6)(2.8,1.5)(2.6,1.5)
    \psbezier(2.6,1.5)(2.3,1.5)(2.0,1.6)(2.0,1.8)
    \psbezier(2.0,1.8)(1.8,1.0)(0.9,1.0)(0.7,1.8)
    \psbezier(0.7,1.8)(0.7,1.6)(0.3,1.5)(0.0,1.5)
  }
  \psline[linestyle=dotted,dotsep=2pt](0.2,1.52)(0.2,0.6)
  \psbezier[linestyle=dotted,dotsep=2pt](0.2,0.6)(0.7,0.4)(2.2,0.6)(2.5,0.6)
  \psbezier[linestyle=dotted,dotsep=2pt](2.5,0.6)(2.7,0.6)(3.3,0.4)(3.3,0.2)
}

\newcommand{\multc}{
      \pscustom[fillstyle=gradient,
    gradbegin=white, gradend=gray,gradmidpoint=0,gradangle=70]{
        \psbezier(1.5,2.5)(1.5,1.1)(.4,1.6)(.5,0)
        \psbezier(.5,0)(.4,-.25)(-.4,-.25)(-.5,0)
        \psbezier(-0.5,0)(-.4,1.6)(-1.5,1.1)(-1.5,2.5)
        \psline(-.5,2.5)
        \psbezier(-.5,2.5)(-.6,1.5)(0.6,1.5)(.5,2.5)
        \psline(1.5,2.5)
    }
    \psellipse[fillcolor=lightgray,fillstyle=gradient,
        gradbegin=lightgray, gradend=gray,gradmidpoint=1,gradangle=110](-1,2.5)(.5,.2)
    \psellipse[fillcolor=lightgray,fillstyle=gradient,
        gradbegin=lightgray, gradend=gray,gradmidpoint=1,gradangle=110](1,2.5)(.5,.2)
     \begin{psclip}{
 \pspolygon[linestyle=none](.5,0)(.5,.3)(-.5,.3)(-.5,0)(.5,0)
 }
 \psellipse[linestyle=dotted](0,0)(.5,0.2)
 \end{psclip}
 }

\newcommand{\comultc}{
  \pscustom[fillstyle=gradient,
    gradbegin=white, gradend=gray,gradmidpoint=0,gradangle=110]{
        \psbezier(1.5,0)(1.5,1.4)(.4,.9)(.5,2.5)
        \psline(-0.5,2.5)
        \psbezier(-0.5,2.5)(-.4,.9)(-1.5,1.4)(-1.5,0)
        \psbezier(-1.5,0)(-1.4,-.25)(-.6,-.25)(-.5,0)
        \psbezier(-.5,0)(-.6,1)(0.6,1)(.5,0)
        \psbezier(.5,0)(.6,-.25)(1.4,-.25)(1.5,0)
    }
  \psellipse[fillcolor=lightgray,fillstyle=gradient,
        gradbegin=lightgray, gradend=gray,gradmidpoint=1,gradangle=110](0,2.5)(.5,.2)
\begin{psclip}{
 \pspolygon[linestyle=none](1.5,0)(1.5,.3)(-1.5,.3)(-1.5,0)(1.5,0)
 }
 \psellipse[linestyle=dotted](1,0)(.5,0.2)
 \psellipse[linestyle=dotted](-1,0)(.5,0.2)
 \end{psclip}
 }

\newcommand{\birthc}{
 \pscustom[fillstyle=gradient,
    gradbegin=white, gradend=gray,gradmidpoint=0,gradangle=110]{
        \psbezier(-.5,0)(-.5,.9)(0.5,.9)(.5,0)
        \psbezier(.5,0)(.4,-.25)(-.4,-.25)(-.5,0)
    }
 \begin{psclip}{
 \pspolygon[linestyle=none](.5,0)(.5,.3)(-.5,.3)(-.5,0)(.5,0)
 }
 \psellipse[linestyle=dotted](0,0)(.5,0.2)
 \end{psclip}
 }

\newcommand{\deathc}{
 \pscustom[fillstyle=gradient,
    gradbegin=white, gradend=gray,gradmidpoint=0,gradangle=70]{
        \psbezier(-.5,1)(-.5,.1)(0.5,.1)(.5,1)
        \psline(-.5,1)
 }
  \psellipse[fillcolor=lightgray,fillstyle=gradient,
        gradbegin=lightgray, gradend=gray,gradmidpoint=1,gradangle=110](0,1)(.5,.2)
 }

\newcommand{\zagc}{
   \pscustom[fillstyle=gradient,
    gradbegin=white, gradend=gray,gradmidpoint=0,gradangle=110]{
        \psbezier(1.5,0)(1.6,2)(-1.6,2)(-1.5,0)
        \psbezier(-1.5,0)(-1.4,-.25)(-.6,-.25)(-.5,0)
        \psbezier(-.5,0)(-.6,.8)(0.6,.8)(.5,0)
        \psbezier(.5,0)(.6,-.25)(1.4,-.25)(1.5,0)
    }
  \begin{psclip}{
 \pspolygon[linestyle=none](1.5,0)(1.5,.3)(-1.5,.3)(-1.5,0)(1.5,0)
 }
 \psellipse[linestyle=dotted](1,0)(.5,0.2)
 \psellipse[linestyle=dotted](-1,0)(.5,0.2)
 \end{psclip}
 }

\newcommand{\zigc}{
       \pscustom[fillstyle=gradient,
    gradbegin=white, gradend=gray,gradmidpoint=0,gradangle=70]{
        \psbezier(1.5,2)(1.6,0)(-1.6,0)(-1.5,2)
        \psline(-.5,2)
        \psbezier(-.5,2)(-.6,1.2)(0.6,1.2)(.5,2)
        \psline(1.5,2)
    }
 \psellipse[fillcolor=lightgray,fillstyle=gradient,
        gradbegin=lightgray, gradend=gray,gradmidpoint=1,gradangle=110](1,2)(.5,.2)
        \psellipse[fillcolor=lightgray,fillstyle=gradient,
        gradbegin=lightgray, gradend=gray,gradmidpoint=1,gradangle=110](-1,2)(.5,.2)
}

\newcommand{\identc}{
 \pscustom[fillcolor=lightgray,fillstyle=gradient,
        gradbegin=white, gradend=gray,gradmidpoint=0,gradangle=88]{
 \psline(.5,0)(.5,2.5)
 \psline(-.5,2.5)
 \psline(-.5,0)
 \psbezier(-.5,0)(-.4,-.25)(.4,-.25)(.5,0)
 }
\psellipse[fillcolor=lightgray,fillstyle=gradient,
        gradbegin=lightgray, gradend=gray,gradmidpoint=1,gradangle=110](0,2.5)(.5,.2)
 \begin{psclip}{
 \pspolygon[linestyle=none](.5,0)(.5,.3)(-.5,.3)(-.5,0)(.5,0)
 }
 \psellipse[linestyle=dotted](0,0)(.5,0.2)
 \end{psclip}
 }

  \newcommand{\medidentc}{
     \pscustom[fillcolor=lightgray,fillstyle=gradient,
        gradbegin=white, gradend=gray,gradmidpoint=0,gradangle=88]{
        \psline(-.5,2)(-.5,0)
        \psbezier(-.5,0)(-.4,-.25)(.4,-.25)(.5,0)
        \psline(.5,2)
        \psline(-.5,2)
    }
\psellipse[fillcolor=lightgray,fillstyle=gradient,
        gradbegin=lightgray, gradend=gray,gradmidpoint=1,gradangle=110](0,2)(.5,.2)
 \begin{psclip}{
 \pspolygon[linestyle=none](.5,0)(.5,.3)(-.5,.3)(-.5,0)(.5,0)
 }
 \psellipse[linestyle=dotted](0,0)(.5,0.2)
 \end{psclip}
}

\newcommand{\crossc}{
   \pscustom[fillcolor=lightgray,fillstyle=gradient,
        gradbegin=white, gradend=gray,gradmidpoint=0,gradangle=125]{
 \psline(-.5,0)(1.5,2.5)
 \psline(.5,2.5)
 \psline(-1.5,0)
 \psbezier(-1.5,0)(-1.4,-.25)(-.6,-.25)(-.5,0)
 }
 \pscustom[fillcolor=lightgray,fillstyle=gradient,
        gradbegin=white, gradend=gray,gradmidpoint=0,gradangle=125]{
 \psline(.5,0)(-1.5,2.5)
 \psline(-.5,2.5)
 \psline(1.5,0)
 \psbezier(1.5,0)(1.4,-.25)(.6,-.25)(.5,0)
 }
\psellipse[fillcolor=lightgray,fillstyle=gradient,
        gradbegin=lightgray, gradend=gray,gradmidpoint=1,gradangle=110](-1,2.5)(.5,.2)
\psellipse[fillcolor=lightgray,fillstyle=gradient,
        gradbegin=lightgray, gradend=gray,gradmidpoint=1,gradangle=70](1,2.5)(.5,.2)
 \psline[linestyle=dotted](1.5,2.5)(-.5,0)
 \psline[linestyle=dotted](.5,2.5)(-1.5,0)
 \begin{psclip}{
 \pspolygon[linestyle=none](1.5,0)(1.5,.3)(-1.5,.3)(-1.5,0)(1.5,0)
 }
 \psellipse[linestyle=dotted](1,0)(.5,0.2)
 \psellipse[linestyle=dotted](-1,0)(.5,0.2)
 \end{psclip}
 }

\newcommand{\curveleftc}{
  \pscustom[fillstyle=gradient,
    gradbegin=white, gradend=gray,gradmidpoint=0,gradangle=115]{
        \psbezier(-1.5,2.5)(-1.5,1.5)(-.4,1.3)(-.5,0)
        \psbezier(-.5,0)(-.4,-.25)(.4,-.25)(.5,0)
        \psbezier(.5,0)(.6,1.3)(-.5,1.5)(-.5,2.5)
        \psline(-1.5,2.5)
    }
    \psellipse[fillcolor=lightgray,fillstyle=gradient,
        gradbegin=lightgray, gradend=gray,gradmidpoint=1,gradangle=110](-1,2.5)(.5,.2)
 \begin{psclip}{
 \pspolygon[linestyle=none](.5,0)(.5,.3)(-.5,.3)(-.5,0)(.5,0)
 }
 \psellipse[linestyle=dotted](0,0)(.5,0.2)
 \end{psclip}
}

\newcommand{\crossmixlc}{
   \pscustom[fillcolor=lightgray,fillstyle=gradient,
        gradbegin=white, gradend=gray,gradmidpoint=0,gradangle=125]{
 \psline(-.5,0)(1.5,2.5)
 \psline(.5,2.5)
 \psline(-1.5,0)
 \psbezier(-1.5,0)(-1.4,-.25)(-.6,-.25)(-.5,0)
 }
 \pscustom[fillcolor=lightgray,fillstyle=solid]{
 \psline(.5,0)(-1.5,2.5)
 \psline(-.5,2.5)
 \psline(1.5,0)
 \psline(.5,0)
 }
\psellipse[fillcolor=lightgray,fillstyle=gradient,
        gradbegin=lightgray, gradend=gray,gradmidpoint=1,gradangle=110](1,2.5)(.5,.2)
 \begin{psclip}{
 \pspolygon[linestyle=none](1.5,0)(1.5,.3)(-1.5,.3)(-1.5,0)(1.5,0)
 }
 \psellipse[linestyle=dotted](-1,0)(.5,0.2)
 \end{psclip}
 }

\newcommand{\crossmixcl}{
   \pscustom[fillcolor=lightgray,fillstyle=solid]{
 \psline(-.5,0)(1.5,2.5)
 \psline(.5,2.5)
 \psline(-1.5,0)
 \psline(-.5,0)
 }
 \pscustom[fillcolor=lightgray,fillstyle=gradient,
        gradbegin=white, gradend=gray,gradmidpoint=0,gradangle=125]{
 \psline(.5,0)(-1.5,2.5)
 \psline(-.5,2.5)
 \psline(1.5,0)
 \psbezier(1.5,0)(1.4,-.25)(.6,-.25)(.5,0)
 }
\psellipse[fillcolor=lightgray,fillstyle=gradient,
        gradbegin=lightgray, gradend=gray,gradmidpoint=1,gradangle=70](-1,2.5)(.5,.2)
 \begin{psclip}{
 \pspolygon[linestyle=none](1.5,0)(1.5,.3)(0,.3)(0,0)(1.5,0)
 }
 \psellipse[linestyle=dotted](1,0)(.5,0.2)
 \end{psclip}
 }

\newcommand{\siidd}[1]{\xybox{%
  (-3,0)*{};
  (3,0)*{};
  (-1,0);(-1,-12) **\dir{-}; ?(.5)*\dir{>};
  (1,0);(1,-12) **\dir{-}; ?(.5)*\dir{<}
}}
\newcommand{\siidds}[1]{\xybox{%
  (-3,0)*{};
  (3,0)*{};
  (-1,0);(-1,-6) **\dir{-}; ?(.5)*\dir{>};
  (1,0);(1,-6) **\dir{-}; ?(.5)*\dir{<}
}}
\newcommand{\siidu}[1]{\xybox{%
  (-3,0)*{};
  (3,0)*{};
  (-1,0);(-1,-12) **\dir{-}; ?(.5)*\dir{<};
  (1,0);(1,-12) **\dir{-}; ?(.5)*\dir{>}
}}
\newcommand{\siidus}[1]{\xybox{%
  (-3,0)*{};
  (3,0)*{};
  (-1,0);(-1,-6) **\dir{-}; ?(.5)*\dir{<};
  (1,0);(1,-6) **\dir{-}; ?(.5)*\dir{>}
}}
\newcommand{\sidlu}[1]{\xybox{%
  (5,0);(-1,-12) **\crv{(5,-6)&(-1,-6)}; ?(.25)*\dir{<};
  (7,0);(1,-12) **\crv{(7,-6)&(1,-6)}; ?(.25)*\dir{>}
}}
\newcommand{\simu}[1]{\xybox{%
  (-7,0);(-1,-12) **\crv{(-7,-6)&(-1,-6)}; ?(.5)*\dir{>};
  (7,0);(1,-12) **\crv{(7,-6)&(1,-6)}; ?(.5)*\dir{<};
  (-5,0);(5,0) **\crv{(-5,-6)&(5,-6)}; ?(.25)*\dir{<}
}}
\newcommand{\sieta}[1]{\xybox{%
  (-3,0)*{};
  (3,0)*{};
  (-1,-12);(1,-12) **\crv{(-1,-2)&(1,-2)}; ?(.15)*\dir{<}
}}
\newcommand{\sidelta}[1]{\xybox{%
  (-7,-12);(-1,0) **\crv{(-7,-6)&(-1,-6)}; ?(.5)*\dir{<};
  (7,-12);(1,0) **\crv{(7,-6)&(1,-6)}; ?(.5)*\dir{>};
  (-5,-12);(5,-12) **\crv{(-5,-9)&(3,-9)&(3,-5)&(-3,-5)&(-3,-9)&(5,-9)}; ?(.1)*\dir{>}
}}
\newcommand{\siepsilon}[1]{\xybox{%
  (-3,0)*{};
  (3,0)*{};
  (-1,0);(1,0) **\crv{(-1,-5)&(1,-5)&(1,-10)&(-1,-10)&(-1,-5)&(1,-5)}; ?(.1)*\dir{>}
}}
\newcommand{\sicup}[1]{\xybox{%
  (-5,0);(5,0) **\crv{(-5,-9)&(5,-9)}; ?(.25)*\dir{<};
  (-7,0);(7,0) **\crv{(-7,-12)&(7,-12)}; ?(.25)*\dir{>}
}}
\newcommand{\sicap}[1]{\xybox{%
  (0,0)*{};
  (-7,-12);(7,-12) **\crv{(-7,0)&(7,0)}; ?(.25)*\dir{<};
  (-5,-12);(5,-12) **\crv{(-5,-3)&(5,-3)}; ?(.25)*\dir{>}
}}
\newcommand{\sipair}[1]{\xybox{%
  (-5,0);(5,0) **\crv{(-5,-5)&(5,-5)}; ?(.25)*\dir{<};
  (-7,0);(7,0) **\crv{(-7,-7)&(3,-7)&(3,-12)&(-3,-12)&(-3,-7)&(7,-7)}; ?(.1)*\dir{>}
}}
\newcommand{\sidual}[1]{\xybox{%
  (-7,-12);(7,-12) **\crv{(-7,0)&(7,0)}; ?(.25)*\dir{<};
  (-5,-12);(5,-12) **\crv{(-5,-9)&(3,-9)&(3,-4)&(-3,-4)&(-3,-9)&(5,-9)}; ?(.1)*\dir{>}
}}
\newcommand{\sibraiddu}[1]{\xybox{%
  (-7,0);(5,-12) **\crv{(-7,-7)&(5,-7)}; ?(.25)*\dir{>};
  (-5,0);(7,-12) **\crv{(-5,-5)&(7,-5)}; ?(.25)*\dir{<};
  (5,0);(-7,-12) **\crv{(5,-5)&(-7,-5)}; ?(.25)*\dir{<};
  (7,0);(-5,-12) **\crv{(7,-7)&(-5,-7)}; ?(.25)*\dir{>}
}}

\newcommand{\bbmu}[1]{\xybox{%
  (-3,0)*{};
  (3,0)*{};
  (0,-4)*{\bullet}="f";
  (-2,0)*{}="t1";
  (2,0)*{}="t2";
  (0,-8)*{}="b";
  "t1";"f" **\crv{(-2,-2)}; ?(.35)*\dir{>};
  "t2";"f" **\crv{(2,-2)}; ?(.35)*\dir{>};
  "f";"b" **\dir{-}; ?(.75)*\dir{>};
}}
\newcommand{\bbmediummu}[1]{\xybox{%
  (-5,0)*{};
  (5,0)*{};
  (0,-4)*{\bullet}="f";
  (-4,0)*{}="t1";
  (4,0)*{}="t2";
  (0,-8)*{}="b";
  "t1";"f" **\crv{(-4,-2)}; ?(.35)*\dir{>};
  "t2";"f" **\crv{(4,-2)}; ?(.35)*\dir{>};
  "f";"b" **\dir{-}; ?(.75)*\dir{>};
}}
\newcommand{\bbeta}[1]{\xybox{%
  (-1,0)*{};
  (1,0)*{};
  (0,-4)*{\bullet}="f";
  (0,-8)*{}="b";
  "f";"b" **\dir{-}; ?(.75)*\dir{>};
}}

\newcommand{\bbcardy}[1]{\xybox{%
  (-2,0)*{};
  (2,0)*{};
  (0,-4)*{\bullet}="g";
  (0,4)*{\bullet}="f";
  (0,8)*{}="t";
  (0,-8)*{}="b";
  "t";"f" **\dir{-}; ?(.25)*\dir{>};
  "g";"b" **\dir{-}; ?(.75)*\dir{>};
  "f";"g" **\crv{(-6,3)&(6,-3)}; ?(.4)*\dir{>};
  "f";"g" **\crv{(6,3)&(-6,-3)}; ?(.4)*\dir{>};
}}

\newcommand{\bbdelta}[1]{\xybox{%
  (-3,0)*{};
  (3,0)*{};
  (0,-4)*{\bullet}="f";
  (-2,-8)*{}="b1";
  (2,-8)*{}="b2";
  (0,-0)*{}="t";
  "f";"b1" **\crv{(-2,-6)}; ?(.75)*\dir{>};
  "f";"b2" **\crv{(2,-6)}; ?(.75)*\dir{>};
  "t";"f" **\dir{-}; ?(.35)*\dir{>};
}}
\newcommand{\bbmeddelta}[1]{\xybox{%
  (-5,0)*{};
  (5,0)*{};
  (0,-4)*{\bullet}="f";
  (-4,-8)*{}="b1";
  (4,-8)*{}="b2";
  (0,-0)*{}="t";
  "f";"b1" **\crv{(-4,-6)}; ?(.75)*\dir{>};
  "f";"b2" **\crv{(4,-6)}; ?(.75)*\dir{>};
  "t";"f" **\dir{-}; ?(.35)*\dir{>};
}}
\newcommand{\bbepsilon}[1]{\xybox{%
  (-1,0)*{};
  (1,0)*{};
  (0,-4)*{\bullet}="f";
  (0,-0)*{}="t";
  "t";"f" **\dir{-}; ?(.35)*\dir{>};
}}
\newcommand{\bbepsilonp}[1]{\xybox{%
  (-1,0)*{};
  (1,0)*{};
  (0,-4)*{\ast}="f";
  (0,-0)*{}="t";
  "t";"f" **\dir{-}; ?(.35)*\dir{>};
}}
\newcommand{\bbid}[1]{\xybox{%
  (-1,0)*{};
  (1,0)*{};
  (0,0);(0,-9) **\dir{-}; ?(.5)*\dir{>}
}}
\newcommand{\bbidu}[1]{\xybox{%
  (-1,0)*{};
  (1,0)*{};
  (0,0);(0,-9) **\dir{-}; ?(.5)*\dir{<}
}}
\newcommand{\bbdl}[1]{\xybox{%
  (2,0);(0,-8) **\crv{(2,-2)&(0,-6)}; ?(.5)*\dir{>}
}}
\newcommand{\bbdlu}[1]{\xybox{%
  (2,0);(0,-8) **\crv{(2,-2)&(0,-6)}; ?(.5)*\dir{<}
}}
\newcommand{\bbdr}[1]{\xybox{%
  (-2,0);(0,-8) **\crv{(-2,-2)&(0,-6)}; ?(.5)*\dir{>}
}}
\newcommand{\bbdru}[1]{\xybox{%
  (-2,0);(0,-8) **\crv{(-2,-2)&(0,-6)}; ?(.5)*\dir{<}
}}
\newcommand{\bbbraid}[1]{\xybox{%
  (-3,0)*{};
  (3,0)*{};
  (-2,0);(2,-8) **\crv{(-2,-2)&(2,-6)}; ?(.25)*\dir{>};
  (2,0);(-2,-8) **\crv{(2,-2)&(-2,-6)}; ?(.25)*\dir{>};
}}
\newcommand{\bbbraiddu}[1]{\xybox{%
  (-3,0)*{};
  (3,0)*{};
  (-2,0);(2,-8) **\crv{(-2,-2)&(2,-6)}; ?(.25)*\dir{>};
  (2,0);(-2,-8) **\crv{(2,-2)&(-2,-6)}; ?(.25)*\dir{<};
}}
\newcommand{\bbbraidr}[1]{\xybox{%
  (-5,0)*{};
  (5,0)*{};
  (-4,0);(4,-8) **\crv{(-4,-2)&(4,-6)}; ?(.25)*\dir{>};
  (0,0);(-4,-8) **\crv{(0,-2)&(-4,-6)}; ?(.25)*\dir{>};
  (4,0);(0,-8) **\crv{(4,-2)&(0,-6)}; ?(.25)*\dir{>};
}}
\newcommand{\bbbraidl}[1]{\xybox{%
  (-5,0)*{};
  (5,0)*{};
  (4,0);(-4,-8) **\crv{(4,-2)&(-4,-6)}; ?(.25)*\dir{>};
  (0,0);(4,-8) **\crv{(0,-2)&(4,-6)}; ?(.25)*\dir{>};
  (-4,0);(0,-8) **\crv{(-4,-2)&(0,-6)}; ?(.25)*\dir{>};
}}
\newcommand{\bbpair}[1]{\xybox{%
  (-3,0)*{};
  (3,0)*{};
  (0,-4)*{\bullet}="f";
  (-2,0)*{}="t1";
  (2,0)*{}="t2";
  "t1";"f" **\crv{(-2,-2)}; ?(.35)*\dir{>};
  "t2";"f" **\crv{(2,-2)}; ?(.35)*\dir{>};
}}
\newcommand{\bbpairp}[1]{\xybox{%
  (-3,0)*{};
  (3,0)*{};
  (0,-4)*{\ast}="f";
  (-2,0)*{}="t1";
  (2,0)*{}="t2";
  "t1";"f" **\crv{(-2,-2)}; ?(.35)*\dir{>};
  "t2";"f" **\crv{(2,-2)}; ?(.35)*\dir{>};
}}
\newcommand{\bbtriple}[1]{\xybox{%
  (-5,0)*{};
  (5,0)*{};
  (0,-4)*{\bullet}="f";
  (-4,0)*{}="t1";
  (0,0)*{}="t2";
  (4,0)*{}="t3";
  "t1";"f" **\crv{(-4,-2)}; ?(.35)*\dir{>};
  "t2";"f" **\dir{-}; ?(.35)*\dir{>};
  "t3";"f" **\crv{(4,-2)}; ?(.35)*\dir{>};
}}
\newcommand{\bbdual}[1]{\xybox{%
  (-3,0)*{};
  (3,0)*{};
  (0,-4)*{\bullet}="f";
  (-2,-8)*{}="b1";
  (2,-8)*{}="b2";
  "f";"b1" **\crv{(-2,-6)}; ?(.75)*\dir{>};
  "f";"b2" **\crv{(2,-6)}; ?(.75)*\dir{>};
}}
\newcommand{\bbdualp}[1]{\xybox{%
  (-3,0)*{};
  (3,0)*{};
  (0,-4)*{\ast}="f";
  (-2,-8)*{}="b1";
  (2,-8)*{}="b2";
  "f";"b1" **\crv{(-2,-6)}; ?(.75)*\dir{>};
  "f";"b2" **\crv{(2,-6)}; ?(.75)*\dir{>};
}}
\newcommand{\bbwidedual}[1]{\xybox{%
  (-5,0)*{};
  (5,0)*{};
  (0,-4)*{\bullet}="f";
  (-4,-8)*{}="b1";
  (4,-8)*{}="b2";
  "f";"b1" **\crv{(-4,-6)}; ?(.75)*\dir{>};
  "f";"b2" **\crv{(4,-6)}; ?(.75)*\dir{>};
}}
\newcommand{\bbhugedual}[1]{\xybox{%
  (-7,0)*{};
  (7,0)*{};
  (0,-4)*{\bullet}="f";
  (-6,-8)*{}="b1";
  (6,-8)*{}="b2";
  "f";"b1" **\crv{(-6,-6)}; ?(.75)*\dir{>};
  "f";"b2" **\crv{(6,-6)}; ?(.75)*\dir{>};
}}
\newcommand{\bbcup}[1]{\xybox{%
  (-3,0)*{};
  (3,0)*{};
  (-2,0);(2,0) **\crv{(-2,-4)&(2,-4)}; ?(.25)*\dir{<};
}}
\newcommand{\bbcap}[1]{\xybox{%
  (-3,0)*{};
  (3,0)*{};
  (-2,-8);(2,-8) **\crv{(-2,-4)&(2,-4)}; ?(.25)*\dir{<};
}}

\newcommand{\bblr}[1]{\xybox{%
  (-3,0)*{};
  (3,0)*{};
  (-2,0);(2,-8) **\crv{(-2,-2)&(2,-6)}; ?(.25)*\dir{>};
}}
\newcommand{\bbrl}[1]{\xybox{%
  (-3,0)*{};
  (3,0)*{};
  (2,0);(-2,-8) **\crv{(2,-2)&(-2,-6)}; ?(.25)*\dir{>};
}}
\newcommand{\bbdouble}[1]{\xybox{%
  (-5,0)*{};
  (5,0)*{};
  (0,-4)*{\bullet}="f";
  (-4,0)*{}="t1";
  (4,0)*{}="t3";
  "t1";"f" **\crv{(-4,-2)}; ?(.35)*\dir{>};
  "t3";"f" **\crv{(4,-2)}; ?(.35)*\dir{>};
}}

\newcommand{\bbrllong}[1]{\xybox{%
  (-5,0)*{};
  (5,0)*{};
  (-4,0);(4,-8) **\crv{(-4,-2)&(4,-6)}; ?(.25)*\dir{>};
 }}

\newcommand{\bblrlong}[1]{
  \xybox{%
  (-5,0)*{};
  (5,0)*{};
  (4,0);(-4,-8) **\crv{(4,-2)&(-4,-6)}; ?(.25)*\dir{>};
  }}
\newcommand{\bbwidemu}[1]{\xybox{%
  (-5,0)*{};
  (5,0)*{};
  (0,-4)*{\bullet}="f";
  (-6,0)*{}="t1";
  (0,-8)*{}="t2";
  (6,0)*{}="t3";
  "t1";"f" **\crv{(-6,-2)}; ?(.45)*\dir{>};
  "t3";"f" **\crv{(6,-2)}; ?(.45)*\dir{>};
  "f";"t2" **\dir{-}; ?(.65)*\dir{>};
}}
\newcommand{\bbmedmu}[1]{\xybox{%
  (-5,0)*{};
  (5,0)*{};
  (0,-4)*{\bullet}="f";
  (-4,0)*{}="t1";
  (0,-8)*{}="t2";
  (4,0)*{}="t3";
  "t1";"f" **\crv{(-4,-2)}; ?(.45)*\dir{>};
  "t3";"f" **\crv{(4,-2)}; ?(.45)*\dir{>};
  "f";"t2" **\dir{-}; ?(.65)*\dir{>};
}}

\newcommand{\bbproject}[1]{\xybox{%
  (-3,0)*{};
  (5,0)*{};
    (2,0)*\bbmu{};
    (0,8)*\bbmu{};
    (4,8)*\bbdl{};
    (0,16)*\bbbraid{};
    (6,16)*\bbid{};
    (-2,24)*\bbid{};
    (4,24)*\bbdual{};
}}

\newcommand{\bbsmallproject}[1]{\xybox{%
  (-4,0)*{};
  (4,0)*{};
  (0,-4)*{\bullet}="f";
  (-2,0)*{\bullet}="t1";
  (2,0)*{}="t2";
  (0,-8)*{}="b";
  "t1";"f" **\crv{(-2,-2)}; ?(.55)*\dir{>};
  "t2";"f" **\crv{(2,-2)}; ?(.55)*\dir{>};
  "f";"b" **\dir{-}; ?(.75)*\dir{>};
  (-2,8)*{}="t";
  (1,5)*{\bullet}="m";
  "t";"t1" **\dir{-}; ?(.35)*\dir{>};
  "t1";"m" **\crv{(-8,2)&(1,4.5)}; ?(.68)*\dir{<};
   "t2";"m" **\crv{(2,3)}; ?(.3)*\dir{<};
}}
\newcommand{\bbP}[1]{\xybox{%
   (0,0)*\xycircle(2.15,2.15){-}="m";
   (0,-6)*{}="b";
   (0,6)*{}="t";
   (0,0)*{p};
         "b";"m" **\dir{-} ?(.3)*\dir{<}  ;
         "t";"m" **\dir{-} ?(.3)*\dir{>}  ;
}}
\newcommand{\bblongid}[1]{\xybox{%
  (-1,0)*{};
  (1,0)*{};
  (0,6);(0,-6) **\dir{-}; ?(.5)*\dir{>}
}}

\newcommand{\hpeprint}[1]{%
  \href{http://arXiv.org/abs/#1}{\texttt{#1}}}%
\newcommand{\hpmathsci}[1]{%
  \href{http://www.ams.org/mathscinet-getitem?mr=#1}{\texttt{MR #1}}}%

\begin{document}
%

\preprint{AEI-2006-012}

\title{Open-closed TQFTs extend Khovanov homology \\from links to tangles}

\author{Aaron D.\ Lauda}
\email{A.Lauda@dpmms.cam.ac.uk}
\affiliation{Department of Pure Mathematics and Mathematical Statistics,\\
  University of Cambridge, Cambridge CB3 0WB, United Kingdom}

\author{Hendryk Pfeiffer}
\email{pfeiffer@aei.mpg.de}
\affiliation{Max Planck Institute for Gravitational Physics,\\
  Am M\"uhlenberg 1, 14476 Potsdam, Germany}

\date{June 14, 2006}

%
\begin{abstract}
%

We use a special kind of $2$-dimensional extended Topological Quantum
Field Theories (TQFTs), so-called open-closed TQFTs, in order to
extend Khovanov homology from links to arbitrary tangles, not
necessarily even. For every plane diagram of an oriented tangle, we
construct a chain complex whose homology is invariant under
Reidemeister moves. The terms of this chain complex are modules of a
suitable algebra $A$ such that there is one action of $A$ or
$A^\mathrm{op}$ for every boundary point of the tangle. We give
examples of such algebras $A$ for which our tangle homology theory
reduces to the link homology theories of Khovanov, Lee, and Bar-Natan
if it is evaluated for links. As a consequence of the Cardy condition,
Khovanov's graded theory can only be extended to tangles if the
underlying field has finite characteristic. In all cases in which the
algebra $A$ is strongly separable, \ie\ for Bar-Natan's theory in any
characteristic and for Lee's theory in characteristic other than $2$,
we also provide the required algebraic operation for the composition
of oriented tangles. Just as Khovanov's theory for links can be
recovered from Lee's or Bar-Natan's by a suitable spectral sequence,
we provide a spectral sequence in order to compute our tangle
extension of Khovanov's theory from that of Bar-Natan's or Lee's
theory. Thus, we provide a tangle homology theory that is locally
computable and still strong enough to recover characteristic $p$
Khovanov homology for links.
\end{abstract}

\noindent
\begin{small}
\end{small}


%
\section{Introduction}
%
\label{sec_intro}

For every plane diagram of an oriented link $L$, Khovanov's link
homology theory~\cite{Kh} yields a chain complex $\tangle{L}$ of
graded vector spaces whose graded Euler characteristic agrees with the
Jones polynomial of the link. This construction can be seen as a
categorification of the unnormalized Jones polynomial, replacing a
polynomial in one indeterminate $q$ by a chain complex of graded
vector spaces. The coefficients of the polynomial arise as the
dimensions of the homogeneous components of the graded homology groups
of the chain complex in such a way that the degree corresponds to the
power of $q$.

If two link diagrams are related by a Reidemeister move, the
corresponding chain complexes of graded vector spaces are homotopy
equivalent, and so their homology groups are isomorphic as graded
vector spaces. This of course implies that their graded Euler
characteristics and thereby their unnormalized Jones polynomials
agree, but in general the homology groups contain more information
about the link than does the Jones polynomial. Indeed, Bar-Natan and
Wehrli~\cite{BN,BN1,Wehrli} have shown that there are knots and
links that have the same Jones polynomial, but which can be
distinguished by their Khovanov homology.

The construction of Khovanov's chain complex heavily relies on a
$2$-dimensional Topological Quantum Field Theory (TQFT). The fact
that the boundary of a cobordism is a closed manifold, naturally
restricts the construction to links rather than arbitrary tangles.
The purpose of the present article is to overcome this limitation.
In the remainder of the introduction, we give a brief overview of
the present article, addressed to the experts. The reader who is
familiar with the details of Khovanov's definition~\cite{Kh} and
with Bar-Natan's work on the extension from links to
tangles~\cite{BN2}, should be able to understand the key ideas of
our work just from reading the introduction. The reader who is not
familiar with all the details, is invited to have a brief look at
the introduction in order to get a first impression, and can then
work though the material step by step starting from
Section~\ref{sect_prelim}.

\subsection{Overview of Khovanov's link homology theory}

\parpic[r]{$\hopfresolutions$} Let $L$ be a plane diagram of an
oriented link with $n_+$ positive crossings $(\overcrossing)$ and
$n_-$ negative ones $(\undercrossing)$, $n:=n_++n_-$. The Kauffman
bracket $\la L \ra$ from which one can compute the unnormalized Jones
polynomial $\hat J(L):={(-1)}^{n_-}q^{n_+-2n_-}\la
L\ra\in\Z[q,q^{-1}]$, can be recursively defined\footnote{We have
adopted the conventions of Khovanov~\cite{Kh}. The usual definition of
the Jones polynomial is obtained by substituting $-\sqrt{t}$ for $q$
in the normalized polynomial $J(L):=\hat J(L)/(q+q^{-1})$.} as
follows:
\begin{equation}
\label{eq_kauffman}
\la\emptyset\ra :=  1; \qquad
 \la \bigcirc L \ra := (q+q^{-1}) \la L \ra ;\qquad
\la\backoverslash\ra  := \la\hsmoothing\ra-q\la\smoothing\ra.
\end{equation}
For every crossing $(\backoverslash)$, each of the two smoothings,
the $0$-smoothing $(\hsmoothing)$ and the $1$-smoothing
$(\smoothing)$, give a contribution to the Kauffman bracket, with a
different sign and a different power of $q$. For the link diagram
$L$ with $n$ crossings, there are $2^n$ smoothings, labeled by
subsets $\alpha\subseteq\cat{n}:=\{1,\ldots,n\}$ where $j\in\alpha$
if the $j$-th crossing is resolved by the $1$-smoothing and
$j\notin\alpha$ if it is resolved by the $0$-smoothing. We denote by
$D_\alpha\subseteq\R^2$, $\alpha\subseteq\cat{n}$, the diagrams
associated with the $2^n$ smoothings of $L$. Each $D_\alpha$ is free
of crossings and therefore forms a disjoint union of a finite number
of circles in $\R^2$.

A crucial ingredient for the construction of Khovanov's chain complex
is the following collection of cobordisms between these diagrams. For
every crossing $j\in\cat{n}$ and for all possible smoothings of the
other crossings, represented by subsets
$\alpha\subseteq\cat{n}\backslash\{j\}$, one specifies a
$2$-dimensional cobordism
\begin{equation}
\label{eq_saddles}
  S_{(\alpha,j)}\colon D_\alpha\to D_{\alpha\sqcup\{j\}}
\end{equation}
which relates the two possible smoothings of the $j$-th crossing.
This cobordism $S_{(\alpha,j)}$ is a cylinder over $D_\alpha$
everywhere except for a neighbourhood of the $j$-th
crossing\footnote{Notice that we have rotated this picture by 90
degrees.} $(\slashoverback)$ where it is a saddle
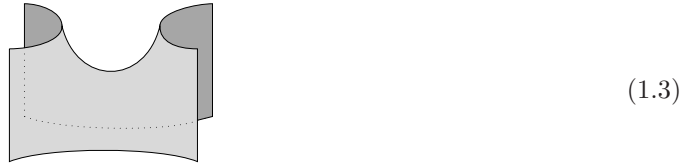
\begin{equation}
\label{eq_saddlepic}
\begin{aligned}
\psset{xunit=1cm,yunit=1cm}
\begin{pspicture}(3.0,2.0)
  \rput(0,0){\bigsaddle}
\end{pspicture}
\end{aligned}
\end{equation}
with the $0$-smoothing $(\smoothing)$ as the source (at the top of
the diagram) and the $1$-smoothing $(\hsmoothing)$ as the target (at
the bottom).

Then one uses a $2$-dimensional TQFT, \ie\ a symmetric monoidal
functor $Z\colon\cat{2Cob}\to\cat{Vect}_k$ from the category
$\cat{2Cob}$ of $2$-dimensional cobordisms to the category
$\cat{Vect}_k$ of vector spaces over some fixed field $k$, in order to
turn the cobordism into a linear map between vector spaces. The chain
groups of Khovanov's chain complex are constructed from the vector
spaces $Z(D_\alpha)$, $\alpha\subseteq\cat{n}$, and the differentials
from the linear maps $Z(S_{(\alpha,j)})$, $j\in\cat{n}$,
$\alpha\subseteq\cat{n}\backslash\{j\}$. For more details, see
Section~\ref{sec_Tangle} below, and for the full explanation of how
Khovanov's theory is related to the Jones polynomial, we refer
to~\cite{Kh,BN}.

The category of $2$-dimensional TQFTs, \ie\ the category
$\cat{Symm-Mon}(\cat{2Cob},\cat{Vect}_k)$ of symmetric monoidal
functors $\cat{2Cob}\to\cat{Vect}_k$ and monoidal natural
transformations, is equivalent as a symmetric monoidal category to the
category of commutative Frobenius algebras over
$k$~\cite{Abrams,Sawin95,Kock}.

The objects of $\cat{2Cob}$ are numbers $k\in\N_0$, representing the
diffeomorphism types of $1$-manifolds that are the disjoint union of
$k$ circles. For the morphisms of $\cat{2Cob}$, one has a description
in terms of generators and relations. The generators are these
cobordisms:
\begin{equation}
\label{eq_closedgen}
\begin{aligned}
\psset{xunit=.4cm,yunit=.4cm}
\begin{pspicture}[.2](12,3.5)
  \rput(0,1){\multc}
  \rput(4,1){\comultc}
  \rput(8,2){\birthc}
  \rput(11,1.6){\deathc}
  \rput(0,0){$\mu$}
  \rput(4,0){$\Delta$}
  \rput(8,0){$\eta$}
  \rput(11,0){$\epsilon$}
\end{pspicture}
\end{aligned}
\end{equation}
We have drawn the cobordisms with their source at the top and their
target at the bottom of the diagram. Composition of morphisms
corresponds to putting the building blocks of~\eqref{eq_closedgen}
on top of each other. The tensor product of morphisms is the
disjoint union of cobordisms. If we denote by $C:=Z(1)$ the vector
space associated with a single circle, the functor $Z$ yields linear
maps $\mu\colon C\otimes C\to C$, $\Delta\colon C\to C\otimes C$,
$\eta\colon k\to C$, and $\epsilon\colon C\to k$ for these
cobordisms. The relations among the morphisms of $\cat{2Cob}$ are
such that $(C,\mu,\eta,\Delta,\epsilon)$ forms a commutative
Frobenius algebra.

Khovanov's original choice of Frobenius algebra, called $A[c]$
in~\cite{Kh}, over the commutative ring $R=\Z[c]$, is such that one
gets a chain complex of graded modules and a categorification of the
Jones polynomial. His TQFT is actually a functor
$\cat{2Cob}\to\cat{Mod}_R$. In this article, we restrict ourselves
to the case $c\equiv 0$ (\cf~\cite{Jac1,Jac2}) and to algebras over
a field $k$. With this restriction, Khovanov's commutative Frobenius
algebra forms the special case $h=0$, $t=0$ of the following
definition:

\begin{defn}[see \cite{Kh2}]
\label{def_khovanov} Let $k$ be a field and $h,t\in k$. $C_{h,t}$
denotes the algebra $C_{h,t}=k[x]/(x^2-hx-t)$ equipped with the
structure of a commutative Frobenius algebra
$(C_{h,t},\mu,\eta,\Delta,\epsilon)$ which is given in the basis
$\{1,x\}$ by $\mu(1\otimes 1)=1$, $\mu(1\otimes x)=x$, $\mu(x\otimes
1)=x$, $\mu(x\otimes x)=hx+t$, $\eta(1)=1$, $\Delta(1)=1\otimes
x+x\otimes 1-h\cdot 1\otimes 1$, $\Delta(x)=x\otimes x + t\cdot
1\otimes 1$, $\epsilon(1)=0$, and $\epsilon(x)=1$.
\end{defn}

Several results below depend on the characteristic of $k$. While
Khovanov originally studied the case $h=0$, $t=0$, Lee~\cite{Lee}
considered $h=0$, $t=1$, and Bar-Natan~\cite{BN2} studied $h=1$,
$t=0$. Below, we refer to $C_{0,0}$ as \emph{Khovanov's}, to $C_{0,1}$
as \emph{Lee's}, and to $C_{1,0}$ as \emph{Bar-Natan's} Frobenius
algebra.

The link homology theory associated with Khovanov's Frobenius algebra
is known to categorify the Jones polynomial, a \emph{quantum
invariant} of links. In some cases, the other two link homology
theories are related to \emph{classical invariants} of links: Both
Lee's theory in characteristic $0$~\cite{Lee} and Bar-Natan's theory
in characteristic $2$~\cite{Turner} categorify combinatorial
expressions involving linking numbers.

Bar-Natan~\cite{BN2} has presented sufficient (but not necessary)
conditions for TQFTs to yield link homology theories. These can be
stated in topological and in algebraic terms as follows:

\begin{defn}\hfill
\label{def_barnatancond}
\begin{enumerate}
\item
  A $2$-dimensional TQFT $Z\colon\cat{2Cob}\to\cat{Vect}_k$ is said to
  \emph{satisfy Bar-Natan's conditions} if the following three
  conditions hold:
\begin{eqnarray}
\label{eq_barnatan1}
Z(
\begin{aligned}
\psset{xunit=.2cm,yunit=.2cm}
\begin{pspicture}[.2](0.9,5)
  \rput(0.3,1.5){\deathc}
  \rput(0.6,2.5){\birthc}
\end{pspicture}
\end{aligned}) &=& 0,\\
\label{eq_barnatan2}
Z\bigl(
\begin{aligned}
\psset{xunit=.2cm,yunit=.2cm}
\begin{pspicture}[.2](3,10)
  \rput(1.3,1.5){\deathc}
  \rput(1.6,2.5){\multc}
  \rput(1.6,5){\comultc}
  \rput(1.6,7.5){\birthc}
\end{pspicture}
\end{aligned}\bigr) &=& 2,\\
\label{eq_barnatan3}
Z\biggl(
\begin{aligned}
\psset{xunit=.2cm,yunit=.2cm}
\begin{pspicture}[.2](3.5,5)
  \rput(0.8,0){\birthc}
  \rput(2.8,0){\birthc}
  \rput(1.5,1.5){\deathc}
  \rput(1.8,2.5){\multc}
\end{pspicture}
\end{aligned}
\biggr)\;+\;Z\biggl(
\begin{aligned}
\psset{xunit=.2cm,yunit=.2cm}
\begin{pspicture}[.2](4,5)
  \rput(2.3,0){\comultc}
  \rput(2.3,2.5){\birthc}
  \rput(1,4){\deathc}
  \rput(3,4){\deathc}
\end{pspicture}
\end{aligned}
\biggr)\;-\;Z\biggl(
\begin{aligned}
\psset{xunit=.2cm,yunit=.2cm}
\begin{pspicture}[.2](4,5)
  \rput(1,0){\identc}
  \rput(1,2.5){\identc}
  \rput(3.3,0){\birthc}
  \rput(3,4){\deathc}
\end{pspicture}
\end{aligned}
\biggr)\;-\;Z\biggl(
\begin{aligned}
\psset{xunit=.2cm,yunit=.2cm}
\begin{pspicture}[.2](4,5)
  \rput(3.3,0){\identc}
  \rput(3.3,2.5){\identc}
  \rput(1,0){\birthc}
  \rput(0.7,4){\deathc}
\end{pspicture}
\end{aligned}
\biggr) &=& 0.
\end{eqnarray}
\item
  A commutative Frobenius algebra $(C,\mu,\eta,\Delta,\epsilon)$ is
  said to \emph{satisfy Bar-Natan's conditions} if
\begin{alignat}{2}
\label{eq_sphere}
  (\epsilon_C\circ\eta_C)(1) &= 0&\qquad&\mbox{(S=`sphere'),}\\
\label{eq_torus}
  (\epsilon_C\circ\mu_C\circ\Delta_C\circ\eta_C)(1) &= 2 &\qquad&\mbox{(T=`torus'),}\\
\label{eq_tubes}
  \Delta_C\circ\eta_C\circ(\epsilon_C\otimes\epsilon_C)
  +(\eta_C\otimes\eta_C)\circ\epsilon_C\circ\mu_C &&&\nn\\
  -(\eta_C\circ\epsilon_C)\otimes\id_C
  -\id_C\otimes(\eta_C\circ\epsilon_C) &=0 &\qquad&\mbox{(4Tu=`four tubes').}
\end{alignat}
\end{enumerate}
\end{defn}

The Frobenius algebra $C_{h,t}$ of Definition~\ref{def_khovanov}
satisfies Bar-Natan's conditions for all $h,t\in k$.
Khovanov~\cite{Kh2} has classified the Frobenius algebras that give
rise to link homology theories. The algebra $C_{h,t}$ of
Definition~\ref{def_khovanov} is a specialization of $A_5$
of~\cite{Kh2} to a field $k$ and to $h,t\in k$. The classification
includes examples, for instance Khovanov's $A[c]$ without evaluation
at $c=0$, called $A_2$ in~\cite{Kh2}, which do not satisfy all of
Bar-Natan's conditions~\eqref{eq_barnatan1} to~\eqref{eq_barnatan3}.
In the present article, we restrict ourselves to Frobenius algebras
that do satisfy these conditions.

\parpic[r]{$\tangleresolutions$} We wish to extend these link homology
theories from links to tangles. For tangles, the diagrams $D_{\alpha}$
corresponding to the smoothings would consist not only of circles, but
rather of circles and arcs. Khovanov has already remarked
in~\cite{Kh3} that one would need an extended $2$-dimensional TQFT in
which the cobordisms are generalized to suitable manifolds with
corners.

Even without an extended TQFT, there have been two workarounds.
Khovanov~\cite{Kh3} considers tangles with an even number of points
both for the source and the target of the tangle, so-called
\emph{even tangles}, and presents a definition in which the tangle
homology is reduced to his link homology by closing the open ends of
the tangles in all possible ways and taking the sum over the
resulting expressions for the links. He thus obtains a tangle
homology theory only for even tangles, but it is still possible to
tell which expression this construction categorifies.

Bar-Natan~\cite{BN2} works with formal linear combinations of
manifolds with corners. As long as one stays in this topological
setting, called the `picture world', one has good composition laws
for tangles, but Bar-Natan still relies on a TQFT and can therefore
translate from topology to algebra only for links, \ie\ after all
open ends have been closed.

\noindent
The following two questions remain to be answered and form the subject
of the present article:
\begin{enumerate}
\item
  How does one define a homology theory for arbitrary oriented
  tangles, using only algebraic as opposed to topological
 `picture world' constructions?
\item
  Which algebraic operation corresponds to the composition of tangles?
\end{enumerate}

In particular, we wish to understand which algebraic data to assign to
a single arc ($|$) so that one can take a disjoint union of two such
arcs ($|\,|$), then glue on a cap ($\cap$) and a cup ($\cup$) in order
to obtain a circle in such a way that both operations, disjoint union
(tensor product) and gluing (composition) have a correspondence on the
algebraic side. If we write $\tangle{T}$ for the algebraic data
associated with a tangle diagram $T$, this means that:
\begin{equation}
\label{eq_compositionrule}
  \tangle{\circ} = \tangle{\cup\circ(|\,\otimes\,|)\circ\cap}
  = \tangle{\cup}\circ(\tangle{|}\otimes\tangle{|})\circ\tangle{\cap}.
\end{equation}
This cannot be done by Khovanov's method~\cite{Kh3} since a single
arc is not an even tangle, and it cannot be done by Bar-Natan's
method~\cite{BN2} since there one can translate from topology to
algebra only for links, but not for a single arc which forms a
proper tangle.

A further limitation of the existing link homology theories is the
following. Recall that the cobordisms $S_{(\alpha,j)}$
of~\eqref{eq_saddles} are the disjoint union of either
$\psset{xunit=.3em,yunit=.3em}\begin{pspicture}(3,2.5)\rput(2,0){\multc}\end{pspicture}$
or
$\psset{xunit=.3em,yunit=.3em}\begin{pspicture}(3,2.5)\rput(2,0){\comultc}\end{pspicture}$
and a number of cylinders over circles. In order to decide which one
of the two cobordisms with critical point of index $1$ is needed, one
has to know the diagrams $D_\alpha$ and $D_{\alpha\sqcup\{j\}}$
\emph{globally}. An extension of Khovanov homology from links to
tangles in contrast would provide us with a \emph{local} description
of such a homology theory. This would not only be conceptually more
transparent, but also form a valuable property in computations.

\subsection{Open-closed Topological Quantum Field Theories}

In a previous article~\cite{LP1}, we have studied a suitable class of
manifolds with corners, called \emph{open-closed cobordisms}, and
classified the corresponding extended TQFTs which we call
\emph{open-closed TQFTs}. They are classified in terms of a
custom-made algebraic structure, so-called \emph{knowledgeable
Frobenius algebras}. The category $\cat{2Cob}^{\mathrm{ext}}$ of
open-closed cobordisms has a description in terms of generators and
relations, too. The objects are sequences
$\underline{m}=(m_1,\ldots,m_k)\in{\{0,1\}}^k$, $k\in\N_0$, which
specify a particular sequence of circles ($m_j=0$) and arcs
($m_j=1$). The generators for the morphisms of
$\cat{2Cob}^{\mathrm{ext}}$ are these open-closed cobordisms:
\begin{equation}
\label{eq_generators}
\begin{aligned}
\psset{xunit=.4cm,yunit=.4cm}
\begin{pspicture}[.2](27,3.5)
  \rput(0,1){\multl}
  \rput(0,0){$\mu_A$}
  \rput(4,1){\comultl}
  \rput(4,0){$\Delta_A$}
  \rput(7,2){\birthl}
  \rput(7,0){$\eta_A$}
  \rput(9,2.4){\deathl}
  \rput(9,0){$\epsilon_A$}
  \rput(12,1){\multc}
  \rput(12,0){$\mu_C$}
  \rput(16,1){\comultc}
  \rput(16,0){$\Delta_C$}
  \rput(19,2){\birthc}
  \rput(19,0){$\eta_C$}
  \rput(21,1.6){\deathc}
  \rput(21,0){$\epsilon_C$}
  \rput(23,1){\ctl}
  \rput(23,0){$\imath$}
  \rput(26,1){\ltc}
  \rput(26,0){$\imath^\ast$}
\end{pspicture}
\end{aligned}
.
\end{equation}
The global structure of such open-closed cobordisms $M$ is described
by specifying two faces $\del_0M$ and $\del_1M$ of $M$ as follows,
called the \emph{black} and the \emph{coloured} boundary,
respectively:
\begin{equation}
\label{eq_faces}
\begin{aligned}
\psset{xunit=.4cm,yunit=.4cm}
\begin{pspicture}[0.5](4,5.4)
  \rput(2, 0){\multl}
  \rput(1, 2.5){\medidentl}
  \rput(2.92, 2.5){\ctl}
  \rput(2, 5.25){$M$}
\end{pspicture}
\qquad\quad
\begin{pspicture}[0.5](4,5.4)
  \rput(2, 5.255){$\partial_0 M$}
  \psline[linewidth=0.5pt](1.5,4.5)(.5,4.5)
  \psline[linewidth=0.5pt](2.5,0)(1.5,0)
  \psellipse[linewidth=0.5pt](3,4.5)(.5,0.2)
\end{pspicture}
\qquad\quad
\begin{pspicture}[0.5](4,5.4)
  \rput(2, 5.25){$\partial_1 M$}
  \rput(2,0){
    \psline[linewidth=0.5pt](-1.5,4.5)(-1.5,2.5)
    \psbezier[linewidth=0.5pt](-1.5,2.5)(-1.5,1.5)(-.4,1.3)(-.5,0)
    \psline[linewidth=0.5pt](-.5,4.5)(-.5,2.5)
  }
  \rput(3,0){
    \psbezier[linewidth=0.5pt](-.5,0)(-.6,1.3)(.5,1.5)(.5,2.5)
  }
  \rput(3.08,2.5){
    \psbezier[linewidth=0.5pt](-.58,0)(-.48,.5)(-.48,.7)(-.08,1)
    \psbezier[linewidth=0.5pt](-.08,1)(.32,.7)(.32,.5)(.42,0)
  }
  \rput(2,.5){
    \psbezier[linewidth=0.5pt](-.5,2)(-.6,1)(0.6,1)(.5,2)
  }
\end{pspicture}
\qquad\quad
\begin{pspicture}[0.5](4,5.4)
  \rput(2, 5.25){$\partial_0 M \cup\partial_1 M$}
  \rput(2,0){
     \psline[linewidth=0.5pt](-1.5,4.5)(-1.5,2.5)
     \psbezier[linewidth=0.5pt](-1.5,2.5)(-1.5,1.5)(-.4,1.3)(-.5,0)
     \psline[linewidth=0.5pt](-.5,4.5)(-.5,2.5)
  }
  \rput(3,0){
    \psbezier[linewidth=0.5pt](-.5,0)(-.6,1.3)(.5,1.5)(.5,2.5)
  }
  \rput(3.08,2.5){
    \psbezier[linewidth=0.5pt](-.58,0)(-.48,.5)(-.48,.7)(-.08,1)
    \psbezier[linewidth=0.5pt](-.08,1)(.32,.7)(.32,.5)(.42,0)
  }
  \rput(2,.5){
    \psbezier[linewidth=0.5pt](-.5,2)(-.6,1)(0.6,1)(.5,2)
  }
  \psline[linewidth=0.5pt](1.5,4.5)(.5,4.5)
  \psline[linewidth=0.5pt](2.5,0)(1.5,0)
  \psellipse[linewidth=0.5pt](3,4.5)(.5,0.2)
\end{pspicture}
\qquad\quad
\begin{pspicture}[0.5](4,5.4)
  \rput(2, 5.25){$\partial_0 M \cap\partial_1 M$}
  \psdots(.5,4.5)(1.5,4.5)(1.5,0)(2.5,0)
\end{pspicture}
\end{aligned}
\end{equation}
The union $\del_0M\cup\del_1M=\del M$ is the boundary of $M$ as a
topological manifold, and their intersection $\del_0M\cap\del_1M$ is
the set of corners of $M$. Open-closed cobordisms are composed by
gluing them along their \emph{black} boundaries.

An open-closed TQFT is a symmetric monoidal functor
$Z\colon\cat{2Cob}^{\mathrm{ext}}\to\cat{Vect}_k$. If we denote by
$A:=Z((1))$ the vector space associated with a single arc and by
$C:=Z((0))$ the vector space for a single circle, then the relations
of $\cat{2Cob}^{\mathrm{ext}}$ are precisely the Moore--Segal
relations~\cite{MS}. They imply that $(A,C,\imath,\imath^\ast)$ forms
what we call a knowledgeable Frobenius algebra, \ie\ that
$(A,\mu_A,\eta_A,\Delta_A,\epsilon_A)$ is a symmetric Frobenius
algebra, that $(C,\mu_C,\Delta_C,\eta_C,\epsilon_C)$ is a commutative
Frobenius algebra, and that some additional conditions on $\imath$ and
$\imath^\ast$ hold. For the details, we refer to
Section~\ref{sect_kfrob} below.

From Bar-Natan's work~\cite{BN2} and from his purely topological
proof of the invariance of his tangle homology under Reidemeister
moves, it is immediately obvious that an open-closed TQFT is just
the appropriate tool for turning Bar-Natan's topological setup into
an algebraic tangle homology theory. What remains to be done is to
find knowledgeable Frobenius algebras $(A,C,\imath,\imath^\ast)$ and
thereby open-closed TQFTs that satisfy the
conditions~\eqref{eq_barnatan1} to~\eqref{eq_barnatan3}. We present
examples of such knowledgeable Frobenius algebras
$(A,C,\imath,\imath^\ast)$ for which $C=C_{h,t}$
(Definition~\ref{def_khovanov}) and for which the resulting tangle
homology theories reduce to the link homology theories of Khovanov,
Lee, and Bar-Natan when they are evaluated for links. This resolves
the first issue, namely to find homology theories for arbitrary
oriented tangles (see Theorem~\ref{thm_reidemeister} below).

If we wish to extend Khovanov's Frobenius algebra $C_{0,0}$
(Definition~\ref{def_khovanov}) with its grading to a knowledgeable
Frobenius algebra that is graded, too, then the Cardy condition
$\mu_A\circ\tau_{A,A}\circ\Delta_A=\imath\circ\imath^\ast$ which is
a consequence of the following diffeomorphism,
\begin{equation}
\label{cardy}
\psset{xunit=.3cm,yunit=.3cm}
\begin{aligned}
\begin{pspicture}(4,8)
  \rput(2,0){\multl}
  \rput(2,2.5){\crossl}
  \rput(2,5){\comultl}
\end{pspicture}
\end{aligned}
\quad\cong\quad
\begin{aligned}
\begin{pspicture}(3,5.5)
  \rput(1.8,1){\ctl}
  \rput(2,3){\ltc}
\end{pspicture}
\end{aligned}
,
\end{equation}
implies that the underlying field has finite characteristic (see
Proposition~\ref{prop_gradedlimitation} below).

\subsection{Composition of tangles}

The answer to the second question, namely to find an algebraic
operation that represents the composition of tangles, is inspired by
our recent work~\cite{LP2} on the state sum construction of
open-closed TQFTs, \ie\ the construction of open-closed TQFTs from
triangulations of the open-closed cobordisms. The mechanism by which
one can subdivide edges in the black boundary of an open-closed
cobordism, provides us with the blue print for an algebraic operation
to compose tangles. For example, triangulate the cobordism $\mu_C$
of~\eqref{eq_generators}. Then the three circles in its black boundary
are composed from several edges of the triangulation. This composition
of the boundary edges to the circle is our guiding example for
defining the composition of arcs.

To make this mechanism work, we impose the additional conditions (see
Theorem~2.22 of~\cite{LP2} or Proposition~\ref{prop_statesum} below)
that hold for those open-closed TQFTs that we can construct from a
state sum. In particular, these \emph{state sum open-closed TQFTs}
$Z\colon\cat{2Cob}^\mathrm{ext}\to\cat{Vect}_k$ have an associated
knowledgeable Frobenius algebra $(A,C,\imath,\imath^\ast)$, \ie\
$A=Z((1))$, $C=Z((0))$, \etc, for which $A$ is strongly
separable\footnote{This implies that $A$ is finite-dimensional over
$k$ and semisimple.} and for which $C$ is the centre of $A$.

If the algebra $A$ is strongly separable, then the element
$a=\mu_A\circ\Delta_A\circ\eta_A\colon k\to A$ associated with the
surface
\begin{equation}
a = Z(
\begin{aligned}
\psset{xunit=.3cm,yunit=.3cm}
\begin{pspicture}[.2](3.5,7.5)
  \rput(1.6,1){\multl}
  \rput(1.6,3.5){\comultl}
  \rput(1.6,6){\birthl}
\end{pspicture}
\end{aligned}
)
\end{equation}
has a convolution inverse $a^{-1}\colon k\to A$, \ie\
$\mu_A\circ(a\otimes a^{-1})=\eta_A=\mu_A\circ(a^{-1}\otimes a)$. We
call $a$ the \emph{window element}. Both $a$ and $a^{-1}$ are
central, \ie\ for any linear map $\phi\colon A\to A$, we have
$\mu_A\circ(a\otimes\phi)=\mu_A\circ(\phi\otimes a)$. In a strongly
separable algebra, multiplication with the inverse of the window
element can be used to remove holes (technical term:
\emph{windows}). By this, we mean that
\begin{equation}
\label{eq_windowkiller}
  \mu_A\circ\biggl(a^{-1}\otimes Z\bigl(
\begin{aligned}
\psset{xunit=.3cm,yunit=.3cm}
\begin{pspicture}[.2](3.5,5)
  \rput(1.6,0){\multl}
  \rput(1.6,2.5){\comultl}
\end{pspicture}
\end{aligned}
\bigr)\biggr) = \id_A = Z(
\begin{aligned}
\psset{xunit=.3cm,yunit=.3cm}
\begin{pspicture}[.2](1.8,2.5)
  \rput(0.8,0){\identl}
\end{pspicture}
\end{aligned}
),
\end{equation}
\ie\ multiplication by $a^{-1}$ on the algebraic side has the same
effect as removing a window on the topological side.

If one wishes to compose tangles, say a crossing $(\slashoverback)$
with an arc $(\rightarc)$ in order to get the tangle $(\reideone)$,
one needs to compose the extended cobordisms that generalize the
$S_{(\alpha,j)}$ of~\eqref{eq_saddles} in a suitable fashion. Let
$Z\colon\cat{2Cob}^\mathrm{ext}\to\cat{Vect}_k$ be an open-closed
TQFT that satisfies Bar-Natan's conditions, and let
$(A,C,\imath,\imath^\ast)$ be its associated knowledgeable Frobenius
algebra.

For a single crossing $(\slashoverback)$, we have the two smoothings
$D_\emptyset=\smoothing$ and $D_{\{1\}}=\hsmoothing$. With such a
crossing $(\slashoverback)$, we associate the $2$-term chain complex
of vector spaces
\begin{equation}
\xymatrix{
  0\ar[rr]&&
  Z(D_\emptyset)\ar[rr]^{Z(S_{(\emptyset,1)})}&&
  Z(D_{\{1\}})\ar[rr]&&
  0
}
\end{equation}
whose differential is obtained from the saddle $S_{(\emptyset,1)}$
depicted in~\eqref{eq_saddlepic}. This saddle is now viewed as an
open-closed cobordism, and so it gives a map
$Z(S_{(\emptyset,1)})\colon A\otimes A\to A\otimes A$. To the arc
($\rightarc$), we associate the $1$-term chain complex $0\rightarrow
A\rightarrow 0$.

The tensor product of the two chain complexes is the following
$2$-term chain complex
\begin{equation}
\label{eq_juxtaposition}
\xymatrix{
  0\ar[rr]&&
  Z(D_\emptyset)\otimes A\ar[rrr]^{Z(S_{(\emptyset,1)})\otimes\id_A}&&&
  Z(D_{\{1\}})\otimes A\ar[rr]&&
  0
}
\end{equation}
which can be associated with the following disjoint union of
open-closed cobordisms:
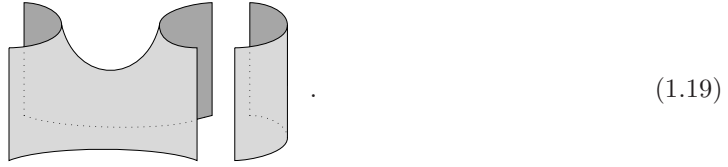
\begin{equation}
\label{eq_constituents}
\begin{aligned}
\psset{xunit=1cm,yunit=1cm}
\begin{pspicture}(4.0,2.0)
  \rput(0,0){\bigsaddle}
  \rput(3,0){\bigarc}
\end{pspicture}
\end{aligned}
.
\end{equation}

We would like to glue the two open-cobordisms along their coloured
boundaries like this,
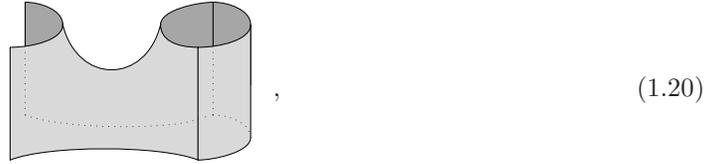
\begin{equation}
\label{eq_composite}
\begin{aligned}
\psset{xunit=1cm,yunit=1cm}
\begin{pspicture}(3.5,2.0)
  \rput(0,0){\bigsaddle}
  \rput(2.5,0){\bigarcc}
\end{pspicture}
\end{aligned}
,
\end{equation}
in order to obtain a cobordism from the $0$-smoothing
$(\reideonezero)$ of $(\reideone)$ to the $1$-smoothing
$(\reideoneone)$ of $(\reideone)$. In general, however, an open-closed
TQFT does not have any operation for such a gluing along coloured
boundaries.

Our work~\cite{LP2} on the state sum construction of open-closed TQFTs
suggests the following solution to this problem. We assume that the
open-closed TQFT is a state sum open-closed TQFT
(Proposition~\ref{prop_statesum} below). We imagine that the composite
surface~\eqref{eq_composite} is triangulated in such a way that the
components of the coloured boundary in~\eqref{eq_composite} along
which we want to glue, coincide with edges of the triangulation. This
implies that the circle in the source of the open-closed
cobordism~\eqref{eq_composite} is triangulated with two edges and the
arc in its target with three edges. Such a triangulation is displayed
here:
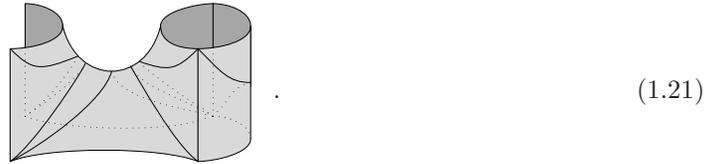
\begin{equation}
\label{eq_triangulated}
\begin{aligned}
\psset{xunit=1cm,yunit=1cm}
\begin{pspicture}(3.5,2.0)
  \rput(0,0){\bigsaddlet}
  \rput(2.5,0){\bigarct}
\end{pspicture}
\end{aligned}
.
\end{equation}
We can now employ the state sum construction presented in~\cite{LP2}
in order to compute the linear map $A\otimes C\to A$ that the
open-closed TQFT associates with the open-closed
cobordism
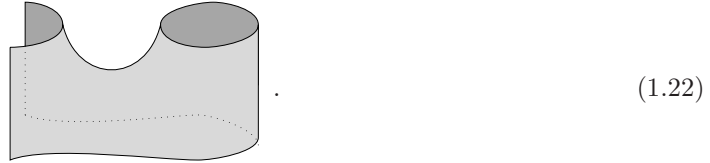
\begin{equation}
\label{eq_glued}
\begin{aligned}
\psset{xunit=1cm,yunit=1cm}
\begin{pspicture}(3.5,2.0)
  \rput(0,0){\bigcomposed}
\end{pspicture}
\end{aligned}
.
\end{equation}
This gives the chain complex $0\rightarrow A\otimes C\rightarrow
A\rightarrow 0$ which relates the two smoothings $(\reideonezero)$
and $(\reideoneone)$ of the composite tangle diagram $(\reideone)$.

We now present an equivalent way of computing the linear map $A\otimes
C\to A$ which makes transparent how this linear map can be computed
from the two constituents of~\eqref{eq_constituents}. This procedure
does not refer to triangulations or to state sums, and so it can be
understood without being familiar with~\cite{LP2}. We nevertheless
encourage the reader to study~\cite{LP2} in greater detail and to
verify that both procedures indeed result in the same linear map.

In an open-closed TQFT, we cannot compose the two constituents
of~\eqref{eq_constituents} by gluing them along their coloured
boundary. We rather pre- and postcompose~\eqref{eq_constituents} with
suitable open-closed cobordisms by gluing along their black boundary
as follows:
\begin{equation}
\label{eq_bigcomposition}
a^{-2}\cdot Z\biggl(
\begin{aligned}
\psset{xunit=6mm,yunit=6mm}
\begin{pspicture}(3.85,7.5)
  \rput(0,2.5){\bigsaddle}
  \rput(3,2.5){\bigarc}
  \rput(2,5){\bigtopright}
  \rput(0,5){\bigtopleft}
  \rput(0,0){\bigbottom}
  \rput(1.92,2.4){$\circ$}
  \rput(1.92,4.8){$\circ$}
\end{pspicture}
\end{aligned}
\biggr) = a^{-2}\cdot Z\biggl(
\begin{aligned}
\begin{pspicture}(3.85,4.5)
  \rput(0,0){\bigbottom}
  \rput(0,1.5){\bigsaddleb}
  \rput(3,1.5){\bigarcb}
  \rput(2,3){\bigtoprightb}
  \rput(0,3){\bigtopleftb}
\end{pspicture}
\end{aligned}
\biggr) = Z\biggl(
\begin{aligned}
\begin{pspicture}(3.5,2.5)
  \rput(0,0){\bigcomposed}
\end{pspicture}
\end{aligned}
\biggr).
\end{equation}
The composition is an open-closed cobordism with two windows, but
multiplication by the appropriate power of the inverse window element
$a^{-1}$ removes these windows and results in the desired composite
open-closed cobordism.

How do we find the appropriate maps
\begin{equation}
\label{eq_universalmaps}
\phi:=a^{-2}\cdot Z\biggl(
\begin{aligned}
\psset{xunit=6mm,yunit=6mm}
\begin{pspicture}(3.85,2.5)
  \rput(2,0){\bigtopright}
  \rput(0,0){\bigtopleft}
\end{pspicture}
\end{aligned}
\biggr)\qquad\mbox{and}\qquad
\psi:=Z\biggl(
\begin{aligned}
\psset{xunit=6mm,yunit=6mm}
\begin{pspicture}(3.85,2.5)
  \rput(0,0){\bigbottom}
\end{pspicture}
\end{aligned}
\end{equation}
by which to pre- and postcompose? The linear map
$Z(S_{(\emptyset,1)})\colon A\otimes A\to A\otimes A$ associated with
the saddle is a morphism of $(A^{\otimes 2},A^{\otimes 2})$-bimodules,
and so the $2$-term chain complex associated with the crossing is a
chain complex of $(A^{\otimes 2},A^{\otimes 2})$-bimodules. Similarly,
the algebra $A$ associated with the arc forms an $(A,A)$-bimodule, and
so the $1$-term chain complex associated with the arc is a chain
complex of $(A,A)$-bimodules. In general, for each open end of a
tangle, we have one action of $A$.

The linear map associated with the composite~\eqref{eq_glued} is the
tensor product over $A\otimes A$ of these chain complexes, using the
appropriate left- and right-actions of $A\otimes A$ that correspond
to the two pairs of open ends of the tangles that are glued. This is
completely analogous to the situation for even tangles~\cite{Kh3}.
The map $\psi$ is the coequalizer that defines the tensor product
over $A\otimes A$ in the target of the saddle, whereas the map
$\phi$ is the unique map obtained from the universal property of a
similar coequalizer that defines the tensor product over $A\otimes
A$ in the source of the saddle because the differential of the chain
complex factors through that coequalizer. We make all this precise
in Section~\ref{sec_composition} below.

\subsection{Summary}

In an open-closed TQFT, we have operations for putting the building
blocks of~\eqref{eq_generators} on top of each other by gluing them
along their black boundaries (composition) and for putting them next
to each other by taking their disjoint unions (tensor product), but
there is no elementary operation for gluing them `sideways' along
their coloured boundaries. If in addition $A$ is strongly separable
and if $C=Z(A)$ with the canonical inclusion $\imath\colon C\to A$,
then the above prescription can be viewed as an additional operation
that allows us to glue open-closed cobordisms along their coloured
boundaries.

If $(A,C,\imath,\imath^\ast)$ is a state sum knowledgeable Frobenius
algebra with $A$ strongly separable, and $C=Z(A)$, then $C$ is
strongly separable, too. The commutative Frobenius algebra $C_{0,1}$
(Definition~\ref{def_khovanov}) of Lee in characteristic other than
$2$ and $C_{1,0}$ of Bar-Natan in any characteristic are strongly
separable, and so our strategy for composing tangles works in these
cases. Khovanov's Frobenius algebra $C_{0,0}$ in any characteristic
and Lee's Frobenius algebra in characteristic $2$ are not semisimple
and both fail to yield a tangle homology theory that is well behaved
under tangle composition defined in this way. This is no surprise.
The non-semisimplicity of Khovanov's Frobenius algebra is closely
related to the fact that it can detect \emph{global} properties of
the tangle that cannot be easily computed \emph{locally}, \ie\ from
the constituents of the tangle.

But does this mean that our method of composing tangles rules out the
most interesting example, namely Khovanov's? No, it does not.  There
is still a way out. It is known, for example, from the work of Lee and
Rasmussen~\cite{Lee,Rasmussen} and of Turner~\cite{Turner} that for
links, the chain complex and homology for Khovanov's Frobenius algebra
can be recovered from a filtration on the chain complex of Lee or
Bar-Natan, respectively, by computing the $E_0$- and $E_1$-pages of a
certain spectral sequence. When we extend this method of computing
Khovanov's chain complex from links to tangles, there are additional
conditions to satisfy, and an analogous spectral sequence only exists
over fields of finite characteristic (see Example~\ref{ex_barnatan},
Theorem~\ref{thm_sskfrob} and Examples~\ref{ex_filtered}
to~\ref{ex_filtered3}). These extensions are based on strongly
separable algebras, and they are well behaved under tangle
composition. Already the simplest example in characteristic $2$
(Example~\ref{ex_barnatan}) carries a filtration sophisticated enough
to recover our tangle extension of characteristic $2$ Khovanov
homology with all its non-trivial properties
(Theorem~\ref{thm_spectral} and Remark~\ref{rem_chartwo} below).

The tangle homology presented in this article then has the following
desirable properties:
\begin{enumerate}
\item
  The tangle homology behaves well with respect to composition.
\item
  It respects the monoidal structure of the category of tangles.
\item
  It is rich enough to recover the finite characteristic Khovanov
  homology for links as well as Rasmussen's $s$-invariant.
\end{enumerate}

Our use of open-closed TQFTs in order to define a tangle homology
theory is by no means the most general conceivable approach. For
example, in an open-closed TQFT, all the open-closed cobordisms
displayed in~\eqref{eq_generators} have algebraic operations
associated with them although only the composites displayed
in~\eqref{eq_crossing} to~\eqref{eq_copants} below are needed in order
to define the chain complexes. This is in analogy to the observation
that in Khovanov's link homology theory, the unit and counit of the
commutative Frobenius algebra are not needed to define the chain
complex. Notice, however, that both Bar-Natan's~\cite{BN2} and our
proof (Section~\ref{sec_Reidemeister} below) of Reidemeister move
invariance still need some of the additional cobordisms. We have
chosen to work with open-closed TQFTs mainly because the extension of
a $2$-dimensional TQFT to an open-closed TQFT is much better
understood than the corresponding question for the Temperley--Lieb
$2$-category. Notice also that in any tangle cobordism theory, the
object associated with the unknot, \ie\ in our case the homology of
the $1$-term complex $0\rightarrow C\rightarrow 0$ which is just $C$,
automatically carries the structure of a commutative Frobenius
algebra. Since we work with open-closed TQFTs, the object associated
with the arc, \ie\ the identity tangle on a single point, which is the
homology of $0\rightarrow A\rightarrow 0$ and thereby just $A$, forms
a symmetric Frobenius algebra. This is sufficient to yield an
algebraic model for the Temperley--Lieb $2$-category, but presumably
not necessary.

%
\section{Preliminaries}
%
\label{sect_prelim}

In this section, we introduce the notion of a knowledgeable
Frobenius algebra. This is the algebraic structure that is needed to
classify open-closed TQFTs. We then summarize the definitions and
the main properties of the commutative and anti-commutative cubes of
graded and filtered vector spaces that are used to define
Khovanov's, Lee's and Bar-Natan's chain complexes.

\subsection{Knowledgeable Frobenius algebras and open-closed TQFTs}
\label{sect_kfrob}

This section is a brief summary of the relevant definitions and
results. For more details, we refer to~\cite{LP1,LP2}.

\begin{defn}
Let $k$ be a field.
\begin{itemize}
\item
  A \emph{Frobenius algebra} $(A,\mu,\eta,\Delta,\epsilon)$ is an
  associative algebra $(A,\mu,\eta)$ with multiplication $\mu\colon
  A\otimes A\to A$ and unit $\eta\colon k\to A$ that forms a
  coassociative coalgebra $(A,\Delta,\epsilon)$ with
  comultiplication $\Delta\colon A\to A\otimes A$ and counit
  $\epsilon\colon A\to k$ such that the following condition holds,
\begin{equation}
  (\id_A\otimes\mu)\circ(\Delta\otimes\id_A)
  = \Delta\circ\mu
  = (\mu\otimes\id_A)\circ (\id_A\otimes\Delta).
\end{equation}
\item
  A Frobenius algebra $(A,\mu,\eta,\Delta,\epsilon)$ is called
  \emph{symmetric} if
  $\epsilon\circ\mu\circ\tau_{A,A}=\epsilon\circ\mu$. It is called
  \emph{commutative} if $\mu\circ\tau_{A,A}=\mu$. We have used the
  notation $\tau_{A,A}\colon A\otimes A\to A\otimes A, a\otimes
  b\mapsto b\otimes a$.
\item
  In a symmetric Frobenius algebra $(A,\mu,\eta,\Delta,\epsilon)$, we
  call the element $a=\mu\circ\Delta\circ\eta\colon k\to A$ the
  \emph{window element}.
\item
  An algebra $(A,\mu,\eta)$ is called \emph{strongly separable} if the
  canonical bilinear form
\begin{equation}
  g_{\mathrm{can}}\colon A\otimes A\to k,\quad
  a\otimes b\mapsto \tr_A(L_a\circ L_b)
\end{equation}
  is non-degenerate. Here we have written $L\colon A\to\End(A)$,
  $a\mapsto L_a$ with $L_a\colon A\to A, b\mapsto\mu(a\otimes b)$ for
  the left-regular representation.
\item
  A \emph{homomorphism of Frobenius algebras} $f\colon A\to A^\prime$
  is a linear map which is both a homomorphism of unital algebras and
  of counital coalgebras.
\end{itemize}
\end{defn}

Every strongly separable algebra is finite-dimensional and
semisimple. The converse implication holds in characteristic $0$ but
not in general in finite characteristic.

A symmetric Frobenius algebra is strongly separable if and only if the
window element is convolution invertible, \ie\ if there exists some
$a^{-1}\colon k\to A$ such that $\mu\circ(a\otimes
a^{-1})=\eta=\mu\circ(a^{-1}\otimes a)$. The window element is
central, \ie\ $\mu\circ(\phi\otimes a)=\mu\circ(a\otimes\phi)$ for all
linear maps $\phi\colon A\to A$, and so is its inverse when it
exists. We write $a\cdot\id_A:=\mu_A\circ(a\otimes\id_A)$, \etc.

The Frobenius algebra $C_{h,t}$ of Definition~\ref{def_khovanov} has
the window element $a=(\mu\circ\Delta\circ\eta)(1)=2x-h$ which
satisfies $a^2=(h^2+4t)\cdot 1$, and so $C_{h,t}$ is strongly
separable if and only if $h^2+4t\neq 0$. In particular, Khovanov's
Frobenius algebra $C_{0,0}$ is not strongly separable; Lee's $C_{0,1}$
is strongly separable if and only if $\chr k\neq 2$; and Bar-Natan's
$C_{1,0}$ is strongly separable in any characteristic.

\begin{defn}
A \emph{knowledgeable Frobenius algebra} $(A,C,\imath,\imath^\ast)$
consists of,
\begin{itemize}
\item
  a symmetric Frobenius algebra $A=(A,\mu_A,\eta_A,\Delta_A,\epsilon_A)$,
\item
  a commutative Frobenius algebra $C=(C,\mu_C,\eta_C,\Delta_C,\epsilon_C)$,
\item
  linear maps $\imath\colon C\to A$ and $\imath^\ast\colon A\to C$,
\end{itemize}
such that $\imath\colon C\to A$ is a homomorphism of algebras and the
following conditions hold,
\begin{alignat}{2}
\label{eq_kfrob1}
  \mu_A\circ(\imath\otimes\id_A)
    &= \mu_A\circ\tau_{A,A}\circ(\imath\otimes\id_A)
    &\qquad& \mbox{(knowledge),}\\
\label{eq_kfrob2}
  \epsilon_C\circ\mu_C\circ(\id_C\otimes\imath^\ast)
    &= \epsilon_A\circ\mu_A\circ(\imath\otimes\id_A)
    &\qquad& \mbox{(duality),}\\
\label{eq_kfrob3}
  \mu_A\circ\tau_{A,A}\circ\Delta_A
    &= \imath\circ\imath^{\ast}
    &\qquad& \mbox{(Cardy condition).}
\end{alignat}
The defining relations of a knowledgeable Frobenius algebra are
precisely the Moore--Segal relations~\cite{MS}. A \emph{homomorphism
of knowledgeable Frobenius algebras} $f\colon
(A,C,\imath,\imath^\ast)\to(A^\prime,C^\prime,\imath^\prime,{\imath^\prime}^\ast)$
is a pair $f=(f_1,f_2)$ such that $f_1\colon A\to A^\prime$ and
$f_2\colon C\to C^\prime$ are homomorphisms of Frobenius algebras and
such that $f_1\circ\imath = \imath^\prime\circ f_2$ and
$f_2\circ\imath^\ast ={\imath^\prime}^\ast\circ f_1$. We denote the
category of knowledgeable Frobenius algebras and their homomorphisms
by $\cat{K-Frob}(\cat{Vect}_k)$.
\end{defn}

The category $\cat{K-Frob}(\cat{Vect}_k)$ is a symmetric monoidal
category. The tensor product of two knowledgeable Frobenius algebras
${\mathbbm A}=(A,C,\imath,\imath^\ast)$ and ${\mathbbm
A^\prime}=(A^\prime,C^\prime,\imath^\prime,{\imath^\prime}^\ast)$ is
${\mathbbm A}\otimes{\mathbbm A^\prime}=
(A\otimes_kA^\prime,C\otimes_kC^\prime,\imath\otimes\imath^\prime,\imath^\ast\otimes{\imath^\prime}^\ast)$.
The unit of the monoidal structure is $(k,k,\id_k,\id_k)$ where $k$ is
equipped with the trivial Frobenius algebra structure
$\Delta_k(1)=1\otimes 1$ and $\epsilon_k(1)=1$. The associativity and
unit constraints and the symmetric braiding are inherited from
$\cat{Vect}_k$.

There is a category $\cat{2Cob}^{\mathrm{ext}}$ of \emph{open-closed
cobordisms} which generalizes the familiar category $\cat{2Cob}$ of
$2$-dimensional cobordisms from manifolds with boundary to manifolds
with corners that have a particular global structure as sketched
in~\eqref{eq_faces}. The objects of $\cat{2Cob}^{\mathrm{ext}}$ are
sequences $(m_1,\ldots,m_k)\in{\{0,1\}}^k$, $k\in\N_0$, which specify
a particular sequence of circles ($m_j=0$) and arcs ($m_j=1$). The
morphisms are the equivalence classes of suitable $2$-manifolds with
corners up to orientation-preserving diffeomorphisms that restrict to
the identity on the black boundary.  There exists a generators and
relations description for the morphisms of
$\cat{2Cob}^{\mathrm{ext}}$. The generators are the equivalence
classes of the open-closed cobordisms of~\eqref{eq_generators}, and
the relations among the morphisms of $\cat{2Cob}^{\mathrm{ext}}$ are
precisely the defining equations of a knowledgeable Frobenius
algebra. We recommend that the reader draw the pictures corresponding
to all equations that appear in this section. The main result
of~\cite{LP1} can be stated as follows.

\begin{thm}[see~\cite{LP1}]
The category $\cat{2Cob}^{\mathrm{ext}}$ of open-closed cobordisms is
the strict symmetric monoidal category freely generated by a
knowledgeable Frobenius algebra object
$(A,C,\imath,\imath^\ast)$. Here, $A=(1)$ represents the
diffeomorphism type of the arc; $C=(0)$ is the one of the circle; and
the structure maps of the knowledgeable Frobenius algebra object are
the equivalence classes of the open-closed cobordisms depicted
in~\eqref{eq_generators}.
\end{thm}

\begin{defn}
An \emph{open-closed Topological Quantum Field Theory (TQFT)} is a
symmetric monoidal functor
$Z\colon\cat{2Cob}^{\mathrm{ext}}\to\cat{Vect}_k$. A
\emph{homomorphism of open-closed TQFTs} is a monoidal natural
transformation between such functors. We denote the category of
open-closed TQFTs and their homomorphisms by
\begin{equation}
\cat{OC-TQFT}(\cat{Vect}_k):=\cat{Symm-Mon}(\cat{2Cob}^\mathrm{ext},\cat{Vect}_k)
\end{equation}
\end{defn}

\begin{cor}
There is an equivalence of symmetric monoidal categories
\begin{equation}
\begin{aligned}
\xymatrix{
  \cat{OC-TQFT}(\cat{Vect}_k)\ar@<.5ex>[rr]^{F}&&\cat{K-Frob}(\cat{Vect}_k)\ar@<.5ex>[ll]^{G}
}
\end{aligned}
\end{equation}
between the category of open-closed TQFTs and that of knowledgeable
Frobenius algebras. Given an open-closed TQFT $Z$, we call $F(Z)$ the
\emph{associated knowledgeable Frobenius algebra}. It is determined up to
isomorphism of knowledgeable Frobenius algebras by
$(A,C,\imath,\imath^\ast)$ where $A=Z((1))$, $C=Z((0))$, and the
structure maps are the images of the open-closed cobordisms
of~\eqref{eq_generators} under $Z$.
\end{cor}

In order to turn Bar-Natan's `picture world'~\cite{BN2} into
algebra, we show below that it is sufficient if the corresponding
knowledgeable Frobenius algebra $(A,C,\imath,\imath^\ast)$ has a
commutative Frobenius algebra $(C,\mu_C,\eta_C,\Delta_C,\epsilon_C)$
that satisfies Bar-Natan's conditions.

\begin{defn}\hfill
\begin{enumerate}
\item
  A knowledgeable Frobenius algebra $(A,C,\imath,\imath^\ast)$ is said
  to \emph{satisfy Bar-Natan's conditions} if $C$ satisfies them
  (Definition~\ref{def_barnatancond}).
\item
  An open-closed TQFT is said to \emph{satisfy Bar-Natan's conditions}
  if its associated knowledgeable Frobenius algebra does so.
\end{enumerate}
\end{defn}

\subsection{Strongly separable algebras and state sum TQFTs}

Let us now focus on strongly separable symmetric Frobenius algebras.
These are used to construct state sum open-closed TQFTs
in~\cite{LP2}. From every such algebra, one can construct an
open-closed TQFT which is described by the knowledgeable Frobenius
algebra defined in the following proposition.

\begin{prop}[see Theorem 2.22 of~\cite{LP2}]
\label{prop_statesum}
Let $(A,\mu_A,\eta_A,\Delta_A,\epsilon_A)$ be a strongly separable
symmetric Frobenius algebra with window element $a$. Then
\begin{equation}
\label{eq_idempotent}
  p:=(a^{-1}\cdot\id_A)\circ\mu_A\circ\tau_{A,A}\circ\Delta_A
  \colon A\to A
\end{equation}
forms an idempotent whose image is the centre of $A$, \ie\
$p(A)=Z(A)$.

In this case, there is a knowledgeable Frobenius algebra
$(A,C,\imath,\imath^\ast)$ with $C:=Z(A)$, equipped with the
commutative Frobenius algebra structure
$(C,\mu_C,\eta_C,\Delta_C,\epsilon_C)$ given by
\begin{eqnarray}
\mu_C      &=& p\circ\mu_A,\\
\eta_C     &=& p\circ\eta_A,\\
\Delta_C   &=& (p\otimes p)\circ\Delta_A\circ(a\cdot\id_A),\\
\epsilon_C &=& \epsilon_A\circ(a^{-1}\cdot\id_A),
\end{eqnarray}
and with the canonical inclusion $\imath\colon C\to A$ and the map
$\imath^\ast=(a\cdot\id_A)\circ p$.
\end{prop}

\begin{defn}
A \emph{state sum knowledgeable Frobenius algebra} is a
knowledgeable Frobenius algebra that arises from a strongly
separable symmetric Frobenius algebra by the construction of
Proposition~\ref{prop_statesum}. A \emph{state sum open-closed TQFT}
$Z\colon\cat{2Cob}^\mathrm{ext}\to\cat{Vect}_k$ is an open-closed
TQFT whose associated knowledgeable Frobenius algebra is a state sum
knowledgeable Frobenius algebra.
\end{defn}

In order to compose tangles, we need some idempotents defined in the
subsequent proposition. Let $(A,\mu_A,\eta_A,\Delta_A,\epsilon_A)$ be
a symmetric Frobenius algebra. For $j\in\N$, we denote by
\begin{equation}
  \mu_A^{(j+1)}:=\mu_A\circ(\mu_A^{(j)}\otimes\id_A),\qquad
  \mu_A^{(2)}:=\mu_A,\qquad
  \mu_A^{(1)}:=\id_A
\end{equation}
and by
\begin{equation}
  \Delta_A^{(j+1)}:=(\Delta_A^{(j)}\otimes\id_A)\circ\Delta_A,\qquad
  \Delta_A^{(2)}:=\Delta_A,\qquad
  \Delta_A^{(1)}:=\id_A
\end{equation}
the iterated multiplication and comultiplication. We also write
$A^{\otimes(j+1)}:=A^{\otimes j}\otimes A$, $A^{\otimes 1}:=A$ and
$A^{\otimes 0}:=\1$, and
$a^{j+1}\cdot\id_A:=(a^j\cdot\id_A)\circ(a\cdot\id_A)$ and
$a^0\cdot\id_A:=\id_A$.

\begin{prop}
\label{prop_projproperty}
Let $(A,\mu_A,\eta_A,\Delta_A,\epsilon_A)$ be a strongly separable
symmetric Frobenius algebra with window element $a$, and let $p$ be
as in Proposition~\ref{prop_statesum}. Then the linear maps,
\begin{alignat}{2}
  P_{j\ell}&:=\Delta^{(j)}\circ(a^{-(j-1)}\cdot\id_A)\circ\mu^{(\ell)}&\colon
    A^{\otimes\ell}\to A^{\otimes j},\\
  Q_{j\ell}&:=\Delta^{(j)}\circ(a^{-(j-1)}\cdot\id_A)\circ p\circ\mu^{(\ell)}&\colon
    A^{\otimes\ell}\to A^{\otimes j},
\end{alignat}
satisfy
\begin{equation}
\label{eq_projproperty}
  P_{j\ell}\circ P_{\ell m} = P_{jm}\qquad\mbox{and}\qquad
  Q_{j\ell}\circ Q_{\ell m} = Q_{jm}
\end{equation}
for all $j,\ell,m\in\N$. In particular, $P_{jj}$ and $Q_{jj}$ are
idempotents, and we have $P_{11}=\id_A$ and $Q_{11}=p$. There are
isomorphisms,
\begin{equation}
  P_{jj}(A^{\otimes j})\cong A\qquad\mbox{and}\qquad
  Q_{jj}(A^{\otimes j})\cong p(A),
\end{equation}
given by $P_{1j}$ and $Q_{1j}$, respectively, whose inverses are
$P_{j1}$ and $Q_{j1}$, respectively.
\end{prop}

\subsection{Tensor products and coequalizers}

The images of these idempotents can be characterized by the following
universal properties that may be more familiar to the reader.

\begin{prop}
\label{prop_colimits}
Let $(A,\mu_A,\eta_A,\Delta_A,\epsilon_A)$ be a strongly separable
symmetric Frobenius algebra, and let $p=Q_{11}$ and $P_{22}$ be as
above.
\begin{enumerate}
\item
  The linear map $p\colon A\to p(A)$ is the coequalizer
\begin{equation}
\xymatrix{
  A\otimes A\ar@<.5ex>[rr]^-{\mu_A}\ar@<-.5ex>[rr]_-{\mu_A\circ\tau_{A,A}}&& A\ar[rr]^-{p}&&p(A).
}
\end{equation}
\item
  The linear map $P_{22}\colon A\otimes A\to P_{22}(A\otimes A)$ is
  the coequalizer
\begin{equation}
\xymatrix{
  A\otimes A\otimes A
    \ar@<.5ex>[rrr]^-{\mu_A\otimes\id_A}
    \ar@<-.5ex>[rrr]_-{\id_A\otimes\mu_A}
  &&& A\otimes A\ar[rr]^-{P_{22}}&&P_{22}(A\otimes A).
}
\end{equation}
\end{enumerate}
This means that $p(A)\cong A/[A,A]$ and that $P_{22}(A\otimes
A)\cong A\otimes_A A$ where $A$ is viewed as an $(A,A)$-bimodule.
\end{prop}

\begin{proof}
We prove this proposition for algebra objects in an arbitrary abelian
symmetric monoidal category $\cal{C}$. Given the
idempotent~\eqref{eq_idempotent}, we denote its image splitting by
$p=\im p\circ\coim p$ with $\im p\colon p(A)\to A$ and $\coim p\colon
A\to p(A)$. In the abelian category $\cal{C}$, the idempotent $p$
splits, and so $\coim p\circ\im p=\id_{p(A)}$.

We assume that $(A,\mu_A,\eta_A,\Delta_A,\epsilon_A)$ is a strongly
separable symmetric Frobenius algebra object in some abelian symmetric
monoidal category $\cal{C}$ and prove the following claim:
\begin{enumerate}
\item
  The coimage $\coim p\colon A\to p(A)$ is the coequalizer
\begin{equation}
\xymatrix{
  A\otimes
A\ar@<.5ex>[rr]^-{\mu_A}\ar@<-.5ex>[rr]_-{\mu_A\circ\tau_{A,A}}&&
A\ar[rr]^-{\coim p}&&p(A).  }
\end{equation}
\item
  The coimage $\coim P_{22}\colon A\otimes A\to P_{22}(A\otimes A)$
  is the coequalizer
\begin{equation}
\xymatrix{
  (A\otimes A)\otimes A \ar@<.5ex>[rrr]^-{\mu_A\otimes\id_A}
    \ar@<-.5ex>[rrr]_-{(\id_A\otimes\mu_A)\circ\alpha_{A,A,A}} &&&
    A\otimes A\ar[rr]^-{\coim P_{22}}&&P_{22}(A\otimes A).  }
\end{equation}
\end{enumerate}

In order to prove (1.), note that since $\coim p=\coim p\circ p$, the
morphism $\coim p$ satisfies $\coim p\circ\mu_A=\coim
p\circ\mu_A\circ\tau_{A,A}$. Given any morphism $f\colon A\to B$ such
that $f\circ\mu_A=f\circ\mu_A\circ\tau_{A,A}$, there is a morphism
$\phi\colon p(A)\to B$ given by $\phi:=f\circ\im p$ such that
\begin{equation}
  \phi\circ\coim p = f\circ p = f\circ\mu_A\circ\tau_{A,A}\circ\Delta_A\circ(a^{-1}\cdot\id_A)
  = f\circ\mu_A\circ\Delta_A\circ(a^{-1}\cdot\id_A) = f.
\end{equation}
If $\psi\colon p(A)\to B$ also satisfies $f=\psi\circ\coim p$, then
$\phi=f\circ\im p=\psi\circ\id_{p(A)} = \psi$, and so $\phi$ is unique
with that property.

In order to prove (2.), note that since $\coim P_{22}=\coim
P_{22}\circ P_{22}$, and because of associativity, we have $\coim
P_{22}\circ(\mu_A\otimes\id_A) = \coim
P_{22}\circ(\id_A\otimes\mu_A)\circ\alpha_{A,A,A}$. Given any morphism
$f\colon A\otimes A\to B$ such that
$f\circ(\mu_A\otimes\id_A)=f\circ(\id_A\otimes\mu_A)\circ\alpha_{A,A,A}$,
there is a morphism $\phi:=f\circ\im P_{22}\colon P_{22}(A\otimes
A)\to B$ such that
\begin{eqnarray}
  \phi\circ\coim P_{22} &=& f\circ P_{22}
    = f\circ(\mu_A\otimes\id_A)\circ\alpha_{A,A,A}^{-1}\circ(\id_A\otimes\Delta_A)\circ(\id_A\otimes(a^{-1}\cdot\id_A))\nn\\
    &=& f\circ(\id_A\otimes (\mu_A\circ\Delta_A\circ(a^{-1}\cdot\id_A)))=f.
\end{eqnarray}
If $\psi\colon P_{22}(A\otimes A)\to B$ also satisfies $\psi\circ\coim
P_{22}=f$, then $\phi=f\circ\im P_{22}=\psi\circ\coim P_{22}\circ\im
P_{22}=\psi$, and so $\phi$ is unique with that property.
\end{proof}

Using the fact that $P_{22}$ is a morphism of $(A,A)$-bimodules and
that $A\otimes_A A\cong A$ as $(A,A)$-bimodules, one can obtain
similar characterizations for the images $P_{jj}(A^{\otimes j})$ and
$Q_{jj}(A^{\otimes j})$ for all $j\in\N$. Note that the map $\phi$
of~\eqref{eq_universalmaps} in the introduction is $\phi=\id_A\otimes
Q_{21}$ whereas $\psi=P_{13}$.

The remainder of the present subsection is merely a `wild' combination
of the above results which introduces a notation relevant to
Section~\ref{sec_composition} below. The reader may wish to skip this
material on first reading.

The $k$-fold tensor product $M:=A^{\otimes k}$, $k\in\N$, forms both a
left $A$-module and a right $A$-module in $k$ different ways, by
acting with $A$ from the left or from the right on the $j$-th tensor
factor, $j\in\{1,\ldots,k\}$. We denote the left actions by
$\lambda_j\colon A\otimes M\to M$ and the right actions by
$\rho_j\colon M\otimes A\to A$, respectively. We introduce the
following notation for the coequalizer in $\cat{Vect}_k$ that
`coequalizes' the $j$-th left action with the $\ell$-th right action,
$j,\ell\in\{1,\ldots,k\}$:
\begin{equation}
\label{eq_tensorover}
\xymatrix{
  A\otimes M\ar@<.5ex>[rrr]^-{\lambda_j}\ar@<-.5ex>[rrr]_-{\rho_\ell\circ\tau_{A,M}}&&&
  M\ar[rr]^-{{}_{\rho_\ell}\otimes_{\lambda_j}}&&
  {}_{\rho_\ell}\otimes_{\lambda_j}(M).
}
\end{equation}
We note that ${}_{\rho_\ell}\otimes_{\lambda_j}(M)$ is an
$(A^{\otimes(k-1)},A^{\otimes(k-1)})$-bimodule whose associated
$A$-actions we denote by
$\lambda_1,\ldots,\lambda_{j-1},\lambda_{j+1},\ldots,\lambda_k$ and
$\rho_1,\ldots,\rho_{\ell-1},\rho_{\ell+1},\ldots,\rho_k$,
respectively. For $j\neq r$, $\ell\neq s$, taking these coequalizers
commutes, and we have canonical isomorphisms
\begin{equation}
  {}_{\rho_\ell}\otimes_{\lambda_j}\bigl({}_{\rho_s}\otimes_{\lambda_r}(M)\bigr)
\cong
  {}_{\rho_s}\otimes_{\lambda_r}\bigl({}_{\rho_\ell}\otimes_{\lambda_j}(M)\bigr)
\end{equation}
of $(A^{\otimes(k-2)},A^{\otimes(k-2)})$-bimodules.

\begin{cor}
\label{cor_standardpkl}
Let $(A,\mu_A,\eta_A,\Delta_A,\epsilon_A)$ be a strongly separable
symmetric Frobenius algebra, $M:=A^{\otimes k}$, $k\in\N$, and denote
the left $A$- and right $A$-actions on $M$ as above.
\begin{enumerate}
\item
  $P_{1k}\colon M\to A$ is the composition of coequalizer maps
\begin{equation}
\label{eq_standardpkl}
  P_{1k} =
    {}_{\rho_1}\otimes_{\lambda_2}\circ
    {}_{\rho_2}\otimes_{\lambda_3}\circ\cdots\circ
    {}_{\rho_{k-1}}\otimes_{\lambda_k}
\end{equation}
  in $\cat{Vect}_k$. 
\item
  $Q_{1k}\colon M\to p(A)$ is the composition of coequalizer maps
\begin{equation}
\label{eq_standardqkl}
  Q_{1k} = {}_{\rho_k}\otimes_{\lambda_1}\circ P_{1k}
\end{equation}
  in $\cat{Vect}_k$.
\end{enumerate}
\end{cor}

This result allows us to characterize the map $\psi$
of~\eqref{eq_universalmaps} in the introduction by a universal
property as soon as we have organized the tensor factors and the
associated left and right actions accordingly.

For permutations $\sigma,\tau\in S_k$, we define
\begin{equation}
  P_{1k}^{(\sigma,\tau)} :=
    {}_{\rho_{\tau(1)}}\otimes_{\lambda_{\sigma(2)}}\circ
    {}_{\rho_{\tau(2)}}\otimes_{\lambda_{\sigma(3)}}\circ\cdots\circ
    {}_{\rho_{\tau(k-1)}}\otimes_{\lambda_{\sigma(k)}}
\end{equation}
and
\begin{equation}
  Q_{1k}^{(\sigma,\tau)} := {}_{\rho_{\tau(k)}}\otimes_{\lambda_{\sigma(1)}}\circ P_{1k}^{(\sigma,\tau)}.
\end{equation}
Consider a linear map $f\colon M\to M$, $M:=A^{\otimes k}$, $k\in\N$,
that satisfies
\begin{equation}
  \lambda_{\sigma(j)}\circ(\id_A\otimes f) = f\circ\lambda_j
  \quad\mbox{and}\quad
  \rho_{\tau(\ell)}\circ(f\otimes\id_A) = f\circ\rho_\ell,
\end{equation}
for all $j,\ell\in\{1,\ldots,k\}$, for some permutations
$\sigma,\tau\in S_k$. We call $f$ a \emph{$(\sigma,\tau)$-morphism
of $(A^{\otimes k},A^{\otimes k})$-bimodules}. In fact, $f$ is an
honest morphism of $(A^{\otimes k},A^{\otimes k})$-bimodules as soon
as one defines the left and right actions of $A^{\otimes k}$ on its
target by taking the actions on the source and pre-composing them
with the permutations $\sigma$ and $\tau$, respectively. Our prime
examples of $(\sigma,\tau)$-morphisms of bimodules are the maps
\begin{equation}
  (\mu_A\otimes\id_A)\circ(\id_A\otimes\tau_{A,A})\circ(\Delta_A\otimes\id_A)\colon
    A\otimes A\to A\otimes A
\end{equation}
with $\sigma=\id$ and $\tau=(12)$ in cycle notation, and
\begin{equation}
  (\id_A\otimes\mu_A)\circ(\tau_{A,A}\otimes\id_A)\circ(\id_A\otimes\Delta_A)\colon
    A\otimes A\to A\otimes A
\end{equation}
with $\sigma=(12)$ and $\tau=\id$. These maps can be obtained from the
open-closed cobordisms depicted in~\eqref{eq_crossing}
and~\eqref{eq_crossing2}.

If $f\colon M\to M$ is a $(\sigma,\tau)$-morphism of $(A^{\otimes
k},A^{\otimes k})$-bimodules, then for $j,\ell\in\{1,\ldots,k\}$, the
composition ${}_{\rho_{\tau(\ell)}}\otimes_{\lambda_{\sigma(j)}}\circ
f$ factors through the coequalizer ${}_{\rho_\ell}\otimes_{\lambda_j}$
and uniquely defines a map $\hat f$ such that the diagram
\begin{equation}
\label{eq_coeqmaps}
\begin{aligned}
\xymatrix{
  M\ar[rr]^f\ar[dd]_{{}_{\rho_\ell}\otimes_{\lambda_j}}&&
    M\ar[dd]^{{}_{\rho_{\tau(\ell)}}\otimes_{\lambda_{\sigma(j)}}}\\
  \\
  {}_{\rho_\ell}\otimes_{\lambda_j}(M)\ar[rr]_{\hat f}&&
    {}_{\rho_{\tau(\ell)}}\otimes_{\lambda_{\sigma(j)}}(M)
}
\end{aligned}
\end{equation}
commutes. Iterating this result, we obtain:

\begin{cor}
Let $(A,\mu_A,\eta_A,\Delta_A,\epsilon_A)$ be a strongly separable
symmetric Frobenius algebra, $M:=A^{\otimes k}$, $k\in\N$, and let
$f\colon M\to M$ be a $(\sigma,\tau)$-morphism of $(A^{\otimes
k},A^{\otimes k})$-bimodules for some permutations $\sigma,\tau\in
S_k$.
\begin{itemize}
\item
  The map $\hat f:=P_{1k}^{(\sigma,\tau)}\circ f\circ P_{k1}$ is the
  unique linear map such that the diagram
\begin{equation}
\begin{aligned}
\xymatrix{
  M\ar[rr]^f\ar[dd]_{P_{1k}}&&
    M\ar[dd]^{P_{1k}^{(\sigma,\tau)}}\\
  \\
  A\ar[rr]_{\hat f}&&A
}
\end{aligned}
\end{equation}
  commutes.
\item
  The map $\hat f:=Q_{1k}^{(\sigma,\tau)}\circ f\circ Q_{k1}$ is the
  unique linear map such that the diagram
\begin{equation}
\begin{aligned}
\xymatrix{
  M\ar[rr]^f\ar[dd]_{Q_{1k}}&&
    M\ar[dd]^{Q_{1k}^{(\sigma,\tau)}}\\
  \\
  p(A)\ar[rr]_{\hat f}&&p(A)
}
\end{aligned}
\end{equation}
  commutes.
\end{itemize}
\end{cor}

\begin{proof}
This result follows from~\eqref{eq_projproperty}.
\end{proof}

This result is the origin of the map $\phi$
in~\eqref{eq_universalmaps} in the introduction. The $Q_{21}$ in
$\phi=\id_A\otimes Q_{21}$ is the first factor in such a composition
$\hat f$.

\subsection{Gradings and filtrations}

A $k$-vector space $V$ is called \emph{graded} if it is of the form
$V\cong\bigoplus_{j\in\Z}V^j$. We say that a vector $v\in V^j$ is
\emph{homogeneous of degree} $\deg(v)=j\in\Z$. If both $V$ and $W$ are graded
vector spaces, then their tensor product $V\otimes W$ is graded with
${(V\otimes W)}^j=\bigoplus_{\ell\in\Z}V^\ell\otimes W^{j-\ell}$. The
underlying field $k$, \ie\ the unit object of the monoidal structure
of $\cat{Vect}_k$, is considered graded with every vector homogeneous
of degree $0$, \ie\ $k=k^0$. A linear map $f\colon V\to W$ is said to
be \emph{graded of degree} $\ell\in\Z$ if $f(V^j)\subseteq W^{j+\ell}$
for all $j\in\Z$.  In this case, we write $\deg(f)=\ell$. For the
composition of graded linear maps $f\colon V\to W$ and $g\colon W\to
U$, we have $\deg(g\circ f)=\deg(g)+\deg(f)$. If $V$ is a graded
vector space, we define for each $j\in\Z$ a graded vector space
$V\{j\}$ by ${(V\{j\})}^\ell:=V^{\ell-j}$. The operator $\{j\}$ is
called the \emph{grading shift by} $j$. If the linear map $f\colon
V\to W$ is graded of degree $\ell$, then the map $f\colon V\to
W\{-\ell\}$ is graded of degree $0$.

A $k$-vector space $V$ is called \emph{filtered} if there is a
family ${\{F^jV\}}_{j\in\Z}$ of subspaces $F^jV\subseteq V$ such
that $F^{j+1}V\subseteq F^jV$ for all $j\in\Z$. The filtration is
said to be \emph{bounded} if there exist some $j_0,j_1\in\Z$ with
$F^{j_0}V=V$ and $F^{j_1}V=\{0\}$. We say that a vector $v\in V$ is
\emph{of degree} $\deg(v):=\max\{j\in\Z\mid\, v\in F^jV\}$. If both
$V$ and $W$ are filtered vector spaces, then their tensor product
$V\otimes W$ is filtered with $F^j(V\otimes W)=\spann\{v\otimes
w\mid\,v\in V, w\in W, \deg(v)+\deg(w)\geq j\}$.  The underlying
field is considered filtered with $F^0k=k$ and $F^1k=\{0\}$. A
linear map $f\colon V\to W$ is said to be \emph{filtered of degree}
$\ell\in Z$ if $f(F^jV)\subseteq F^{j+\ell}W$ for all $j\in\Z$. We
write $\deg(f)=\ell$. Again, for the composition of filtered linear
maps $f\colon V\to W$ and $g\colon W\to U$, we have $\deg(g\circ
f)=\deg(g)+\deg(f)$. If $V$ is a filtered vector space, we define
for each $j\in\Z$ a filtered vector space $V\{j\}$ by
$F^\ell(V\{j\}):=F^{\ell-j}V$. If the linear map $f\colon V\to W$ is
filtered of degree $\ell$, then the map $f\colon V\to W\{-\ell\}$ is
filtered of degree $0$.

If $V$ is a graded vector space with $V\cong\bigoplus_{j\in\Z}V^j$,
then it is filtered with $F^jV=\bigoplus_{\ell\geq j}V^\ell$. If $V$
is a filtered vector space, there is an \emph{associated graded
vector space} $V\cong\bigoplus_{j\in\Z}V^j$ with
$V^j:=F^jV/F^{j+1}V$. If $f\colon V\to W$ is a filtered map of
degree $\ell$ between filtered vector spaces, then $f(F^jV)\subseteq
F^{j+\ell}W$ for all $j$, and so $f$ descends to the quotients to
give maps $f^j\colon V^j\to W^{j+\ell}$. It thereby yields the
\emph{associated graded map} $\bar f\colon V\to W$ which is graded
of degree $\ell$. If $f\colon V\to W$ and $g\colon W\to U$ are
filtered maps of some degree with associated graded maps $\bar f$
and $\bar g$, then $\bar g\circ\bar f$ is the associated graded map
of $g\circ f$.

The category $\cat{grdVect}_k$ of \emph{graded $k$-vector spaces}
has graded $k$-vector spaces as objects and graded $k$-linear maps
of degree $0$ as morphisms. The category $\cat{fltVect}_k$ of
\emph{filtered $k$-vector spaces} has $k$-vector spaces with bounded
filtration as objects and filtered $k$-linear maps of degree $0$ as
morphisms. Both categories are abelian symmetric monoidal
categories.

\begin{defn}
An open-closed TQFT $Z\colon\cat{2Cob}^\mathrm{ext}\to\cat{Vect}_k$
is called \emph{Euler-graded} [\emph{Euler-filtered}] if the vector
spaces $C:=Z((0))$ and $A:=Z((1))$ are both graded [filtered] and if
for every morphism $M$ of $\cat{2Cob}^\mathrm{ext}$, the associated
linear map $Z(M)$ is graded [filtered] of degree $\deg(Z(M))=d(M)$
where
\begin{equation}
\label{eq_degree}
  d(M) := 2\chi(M)-|\Pi_0(\del_1M)| + \omega(M),
\end{equation}
\ie\ twice the Euler characteristic minus the number of connected
components of the coloured boundary $\del_1M$ of $M$ plus the window
number $\omega(M)$ of $M$. The window number $\omega(M)$ is the
number of components of $\del_1M$ that are diffeomorphic to $S^1$.
\end{defn}

\begin{prop}
\label{prop_degree} The degree $d(M)$ of~\eqref{eq_degree} has the
following properties.
\begin{enumerate}
\item
  It forms a topological invariant, and so it is well defined for the
  morphisms of $\cat{2Cob}^\mathrm{ext}$ which are equivalence classes
  of open-closed cobordisms modulo orientation-preserving
  diffeomorphisms that restrict to the identity on the black boundary.
\item
  If $M$ is a cylinder over an arc or a circle, then $d(M)=0$.
\item
  If $M=M_1\otimes M_2$ is a tensor product (disjoint union) of
  open-closed cobordisms, then $d(M)=d(M_1)+d(M_2)$.
\item
  If $M=M_1\circ M_2$ is a composition of open-closed cobordisms, then
  $d(M)=d(M_1)+d(M_2)$.
\end{enumerate}
\end{prop}

\begin{proof}\hfill
\begin{enumerate}
\item
  This is obvious.
\item
  In order to prove (2.) to (4.), it is easiest to consider a
  simplicial model for $M$. For (2.), one can read off the degree from
  equations~(4.10) and~(4.11) of~\cite{LP2}.
\item
  This is obvious.
\item
  First, gluing two open-closed cobordisms along a circle has the
  following effect. For simplicity suppose that as a cell complex the
  circle is constructed from a single vertex and a single
  edge. Then gluing along this circle reduces the number of vertices
  by one and the number of edges by one, and so it leaves the Euler
  characteristic unchanged. Second, gluing two open-closed cobordisms
  along an arc with the simplest cell structure has the following
  effect. It reduces the number of vertices by two and the number of
  edges by one, and so it reduces the Euler characteristic by one. The
  claim holds because at the same time the gluing reduces
  $|\Pi_0(\del_1M)|-\omega(M)$ by two. If two pairs of disjoint arcs
  of $\del_1M$ are joined by the gluing, $|\Pi_0(\del_1M)|$ is reduced
  by two. If the gluing does not join two disjoint arcs, but rather
  the two ends of a single arc of $\del_1M$, then $\Pi_0(\del_1M)$ is
  unchanged, but $\omega(M)$ increases.
\end{enumerate}
\end{proof}

\begin{defn}
A knowledgeable Frobenius algebra $(A,C,\imath,\imath^\ast)$ is called
\emph{Euler-graded} [\emph{Euler-filtered}] if both $A$ and $C$ are
graded [filtered] vector spaces and if the structure maps are graded
[filtered] of the following degrees:
\begin{alignat}{4}
\label{eq_defdegree}
  \deg(\mu_A)&=&-1,\qquad\deg(\eta_A)&=&1,\qquad
    \deg(\Delta_A)&=&-1,\qquad\deg(\epsilon_A)&=1,\nn\\
  \deg(\mu_C)&=&-2,\qquad\deg(\eta_C)&=&2,\qquad
    \deg(\Delta_C)&=&-2,\qquad\deg(\epsilon_C)&=2,\\
  &&\deg(\imath)&=&-1,\qquad\deg(\imath^\ast)&=&-1.\qquad\qquad\quad&\nn
\end{alignat}
\end{defn}

\begin{prop}
An open-closed TQFT $Z\colon\cat{2Cob}^\mathrm{ext}\to\cat{Vect}_k$ is
Euler-graded [Euler-filtered] if and only if its associated
knowledgeable Frobenius algebra is Euler-graded [Euler-filtered].
\end{prop}

\begin{proof}
Given an Euler-graded [Euler-filtered] open-closed TQFT $Z$, one
simply computes~\eqref{eq_degree} for the generators
of~\eqref{eq_generators} in order to obtain~\eqref{eq_defdegree}.
Triangulations of the generators are displayed, for example, in
equation~(4.15) of~\cite{LP2}.

Let now $(A,C,\imath,\imath^\ast)$ be an Euler-graded [Euler-filtered]
knowledgeable Frobenius algebra and $Z$ be an open-closed TQFT whose
associated knowledgeable Frobenius algebra is given by
$(A,C,\imath,\imath^\ast)$. Every open-closed cobordism $M$ is
equivalent to a composition and tensor product of the
generators~\eqref{eq_generators} (Proposition~3.9 of~\cite{LP1}).  The
degree is defined in~\eqref{eq_defdegree} in such a way that the
claim~\eqref{eq_degree} holds for each generator. By
Proposition~\ref{prop_degree}, it holds for tensor products with other
generators, for tensor products with cylinders, and for composites of
these. This implies that $\deg(Z(M))=d(M)$ as in~\eqref{eq_degree} for
all morphisms $M$.
\end{proof}

\begin{rem}\hfill
\label{rem_degree}
\begin{enumerate}
\item
  In~\eqref{eq_degree}, we employ twice the degree used by
  Bar-Natan~\cite{BN2} in order to have integer degrees for all
  open-closed cobordisms. This means in particular that on $C$, our
  degrees are twice Khovanov's~\cite{Kh} and that the Jones polynomial
  arises in our case as a polynomial in an indeterminate $A$ with
  $q=A^2$ rather than in terms of $q$ itself.
\item
  The right hand side of~\eqref{eq_degree} is also additive under
  gluing open-closed cobordisms along a component of their coloured
  boundary.
\item
  A TQFT $Z\colon\cat{2Cob}\to\cat{Vect}_k$ is Euler-graded
  [Euler-filtered] if and only if its associated commutative Frobenius
  algebra $C=Z(1)$ is Euler-graded [Euler-filtered], \ie\ satisfies
  the conditions of~\eqref{eq_defdegree} for $\mu_C$, $\eta_C$,
  $\Delta_C$, and $\epsilon_C$.
\end{enumerate}
\end{rem}

\begin{example}
\label{ex_grading}
The commutative Frobenius algebra $C_{0,0}$ of Khovanov
(Definition~\ref{def_khovanov}) is Euler-graded with $\deg(1)=2$ and
$\deg(x)=-2$. This grading is uniquely determined by the conditions
$\deg(\mu_C)=-2$ and $\deg(\Delta_C)=-2$, exploiting that
$\mu_C(1\otimes 1)=1$ and $\Delta_C(x)=x\otimes x$. For arbitrary
$h,t\in k$, the commutative Frobenius algebra $C_{h,t}$ is
Euler-filtered with $\deg_C(1)=2$ and $\deg_C(x)=-2$.
\end{example}

\subsection{Commutative and anti-commutative cubes}

In this and in the next section, we recall the definition of the
commutative and anti-commutative cubes of~\cite{Kh}. We write this
section for commutative cubes in arbitrary abelian monoidal
categories $\cal{C}$, but for the remainder of this article, it
suffices if the reader specializes everything to the categories
$\cat{Vect}_k$, $\cat{grdVect}_k$, or $\cat{fltVect}_k$. These are
all abelian, \cat{Ab}-enriched, and symmetric monoidal (see, for
example~\cite{maclane}).

Given a finite set $\cal{I}$, we denote by
$r(\cal{I}):=\{\,(\alpha,j)\in\cal{P}(\cal{I})\times\cal{I}\mid
j\notin\alpha\,\}$ the set of all pairs $(\alpha,j)$ where
$\alpha\subseteq\cal{I}$ is a subset that does not contain
$j\in\cal{I}$. We denote by $\sqcup$ the disjoint union of sets.

\begin{defn}
For a category $\cal{C}$ and a finite set $\cal{I}$, an
$\cal{I}$-\emph{cube}
$X=({\{X_\alpha\}}_{\alpha\subseteq\cal{I}},{\{X_{(\alpha,j)}\}}_{(\alpha,j)\in
r(\cal{I})})$ \emph{in} $\cal{C}$ consists of a family of objects
$X_\alpha$ and of a family of morphisms $X_{(\alpha,j)}\colon
X_\alpha\to X_{\alpha\sqcup\{j\}}$ of $\cal{C}$. The objects
$X_\alpha$ are called the \emph{vertices} of the cube, and the
morphisms $X_{(\alpha,j)}$ the \emph{edges}. For a given vertex
$X_\alpha$, we call $|\alpha|\in\{0,\ldots,|\cal{I}|\}$ its
\emph{height}.
\begin{enumerate}
\item
  The cube $X$ is called \emph{commutative} if for all
  $j,\ell\in\cal{I}$, $j\neq\ell$, and for all
  $\alpha\subseteq\cal{I}\backslash\{j,\ell\}$, the diagram
\begin{equation}
\label{eq_commutes}
\begin{aligned}
\xymatrix@C=1.8pc@R=1pc{
  X_\alpha\ar[rr]^-{X_{(\alpha,j)}}\ar[dd]_-{X_{(\alpha,\ell)}}&&
    X_{\alpha\sqcup\{j\}}\ar[dd]^-{X_{(\alpha\sqcup\{j\},\ell)}}\\
  \\
  X_{\alpha\sqcup\{\ell\}}\ar[rr]_-{X_{(\alpha\sqcup\{\ell\},j)}}&&
    X_{\alpha\sqcup\{j,\ell\}}
}
\end{aligned}
\end{equation}
  commutes.
\item
  If $\cal{C}$ is an $\cat{Ab}$-category, the cube $X$ is called
  \emph{anti-commutative} if for all $j,\ell\in\cal{I}$, $j\neq\ell$,
  and for all $\alpha\subseteq\cal{I}\backslash\{j,\ell\}$, the above
  diagram anti-commutes, \ie
\begin{equation}
  X_{(\alpha\sqcup\{\ell\},j)}\circ X_{(\alpha,\ell)} =
  - X_{(\alpha\sqcup\{j\},\ell)}\circ X_{(\alpha,j)}.
\end{equation}
\end{enumerate}
\end{defn}

\begin{defn}
Let $\cal{C}$ be a monoidal category.
\begin{enumerate}
\item
  Given a finite set $\cal{I}$ and $\cal{I}$-cubes
  $X=({\{X_\alpha\}}_\alpha,{\{X_{(\alpha,j)}\}}_{(\alpha,j)})$ and
  $Y=({\{Y_\alpha\}}_\alpha,{\{Y_{(\alpha,j)}\}}_{(\alpha,j)})$, the
  \emph{internal tensor product} $X\otimes Y$ is the $\cal{I}$-cube
  with the vertices ${(X\otimes Y)}_\alpha:=X_\alpha\otimes Y_\alpha$
  and with the edges ${(X\otimes
  Y)}_{(\alpha,j)}:=X_{(\alpha,j)}\otimes Y_{(\alpha,j)}$.
\item
  Given finite sets $\cal{I}$ and $\cal{J}$, an $\cal{I}$-cube
  $X=({\{X_\alpha\}}_{\alpha\subseteq\cal{I}},{\{X_{(\alpha,j)}\}}_{(\alpha,j)\in
  r(\cal{I})})$ and a $\cal{J}$-cube
  $Y=({\{Y_\beta\}}_{\beta\subseteq\cal{J}},{\{Y_{(\beta,\ell)}\}}_{(\beta,\ell)\in
  r(\cal{J})})$, the \emph{external tensor product} $X\boxtimes Y$ is
  the $\cal{I}\sqcup\cal{J}$-cube with the vertices ${(X\boxtimes
  Y)}_{\alpha\sqcup\beta}:=X_\alpha\otimes Y_\beta$ for
  $\alpha\subseteq\cal{I}$, $\beta\subseteq\cal{J}$, and whose edges
  are defined as follows. If $j\in\cal{I}\sqcup\cal{J}$, then either
  $j\in\cal{I}$ or $j\in\cal{J}$. In the former case, for any
  $\alpha\subseteq\cal{I}\backslash\{j\}$, $\beta\subseteq\cal{J}$, we
  define ${(X\boxtimes
  Y)}_{(\alpha\sqcup\beta,j)}:=X_{(\alpha,j)}\otimes\id_{Y_\beta}$. In
  the latter case, for any $\alpha\subseteq\cal{I}$ and
  $\beta\subseteq\cal{J}\backslash\{j\}$, we define ${(X\boxtimes
  Y)}_{(\alpha\sqcup\beta,j)}:=\id_{X_\alpha}\otimes
  Y_{(\beta,j)}$.
\end{enumerate}
\end{defn}

Both the internal and the external tensor product of two commutative
cubes form commutative cubes, too.

\begin{defn}
Let $\cal{C}$ be a monoidal $\cat{Ab}$-category and $\cal{I}$ be a
finite set equipped with a linear order `$\leq$'. We define the
anti-commutative $\cal{I}$-cube
$E_\cal{I}=({\{E^\cal{I}_\alpha\}}_\alpha,{\{E^\cal{I}_{(\alpha,j)}\}}_{(\alpha,j)})$,
as follows. The vertices are given by the unit object of the monoidal
structure, $E^\cal{I}_\alpha:=\1$ for all
$\alpha\subseteq\cal{I}$. The edges are defined by
$E^\cal{I}_{(\alpha,j)}:={(-1)}^{|\{\,k\in\alpha\mid k<j\,\}|}\id_\1$.
\end{defn}

If $\cal{C}$ is a monoidal $\cat{Ab}$-category, $\cal{I}$ a finite
set, and $X$ a commutative $\cal{I}$-cube, then the internal tensor
product $X\otimes E_\cal{I}$ forms an anti-commutative $\cal{I}$-cube.

\begin{defn}
Let $\cal{C}$ be a category, $\cal{I}$ and $\cal{J}$ be finite sets,
$X$ be a commutative $\cal{I}$-cube and $Y$ be a commutative
$\cal{J}$-cube. A \emph{homomorphism} $f\colon X\to Y$ \emph{of
commutative cubes} is a pair
$f=(f^0,{\{f_\alpha\}}_{\alpha\subseteq\cal{I}})$ consisting of a
bijection $f^0\colon\cal{I}\to\cal{J}$ and a family of morphisms
$f_\alpha\colon X_\alpha\to Y_{f^0(\alpha)}$ such that for all
$j\in\cal{I}$ and $\alpha\subseteq\cal{I}\backslash\{j\}$, the diagram
\begin{equation}
\begin{aligned}
\xymatrix@C=1.8pc@R=1pc{
  X_\alpha\ar[rr]^{f_\alpha}\ar[dd]_{X_{(\alpha,j)}}&&
    Y_{f^0(\alpha)}\ar[dd]^{Y_{(f^0(\alpha),f^0(j))}}\\
  \\
  X_{\alpha\sqcup\{j\}}\ar[rr]_{f_{\alpha\sqcup\{j\}}}&&
    Y_{f^0(\alpha\sqcup\{j\})}
}
\end{aligned}
\end{equation}
commutes.
\end{defn}

The category $\cat{Cube}(\cal{C})$ of \emph{commutative cubes in}
$\cal{C}$ has as objects pairs $(\cal{I},X)$ consisting of a finite
set $\cal{I}$ equipped with a linear order and a commutative
$\cal{I}$-cube $X$. The morphisms $f\colon(\cal{I},X)\to(\cal{J},Y)$
are homomorphisms $f=(f^0,{\{f_\alpha\}}_\alpha)$ of commutative cubes
such that for all $j_1,j_2\in\cal{I}$, we have $f^0(j_1)<f^0(j_2)$ if
and only if $j_1<j_2$. In particular, if $\cal{I}$ and $\cal{J}$ are
not in bijection, then there is no morphism
$(\cal{I},X)\to(\cal{J},Y)$.

If $\cal{C}$ is a [symmetric] monoidal category, then so is
$\cat{Cube}(\cal{C})$. The tensor product is given by the external
tensor product of cubes,
$(\cal{I},X)\otimes(\cal{J},Y):=(\cal{I}\sqcup\cal{J},X\boxtimes
Y)$. Here $\cal{I}\sqcup\cal{J}$ is equipped with the linear order
that restricts to the linear order of $\cal{I}$ and $\cal{J}$,
respectively, and for which $j<\ell$ for all $j\in\cal{I}$ and
$\ell\in\cal{J}$. The unit object $\1:=(\emptyset,\{\1_\emptyset\})$
is given by the empty set and the commutative $\emptyset$-cube whose
only vertex is $\1_\emptyset:=\1$. The associativity and unit
constraints [and the symmetric braiding] are inherited from
$\cal{C}$.

\subsection{Complexes}

\begin{defn}
\label{def_totalcomplex}
Let $\cal{C}$ be an abelian monoidal category.
\begin{enumerate}
\item
  Let $\cal{I}$ be a finite set and
  $X=({\{X_\alpha\}}_\alpha,{\{X_{(\alpha,j)}\}}_{(\alpha,j)})$ be an
  anti-commutative $\cal{I}$-cube. Then we define a bounded complex
  $C(X)=(C^i(X),d^i)$ as follows. The $i$-th term of this complex is
  given by
\begin{equation}
  C^i(X):=\bigoplus_{\alpha\subseteq\cal{I},|\alpha|=i}X_\alpha,
\end{equation}
  for $0\leq i\leq|\cal{I}|$, and the differential $d^i\colon
  C^i(X)\to C^{i+1}(X)$ is given on the $X_\alpha$,
  $\alpha\subseteq\cal{I}$, $|\alpha|=i$, by
\begin{equation}
  d^i|_{X_\alpha}:=\sum_{j\in\cal{I}\backslash\alpha}X_{(\alpha,j)}.
\end{equation}
\item
  For a pair $(\cal{I},X)$ of a finite set $\cal{I}$, equipped with a
  linear order, and a commutative $\cal{I}$-cube
  $X=({\{X_\alpha\}}_\alpha,{\{X_{(\alpha,j)}\}}_{(\alpha,j)})$, we
  define the \emph{total complex} as
\begin{equation}
  C_\mathrm{tot}(\cal{I},X) := C(X\otimes E_\cal{I}).
\end{equation}
\end{enumerate}
\end{defn}

\begin{prop}[Lemma 3 of~\cite{Kh3}]
\label{prop_total} Let $\cal{C}$ be an abelian monoidal category and
$\cal{I}$, $\cal{J}$ be finite sets each equipped with a linear
order. Let $X$ be a commutative $\cal{I}$-cube and $Y$ be a
commutative $\cal{J}$-cube. Then the following two complexes are
isomorphic:
\begin{equation}
\label{eq_total}
  C_\mathrm{tot}(\cal{I}\sqcup\cal{J},X\boxtimes Y) \cong C_\mathrm{tot}(\cal{I},X)\otimes
  C_\mathrm{tot}(\cal{J},Y).
\end{equation}
Here we choose the linear order on $\cal{I}\sqcup\cal{J}$ so that it
restricts to the linear orders of both $\cal{I}$ and $\cal{J}$, and so
that $j<\ell$ for all $j\in\cal{I}$ and $\ell\in\cal{J}$.
\end{prop}

\begin{prop}
\label{prop_total2}
Let $\cal{C}$ be an abelian monoidal category, $\cal{I}$ and $\cal{J}$
be finite sets each equipped with a linear order, $X$ be a commutative
$\cal{I}$-cube and $Y$ be a commutative $\cal{J}$-cube. For every
homomorphism of commutative cubes
\begin{equation}
  f=(f^0,{\{f_\alpha\}}_\alpha)\colon(\cal{I},X)\to(\cal{J},Y),
\end{equation}
there is a morphism of total complexes $C_\mathrm{tot}(f) \colon
C_\mathrm{tot}(\cal{I},X) \to C_\mathrm{tot}(\cal{J},Y)$ whose
components $C_\mathrm{tot}^i(f)\colon C_\mathrm{tot}^i(\cal{I},X)\to
C_\mathrm{tot}^i(\cal{J},Y)$ are given on the $X_\alpha\otimes\1$,
$\alpha\subseteq\cal{I}$, $|\alpha|=i$, by
\begin{equation}
  C_\mathrm{tot}^i(f)|_{X_\alpha\otimes\1}:=f_\alpha\otimes\id_\1.
\end{equation}
\end{prop}

If $\cal{C}$ is an abelian category, we denote by $\cat{Kom}(\cal{C})$
the category whose objects are complexes in $\cal{C}$ and whose
morphisms are morphisms of complexes. If $\cal{C}$ is [symmetric]
monoidal, then so is $\cat{Kom}(\cal{C})$.  $\cat{Kom}(\cal{C})$ is
actually a strict $2$-category whose $2$-morphisms are homotopies. We
denote by $\cal{K}(\cal{C})$ the
\emph{decategorification} of the $2$-category $\cat{Kom}(\cal{C})$,
\ie\ the homotopy category whose objects are the objects of
$\cat{Kom}(\cal{C})$ and whose morphisms are homotopy classes of
morphisms of $\cat{Kom}(\cal{C})$. As usual, we call complexes
\emph{isomorphic} if they are isomorphic in $\cat{Kom}(\cal{C})$ and
\emph{homotopy equivalent} if they are equivalent in
$\cat{Kom}(\cal{C})$, \ie\ isomorphic in $\cal{K}(\cal{C})$.

The point of Proposition~\ref{prop_total} and~\ref{prop_total2} is
that for an abelian [symmetric] monoidal category $\cal{C}$, there is
a functor
$C_\mathrm{tot}\colon\cat{Cube}(\cal{C})\to\cat{Kom}(\cal{C})$ sending
each commutative cube $(\cal{I},X)$ in $\cal{C}$ to its total complex
$C_\mathrm{tot}(\cal{I},X)$ and each morphism of commutative cubes to
the morphism of Proposition~\ref{prop_total2}. This functor has a
strong monoidal structure
$(C_\mathrm{tot},{C_\mathrm{tot}}_0,{C_\mathrm{tot}}_{(\cal{I},X),(\cal{J},Y)})$. The
isomorphism ${C_\mathrm{tot}}_0\colon\1_{\cat{Kom}(\cal{C})}\to
C_\mathrm{tot}(\1_{\cat{Cube}(\cal{C})})$ is the identity morphism on
the $1$-term complex $\1$, the monoidal unit of $\cal{C}$. The natural
isomorphisms ${C_\mathrm{tot}}_{(\cal{I},X),(\cal{J},Y)}\colon
C_\mathrm{tot}(\cal{I},X)\otimes C_\mathrm{tot}(\cal{J},Y)\to
C_\mathrm{tot}((\cal{I},X)\otimes (\cal{J},Y))$ are the isomorphisms
of~\eqref{eq_total}.

A \emph{graded} [\emph{filtered}] \emph{complex} is an object of
$\cat{Kom}(\cat{grdVect}_k)$ [of $\cat{Kom}(\cat{fltVect}_k)$]. If
$C=(C^i,d^i)$ is a graded [filtered] complex, then we denote by
$C\{j\}$ the complex whose terms are $C^i\{j\}$ (\emph{grading
shift}) and which has the same differentials. By $C[\ell]$ we denote
the complex whose terms are $(C[\ell])^i:=C^{i-\ell}$ and whose
differentials are $d^i_{C[\ell]}:={(-1)}^\ell d^i$. We call the
operator $[\ell]$ the \emph{cohomological degree shift}.

\subsection{Spectral sequences}
\label{sect_spectral}

It turns out that under suitable conditions, one can compute
Khovanov's chain complex and homology from the chain complexes of
Lee or of Bar-Natan by exploiting their filtration. This
construction uses a spectral sequence. In this section, we summarize
the key definitions and fix our notation. We adopt the conventions
of~\cite{McCleary} for spectral sequences except for the indexing of
the associated graded complex of a filtered complex which we
describe in detail below.

Let $C=(C^i,d^i)$ be a filtered complex, \ie\ each $C^i$ has a
bounded filtration ${\{F^kC^i\}}_{k\in\Z}$, and the $d^i\colon
C^i\to C^{i+1}$ are filtered linear maps of degree $0$. The
\emph{associated graded complex} is given by $E_0=(E_0^i,d_0^i)$
with $E_0^i=\bigoplus_{k\in\Z}E_0^{k,i}$ where
$E_0^{k,i}:=F^kC^{k+i}/F^{k+1}C^{k+i}$ and
$d^i_0|_{E^{k,i}_0}:=d^{k+i}|_{E^{k,i}_0}$. Every $v\in E_0^{k,i}$
has a representative of degree $\deg(v)=k$ in the filtered vector
space $C^{k+i}$. Here we have shifted the cohomological degree of
$E_0^{k,i}$ compared with~\cite{McCleary}. Since
$d^i(F^kC^{i})\subseteq F^kC^{i+1}$ for all $i,k$, the differential
descends to the quotient and gives rise to graded linear maps
$d_0^i\colon E_0^i\to E_0^{i+1}$ of degree $0$ which satisfy
$d_0^{i+1}\circ d_0^i=0$.

The bounded filtration on $C$ gives rise to a spectral sequence
${\{(E_r,d_r)\}}_{r\in\N_0}$ of cohomological type whose $E_0$-page,
$E_0\cong\bigoplus_{k,i\in\Z}E_0^{k,i}$ coincides with the
associated graded complex, indexed as above, and which converges to
the associated graded space of the homology of $C$:
\begin{equation}
  E_\infty^{k,i}\cong F^kH^{k+i}(C)/F^{k+1}H^{k+i}(C).
\end{equation}
In particular, each differential $d_r\colon E_r\to E_r$ is of bidegree
$(r,1-r)$, \ie\
\begin{equation}
  d_r(E_r^{k,i})\subseteq E_r^{k+r,i+1-r},
\end{equation}
and each page is determined by the homology of the previous page as
\begin{equation}
\label{eq_specseq}
  E_{r+1}^{k,i}\cong H^{k,i}(E_r,d_r) := \ker d_r^{k,i}/\im d_r^{k-r,i+r-1},
\end{equation}
where we have written $d_r^{k,i}:=d_r|_{E_r^{k,i}}$. The
differentials $d_r$ are obtained from the differential of the
original filtered complex $C$. Note that due to our shifting of the
cohomological degree compared with~\cite{McCleary}, the step from
$E_0$ to $E_1$ already follows the general
pattern~\eqref{eq_specseq}, and that, in our notation, the homology
of the associated graded complex is encoded in the $E_1$-page as
follows:
\begin{equation}
\label{eq_gradedhomology}
  H^i(E_0) \cong \bigoplus_{k\in\Z}E_1^{k,i-k}.
\end{equation}

%
\section{Tangle homology}
%

In this section, we take Bar-Natan's `picture world'
construction~\cite{BN2}, \ie\ chain complexes in a category whose
morphisms are formal linear combinations of certain surfaces, and
employ an open-closed TQFT in order to translate this into algebra and
to turn it into a chain complex of vector spaces. This already yields
a generalization of Khovanov homology from links to tangles. Whereas
Bar-Natan's surfaces can be glued so as to represent the composition
of tangles, the resulting chain complexes of vector spaces are not
necessarily equipped with an operation for the composition of
tangles. We consider this question in Section~\ref{sec_composition}
below.

\subsection{Examples of knowledgeable Frobenius algebras}

We first present examples of knowledgeable Frobenius algebras and
thereby of open-closed TQFTs which extend the commutative Frobenius
algebras of Khovanov, Lee, and Bar-Natan.

\subsubsection{Examples with trivial $C$}
\label{sect_kfrobtriv}

We start with some knowledgeable Frobenius algebras
$(A,C,\imath,\imath^\ast)$ whose $C$ is the underlying field equipped
with the trivial Frobenius algebra structure. Given a knowledgeable
Frobenius algebra $(A',C',\imath',{\imath'}^\ast)$ we can tensor with
one of these $(A,C,\imath,\imath^\ast)$ to form a knowledgeable
Frobenius algebra with isomorphic $C'\cong C'
\otimes k$, but with a more interesting symmetric Frobenius algebra
$A\otimes A'$.

\begin{example}
\label{ex_matrixalgebra} Let $k$ be a field, $m\in\N$, and
$(A,\mu_A,\eta_A)$ be the $m\times m$-matrix algebra $A:=M_m(k)$. We
write ${\{e_{pq}\}}_{1\leq p,q\leq m}$ for the standard basis for
which the multiplication reads $\mu_A(e_{pq}\otimes
e_{rs})=\delta_{qr}e_{ps}$ and the unit
$\eta_A(1)=\sum_{p=1}^me_{pp}$. The centre is $Z(A)\cong k$. Then
$(A,\mu_A,\eta_A,\Delta_A,\epsilon_A)$ forms a symmetric Frobenius
algebra with $\Delta_A(e_{pq})=\alpha\sum_{r=1}^me_{pr}\otimes
e_{rq}$ and $\epsilon_A(e_{pq})=\delta_{pq}\alpha^{-1}$ for any
$\alpha\in k\backslash\{0\}$. The window element reads
$a=\mu_A\circ\Delta_A\circ\eta_A=\alpha m\cdot\eta_A$, and so $A$ is
strongly separable if and only if the characteristic of $k$ does not
divide $m$.

There is a state sum knowledgeable Frobenius algebra
$(A,C,\imath,\imath^\ast)$ with $C=k$, $\mu_C(1\otimes 1)=1$,
$\eta_C(1)=1$, $\Delta_C(1)=\alpha^21\otimes 1$,
$\epsilon_C(1)=\alpha^{-2}$, $\imath(1)=\sum_{p=1}^me_{pp}$, and
$\imath^\ast(e_{pq})=\alpha\delta_{pq}$ constructed as in
Proposition~\ref{prop_statesum}. For $\alpha\in\{-1,1\}$, the
Frobenius algebra structure on $C=k$ is trivial.
\end{example}

\begin{example}
\label{ex_quaternion}
Let $k$ be a field. The algebra $A=\H_k$ of quaternions, \ie\ the
free associative unital $k$-algebra generated by $I,J,K$ subject to
the relations $I^2=-1$, $J^2=-1$, $K^2=-1$, $IJ=K$, $JI=-K$, $JK=I$,
$KJ=-I$, $KI=J$, and $IK=-J$, forms a symmetric Frobenius algebra
$(A,\mu_A,\eta_A,\Delta_A,\epsilon_A)$ with
$\Delta_A(1)=\alpha(1\otimes 1-I\otimes I-J\otimes J-K\otimes K)$,
$\Delta_A(I)=\alpha(1\otimes I+I\otimes 1+J\otimes K-K\otimes J)$,
$\Delta_A(J)=\alpha(1\otimes J+J\otimes 1+K\otimes I-I\otimes K)$,
$\Delta_A(K)=\alpha(1\otimes K+K\otimes 1+I\otimes J-J\otimes I)$,
$\epsilon_A(1)=\alpha^{-1}$, and
$\epsilon_A(I)=\epsilon_A(J)=\epsilon_A(K)=0$ for any $\alpha\in
k\backslash\{0\}$. Its window element is $a=4\alpha\cdot\eta_A$, and
so $A$ is strongly separable if and only if $\chr k\neq 2$.

In this case, there is a state sum knowledgeable Frobenius algebra
$(A,C,\imath,\imath^\ast)$ with $C\cong k$, $\mu_C(1\otimes 1)=1$,
$\eta_C(1)=1$, $\Delta_C(1)=4\alpha^2\cdot 1\otimes 1$,
$\epsilon_C(1)=1/4\alpha^2$, $\imath(1)=1$,
$\imath^\ast(1)=4\alpha$, and
$\imath^\ast(I)=\imath^\ast(J)=\imath^\ast(K)=0$ constructed as in
Proposition~\ref{prop_statesum}. For $\alpha\in\{-1/2,1/2\}$, the
Frobenius structure on $C=k$ is the trivial one.
\end{example}

Similar examples are available for all strongly separable algebras
$A$ that are central over $k$, \ie\ $Z(A)\cong k$. In these
examples, the state sum knowledgeable Frobenius algebra
$(A,C,\imath,\imath^\ast)$ of Proposition~\ref{prop_statesum} has
$C\cong k$, and one can choose the window element in such a way that
the Frobenius algebra structure on $k$ is the trivial one.

\begin{prop}
Let $k$ be a field. Every strongly separable and central $k$-algebra
is a Brauer algebra.
\end{prop}

\begin{proof}
Every strongly separable $k$-algebra $A$ is finite-dimensional over
$k$ and semi-simple. By Wedderburn's theorem, $A$ is a direct
product of $m_j\times m_j$-matrix algebras, $A\cong
\bigoplus_{j=1}^nM_{m_j}(L_j)$, $n\in\N$, whose coefficients are in
finite-dimensional skew field extensions $L_j/k$. Its centre is
$Z(A)\cong\bigoplus_{j=1}^nZ(L_j)$. If $A$ is central over $k$, then
$n=1$ and $Z(L_1)\cong k$. In particular, $A$ is finite-dimensional,
simple, and central and therefore a Brauer algebra.
\end{proof}

\noindent
Example~\ref{ex_quaternion} now fits into the following general
framework.

\begin{prop}
\label{prop_skewfield}
Let $k$ be a field and $L/k$ be a finite-dimensional skew field
extension such that $L$ is a strongly separable and central
$k$-algebra. We denote by $g_\mathrm{can}\colon L\otimes L\to k$ the
canonical bilinear form of $L$ and choose an orthogonal $k$-basis
${(e_j)}_j$ for $L$, \ie
\begin{equation}
  g_\mathrm{can}(e_j\otimes e_\ell) =\delta_{j\ell}\beta_j
\end{equation}
for some $\beta_j\in k\backslash\{0\}$. In particular, we can choose
$e_1:=\eta_L(1)$ to be the unit of $L$, and so $\beta_1=\dim_k L$ and
because of strong separability the characteristic of $k$ does not
divide $\dim_k L$. The symmetric Frobenius algebra structures
$(L,\mu_L,\eta_L,\Delta_L,\epsilon_L)$ are parameterized by their
window elements which are just scalars $\zeta\in k\backslash\{0\}$
because of centrality. The Frobenius algebra structures then read:
\begin{eqnarray}
  \Delta_L(e_j) &=& \zeta\sum_{\ell=1}^{\dim_k
  L}\beta_\ell^{-1}\mu_L(e_j\otimes e_\ell)\otimes e_\ell,\\
  \epsilon_L(e_j) &=& \zeta^{-1}\beta_j.
\end{eqnarray}
The state sum knowledgeable Frobenius algebra of
Proposition~\ref{prop_statesum} is then given by
$(L,k,\imath,\imath^\ast)$ with the following Frobenius algebra
structure for $k$:
\begin{eqnarray}
  \mu_k(1\otimes 1)&=&1,\\
  \eta_k(1) &=& 1,\\
  \Delta_k(1) &=& \zeta^2/\dim_k L\cdot 1\otimes 1,\\
  \epsilon_k(1) &=& \dim_k L/\zeta^2
\end{eqnarray}
and with $\imath(1)=e_1$, $\imath^\ast(e_1)=\zeta$, and
$\imath^\ast(e_j)=0$ for all $j\neq 1$. The Frobenius algebra
structure of $Z(L)\cong k$ is therefore the trivial one if and only if
$\zeta^2=\dim_k L$.
\end{prop}

\subsubsection{Euler-graded examples}

Khovanov's Frobenius algebra $C_{0,0}$
(Definition~\ref{def_khovanov}) is Euler-graded, and the degrees of
its basis vectors $\{1,x\}$ are determined by this requirement
(Example~\ref{ex_grading}). Let us try to extend Khovanov's
Frobenius algebra to an Euler-graded knowledgeable Frobenius
algebra.

\begin{prop}
\label{prop_gradedlimitation}
If $(A,C_{0,0},\imath,\imath^\ast)$ is an Euler-graded knowledgeable
Frobenius algebra such that $C_{0,0}$ is Khovanov's Frobenius algebra
(Definition~\ref{def_khovanov}), then $A$ is not strongly separable,
the underlying field $k$ has finite characteristic, and $\chr k$
divides $\dim_k A$. The value of the disc is,
\begin{equation}
\label{eq_disc}
  \epsilon_A\circ\eta_A = 0,
\end{equation}
and the value of the annulus,
\begin{equation}
\label{eq_annulus}
  \epsilon_A\circ\mu_A\circ\Delta_A\circ\epsilon_A = 0.
\end{equation}
\end{prop}

\begin{proof}
Since $\mu_A(\eta_A(1)\otimes\eta_A(1))=\eta_A(1)$, we have
$\deg_A(\eta_A(1))=1$ and therefore
$\deg_C(\imath^\ast(\eta_A(1)))=0$. But the basis vectors $\{1,x\}$ of
$C_{0,0}$ have degrees $\deg_C(1)=2$ and $\deg_C(x)=-2$, and so
$\imath^\ast\circ\eta_A=0$. The window element is
\begin{equation}
  a = \mu_A\circ\Delta_A\circ\eta_A = \mu_A\circ\tau_{A,A}\circ\Delta_A\circ\eta_A
    = \imath\circ\imath^\ast\circ\eta_A = 0,
\end{equation}
and so $A$ is not strongly separable. Finally, the dimension of $A$ as
an element of $\End(k)$ is the value of the annulus,
\begin{equation}
  \dim_k A = \epsilon_A\circ a = 0,
\end{equation}
and so $\chr k$ divides $\dim_k A$. Recall that $A$ is a rigid
object in $\cat{Vect}_k$, and so $A$ is finite-dimensional over $k$.

Since $\deg_A(\eta_A(1))=1$, we have
$\deg_k(\epsilon_A(\eta_A(1)))=2$, but all non-zero vectors of $k$ are
in degree $0$, and so $\epsilon_A\circ\eta_A=0$.
\end{proof}

\begin{example}
\label{ex_khovanov}
Let $k$ be a field and $C_{0,0}$ be Khovanov's commutative Frobenius
algebra over $k$ (Definition~\ref{def_khovanov}). Consider the algebra
$A:=k[y]/(y^2)$ with the commutative Frobenius algebra structure
$(A,\mu_A,\eta_A,\Delta_A,\epsilon_A)$ with $\Delta_A(1)=1\otimes
y+y\otimes 1$, $\Delta_A(y)=y\otimes y$, $\epsilon_A(1)=0$, and
$\epsilon_A(y)=1$. The window element is $a=2y$ which is a zero
divisor, and so $A$ is not strongly separable. If $\chr k=2$, then
$(A,C_{0,0},\imath,\imath^\ast)$ forms a knowledgeable Frobenius
algebra with $\imath(1)=1$, $\imath(x)=0$, $\imath^\ast(1)=0$, and
$\imath^\ast(y)=x$. It is Euler-graded with $\deg_A(1)=1$,
$\deg_A(y)=-1$, $\deg_C(1)=2$, and $\deg_C(x)=-2$ where the subscripts
indicate which vector space we are referring to.
\end{example}

For a field $k$ of characteristic $\chr k\neq 2$, the above example
fails to satisfy the Cardy condition whose left hand side reads
$(\mu_A\circ\tau_{A,A}\circ\Delta_A)(1)=2y$ and
$(\mu_A\circ\tau_{A,A}\circ\Delta_A)(y)=0$ while its right hand side
is $(\imath\circ\imath^\ast)(1)=0$ and
$(\imath\circ\imath^\ast)(y)=0$.

The previous example is the first one of the following series.

\begin{example}
Let $k$ be a field and $C_{0,0}$ be Khovanov's commutative Frobenius
algebra over $k$ (Definition~\ref{def_khovanov}). The truncated
polynomial algebra $A=k[y]/(y^p)$, $p\geq 2$, forms a commutative and
therefore symmetric Frobenius algebra
$(A,\mu_A,\eta_A,\Delta_A,\epsilon_A)$ with
$\Delta_A(y^\ell)=\sum_{j=0}^{p-1-\ell}y^{j+\ell}\otimes y^{p-1-j}$
for all $\ell\in\{0,\ldots,p-1\}$, $\epsilon_A(y^{p-1})=1$, and
$\epsilon_A(y^\ell)=0$ for all $\ell\in\{0,\ldots,p-2\}$. The window
element is $a=py^{p-1}$ which is a zero divisor, and so $A$ is not
strongly separable.

If $\chr k=p$, then $(A,C_{0,0},\imath,\imath^\ast)$ forms a
knowledgeable Frobenius algebra with $\imath(1)=1$, $\imath(x)=0$,
$\imath^\ast(y^{p-1})=1$, and $\imath^\ast(y^\ell)=0$ for all
$\ell\in\{0,\ldots,p-2\}$.

If this knowledgeable Frobenius algebra is Euler-graded, we have
$\deg(y^{p-1})-1=-2$ because $\imath^\ast(y^{p-1})=x$. But
$-1=\deg(y^{p-1})=(p-1)\deg(y)-(p-2)$ implies that
$\deg(y)=(p-3)/(p-1)$. This is an integer only if
$p\in\{2,3\}$. Indeed, in these two cases, the knowledgeable Frobenius
algebra is Euler-graded. If $p=2$, it coincides with
Example~\ref{ex_khovanov}, and if $p=3$, we have $\deg_A(1)=1$,
$\deg_A(y)=0$, and $\deg_A(y^2)=-1$ as well as $\deg_C(1)=2$, and
$\deg_C(x)=-2$.
\end{example}

For a field $k$ of characteristic $\chr k\neq p$, the above example
fails to satisfy the Cardy condition whose left hand side reads
$(\mu_A\circ\tau_{A,A}\circ\Delta_A)(1)=py^{p-1}$ and
$(\mu_A\circ\tau_{A,A}\circ\Delta_A)(y^\ell)=0$ for all
$\ell\in\{0,\ldots,p-2\}$ whereas its right hand side is
$(\imath\circ\imath^\ast)(y^\ell)=0$ for all
$\ell\in\{0,\ldots,p-1\}$.

The truncated polynomial algebras for $p\geq 4$ do not yield
Euler-graded knowledgeable Frobenius algebras. The truncated
polynomial algebra for $p=3$ is, however, at the same time the first
example of the following series of examples all of which are
Euler-graded.

\begin{example}
\label{ex_modp} Let $k$ be a field and $p=2n+1$, $n\in\N$. Consider
the $p$-dimensional vector space $A$ with basis
$\{X_{-n},X_{-n+1},\ldots,X_n\}$. It forms a symmetric Frobenius
algebra $(A,\mu_A,\eta_A,\Delta_A,\epsilon_A)$ with
$\mu_A(X_1\otimes X_j)=X_j=\mu_A(X_j\otimes X_1)$, $\mu_A(X_j\otimes
X_{-j})=X_{-1}$ for all $-n\leq j\leq n$, $\eta_A(1)=X_1$,
$\Delta_A(1)=\sum_{\ell=-n}^nX_\ell\otimes X_{-\ell}$,
$\Delta_A(X_{-1})=X_{-1} \otimes X_{-1}$,
$\Delta_A(X_j)=X_{-1}\otimes X_j+X_j\otimes X_{-1}$ for all
$j\notin\{1,-1\}$, and $\epsilon_A(X_{-1})=1$. The operations
$\mu_A$, $\eta_A$, $\Delta_A$, and $\epsilon_A$ are $0$ on all other
basis vectors. The window element reads $a=pX_{-1}$. It is a zero
divisor, and so $A$ is not strongly separable.

If $\chr k=p$, then there is a knowledgeable Frobenius algebra
$(A,C_{0,0},\imath,\imath^\ast)$ with $\imath(1)=1$, $\imath(x)=0$,
$\imath^\ast(X_{-1})=x$, and $\imath^\ast(X_j)=0$ for all $j\neq
-1$. If we set $\deg_A( X_i)=i$, then it is Euler graded.
\end{example}

Again, for a field $k$ with $\chr k\neq p$, the example fails to
satisfy the Cardy condition whose left hand side gives
$(\mu_A\circ\tau_{A,A}\circ\Delta_A)(1)=pX_{-1}$ and
$(\mu_A\circ\tau_{A,A}\circ\Delta_A)(X_j)=0$ for all $j\neq 1$
whereas its right hand side is $(\imath\circ\imath^\ast)(X_j)=0$ for
all $j$.

\begin{rem}
For any of these examples $(A,C,\imath,\imath^\ast)$, we can take the
tensor product with a knowledgeable Frobenius algebra
$(A^\prime,C^\prime,\imath^\prime,{\imath^\prime}^\ast)$ of
Section~\ref{sect_kfrobtriv} that has $C^\prime=k$ with the trivial
Frobenius algebra structure. The tensor product will therefore have
the same algebra $C\otimes k\cong C$, but a bigger and possibly
non-abelian $A\otimes A^\prime$. When one takes the tensor product of
knowledgeable Frobenius algebras, the window elements are
multiplicative; if both $A$ and $A^\prime$ are strongly separable,
then so is $A\otimes A^\prime$; and the entire construction of
Proposition~\ref{prop_statesum} is compatible with the tensor
product. If the algebra $A^\prime$ is entirely in degree $0$, then the
tensor product is Euler-graded [Euler-filtered] as soon as
$(A,C,\imath,\imath^\ast)$ is.

In Section~\ref{sect_strongsep} below, we present examples that are
not just tensor products of a commutative algebra $A$ with some
non-commutative central algebra, but whose non-commutativity is
genuine. Also note that if one is interested in modules over a
commutative ring rather than in vector spaces over a field,
Example~\ref{ex_modp} works for any commutative ring $k$ of odd
characteristic, not necessarily prime.
\end{rem}

\subsubsection{Euler-filtered examples}

Both Lee's and Bar-Natan's Frobenius algebra $C_{0,1}$ and $C_{1,0}$
are Euler-filtered (Example~\ref{ex_grading}). Lee and
Rasmussen~\cite{Lee,Rasmussen} and Turner~\cite{Turner} have shown
that one can start from the chain complexes of Lee or Bar-Natan and
then recover the chain complex and homology of Khovanov from a
spectral sequence. In this construction, Khovanov's chain complex
appears as the associated graded chain complex of Lee's or Bar-Natan's
filtered complex. We are therefore interested in examples in which the
basis vectors $\{1,x\}$ of $C_{h,t}$ have matching degrees, \ie\
$\deg_C(1)=2$ and $\deg_C(x)=-2$.

\begin{prop}
\label{prop_filtered}
If $(A,C,\imath,\imath^\ast)$ is an Euler-filtered knowledgeable
Frobenius algebra such that $C=C_{h,t}$
(Definition~\ref{def_khovanov}) with the basis vectors in degree
$\deg_C(1)=2$ and $\deg_C(x)=-2$, then either $A$ is not strongly
separable or the window element is $a=\zeta\cdot\eta_A$ with $\zeta\in
k\backslash\{0\}$. In both cases, the underlying field $k$ has finite
characteristic and $\chr k$ divides $\dim_k A$.  Both the disc and the
annulus are zero, \ie~\eqref{eq_disc} and~\eqref{eq_annulus} hold.
\end{prop}

\begin{proof}
Since $\mu_A(\eta_A(1)\otimes\eta_A(1))=\eta_A(1)$, we have
$\deg_A(\eta_A(1))\leq 1$. Since $\imath(\eta_C(1))=\eta_A(1)$, we
have $1\leq\deg_A(\eta_A(1))$, and so $\deg_A(\eta_A(1))=1$.

This implies that
$0=\deg_A(\eta_A(1))-1\leq\deg_C(\imath^\ast(\eta_A(1)))$. Due to the
degrees of the basis vectors of $C$, we have either
$\imath^\ast(\eta_A(1))=0$ or
$\imath^\ast(\eta_A(1))=\zeta\cdot\eta_C(1)$ with $\zeta\in
k\backslash\{0\}$.

In the first case, the window element is $a=0$ and $A$ therefore not
strongly separable whereas in the second case, the window element is
\begin{equation}
  a=\imath\circ\imath^\ast\circ\eta_A=\zeta\cdot\imath\circ\eta_C=\zeta\cdot\eta_A.
\end{equation}

As $\deg_A(\eta_A(1))=1$, we have
$\deg_k(\epsilon_A(\eta_A(1)))=1+1\leq 0$ or
$\epsilon_A(\eta_A(1))=0$. Only the latter is possible, and so
$\epsilon_A\circ a=\zeta\cdot\epsilon_A\circ\eta_A=0$ and both the
disc and the annulus are zero. The remainder of the argument is as
in Proposition~\ref{prop_gradedlimitation}.
\end{proof}

\begin{rem}
Notice that because of $\dim_k A=0$ in $k$, Proposition~2.15
of~\cite{LP2} does not apply to the algebra $A$ of
Proposition~\ref{prop_filtered} in the strongly separable case. The
algebra $A$ is not \emph{special} in the technical sense.
\end{rem}

\begin{example}
Let $k$ be a field, $h,t\in k$, and $A_{h,t}=k[y]/(y^2-hy-t)$. Then
$(A_{h,t},\mu_A,\eta_A,\Delta_A,\epsilon_A)$ forms a commutative
Frobenius algebra with $\Delta_A(1)=1\otimes y+y\otimes 1-h\cdot
1\otimes 1$, $\Delta_A(y)=y\otimes y+t\cdot 1\otimes 1$,
$\epsilon_A(1)=0$, and $\epsilon_A(y)=1$.

Its window element is $a=(\mu_A\circ\Delta_A\circ\eta_A)(1)=2y-h$
which satisfies $a^2=(h^2+4t)\cdot 1$, and so $A_{h,t}$ is strongly
separable if and only if $h^2+4t\neq 0$.

In this case, Proposition~\ref{prop_statesum} yields a knowledgeable
Frobenius algebra $(A_{h,t},C,\imath,\imath^\ast)$ with
$C=Z(A_{h,t})=A_{h,t}$ equipped with the commutative Frobenius algebra
structure $(C,\mu_C,\eta_C,\Delta_C,\epsilon_C)$ with $\mu_C(1\otimes
1)=1$, $\mu_C(1\otimes y)=y$, $\mu_C(y\otimes 1)=y$, $\mu_C(y\otimes
y)=hy+t$, $\eta_C(1)=1$, $\Delta_C(1)=(h^2+2t)\cdot 1\otimes1 -h\cdot
(1\otimes y+y\otimes 1)+2\cdot y\otimes y$, $\Delta_C(y)=-ht\cdot
1\otimes 1+2t\cdot (1\otimes y+y\otimes 1)+h\cdot y\otimes y$,
$\epsilon_C(1)=2/(h^2+4t)$, and $\epsilon_C(y)=h/(h^2+4t)$, with
$\imath(1)=1$, $\imath(y)=y$, $\imath^\ast(1)=2y-h$, and
$\imath^\ast(y)=hy+2t$.
\end{example}

If $\chr k=2$ and $h=1$, this $C$ is isomorphic as a Frobenius algebra
to $C_{1,t}$ of Definition~\ref{def_khovanov}. We summarize this as
follows.

\begin{example}
\label{ex_barnatan} Let $k$ be a field of characteristic $2$ and
$t\in k$. The algebra $A_{1,t}=k[y]/(y^2-y-t)$ forms a strongly
separable symmetric Frobenius algebra
$(A_{1,t},\mu_A,\eta_A,\Delta_A,\epsilon_A)$ with
$\Delta_A(1)=1\otimes y+y\otimes 1-1\otimes 1$,
$\Delta_A(y)=y\otimes y+t\cdot 1\otimes 1$, $\epsilon_A(1)=0$, and
$\epsilon_A(y)=1$. There is a state sum knowledgeable Frobenius
algebra $(A_{1,t},C_{1,t},\imath,\imath^\ast)$ with $C_{1,t}$ as in
Definition~\ref{def_khovanov} and $\imath(1)=1$, $\imath(x)=y$,
$\imath^\ast(1)=1$, and $\imath^\ast(y)=x$. It is Euler-filtered
with $\deg_A(1)=1$, $\deg_A(y)=-1$, $\deg_C(1)=2$, and
$\deg_C(x)=-2$.
\end{example}

While the above example contains Bar-Natan's case $C_{1,0}$ for
$t=0$, the following example contains Lee's case $C_{0,1}$.

\begin{example}
\label{ex_lee}
Let $k$ be a field of characteristic $2$ and $t\in k$. The algebra
$A_{0,t}=k[y]/(y^2-t)$ forms a commutative and therefore symmetric
Frobenius algebra $(A_{0,t},\mu_A,\eta_A,\Delta_A,\epsilon_A)$ with
$\Delta_A(1)=1\otimes y+y\otimes 1$, $\Delta_A(y)=y\otimes y+t\cdot
1\otimes 1$, $\epsilon_A(1)=0$, and $\epsilon_A(y)=1$. Its window
element is zero, and so $A_{0,t}$ is not strongly separable.

There is a knowledgeable Frobenius algebra
$(A_{0,t},C_{0,t^2},\imath,\imath^\ast)$ with
$C_{0,t^2}=k[x]/(x^2-t^2)$ as in Definition~\ref{def_khovanov}, but
with $t^2$, and $\imath(1)=1$, $\imath(x)=t$, $\imath^\ast(1)=0$, and
$\imath^\ast(y)=t+x$. It is Euler-filtered with $\deg_A(1)=1$,
$\deg_A(y)=-1$, $\deg_C(1)=2$, and $\deg_C(x)=-2$.
\end{example}

For a field of characteristic $\chr k\neq 2$ in this example, the left
hand side of the Cardy condition gives
$(\mu_A\circ\tau_{A,A}\circ\Delta_A)(1)=2y$ and
$(\mu_A\circ\tau_{A,A}\circ\Delta_A)(y)=2t$ whereas the right hand
side reads $(\imath\circ\imath^\ast)(1)=0$ and
$(\imath\circ\imath^\ast)(y)=2t$.

\begin{rem}
\label{rem_spectral}
Consider the knowledgeable Frobenius algebra of
Example~\ref{ex_lee}. Consider its structure maps $\mu_A$, $\eta_A$,
$\Delta_A$, $\epsilon_A$, $\mu_C$, $\eta_C$, $\Delta_C$, $\epsilon_C$,
$\imath$, and $\imath^\ast$ which are filtered of degree $-1$, $1$,
$-1$, $1$, $-2$, $2$, $-2$, $2$, $-1$, and $-1$, respectively. Their
associated graded maps which are graded of the same degrees, are
precisely the structure maps of the knowledgeable Frobenius algebra of
Example~\ref{ex_khovanov}. Similarly, the structure maps of the
knowledgeable Frobenius algebra of Example~\ref{ex_barnatan} have as
their associated graded maps the structure maps of that of
Example~\ref{ex_khovanov} as well.
\end{rem}

\subsubsection{Strongly separable examples}
\label{sect_strongsep}

We have not yet studied examples without filtration in greater detail
except for the strongly separable case which is needed in state sum
open-closed TQFTs (Proposition~\ref{prop_statesum}).

\begin{rem}
Let $k$ be a field and $h,t\in k$. The algebra $C_{h,t}=k[x]/p(x)$ for
$p(x)=x^2-hx-t$ of Definition~\ref{def_khovanov} is a field if and
only if $p(x)$ is irreducible in $k[x]$, \ie\ if and only if there
exist no $\alpha,\beta\in k$ such that $\alpha+\beta=h$ and
$\alpha\beta=-t$. Otherwise, $p(x)=(x-\alpha)(x-\beta)$, and $p(x)$ is
therefore reducible.
\end{rem}

Note that for any field $k$ and $h,t\in k$, neither $C_{h,0}$
including Bar-Natan's example nor $C_{0,t^2}$ including Lee's example
are fields. Over $k=\F_2$, however, $C_{1,1}$ is a field. Let us first
consider the case in which $C_{h,t}$ is not strongly separable, \ie\
$h^2+4t=0$.

\begin{prop}
Let $k$ be a field and $h,t\in k$ such that $h^2+4t=0$, and let
$\alpha,\beta\in k$ such that $\alpha+\beta=h$ and $\alpha\beta=-t$,
\ie\ $C_{h,t}$ is not a field. Then there exists no state sum knowledgeable
Frobenius algebra $(A,C_{h,t},\imath,\imath^\ast)$.
\end{prop}

\begin{proof}
Assume there is such a state sum knowledgeable Frobenius algebra. Then
$A$ is strongly separable over $k$ and therefore finite-dimensional
and semi-simple. By Wedderburn's theorem,
$A\cong\bigoplus_{j=1}^nM_{m_j}(L_j)$ is a direct product of matrix
algebras with coefficients in finite-dimensional skew field extensions
$L_j/k$. Its centre is $C_{h,t}=Z(A)\cong\bigoplus_{j=1}^nZ(L_j)$. As
$C_{h,t}$ is $2$-dimensional over $k$, either $n=1$ and $\dim
Z(L_1)=2$ or $n=2$ and $Z(L_1)\cong k$, $Z(L_2)\cong k$.

Since by assumption, $C_{h,t}$ is not a field, the first case
$C_{h,t}\cong Z(L_1)$ is ruled out and therefore $A\cong
M_{m_1}(L_1)\oplus M_{m_2}(L_2)$ with $C_{h,t}=Z(A)\cong k\oplus
k$. But this algebra is strongly separable which contradicts the
assumption that $h^2+4t=0$.
\end{proof}

In particular, the above proposition rules out the existence of state
sum knowledgeable Frobenius algebras $(A,C,\imath,\imath^\ast)$ with
$C=C_{0,0}$ (Khovanov's example) over any field and with $C=C_{0,t^2}$
(including Lee's example) over fields of characteristic $2$. The
above proposition does not discuss the situation in which $C_{h,t}$
is a field which we have not pursued any further. The following lemma
illustrates that this happens only in very special circumstances.

\begin{lem}
Let $k$ be a field and $h,t\in k$ such that $h^2+4t=0$ and such that
$p(x)=x^2-hx-t$ is irreducible in $k[x]$ and therefore $C_{h,t}$ a
field. Then $\chr k=2$, $h=0$, $k$ is not perfect, and in particular
$k$ is not finite.
\end{lem}

\begin{proof}
Consider the algebraic derivative $p^\prime(x)=2x-h$. If $\chr k\neq
2$ or $h\neq 0$, then $p^\prime(x)\neq 0$ and therefore $p(x)$ is
separable. If $k$ is perfect, $p(x)$ is separable anyway. Since $p(x)$
is irreducible by assumption, the field extension $C_{h,t}/k$ is
separable and therefore the algebra $C_{h,t}$ separable over $k$. As
$C_{h,t}$ is commutative, it is even strongly separable, but this
contradicts the assumption that $h^2+4t=0$.
\end{proof}

For any field $k$ and $h,t\in k$ such that $h^2+4t=0$, the above lemma
implies that if $\chr k\neq 2$ or $h\neq 0$ or $k$ is perfect or $k$
is finite, then there exist $\alpha,\beta\in k$ such that
$\alpha+\beta=h$ and $\alpha\beta=-t$. 

Apart from the case which we have not yet studied in greater detail,
we expect state sum TQFTs only when we extend a strongly separable
commutative Frobenius algebra. We therefore consider the algebra
$C_{h,t}$ in the strongly separable case, and again we restrict
ourselves to the case in which it is not a field.

\begin{lem}
\label{lem_onidempotent}
Let $k$ be a field, $h,t\in k$, and $C_{h,t}=k[x]/(x^2-hx-t)$ be
strongly separable and not a field, \ie\ $h^2+4t\neq 0$ and there
exist $\alpha,\beta\in k$ such that $\alpha+\beta=h$ and
$\alpha\beta=-t$. Then $C_{h,t}\cong k\oplus k$ as a Frobenius
algebra, and the orthogonal idempotent basis $\{z_1,z_2\}$ of
$C_{h,t}$ associated with this decomposition is given by
\begin{equation}
  z_1 = \frac{x-\alpha}{\beta-\alpha},\qquad
  z_2 = \frac{x-\beta}{\alpha-\beta}.
\end{equation}
The Frobenius algebra structure of $C_{h,t}$
(Definition~\ref{def_khovanov}) reads in this basis as follows:
\begin{eqnarray}
  \epsilon(z_1) &=& -1/c,\\
  \epsilon(z_2) &=& 1/c,\\
  \Delta(z_1) &=& -c\cdot z_1\otimes z_1,\\
  \Delta(z_2) &=& c\cdot z_2\otimes z_2,
\end{eqnarray}
where $c=\alpha-\beta$. We have $1=z_1+z_2$ and $x=\beta z_1+\alpha
z_2$.
\end{lem}

\begin{proof}
We compute $c^2={(\alpha-\beta)}^2=h^2+4t\neq 0$, and so $c\neq
0$. Then $z_1^2=z_1$, $z_2^2=z_2$, $z_1z_2=0$, and the other claims
follow from straightforward computations.
\end{proof}

\begin{thm}
\label{thm_sskfrob}
Let $k$ be a field, $h,t\in k$, and let $C_{h,t}=k[x]/(x^2-hx-t)$ with
the Frobenius algebra structure of Definition~\ref{def_khovanov} be
strongly separable and not a field.
\begin{enumerate}
\item
  Every state sum knowledgeable Frobenius algebra
  $(A,C_{h,t},\imath,\imath^\ast)$ is of the form
\begin{equation}
\label{eq_sskfrob}
  A\cong (M_{m_1}(k)\otimes L_1)\oplus(M_{m_2}(k)\otimes L_2)
\end{equation}
  where $L_1/k$ and $L_2/k$ are finite-dimensional strongly separable
  and central skew field extensions such that the characteristic of
  $k$ does not divide $m_1$, $m_2$, $\dim_k L_1$, nor $\dim_k
  L_2$. The window element of $A$ is of the form
\begin{equation}
  a = \xi_1\cdot z_1+\xi_2\cdot z_2
\end{equation}
  for some $\xi_1,\xi_2\in k\backslash\{0\}$ such that
\begin{equation}
  \xi_1^2=(\alpha-\beta)d_1,\qquad\mbox{and}\qquad
  \xi_2^2=(\beta-\alpha)d_2,
\end{equation}
  where we write $z_1=(I_{m_1}\otimes 1,0)$ and $z_2=(0,I_{m_2}\otimes
  1)$ for the basis vectors of the centre $Z(A)$, and $d_j:=m_j^2\dim_k
  L_j=\dim_k(M_{m_j}(k)\otimes L_j)$.
\item
  The knowledgeable Frobenius algebra $(A,C_{h,t},\imath,\imath^\ast)$
  of part (1.) is Euler-filtered if and only if the window element of
  $A$ is of the form
\begin{equation}
  a = \xi\cdot(z_1+z_2)=\xi\cdot\eta_A(1)
\end{equation}
  for some $\xi\in k\backslash\{0\}$ and if the characteristic of $k$
  divides $d_1+d_2$ and if
\begin{equation}
  \xi^2=(\alpha-\beta)d_1.
\end{equation}
\end{enumerate}
\end{thm}

\begin{proof}\hfill
\begin{enumerate}
\item
  If $(A,C_{h,t},\imath,\imath^\ast)$ is a state sum knowledgeable
  Frobenius algebra, then $A$ is strongly separable and therefore
  finite-dimensional over $k$ and semi-simple. By Wedderburn's
  theorem, $A$ is of the form $A\cong\bigoplus_{j=1}^nM_{m_j}(L_j)$
  where $m_j\in\N$ and $L_j/k$ are finite-dimensional skew field
  extensions. Its centre is $Z(A)\cong\bigoplus_{j=1}^nZ(L_j)$. As
  $C_{h,t}$ is $2$-dimensional over $k$, we have either $n=2$ and
  $\dim_kZ(L_1)=2$ or $n=2$ and $Z(L_1)\cong k\cong Z(L_2)$. Since by
  assumption, $C_{h,t}$ is not a field, only the second case is
  possible, and therefore $A$ is of the form~\eqref{eq_sskfrob} with
  skew field extensions $L_1/k$ and $L_2/k$ that are central. Since
  $A$ is strongly separable, so are $M_{m_j}(k)\otimes L_j$ as well as
  $M_{m_j}(k)$ and $L_j$, and so the characteristic of $k$ does not
  divide $m_1$, $m_2$, $\dim_kL_1$, nor $\dim_jL_2$.

  Any symmetric Frobenius algebra structure on the strongly separable
  $A$ is characterized by its invertible window element $a=\xi_1\cdot
  z_1+\xi_2\cdot z_2$ where $\{z_1,z_2\}$ is the orthogonal idempotent
  basis of the centre $Z(A)$ as specified above. We can now compute
  the state sum knowledgeable Frobenius algebra
  $(A,C=Z(A),\imath,\imath^\ast)$ of Proposition~\ref{prop_statesum}
  for this choice of $A$. This computation makes use of
  Example~\ref{ex_matrixalgebra} (don't forget to rescale the window
  element appropriately) and of Proposition~\ref{prop_skewfield}. If
  we denote by ${(e_{pq})}_{p,q}$ the standard basis of the matrix
  algebra $M_m(k)$ and by ${(e_\ell)}_\ell$ the basis of the skew
  field $L$ as in Proposition~\ref{prop_skewfield}, we have the basis
  vectors $e_{pq}^{(1)}\otimes e_\ell^{(1)}:=(e_{pq}\otimes
  e_\ell,0)$, $1\leq p,q\leq m_1$, $1\leq \ell\leq\dim_k L_1$ and
  $e_{pq}^{(2)}\otimes e_\ell^{(2)}:=(0,e_{pq}\otimes e_\ell)$, $1\leq
  p,q\leq m_2$, $1\leq \ell\leq\dim_kL_2$ for $A$ as
  in~\eqref{eq_sskfrob}. The knowledgeable Frobenius algebra structure
  is then given by
\begin{eqnarray}
  \mu_A((e_{pq}^{(j)}\otimes e_\ell^{(j)})\otimes
        (e_{rs}^{(j^\prime)}\otimes e_i^{(j^\prime)}))
    &=& \delta_{jj^\prime}\delta_{qr}e_{ps}^{(j)}\otimes{(\mu_{L_j}(e_\ell\otimes e_i))}^{(j)},\\
  \eta_A(1)
    &=& \sum_{j=1}^2\sum_{p=1}^{m_j}e_{pp}^{(j)}\otimes e_1^{(j)},
\end{eqnarray}
\begin{eqnarray}
  \Delta_A(e_{pq}^{(j)}\otimes e_\ell^{(j)})
    &=& \frac{\xi_j}{m_j}\sum_{i=1}^{\dim_kL_j}\sum_{r=1}^{m_j}\beta_i^{-1}
        (e_{pr}^{(j)}\otimes{(\mu_{L_j}(e_\ell\otimes e_i))}^{(j)})\otimes
        (e_{rq}^{(j)}\otimes e_i^{(j)}),\\
  \epsilon_A(e_{pq}^{(j)}\otimes e_\ell^{(j)})
    &=& \frac{m_j}{\xi_j}\delta_{pq}\beta_\ell,
\end{eqnarray}
\begin{eqnarray}
  \mu_C(z_j\otimes z_{j^\prime})
    &=& \delta_{jj^\prime}z_j,\\
  \eta_C(1)
    &=& z_1+z_2,\\
  \Delta_C(z_j)
    &=& \frac{\xi_j^2}{m_j^2\dim_kL_j}z_j\otimes z_j,\\
  \epsilon_C(z_j)
    &=& \frac{m_j^2\dim_kL_j}{\xi_j^2},
\end{eqnarray}
\begin{eqnarray}
  \imath(z_j)
    &=& \sum_{p=1}^{m_j}e_{pp}^{(j)}\otimes e_1^{(j)},\\
  \imath^\ast(e_{pq}^{(j)}\otimes e_\ell^{(j)})
    &=& \frac{\xi_j}{m_j}\delta_{pq}\delta_{\ell1}.
\end{eqnarray}
  We have denoted by $\mu_{L_j}\colon L_j\otimes L_j\to L_j$ the
  multiplication of $L_j$. The Frobenius algebra structure of the
  centre $C=Z(A)$ is now supposed to agree with that of
  $C_{h,t}$. Comparing the coefficients in the orthogonal idempotent
  basis with Lemma~\ref{lem_onidempotent}, we obtain the claim.
\item
  If $(A,C_{h,t},\imath,\imath^\ast)$ is Euler-filtered,
  Proposition~\ref{prop_filtered} implies that the window element is
  of the form $a=\xi\cdot\eta_A(1)$ for some $\xi\in
  k\backslash\{0\}$, \ie\ $\xi_1=\xi_2=\xi$, because $A$ is strongly
  separable, and that the characteristic of $k$ divides $\dim_k
  A=d_1+d_2$. In this case, however, $\xi^2=cd_1$ already implies
  $\xi^2=-cd_2$.

  For the converse implication, we assume that the characteristic of
  $k$ divides $d_1+d_2$ and that $a=\xi\cdot\eta_A(1)$ with
  $\xi^2=cd_1$. An Euler-filtration for
  $(A,C_{h,t},\imath,\imath^\ast)$ is now given by the following
  filtrations for $A$:
\begin{equation}
  A=A^{-1}\supseteq A^0=\ker\epsilon_A\supseteq A^1=\spann_k\{\eta_A(1)\}
   \supseteq A^2=\{0\},
\end{equation}
  and for $C_{h,t}$:
\begin{equation}
  C^{-2}=C_{h,t}\supseteq C^{-1}=C^0=C^1=C^2=\spann_k\{\eta_C(1)\}
    \supseteq C^3=\{0\}.
\end{equation}
  Since $(\epsilon_A\circ\eta_A)(1)=0$ by
  Proposition~\ref{prop_filtered}, the above inclusions hold. We now
  verify that the structure maps of the knowledgeable Frobenius
  algebra are filtered of the required degrees~\eqref{eq_defdegree}.
\begin{enumerate}
\item
  $\eta_A\colon k\to A$ is of degree $+1$ because $\deg_k(1)=0$ and
  $\eta_A(1)\in A^1$.
\item
  $\mu_A\colon A\otimes A\to A$ is of degree $-1$ because
  $\mu_A(\eta_A(1)\otimes\eta_A(1))=\eta_A(1)\in A^2\subseteq A^0$,
  because $\mu_A(\eta_A(1)\otimes a)=a\in A^j$ for all $a\in A^j$, and
  because $\mu_A(a\otimes b)\in A=A^{-1}$ for all $a,b\in A^0$.
\item
  $\epsilon_A\colon A\to k$ is of degree $+1$ because
  $\epsilon_A(a)=0$ for all $a\in A^1$ or $a\in A^0$ and
  $\epsilon_A(a)\in k^0=k$ for all $a\in A^{-1}$.
\item
  In order to see that $\Delta_A\colon A\to A\otimes A$ is of degree
  $-1$, we consider three cases. First, by a direct computation, we
  show that
\begin{equation}
  \Delta_A(\eta_A(1)) \in A\otimes\spann_k\{\eta_A(1)\}
    + \spann_k\{\eta_A(1)\}\otimes A
    + \ker\epsilon_A\otimes\ker\epsilon_A,
\end{equation}
  and so $\deg_A(\Delta_A(\eta_A(1)))\geq 0$ as required. Second,
  $\epsilon_A\colon A\to k$ is a homomorphism of coalgebras and
  therefore $\ker\epsilon_A$ a coideal, \ie\
\begin{equation}
  \Delta_A(\ker\epsilon_A)\subseteq A\otimes\ker\epsilon_A+\ker\epsilon_A\otimes A,
\end{equation}
  and therefore $\deg_A(\Delta_A(\ker\epsilon_A))\geq -1$. Third, for
  any $a\in A$, $\Delta_A(a)\in A\otimes A$ and so
  $\deg_A(\Delta_A(a))\geq -2$ anyway.
\item
  $\imath\colon C_{h,t}\to A$ is of degree $-1$ because
  $\imath(1)=\eta_A(1)\in A^1$ and $\imath(x)\in A\subseteq A^{-3}$.
\item
  Finally, in order to see that $\imath^\ast\colon A\to C_{h,t}$ is of
  degree $-1$, we consider three cases. First, if $a\in A^{-1}$ then
  $\imath^\ast(a)\in C_{h,t}=C^{-2}$. Second, if $a\in
  A^0=\ker\epsilon_A$, the duality axiom implies that
\begin{equation}
  \epsilon_C(\imath^\ast(a))=\epsilon_C(\eta_C(1)\cdot\imath^\ast(a))
  =\epsilon_A(\imath(1)\cdot a)=\epsilon_A(a)=0
\end{equation}
  and therefore $\imath^\ast(a)\in\ker\epsilon_C$ and so
  $\imath^\ast(a)\in\spann_k\{\eta_C(1)\}=C^2$. This implies that
  $\deg_C(\imath^\ast(a))=2\geq -1$. Third, if $a\in A^1$, then $a\in
  A^0$ and $\deg_C(\imath^\ast(a))=2\geq 0$.
\end{enumerate}
\end{enumerate}
\end{proof}

\begin{rem}
The above theorem establishes conditions under which there exists an
Euler-filtration. In general, this filtration is not unique as we
demonstrate by Example~\ref{ex_filtered3} below.
\end{rem}

\begin{example}
\label{ex_filtered}
If $k$ is a field whose characteristic $p=\chr k$ is prime of the form
$p=4n+1$, $n\in\N$, then there exist $m_1,m_2\in\N$ such that
$m_1<m_2$ and $m_1^2+m_2^2=p$. The above theorem then yields an
example of an Euler-filtered state sum knowledgeable Frobenius algebra
$(A,C_{1,0},\imath,\imath^\ast)$, \ie\ $h=1$, $t=0$, $\alpha=1$,
$\beta=0$, and $c=1$, with $\xi=m_1$, $L_1=k$, and $L_2=k$. Its
associated Euler-graded algebra is a knowledgeable Frobenius algebra
whose $C$ is Khovanov's $C_{0,0}$.
\end{example}

\begin{example}
\label{ex_filtered1}
Let $k$ be a field of characteristic $5$, and let $A=M_2(k)\oplus k$
be equipped with the symmetric Frobenius algebra structure that is
characterized by the window element $a=\eta_A(1)$. This is a special
case of Theorem~\ref{thm_sskfrob}(2.), and the state sum knowledgeable
Frobenius algebra associated with $A$ is given by
$(A,C,\imath,\imath^\ast)$ with $C=C_{h,t}$ for $h=1$, $t=0$,
$\alpha=1$, $\beta=0$, and $c=1$. The Euler-filtration becomes
manifest in the basis $(1,A,B,C,D)$ for $A$ where
\begin{equation}
  1=\left(\begin{pmatrix}1&0\\0&1\end{pmatrix},1\right),\quad
  A=\left(\begin{pmatrix}0&1\\0&0\end{pmatrix},0\right),\quad
  B=\left(\begin{pmatrix}0&0\\1&0\end{pmatrix},0\right),
\end{equation}
\begin{equation}
  C=\left(\begin{pmatrix}1&0\\0&-1\end{pmatrix},0\right),\quad
  D=\left(\begin{pmatrix}0&0\\0&0\end{pmatrix},1\right).\nn
\end{equation}
with $\deg_A(1)=1$, $\deg_A(A)=\deg_A(B)=\deg_A(C)=0$, and
$\deg_A(D)=-1$. We have the basis $\{z_1,z_2\}$ with $z_1=1-D$,
$z_2=D$ for $C=Z(A)$, and the knowledgeable Frobenius algebra
structure given by $\eta_A(1)=1$, $\mu_A(1\otimes 1)=1$,
$\mu_A(1\otimes X)=X$ for $X\in\{A,B,C,D\}$, $\mu_A(C\otimes C)=1-D$,
$\mu_A(D\otimes D)=D$,
$\mu_A(A\otimes B)=(1+C-D)/2$, $\mu_A(B\otimes A)=(1-C-D)/2$,
$\mu_A(A\otimes C)=-A$, $\mu_A(C\otimes A)=A$, $\mu_A(B\otimes C)=B$,
$\mu_A(C\otimes B)=-B$, $\epsilon_A(D)=1$, $\Delta_A(1)=-1\otimes
1-C\otimes C-(1\otimes D+D\otimes 1)+(A\otimes B+B\otimes A)$,
$\Delta_A(A)=-(1\otimes A+A\otimes 1)+(D\otimes A+A\otimes D)-C\otimes
A+A\otimes C$, $\Delta_A(B)=-(1\otimes B+B\otimes 1)+(B\otimes
D+D\otimes B)-B\otimes C+C\otimes B$, $\Delta_A(C)=-(1\otimes
C+C\otimes 1)+(C\otimes D+D\otimes C)+A\otimes B-B\otimes A$,
$\Delta_A(D)=D\otimes D$, $\imath_A(1)=1$, $\imath(x)=D$,
$\imath^\ast(1)=1$, and $\imath^\ast(D)=x$. The structure maps are
zero on all basis vectors not listed here.
\end{example}

\begin{example}
\label{ex_filtered2}
Let $k$ be a field of characteristic $5$, and let $A=\H_k\oplus k$
where $\H_k$ is the quaternion algebra (Example~\ref{ex_quaternion})
over $k$. Equip $A$ with the symmetric Frobenius algebra structure
characterized by the window element $a=\eta_A(1)$. This is a special
case of Theorem~\ref{thm_sskfrob}(2.), and the state sum knowledgeable
Frobenius algebra $(A,C,\imath,\imath^\ast)$ associated with $A$ has
$C=C_{h,t}$ with $h=1$, $t=0$, $\alpha=1$, $\beta=0$, and $c=1$. Its
Euler-filtration becomes obvious in the basis $\{1=(1,1), I=(I,0),
J=(J,0), K=(K,0), L=(0,1)\}$ with $\deg_A(1)=1$,
$\deg_A(I)=\deg_A(J)=\deg_A(K)=0$, and $\deg_A(L)=-1$. The
knowledgeable Frobenius algebra structure is given by $\eta_A(1)=1$,
$\mu_A(1\otimes X)=X=\mu_A(X\otimes 1)$ for $X\in\{1,I,J,K,L\}$,
$\mu_A(I\otimes I)=\mu_A(J\otimes J)=\mu_A(K\otimes K)=L-1$,
$\mu_A(I\otimes J)=-\mu_A(J\otimes I)=K$, $\mu_A(J\otimes
K)=-\mu_A(K\otimes J)=I$, $\mu_A(K\otimes I)-\mu_A(I\otimes K)=J$,
$\mu_A(L\otimes L)=L$, $\epsilon_A(L)=1$, $\Delta_A(1)=-1\otimes
1+(L\otimes 1+1\otimes L)+I\otimes I+J\otimes J+K\otimes K$,
$\Delta_A(I)=-(1\otimes I+I\otimes 1)+(L\otimes I+I\otimes L)+J\otimes
K-K\otimes J$, $\Delta_A(J)=-(1\otimes J+J\otimes 1)+(L\otimes
J+J\otimes L)+K\otimes I-I\otimes K$, $\Delta_A(K)=-(1\otimes
K+K\otimes 1)+(L\otimes K+K\otimes L)+I\otimes J-J\otimes I$,
$\Delta_A(L)=L\otimes L$, $\imath(1)=1$, $\imath(x)=L$,
$\imath^\ast(1)=1$, $\imath^\ast(L)=x$. The structure maps are zero on
all basis vectors not listed here.
\end{example}

The following example finally shows that the state sum knowledgeable
Frobenius algebra of Example~\ref{ex_filtered1} can be equipped with
an Euler-filtration other than that of Theorem~\ref{thm_sskfrob}.

\begin{example}
\label{ex_filtered3}
The state sum knowledgeable Frobenius algebra of
Example~\ref{ex_filtered1} is Euler-filtered with $\deg_A(1)=1$,
$\deg_A(A)=-2$, $\deg_A(B)=2$, $\deg_A(C)=0$, and $\deg_A(D)=-1$. With
this filtration, the associated Euler-graded knowledgeable Frobenius
algebra is the one of Example~\ref{ex_modp} in characteristic
$p=5$. In order to see this, we set $X_{-2}:=3A$, $X_{-1}:=D$,
$X_0:=2C$, $X_1:=D$, and $X_2:=B$.
\end{example}

\subsection{Orienting the smoothings of tangle diagrams}

In this section, we define orientations for the smoothings of every
plane diagram of a tangle.

Let $T$ be a plane diagram of a tangle with $n$ crossings, numbered
$1,\ldots,n$. Let $D_\alpha\subseteq\R^2$ for subsets
$\alpha\subseteq\cat{n}$ denote the plane diagrams associated with the
smoothings of all crossings such that $j\in\alpha$ if the $j$-th
crossing is resolved by the $1$-smoothing and $j\notin\alpha$ if it is
resolved by the $0$-smoothing. Each $D_\alpha$ is the disjoint union
of a finite number of circles and arcs in $\R^2$.

For each tangle diagram $T$, for example for
\begin{equation}
  T =
\begin{aligned}
\begin{pspicture}(4,2)
  \psbezier[linewidth=1pt](2.0,2.0)(2.0,1.5)(1.0,1.5)(1.0,1.0)
  \psline[linewidth=1pt](1.0,1.0)(1.0,0.0)
  \psline[linewidth=3pt,linecolor=white](1.3,1.6)(1.7,1.4)
  \psbezier[linewidth=1pt](1.0,2.0)(1.0,1.5)(2.0,1.5)(2.0,1.0)
  \psbezier[linewidth=1pt](2.0,1.0)(2.0,0.5)(2.3,0.5)(2.5,0.5)
  \psbezier[linewidth=1pt](2.5,0.5)(2.7,0.5)(3.0,0.5)(3.0,1.0)
  \psline[linewidth=1pt](3.0,1.0)(3.0,2.0)
\end{pspicture}
\end{aligned}
\end{equation}
we define two checkerboard colourings. The positive colouring is
obtained by shading the region left of the first strand of the tangle
and by continuing in a checkerboard fashion. If $T$ is empty, all of
$\R^2$ is shaded. The negative colouring is obtained from the positive
one by replacing shaded with non-shared regions and vice versa. We
denote the positively and negatively checkerboard coloured tangle
diagrams by $T^{(\epsilon)}$ for $\epsilon\in\{+1,-1\}$,
respectively. For instance, for the above example,
\begin{equation}
  T^{(+1)} =
\begin{aligned}
\begin{pspicture}(4,2)
  \pscustom[fillstyle=solid,fillcolor=lightgray,linecolor=lightgray]{
    \psbezier(1.0,2.0)(1.0,1.5)(2.0,1.5)(2.0,1.0)
    \psbezier(2.0,1.0)(2.0,0.5)(2.3,0.5)(2.5,0.5)
    \psbezier(2.5,0.5)(2.7,0.5)(3.0,0.5)(3.0,1.0)
    \psline(3.0,1.0)(3.0,2.0)(2.0,2.0)
    \psbezier(2.0,2.0)(2.0,1.5)(1.0,1.5)(1.0,1.0)
    \psline(1.0,1.0)(1.0,0.0)(0.5,0.0)(0.5,2.0)(1.0,2.0)
  }
  \psbezier[linewidth=1pt](2.0,2.0)(2.0,1.5)(1.0,1.5)(1.0,1.0)
  \psline[linewidth=1pt](1.0,1.0)(1.0,0.0)
  \psline[linewidth=3pt,linecolor=white](1.3,1.6)(1.7,1.4)
  \psbezier[linewidth=1pt](1.0,2.0)(1.0,1.5)(2.0,1.5)(2.0,1.0)
  \psbezier[linewidth=1pt](2.0,1.0)(2.0,0.5)(2.3,0.5)(2.5,0.5)
  \psbezier[linewidth=1pt](2.5,0.5)(2.7,0.5)(3.0,0.5)(3.0,1.0)
  \psline[linewidth=1pt](3.0,1.0)(3.0,2.0)
\end{pspicture}
\end{aligned}
\qquad T^{(-1)} =
\begin{aligned}
\begin{pspicture}(4,2)
  \pscustom[fillstyle=solid,fillcolor=lightgray,linecolor=lightgray]{
    \psbezier(1.0,2.0)(1.0,1.5)(2.0,1.5)(2.0,1.0)
    \psbezier(2.0,1.0)(2.0,0.5)(2.2,0.5)(2.3,0.5)
    \psbezier(2.5,0.5)(2.7,0.5)(3.0,0.5)(3.0,1.0)
    \psline(3.0,1.0)(3.0,2.0)(3.5,2.0)(3.5,0.0)(1.0,0.0)(1.0,1.0)
    \psbezier(1.0,1.0)(1.0,1.5)(2.0,1.5)(2.0,2.0)
    \psline(2.0,2.0)(1.0,2.0)
  }
  \psbezier[linewidth=1pt](2.0,2.0)(2.0,1.5)(1.0,1.5)(1.0,1.0)
  \psline[linewidth=1pt](1.0,1.0)(1.0,0.0)
  \psline[linewidth=3pt,linecolor=white](1.3,1.6)(1.7,1.4)
  \psbezier[linewidth=1pt](1.0,2.0)(1.0,1.5)(2.0,1.5)(2.0,1.0)
  \psbezier[linewidth=1pt](2.0,1.0)(2.0,0.5)(2.3,0.5)(2.5,0.5)
  \psbezier[linewidth=1pt](2.5,0.5)(2.7,0.5)(3.0,0.5)(3.0,1.0)
  \psline[linewidth=1pt](3.0,1.0)(3.0,2.0)
\end{pspicture}
\end{aligned}
\end{equation}
The checkerboard colouring of the tangle diagram $T$ induces a
checkerboard colouring on all smoothings $D_\alpha$,
$\alpha\subseteq\cat{n}$, which we denote by
$D^{(\epsilon)}_\alpha$. For the above example $T^{(+1)}$ we have the
two smoothings,
\begin{equation}
  D_\emptyset^{(+1)} =
\begin{aligned}
\begin{pspicture}(4,2)
  \pscustom[fillstyle=solid,fillcolor=lightgray,linecolor=lightgray]{
    \psbezier(1.0,2.0)(1.0,1.8)(1.3,1.7)(1.5,1.7)
    \psbezier(1.5,1.7)(1.7,1.7)(2.0,1.8)(2.0,2.0)
    \psline(2.0,2.0)(3.0,2.0)(3.0,1.0)
    \psbezier(3.0,1.0)(3.0,0.8)(2.7,0.7)(2.5,0.7)
    \psbezier(2.5,0.7)(2.3,0.7)(2.0,0.8)(2.0,1.0)
    \psbezier(2.0,1.0)(2.0,1.2)(1.7,1.3)(1.5,1.3)
    \psbezier(1.5,1.3)(1.3,1.3)(1.0,1.2)(1.0,1.0)
    \psline(1.0,1.0)(1.0,0.0)(0.5,0.0)(0.5,2.0)(1.0,2.0)
  }
  \psbezier[linewidth=1pt](1.0,2.0)(1.0,1.8)(1.3,1.7)(1.5,1.7)
  \psbezier[linewidth=1pt](1.5,1.7)(1.7,1.7)(2.0,1.8)(2.0,2.0)
  \psline[linewidth=1pt](3.0,2.0)(3.0,1.0)
  \psbezier[linewidth=1pt](3.0,1.0)(3.0,0.8)(2.7,0.7)(2.5,0.7)
  \psbezier[linewidth=1pt](2.5,0.7)(2.3,0.7)(2.0,0.8)(2.0,1.0)
  \psbezier[linewidth=1pt](2.0,1.0)(2.0,1.2)(1.7,1.3)(1.5,1.3)
  \psbezier[linewidth=1pt](1.5,1.3)(1.3,1.3)(1.0,1.2)(1.0,1.0)
  \psline[linewidth=1pt](1.0,1.0)(1.0,0.0)
\end{pspicture}
\end{aligned}
\qquad D_{\{1\}}^{(+1)} =
\begin{aligned}
\begin{pspicture}(4,2)
  \pscustom[fillstyle=solid,fillcolor=lightgray,linecolor=lightgray]{
    \psline(1.0,2.0)(1.0,0.0)(0.5,0.0)(0.5,2.0)(1.0,2.0)
  }
  \pscustom[fillstyle=solid,fillcolor=lightgray,linecolor=lightgray]{
    \psline(2.0,2.0)(2.0,1.0)
    \psbezier(2.0,1.0)(2.0,0.8)(2.3,0.7)(2.5,0.7)
    \psbezier(2.5,0.7)(2.7,0.7)(3.0,0.8)(3.0,1.0)
    \psline(3.0,1.0)(3.0,2.0)(2.0,2.0)
  }
  \psline[linewidth=1pt](1.0,2.0)(1.0,0.0)
  \psline[linewidth=1pt](2.0,2.0)(2.0,1.0)
  \psbezier[linewidth=1pt](2.0,1.0)(2.0,0.8)(2.3,0.7)(2.5,0.7)
  \psbezier[linewidth=1pt](2.5,0.7)(2.7,0.7)(3.0,0.8)(3.0,1.0)
  \psline[linewidth=1pt](3.0,1.0)(3.0,2.0)
\end{pspicture}
\end{aligned}
\end{equation}
The checkerboard colouring of a smoothing $D_\alpha^{(\epsilon)}$
defines an orientation of the $1$-manifold
$D_\alpha^{(\epsilon)}\subseteq\R^2$. One therefore draws arrows on
the components of $D_\alpha$ in such a way that the arrows go
counter-clockwise around the shaded regions. Often, we omit the
shading and indicate only the orientation of the smoothings. For the
above examples, we have
\begin{equation}
\label{eq_orientedsmoothing}
  D_\emptyset^{(+1)} =
\begin{aligned}
\begin{pspicture}(4,2)
  \psbezier[linewidth=1pt](1.0,2.0)(1.0,1.8)(1.3,1.7)(1.5,1.7)
  \psbezier[linewidth=1pt]{<-}(1.5,1.7)(1.7,1.7)(2.0,1.8)(2.0,2.0)
  \psline[linewidth=1pt](3.0,2.0)(3.0,1.0)
  \psbezier[linewidth=1pt](3.0,1.0)(3.0,0.8)(2.7,0.7)(2.5,0.7)
  \psbezier[linewidth=1pt](2.5,0.7)(2.3,0.7)(2.0,0.8)(2.0,1.0)
  \psbezier[linewidth=1pt](2.0,1.0)(2.0,1.2)(1.7,1.3)(1.5,1.3)
  \psbezier[linewidth=1pt]{<-}(1.5,1.3)(1.3,1.3)(1.0,1.2)(1.0,1.0)
  \psline[linewidth=1pt](1.0,1.0)(1.0,0.0)
\end{pspicture}
\end{aligned}
,\qquad D_{\{1\}}^{(+1)} =
\begin{aligned}
\begin{pspicture}(4,2)
  \psline[linewidth=1pt](1.0,2.0)(1.0,1.0)
  \psline[linewidth=1pt]{<-}(1.0,1.0)(1.0,0.0)
  \psline[linewidth=1pt](2.0,2.0)(2.0,1.0)
  \psbezier[linewidth=1pt]{->}(2.0,1.0)(2.0,0.8)(2.3,0.7)(2.5,0.7)
  \psbezier[linewidth=1pt](2.5,0.7)(2.7,0.7)(3.0,0.8)(3.0,1.0)
  \psline[linewidth=1pt](3.0,1.0)(3.0,2.0)
\end{pspicture}
\end{aligned}
\end{equation}
Changing the checkerboard colouring from $\epsilon$ to $-\epsilon$
reverses the orientation of the diagram $D_\alpha$. Notice that the
orientation of the smoothing $D_\alpha$ depends only on the
checkerboard colouring, but it has nothing to do with a possible
orientation of the tangle.

\subsection{Constructing the tangle chain complex}
\label{sec_Tangle}

In the following definition and theorem, we translate the
$2$-manifolds with corners used by Bar-Natan in his `picture world'
construction~\cite{BN2} into the language of our open-closed
cobordisms.

\begin{defn}
\label{def_tanglecube}
Let $T$ be a plane diagram of an oriented tangle with $n_+$ positive
crossings $(\overcrossing)$ and $n_-$ negative ones
$(\undercrossing)$, $n:=n_++n_-$. Number the crossings
$1,\ldots,n$. Choose a checkerboard colouring $\epsilon\in\{-1,+1\}$.

Let $D^{(\epsilon)}_\alpha\subseteq\R^2$, $\alpha\subseteq\cat{n}$,
denote the plane diagrams associated with the smoothings of all
crossings, equipped with the orientation induced from the
checkerboard colouring. Number the components of each
$D^{(\epsilon)}_\alpha$ by $1,\ldots,k_\alpha$, $k_\alpha\in\N_0$,
and define a sequence
$m^{(\alpha)}=(m_1^{(\alpha)},\ldots,m_{k_\alpha}^{(\alpha)})\in{\{0,1\}}^{k_\alpha}$
such that $m_\ell^{(\alpha)}=0$ if the $\ell$-th component is a
circle and $m_\ell^{(\alpha)}=1$ if it is an arc.

For all $j\in\cat{n}$ and $\alpha\subseteq\cat{n}\backslash\{j\}$,
consider the open-closed cobordism
\begin{equation}
  S_{(\alpha,j)}\colon D_\alpha\to D_{\alpha\sqcup\{j\}},
\end{equation}
which is a cylinder over $D_\alpha^{(\epsilon)}$ except for a sufficiently small
neighbourhood of the $j$-th crossing $(\slashoverback)$ where it is
the saddle
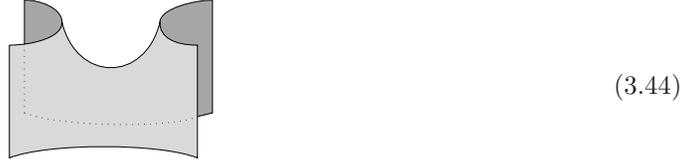
\begin{equation}
\label{eq_saddleagain}
\begin{aligned}
\psset{xunit=1cm,yunit=1cm}
\begin{pspicture}(3.0,2.0)
  \rput(0,0){\bigsaddle}
\end{pspicture}
\end{aligned}
\end{equation}
with the $0$-smoothing $(\smoothing)$ as the source and the
$1$-smoothing $(\hsmoothing)$ as the target.

We define the $\cat{n}$-cube $S_{T,\epsilon}$ in
$\cat{2Cob}^\mathrm{ext}$ to have as vertices the sequences
$m^{(\alpha)}$ and as edges the (equivalence classes of) the
open-closed cobordisms $S_{(\alpha,j)}$.
\end{defn}

\begin{rem}
For $j\in\cat{n}$, $\alpha\in\cat{n}\backslash\{j\}$, both
$D_\alpha^{(\epsilon)}$ and $D_{\alpha\sqcup\{j\}}^{(\epsilon)}$ are
equipped with the orientation induced from the checkerboard
colouring $\epsilon$. This determines the orientation of the
open-closed cobordism $S_{(\alpha,j)}$. In order to see this, we
display below the relevant saddle which is a smooth oriented
$2$-manifold with corners $M$ representing an open-closed cobordism,
together with its source $D_\alpha^{(\epsilon)}$ and target
$D_{\alpha\sqcup\{j\}}^{(\epsilon)}$ as oriented manifolds.
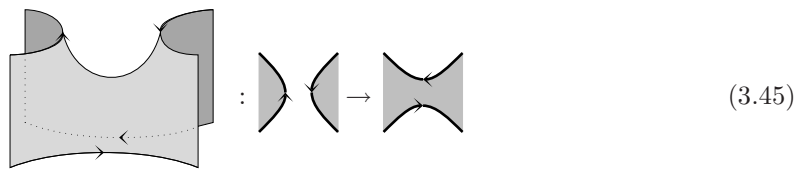
\begin{equation}
\label{eq_saddleorient}
\begin{aligned}
\psset{xunit=1cm,yunit=1cm}
\begin{pspicture}(3.0,2.0)
  \rput(0,0){\bigsaddleorient}
\end{pspicture}
\end{aligned}
\colon
\begin{aligned}
\begin{pspicture}(1.5,1.5)
  \pscustom[fillstyle=solid,fillcolor=lightgray,linecolor=lightgray]{
    \psbezier(0.0,0.0)(0.25,0.25)(0.5,0.5)(0.5,0.75)
    \psbezier(0.5,0.75)(0.5,1.0)(0.25,1.25)(0.0,1.5)
    \psline(0.0,1.5)(0.0,0.0)
  }
  \pscustom[fillstyle=solid,fillcolor=lightgray,linecolor=lightgray]{
    \psbezier(1.5,1.5)(1.25,1.25)(1.0,1.0)(1.0,0.75)
    \psbezier(1.0,0.75)(1.0,0.5)(1.25,0.25)(1.5,0.0)
    \psline(1.5,0.0)(1.5,1.5)
  }
  \psbezier[linewidth=1pt]{->}(0.0,0.0)(0.25,0.25)(0.5,0.5)(0.5,0.75)
  \psbezier[linewidth=1pt](0.5,0.75)(0.5,1.0)(0.25,1.25)(0.0,1.5)
  \psbezier[linewidth=1pt]{->}(1.5,1.5)(1.25,1.25)(1.0,1.0)(1.0,0.75)
  \psbezier[linewidth=1pt](1.0,0.75)(1.0,0.5)(1.25,0.25)(1.5,0.0)
\end{pspicture}
\end{aligned}
\to
\begin{aligned}
\begin{pspicture}(1.5,1.5)
  \pscustom[fillstyle=solid,fillcolor=lightgray,linecolor=lightgray]{
    \psbezier(0.0,0.0)(0.25,0.25)(0.5,0.5)(0.75,0.5)
    \psbezier(0.75,0.5)(1.0,0.5)(1.25,0.25)(1.5,0.0)
    \psline(1.5,0.0)(1.5,1.5)
    \psbezier(1.5,1.5)(1.25,1.25)(1.0,1.0)(0.75,1.0)
    \psbezier(0.75,1.0)(0.5,1.0)(0.25,1.25)(0.0,1.5)
    \psline(0.0,1.5)(0.0,0.0)
  }
  \psbezier[linewidth=1pt]{->}(0.0,0.0)(0.25,0.25)(0.5,0.5)(0.75,0.5)
  \psbezier[linewidth=1pt](0.75,0.5)(1.0,0.5)(1.25,0.25)(1.5,0.0)
  \psbezier[linewidth=1pt]{->}(1.5,1.5)(1.25,1.25)(1.0,1.0)(0.75,1.0)
  \psbezier[linewidth=1pt](0.75,1.0)(0.5,1.0)(0.25,1.25)(0.0,1.5)
\end{pspicture}
\end{aligned}
\end{equation}
In our pictures, we equip the target of $M$ with the orientation
induced from the orientation of $M$ and the source with the opposite
of the induced orientation. In this way, the orientation of the
smoothing $D_\alpha^{(\epsilon)}$ matches that of the source of $M$
and the orientation of $D_{\alpha\sqcup\{j\}}^{(\epsilon)}$ matches
that of the target.
\end{rem}

\begin{thm}[see~\cite{BN2}]
Given a plane diagram $T$ of an oriented tangle with $n$ crossings and
a checkerboard colouring $\epsilon\in\{-1,+1\}$, the $\cat{n}$-cube
$S_{T,\epsilon}$ is commutative.
\end{thm}

\begin{proof}
The coordinate of the vertical axis of our
diagrams~\eqref{eq_saddleagain} is a Morse function for the
open-closed cobordism. The saddle has a single critical point of index
$1$. Let $j,\ell\in\cat{n}$, $j\neq\ell$, and
$\alpha\subseteq\cat{n}\backslash\{j,\ell\}$. The diagram
\begin{equation}
\begin{aligned}
\xymatrix@C=1.8pc@R=1pc{
  D^{(\epsilon)}_\alpha\ar[rr]^-{S_{(\alpha,j)}}\ar[dd]_-{S_{(\alpha,\ell)}}&&
    D^{(\epsilon)}_{\alpha\sqcup\{j\}}\ar[dd]^-{S_{(\alpha\sqcup\{j\},\ell)}}\\
  \\
  D^{(\epsilon)}_{\alpha\sqcup\{\ell\}}\ar[rr]_-{S_{(\alpha\sqcup\{\ell\},j)}}&&
    D^{(\epsilon)}_{\alpha\sqcup\{j,\ell\}}
}
\end{aligned}
\end{equation}
commutes because the two open-closed cobordisms
$S_{(\alpha\sqcup\{j\},\ell)}\circ S_{(\alpha,j)}$ and
$S_{(\alpha\sqcup\{\ell\},j)}\circ S_{(\alpha,\ell)}$ have handle
decompositions related by interchanging the critical values of the two
index $1$ critical points.
\end{proof}

Departing from Bar-Natan's work, we now apply an open-closed TQFT
$Z\colon\cat{2Cob}^\mathrm{ext}\to\cat{Vect}_k$ and turn the
commutative $\cat{n}$-cube $S_{T,\epsilon}$ in
$\cat{2Cob}^\mathrm{ext}$ (topology) into a commutative
$\cat{n}$-cube $Z(S_{T,\epsilon})$ in $\cat{Vect}_k$ (algebra).

\begin{defn}
\label{def_vectcube}
Let $T$ be a plane diagram of an oriented tangle with $n$ crossings,
let $\epsilon\in\{-1,+1\}$ denote a checkerboard colouring of $T$, and
let $Z\colon\cat{2Cob}^\mathrm{ext}\to\cat{Vect}_k$ be an Euler-graded
[Euler-filtered] open-closed TQFT. We define the $\cat{n}$-cube
$Z(S_{T,\epsilon})$ in $\cat{Vect}_k$ to have as vertices the vector
spaces $Z(m^{(\alpha)})\{2|\alpha|\}$ and as edges the linear maps
$Z(S_{(\alpha,j)})$.
\end{defn}

We use the notation $Z(S_{T,\epsilon})$ because up to a grading shift
by $2|\alpha|$ this cube is just the image of $S_{T,\epsilon}$ under
the functor $Z$. We will sometimes just write $S$ and $Z(S)$ rather
than $S_{T,\epsilon}$ and $Z(S_{T,\epsilon})$ if the plane tangle
diagram $T$ and the checkerboard colouring $\epsilon$ are evident from
the context.

\begin{prop}
\label{prop_degreezero} The $\cat{n}$-cube $Z(S_{T,\epsilon})$ is
commutative. All edges are graded [filtered] $k$-linear maps of
degree $0$, and so $Z(S_{T,\epsilon})$ forms an object of
$\cat{Cube}(\cat{grdVect}_k)$ [or $\cat{Cube}(\cat{fltVect}_k)$].
\end{prop}

\begin{proof}
The cube $Z(S_{T,\epsilon})$ is commutative because $S_{T,\epsilon}$
is commutative and $Z$ is a functor. Let $M$ be the
saddle~\eqref{eq_saddleagain} as an open-closed cobordism. It has
$\chi(M)=1$, $|\Pi_0(\del_1M)|=4$, and $\omega(M)=0$. By
Proposition~\ref{prop_degree} and Remark~\ref{rem_degree}(2.), we
have the degree
\begin{equation}
  \deg S_{(\alpha,j)} = -2
\end{equation}
for all $j\in\cat{n}$, $\alpha\subseteq\cat{n}\backslash\{j\}$, \ie\
for all edges of the cube. The degree shift by twice the height of
each vertex $\alpha\subseteq\cat{n}$ in $Z(m^{(\alpha)})\{2|\alpha|\}$
makes sure that all edges are graded [filtered] of degree $0$.
\end{proof}

\begin{rem}
\label{rem_intersect}\hfill
\begin{enumerate}
\item
  All the cobordisms $S_{(\alpha,j)}$ have the following property:
  Each connected component of $\del_1S_{(\alpha,j)}$ has one boundary
  point belonging to the source and the other to the target of the
  cobordism. In particular, $S_{(\alpha,j)}$ has window number zero.
\item
  Inverting the checkerboard colouring does not affect the tensor
  factors $C$ of the $Z(D_\alpha)$, but it does affect $A$ by
  replacing it with $A^\mathrm{op}$. This is the Frobenius
  algebra opposite to $A$ whose multiplication and comultiplication
  are given by $\mu_{A^\mathrm{op}}=\mu_A\circ\tau_{A,A}$ and
  $\Delta_{A^\mathrm{op}}=\tau_{A,A}\circ\Delta_A$.
\end{enumerate}
\end{rem}

\begin{rem}
\label{rem_crossings}
If $T$ is a plane diagram of an oriented link, each $2$-dimensional
cobordism $S_{(\alpha,j)}$ is the disjoint union of either
$\psset{xunit=.3em,yunit=.3em}\begin{pspicture}(3,2.5)\rput(2,0){\multc}\end{pspicture}$
or
$\psset{xunit=.3em,yunit=.3em}\begin{pspicture}(3,2.5)\rput(2,0){\comultc}\end{pspicture}$
(pairs of pants) with a finite number of cylinders over the circle.

Now let $T$ be a plane diagram of an oriented tangle with $n$
crossings, $\epsilon\in\{-1,+1\}$ be a checkerboard colouring and
$j\in\cat{n}$ denote a crossing. The open-closed cobordism
$S_{(\alpha,j)}\colon D_\alpha^{(\epsilon)}\to
D_{\alpha\sqcup\{j\}}^{(\epsilon)}$ is the disjoint union of one
connected component that contains the saddle with a finite number of
cylinders over the circle and over the arc. Depending on the
smoothing $\alpha\subseteq\cat{n}\backslash\{j\}$ of the other
crossings, the component with the saddle can have one of the
following forms. On the left, we display the $j$-th crossing with
the checkerboard colouring and with the connected components of the
diagram $T$ that are attached to that crossing. Then we show the two
smoothings $D_\alpha^{(\epsilon)}$ and
$D_{\alpha\sqcup\{j\}}^{(\epsilon)}$ of the $j$-th crossing, and
finally on the right we show the associated open-closed cobordism as
a composition of the generators, read from top to bottom.
\begin{eqnarray}
\label{eq_crossing}
\begin{aligned}
\begin{pspicture}(1.5,1.5)
  \pscustom[fillstyle=solid,fillcolor=lightgray,linecolor=lightgray]{
    \psline(0.0,0.0)(1.5,1.5)(1.5,0.0)(0.0,1.5)(0.0,0.0)
  }
  \psline[linewidth=1pt](0.0,1.5)(1.5,0.0)
  \psline[linewidth=3pt,linecolor=white](0.7,0.8)(0.8,0.7)
  \psline[linewidth=1pt](0.0,0.0)(1.5,1.5)
\end{pspicture}
\end{aligned}
&\colon&
\begin{aligned}
\begin{pspicture}(1.5,1.5)
  \psbezier[linewidth=1pt]{->}(0.0,0.0)(0.25,0.25)(0.5,0.25)(0.5,0.75)
  \psbezier[linewidth=1pt](0.5,0.75)(0.5,1.25)(0.25,1.25)(0.0,1.5)
  \psbezier[linewidth=1pt]{->}(1.5,1.5)(1.25,1.25)(1.0,1.25)(1.0,0.75)
  \psbezier[linewidth=1pt](1.0,0.75)(1.0,0.25)(1.25,0.25)(1.5,0.0)
  \rput(0.2,0.75){$1$}
  \rput(1.3,0.75){$2$}
\end{pspicture}
\end{aligned}
\to
\begin{aligned}
\begin{pspicture}(1.5,1.5)
  \psbezier[linewidth=1pt]{->}(0.0,0.0)(0.25,0.25)(0.25,0.5)(0.75,0.5)
  \psbezier[linewidth=1pt](0.75,0.5)(1.25,0.5)(1.25,0.25)(1.5,0.0)
  \psbezier[linewidth=1pt]{->}(1.5,1.5)(1.25,1.25)(1.25,1.0)(0.75,1.0)
  \psbezier[linewidth=1pt](0.75,1.0)(0.25,1.0)(0.25,1.25)(0.0,1.5)
  \rput(0.75,0.15){$1$}
  \rput(0.75,1.35){$2$}
\end{pspicture}
\end{aligned}
\qquad\qquad
\begin{aligned}
\psset{xunit=.3cm,yunit=.3cm}
\begin{pspicture}(4.5,5.5)
  \rput(0,0){\saddler}
  \psline{->}(3.5,5)(4.2,5)
  \psline{->}(-.5,5)(0.2,5)
  \psline{->}(3.5,0)(4.2,0)
  \psline{->}(-.5,0)(0.2,0)
\end{pspicture}
\end{aligned}\colon A\otimes A\to A\otimes A\\
\label{eq_crossing2}
\begin{aligned}
\begin{pspicture}(1.5,1.5)
  \pscustom[fillstyle=solid,fillcolor=lightgray,linecolor=lightgray]{
    \psline(0.0,0.0)(1.5,0.0)(0.0,1.5)(1.5,1.5)(0.0,0.0)
  }
  \psline[linewidth=1pt](0.0,0.0)(1.5,1.5)
  \psline[linewidth=3pt,linecolor=white](0.7,0.7)(0.8,0.8)
  \psline[linewidth=1pt](0.0,1.5)(1.5,0.0)
\end{pspicture}
\end{aligned}
&\colon&
\begin{aligned}
\begin{pspicture}(1.5,1.5)
  \psbezier[linewidth=1pt](0.0,0.0)(0.25,0.25)(0.25,0.5)(0.75,0.5)
  \psbezier[linewidth=1pt]{<-}(0.75,0.5)(1.25,0.5)(1.25,0.25)(1.5,0.0)
  \psbezier[linewidth=1pt](1.5,1.5)(1.25,1.25)(1.25,1.0)(0.75,1.0)
  \psbezier[linewidth=1pt]{<-}(0.75,1.0)(0.25,1.0)(0.25,1.25)(0.0,1.5)
  \rput(0.75,0.15){$1$}
  \rput(0.75,1.35){$2$}
\end{pspicture}
\end{aligned}
\to
\begin{aligned}
\begin{pspicture}(1.5,1.5)
  \psbezier[linewidth=1pt](0.0,0.0)(0.25,0.25)(0.5,0.25)(0.5,0.75)
  \psbezier[linewidth=1pt]{<-}(0.5,0.75)(0.5,1.25)(0.25,1.25)(0.0,1.5)
  \psbezier[linewidth=1pt](1.5,1.5)(1.25,1.25)(1.0,1.25)(1.0,0.75)
  \psbezier[linewidth=1pt]{<-}(1.0,0.75)(1.0,0.25)(1.25,0.25)(1.5,0.0)
  \rput(0.2,0.75){$1$}
  \rput(1.3,0.75){$2$}
\end{pspicture}
\end{aligned}
\qquad\qquad
\begin{aligned}
\psset{xunit=.3cm,yunit=.3cm}
\begin{pspicture}[.4](4.5,5.5)
  \rput(4,0){\saddlel}
  \psline{->}(3.5,5)(4.2,5)
  \psline{->}(-.5,5)(0.2,5)
  \psline{->}(3.5,0)(4.2,0)
  \psline{->}(-.5,0)(0.2,0)
\end{pspicture}
\end{aligned}\colon A\otimes A\to A\otimes A\\
\begin{aligned}
\begin{pspicture}(1.5,1.5)
  \pscustom[fillstyle=solid,fillcolor=lightgray,linecolor=lightgray]{
    \psline(0.0,0.0)(1.5,0.0)(1.5,1.5)(0.0,1.5)(0.0,0.0)
  }
  \pscustom[fillstyle=solid,fillcolor=white,linecolor=white]{
    \psline(0.0,0.0)(1.0,0.5)
    \psbezier(1.0,0.5)(1.5,1.0)(1.25,1.5)(0.75,1.5)
    \psbezier(0.75,1.5)(0.25,1.5)(0.0,1.0)(0.5,0.5)
    \psline(0.5,0.5)(1.5,0.0)(0.0,0.0)
  }
  \psbezier[linewidth=1pt](0.75,1.5)(0.25,1.5)(0.0,1.0)(0.5,0.5)
  \psline[linewidth=1pt](0.5,0.5)(1.5,0.0)
  \psline[linewidth=3pt,linecolor=white](0.7,0.45)(0.8,0.3)
  \psline[linewidth=1pt](0.0,0.0)(1.0,0.5)
  \psbezier[linewidth=1pt](1.0,0.5)(1.5,1.0)(1.25,1.5)(0.75,1.5)
\end{pspicture}
\end{aligned}
&\colon&
\begin{aligned}
\begin{pspicture}(1.5,1.5)
  \psbezier[linewidth=1pt](0.0,0.0)(0.25,0.25)(0.5,0.5)(0.25,0.75)
  \psbezier[linewidth=1pt]{->}(0.25,0.75)(0.0,1.0)(0.25,1.5)(0.75,1.5)
  \psbezier[linewidth=1pt](0.75,1.5)(1.25,1.5)(1.5,1.0)(1.25,0.75)
  \psbezier[linewidth=1pt](1.25,0.75)(1.0,0.5)(1.25,0.25)(1.5,0.0)
\end{pspicture}
\end{aligned}
\to
\begin{aligned}
\begin{pspicture}(1.5,1.5)
  \psbezier[linewidth=1pt]{->}(0.0,0.0)(0.25,0.25)(0.65,0.2)(0.75,0.2)
  \psbezier[linewidth=1pt](0.75,0.2)(0.85,0.2)(1.25,0.25)(1.5,0.0)
  \psbezier[linewidth=1pt](0.25,1.0)(0.25,1.25)(0.5,1.5)(0.75,1.5)
  \psbezier[linewidth=1pt]{->}(0.75,1.5)(1.0,1.5)(1.25,1.25)(1.25,1.0)
  \psbezier[linewidth=1pt](1.25,1.0)(1.25,0.75)(1.0,0.5)(0.75,0.5)
  \psbezier[linewidth=1pt](0.75,0.5)(0.5,0.5)(0.25,0.75)(0.25,1.0)
  \rput(0.75,-0.1){$2$}
  \rput(0.75,1.0){$1$}
\end{pspicture}
\end{aligned}
\qquad\qquad
\begin{aligned}
\psset{xunit=.3cm,yunit=.3cm}
\begin{pspicture}(4,5)
  \rput(2,2){\comultl}
  \rput(1.1,0){\ltc}
  \rput(3,0){\medidentl}
\end{pspicture}
\end{aligned}\colon A\to C\otimes A\\
\begin{aligned}
\begin{pspicture}(1.5,1.5)
  \pscustom[fillstyle=solid,fillcolor=lightgray,linecolor=lightgray]{
    \psline(0.0,0.0)(1.5,0.0)(1.5,1.5)(0.0,1.5)(0.0,0.0)
  }
  \pscustom[fillstyle=solid,fillcolor=white,linecolor=white]{
    \psline(0.0,0.0)(1.0,0.5)
    \psbezier(1.0,0.5)(1.5,1.0)(1.25,1.5)(0.75,1.5)
    \psbezier(0.75,1.5)(0.25,1.5)(0.0,1.0)(0.5,0.5)
    \psline(0.5,0.5)(1.5,0.0)(0.0,0.0)
  }
  \psline[linewidth=1pt](0.0,0.0)(1.0,0.5)
  \psbezier[linewidth=1pt](1.0,0.5)(1.5,1.0)(1.25,1.5)(0.75,1.5)
  \psline[linewidth=3pt,linecolor=white](0.7,0.3)(0.8,0.45)
  \psbezier[linewidth=1pt](0.75,1.5)(0.25,1.5)(0.0,1.0)(0.5,0.5)
  \psline[linewidth=1pt](0.5,0.5,)(1.5,0.0)
\end{pspicture}
\end{aligned}
&\colon&
\begin{aligned}
\begin{pspicture}(1.5,1.5)
  \psbezier[linewidth=1pt]{->}(0.0,0.0)(0.25,0.25)(0.65,0.2)(0.75,0.2)
  \psbezier[linewidth=1pt](0.75,0.2)(0.85,0.2)(1.25,0.25)(1.5,0.0)
  \psbezier[linewidth=1pt](0.25,1.0)(0.25,1.25)(0.5,1.5)(0.75,1.5)
  \psbezier[linewidth=1pt]{->}(0.75,1.5)(1.0,1.5)(1.25,1.25)(1.25,1.0)
  \psbezier[linewidth=1pt](1.25,1.0)(1.25,0.75)(1.0,0.5)(0.75,0.5)
  \psbezier[linewidth=1pt](0.75,0.5)(0.5,0.5)(0.25,0.75)(0.25,1.0)
  \rput(0.75,-0.1){$2$}
  \rput(0.75,1.0){$1$}
\end{pspicture}
\end{aligned}
\to
\begin{aligned}
\begin{pspicture}(1.5,1.5)
  \psbezier[linewidth=1pt](0.0,0.0)(0.25,0.25)(0.5,0.5)(0.25,0.75)
  \psbezier[linewidth=1pt]{->}(0.25,0.75)(0.0,1.0)(0.25,1.5)(0.75,1.5)
  \psbezier[linewidth=1pt](0.75,1.5)(1.25,1.5)(1.5,1.0)(1.25,0.75)
  \psbezier[linewidth=1pt](1.25,0.75)(1.0,0.5)(1.25,0.25)(1.5,0.0)
\end{pspicture}
\end{aligned}
\qquad\qquad
\begin{aligned}
\psset{xunit=.3cm,yunit=.3cm}
\begin{pspicture}(4,5)
  \rput(2,0){\multl}
  \rput(.92,2.5){\ctl}
  \rput(3,2.5){\medidentl}
\end{pspicture}
\end{aligned}\colon C\otimes A\to A\\
\begin{aligned}
\begin{pspicture}(1.5,1.5)
  \pscustom[fillstyle=solid,fillcolor=lightgray,linecolor=lightgray]{
    \psline(0.0,0.0)(1.5,0.0)(1.5,1.5)(0.0,1.5)(0.0,0.0)
  }
  \pscustom[fillstyle=solid,fillcolor=white,linecolor=white]{
    \psbezier(1.5,0.75)(1.5,0.25)(1.25,0.25)(0.75,0.75)
    \psbezier(0.75,0.75)(0.25,1.25)(0.0,1.25)(0.0,0.75)
    \psbezier(0.0,0.75)(0.0,0.25)(0.25,0.25)(0.75,0.75)
    \psbezier(0.75,0.75)(1.25,1.25)(1.5,1.25)(1.5,0.75)
  }
  \psbezier[linewidth=1pt](1.5,0.75)(1.5,0.25)(1.25,0.25)(0.75,0.75)
  \psbezier[linewidth=1pt](0.75,0.75)(0.25,1.25)(0.0,1.25)(0.0,0.75)
  \psline[linewidth=3pt,linecolor=white](0.7,0.8)(0.8,0.7)
  \psbezier[linewidth=1pt](0.0,0.75)(0.0,0.25)(0.25,0.25)(0.75,0.75)
  \psbezier[linewidth=1pt](0.75,0.75)(1.25,1.25)(1.5,1.25)(1.5,0.75)
\end{pspicture}
\end{aligned}
&\colon&
\begin{aligned}
\begin{pspicture}(1.5,1.5)
  \pscircle[linewidth=1pt](0.4,0.75){0.2}
  \pscircle[linewidth=1pt](1.1,0.75){0.2}
  \psline[linewidth=1pt]{->}(0.41,1.05)(0.43,1.05)
  \psline[linewidth=1pt]{->}(1.11,1.05)(1.13,1.05)
\end{pspicture}
\end{aligned}
\to
\begin{aligned}
\begin{pspicture}(1.5,1.5)
  \pscircle[linewidth=1pt](0.75,0.75){0.4}
  \psline[linewidth=1pt]{->}(0.76,1.35)(0.78,1.35)
\end{pspicture}
\end{aligned}
\qquad\qquad
\begin{aligned}
\psset{xunit=.3cm,yunit=.3cm}
\begin{pspicture}(4,5)
  \rput(2,1){\multc}
\end{pspicture}
\end{aligned}\colon C\otimes C\to C\\
\label{eq_copants}
\begin{aligned}
\begin{pspicture}(1.5,1.5)
  \pscustom[fillstyle=solid,fillcolor=lightgray,linecolor=lightgray]{
    \psline(0.0,0.0)(1.5,0.0)(1.5,1.5)(0.0,1.5)(0.0,0.0)
  }
  \pscustom[fillstyle=solid,fillcolor=white,linecolor=white]{
    \psbezier(0.75,0.0)(1.25,0.0)(1.25,0.25)(0.75,0.75)
    \psbezier(0.75,0.75)(0.25,1.25)(0.25,1.5)(0.75,1.5)
    \psbezier(0.75,1.5)(1.25,1.5)(1.25,1.25)(0.75,0.75)
    \psbezier(0.75,0.75)(0.25,0.25)(0.25,0.0)(0.75,0.0)
  }
  \psbezier[linewidth=1pt](0.75,0.0)(1.25,0.0)(1.25,0.25)(0.75,0.75)
  \psbezier[linewidth=1pt](0.75,0.75)(0.25,1.25)(0.25,1.5)(0.75,1.5)
  \psline[linewidth=3pt,linecolor=white](0.7,0.8)(0.8,0.7)
  \psbezier[linewidth=1pt](0.75,1.5)(1.25,1.5)(1.25,1.25)(0.75,0.75)
  \psbezier[linewidth=1pt](0.75,0.75)(0.25,0.25)(0.25,0.0)(0.75,0.0)
\end{pspicture}
\end{aligned}
&\colon&
\begin{aligned}
\begin{pspicture}(1.5,1.5)
  \pscircle[linewidth=1pt](0.75,0.75){0.4}
  \psline[linewidth=1pt]{->}(0.15,0.76)(0.15,0.78)
\end{pspicture}
\end{aligned}
\to
\begin{aligned}
\begin{pspicture}(1.5,1.5)
  \pscircle[linewidth=1pt](0.75,0.4){0.2}
  \pscircle[linewidth=1pt](0.75,1.1){0.2}
  \psline[linewidth=1pt]{->}(0.45,0.41)(0.45,0.43)
  \psline[linewidth=1pt]{->}(0.45,1.11)(0.45,1.13)
\end{pspicture}
\end{aligned}
\qquad\qquad
\begin{aligned}
\psset{xunit=.3cm,yunit=.3cm}
\begin{pspicture}(4,5)
  \rput(2,1){\comultc}
\end{pspicture}
\end{aligned}\colon C\to C\otimes C
\end{eqnarray}
The components of the diagrams $D_\alpha^{(\epsilon)}$ and
$D_{\alpha\sqcup\{j\}}^{(\epsilon)}$ are numbered $1,2,\ldots$
corresponding to the tensor factors in the source and target of the
open-closed cobordism on the right. If necessary, we indicate with a
little arrow the orientation of the arc in the black boundary of the
open-closed cobordism.

All other possibilities for the component of the open-closed
cobordism that contains the saddle and the remaining checkerboard
colourings can be obtained from these by rotations in the plane.
Changing the numbering of the components results in permutations of
the tensor factors.

We encourage the reader to visualize how the open-closed cobordism
displayed in~\eqref{eq_crossing} can be obtained from the
saddle~\eqref{eq_saddleorient} by `smashing' the saddle flat on the
drawing plane and then pulling the source to the top and the target to
the bottom of the diagram.
\end{rem}

\subsection{The tangle complex is well defined}

In this section, we show that the tangle complex
${\tangle{T}}_{\epsilon,Z}$ of~\eqref{eq_tanglecomplex} is well
defined up to isomorphism regardless of the numbering of the
crossings of $T$ and regardless of the numbering of the components
of each smoothing $D^{(\epsilon)}_\alpha\subseteq\R^2$.

\begin{prop}
\label{prop_indep} Let $T$ be a plane diagram of an oriented tangle
with $n$ crossings, $\epsilon\in\{-1,+1\}$ be a checkerboard
colouring, and let $Z\colon\cat{2Cob}^\mathrm{ext}\to\cat{Vect}_k$
be an open-closed TQFT. The tangle complex
${\tangle{T}}_{\epsilon,Z}$ is
\begin{enumerate}
\item
  independent of the numbering of the crossings up to isomorphism.
\item
  independent of the numbering of components of any smoothing
  $D_\alpha^{(\epsilon)}$.
\end{enumerate}
\end{prop}

\begin{proof}\hfill
\begin{enumerate}
\item
  Given some numbering $1,\ldots,n$ of the crossings of $T$, we
  represent a change of this numbering by a bijection
  $f^0\colon\cat{n}\to\cat{n}$. Consider the commutative cube
  $S_{T,\epsilon}$ of Definition~\ref{def_tanglecube}. There is a
  homomorphism of commutative cubes
  $f=(f^0,{\{f_\alpha\}}_{\alpha\subseteq\cat{n}})\colon
  S_{T,\epsilon}\to S_{T,\epsilon}$ with
  $f_\alpha:=\id_{m^{(\alpha)}}$ whose inverse is given similarly by
  the inverse bijection ${(f^0)}^{-1}$.

  This isomorphism of commutative cubes yields an isomorphism
  $Z(f)\colon Z(S_{T,\epsilon})\to Z(S_{T,\epsilon})$ for the
  commutative cubes of Definition~\ref{def_vectcube} and by
  Proposition~\ref{prop_total2} an isomorphism of complexes
  $C_\mathrm{tot}(Z(f))\colon{\tangle{T}}_{\epsilon,Z}\to{\tangle{T}}_{\epsilon,Z}$. Note
  that permuting the index set $\cat{n}$ of the commutative cube
  induces a permutation on the set of vertices that leaves the height
  of each vertex fixed.
\item
  Recall that a linear order of the components of
  $D_\alpha^{(\epsilon)}$ is needed in order to define the object
  $m^{(\alpha)}=(m_1^{(\alpha)},\ldots,m_{k_\alpha}^{(\alpha)})\in{\{0,1\}}^{k_\alpha}$
  of $\cat{2Cob}^\mathrm{ext}$ associated with $D_\alpha^{(\epsilon)}$
  in Definition~\ref{def_tanglecube}.

  Let $m^{(\alpha)}$ be this object of $\cat{2Cob}^\mathrm{ext}$, and
  let $\sigma\colon\{1,\ldots,k_\alpha\}\to\{1,\ldots,k_\alpha\}$ be a
  bijection. Then $\sigma$ induces a bijection
  $\bar\sigma\colon{\{0,1\}}^{k_\alpha}\to{\{0,1\}}^{k_\alpha}$,
  $m^{(\alpha)}\mapsto
  \sigma(m^{(\alpha)}):=(m_{\sigma^{-1}(1)}^{(\alpha)},\ldots,m_{\sigma^{-1}(k_\alpha)}^{(\alpha)})$.

  Let $S_{T,\epsilon}$ denote the commutative cube constructed using
  the ordering of the components $(1,\ldots,k_\alpha)$ for some vertex
  $\alpha\subseteq\cat{n}$, and let $S^\prime_{T,\epsilon}$ denote the
  commutative cube constructed using the same ordering of components as
  in $S_{T,\epsilon}$ except at the vertex $\alpha$, where we use the
  ordering $(\sigma^{-1}(1),\ldots,\sigma^{-1}(k_\alpha))$. Note that
  the edges of $S^\prime_{T,\epsilon}$ will in general differ from
  those of $S_{T,\epsilon}$ due to the different ordering of the
  components.

  There is a homomorphism
  $f=(f^0,{\{f_\beta\}}_{\beta\subseteq\cat{n}})\colon
  S_{T,\epsilon}\to S^\prime_{T,\epsilon}$ of commutative cubes with
  $f^0=\id_{\cat{n}}$, $f_\alpha=\bar\sigma$, and
  $f_\beta=\id_{{\{0,1\}}^{k_\alpha}}$ for all
  $\beta\subseteq\cat{n}$, $\beta\neq\alpha$. Its inverse is given
  similarly using ${\bar\sigma}^{-1}$ instead of $\bar\sigma$.

  Just as in part (1.), this isomorphism of commutative cubes induces
  an isomorphism of tangle complexes.
\end{enumerate}
\end{proof}

\begin{rem}\hfill
\begin{enumerate}
\item
  Part (2.) of the proposition implies in particular that the two
  saddles~\eqref{eq_crossing} and~\eqref{eq_crossing2} are related by
  permutations as follows:
\begin{equation}
\label{eq_saddleinversion} \psset{xunit=.25cm,yunit=.25cm}
\begin{pspicture}[.5](5,7.5)
  \rput(4,0){\saddlel}
  \rput(2,5){\widecrossl}
\end{pspicture}
\quad = \quad
\begin{pspicture}[.5](5,5)
  \rput(1,0){\saddler}
\end{pspicture}
\quad = \quad
\begin{pspicture}[.5](5,7.5)
  \rput(4,2.5){\saddlel}
  \rput(2,0){\widecrossl}
\end{pspicture}
\end{equation}
\item
  Inverting the checkerboard colouring has the effect of
  systematically replacing $A$ by the opposite Frobenius algebra
  $A^\mathrm{op}$, and of replacing $A^\mathrm{op}$-left actions by
  $A$-right actions and vice versa.
\end{enumerate}
\end{rem}

\subsection{Main results}

The following definition and theorem are the main results of this
section. They improve Bar-Natan's result of~\cite{BN2} and yield for
every tangle (not just for every link) a complex of \emph{vector
spaces} that is invariant under Reidemeister moves up to homotopy
equivalence.

\begin{defn}
\label{def_tanglecomplex}
Let $T$ be a plane diagram of an oriented tangle with $n_+$ positive
and $n_-$ negative crossings, $n:=n_++n_-$, let $\epsilon\in\{-1,+1\}$
denote a checkerboard colouring of $T$, and let
$Z\colon\cat{2Cob}^\mathrm{ext}\to\cat{Vect}_k$ be an Euler-graded
[Euler-filtered] open-closed TQFT.  We define the \emph{tangle
complex} as the total complex
\begin{equation}
\label{eq_tanglecomplex}
  {\tangle{T}}_{\epsilon,Z}:=
  C_\mathrm{tot}(\cat{n},Z(S_{T,\epsilon}))\{2n_+-4n_-\}[-n_-].
\end{equation}
We will often denote the complex ${\tangle{T}}_{\epsilon,Z}$ as
$\tangle{T}$ when no confusion is likely to arise.
\end{defn}

\begin{prop}
${\tangle{T}}_{\epsilon,Z}$ is an object of
$\cat{Kom}(\cat{grdVect}_k)$ or of $\cat{Kom}(\cat{fltVect}_k)$,
respectively.
\end{prop}

\begin{proof}
By Proposition~\ref{prop_degreezero}.
\end{proof}

\begin{thm}
\label{thm_reidemeister}
Let $T,T^\prime$ be plane diagrams of oriented tangles,
$\epsilon\in\{-1,+1\}$, and
$Z\colon\cat{2Cob}^\mathrm{ext}\to\cat{Vect}_k$ be an Euler-graded
[Euler-filtered] open-closed TQFT that satisfies Bar-Natan's
conditions. If $T$ and $T^\prime$ are related by a Reidemeister move
or by a plane isotopy, then the complexes ${\tangle{T}}_{\epsilon,Z}$
and ${\tangle{T^\prime}}_{\epsilon,Z}$ are homotopy equivalent as
graded [filtered] complexes.
\end{thm}

\begin{proof}
The proof is almost the same as Bar-Natan's~\cite{BN2}, except that
whenever Bar-Natan writes down a $2$-manifold with corners, we apply
the open-closed TQFT and turn it into a linear map between vector
space. We do not reproduce the proof here. We do, however, present the
proofs for the invariance under Reidemeister moves up to homotopy
equivalence when we develop the composition operation for tangles in
Section~\ref{sec_composition} below.
\end{proof}

\begin{rem}\hfill
\begin{enumerate}
\item
  While the commutative cube $Z(S_{T,\epsilon})$ did not depend on the
  orientation of the tangle, the tangle complex
  ${\tangle{T}}_{\epsilon,Z}$ does depend on it via the grading shift
  by $2n_+-4n_-$ and the homological degree shift by $-n_-$. These
  shifts correspond to the prefactor $(-1)^{n_-}A^{2n_+-4n_-}$,
  $q=A^2$, which turns the Kauffman bracket into an invariant under
  ambient isotopy, and they are needed to prove
  Theorem~\ref{thm_reidemeister}.
\item
  Since we have been able to translate from topology to algebra right
  in the beginning and since $\cat{grdVect}_k$ and
  $\cat{fltVect}_k$ are abelian categories, there is no
  need to consider formal linear combinations of surfaces as
  in~\cite{BN2}.
\item
  In order to obtain the case in which there is no filtration, one can
  simply define the grading shift to be the identity, use the trivial
  filtration $V=F^0V\supseteq F^1V=\{0\}$ for all vector spaces $V$,
  and omit the assumption of an Euler-filtration for the TQFT.
\end{enumerate}
\end{rem}

Denote by $\cal{H}^r(C)$ the $r$-th cohomology of the tangle complex
$C:=\llbracket T \rrbracket_{\varepsilon, Z}$. If $Z$ is an
Euler-filtered open-closed TQFT, then the filtration $F^*$ on the
complex $C$ induces a filtration on each homology group given by the
image of the inclusion
\begin{equation}
F^k\cal{H}^r(C):= {\rm image}\big(\cal{H}^r(F^kC) \hookrightarrow
\cal{H}^r(C) \big)
\end{equation}
where $F^kC:=\bigoplus_{r\in\Z}F^kC^r$. The filtration defines the
following bigrading on homology:
\begin{equation}
  \cal{H}^{k,r}(C):=F^k\cal{H}^r(C)/F^{k+1}\cal{H}^r(C).
\end{equation}

\begin{cor} \label{cor_poly}
If $Z$ is an Euler-filtered open-closed TQFT that satisfies
Bar-Natan's conditions, then the 2-variable tangle polynomial
${\cal{P}(T)}_{\epsilon,Z}\in\Z[t,t^{-1},A,A^{-1}]$ defined as the
filtered Poincar\'{e} polynomial of the complex $C:=\llbracket
T\rrbracket_{\varepsilon, Z}$, \ie
\begin{equation}
  {\cal{P}(T)}_{\varepsilon, Z} :=
  \sum_{r,k\in\Z} t^{r}A^{k}\rk \cal{H}^{k,r}(C),
\end{equation}
is an invariant of the tangle $T$.
\end{cor}

\begin{proof}
This follows from Theorem~\ref{thm_reidemeister}.
\end{proof}

\begin{prop}
\label{prop_actions}
Let $T$ be a plane diagram of an oriented $(p,q)$-tangle, \ie\ its
source and its target are a disjoint union of $p$ and $q$ points,
respectively. We define $r\in\N_0$ in such a way that $2r=p+q$. Let
$\epsilon\in\{-1,+1\}$ denote a checkerboard colouring and let
$Z\colon\cat{2Cob}^\mathrm{ext}\to\cat{Vect}_k$ be an open-closed TQFT
with associated knowledgeable Frobenius algebra
$(A,C,\imath,\imath^\ast)$. Then the tangle complex
${\tangle{T}}_{\epsilon,Z}$ is a complex of $(A^{\otimes r},A^{\otimes
r})$-bimodules.
\end{prop}

\begin{proof}
Let $T$ have $n$ crossings. We show that the vertices of the
$\cat{n}$-cube $Z(S_{T,\epsilon})$ of Definition~\ref{def_vectcube}
are $(A^{\otimes r},A^{\otimes r})$-bimodules and that its edges are
bimodule morphisms. The definition of the total complex which is used
in~\eqref{eq_tanglecomplex} then implies that the terms of
${\tangle{T}}_{\epsilon,Z}$ are $(A^{\otimes r},A^{\otimes
r})$-bimodules and that the differential consists of bimodule
morphisms.

We use the notation of Definition~\ref{def_tanglecube}. The vertex
$\alpha\subseteq\cat{n}$ of $Z(S_{T,\epsilon})$ is the vector space
\begin{equation}
  Z(m^{(\alpha)}) = A^{(m_1^{(\alpha)})}\otimes\cdots\otimes
  A^{(m_{k_\alpha}^{(\alpha)})},
\end{equation}
where $A^{(0)}:=C$ and $A^{(1)}:=A$. Obviously
$r=|\{\,j\in\{1,\ldots,k_\alpha\}\mid\,m_j^{(\alpha)}=1\,\}|$, \ie\
for every tensor factor of $Z(m^{(\alpha)})$ that is an $A$ as opposed
to a $C$, the tangle has two boundary points.

By the checkerboard colouring, $r$ out of the $p+q$ boundary points
are '$-$'-boundaries of the oriented smoothing
$D_\alpha^{(\epsilon)}\subseteq\R^2$, \ie\ they coincide with the tail
of an arrow in the example~\eqref{eq_orientedsmoothing}. The other $r$
boundary points are `$+$'-boundaries and coincide with the head of an
arrow. Note that this partitioning of the boundary points of $T$ into
`$-$'- versus `$+$'-boundaries depends only on the checkerboard
colouring, but that it is independent of the vertex
$\alpha\subseteq\cat{n}$.

We number the `$-$'-boundary points by $1,\ldots,r$. For each diagram
$D_\alpha^{(\epsilon)}$, $\alpha\subseteq\cat{n}$, the boundary point
number $j\in\{1,\ldots,r\}$ is contained in an arc that forms a
component of $D_\alpha^{(\epsilon)}$. We denote by
$s_j^{(\alpha,\epsilon)}\in\{1,\ldots,k_\alpha\}$ the number
corresponding to the tensor factor of this component. We have in
particular $m_{s_j^{(\alpha,\epsilon)}}^{(\alpha)}=1$. Similarly, we
number the `$+$'-boundaries by $1,\ldots,r$. For each
$\alpha\subseteq\cat{n}$, the `$+'$-boundary number
$j\in\{1,\ldots,r\}$ is contained in a component of
$D_\alpha^{(\epsilon)}$. We denote its number by
$t_j^{(\alpha,\epsilon)}\in\{1,\ldots,k_\alpha\}$.

Then $Z(m^{(\alpha)})$ forms an left $A^{\otimes r}$-module where the
$j$-th tensor factor of $A^{\otimes r}$ acts by left-multiplication on
the $s_j^{(\alpha,\epsilon)}$-th tensor factor of
$Z(m^{(\alpha)})$. Similarly, $Z(m^{(\alpha)})$ forms a right
$A^{\otimes r}$-module where the $j$-th tensor factor of $A^{\otimes
r}$ acts by right-multiplication on the $t_j^{(\alpha,\epsilon)}$-th
tensor factor of $Z(m^{(\alpha)})$. By associativity of $A$,
$Z(m^{(\alpha)})$ forms an $(A^{\otimes r},A^{\otimes r})$-bimodule.

For instance, the crossing of~\eqref{eq_crossing} has two
`$+$'-boundary points. Say, number $1$ at the top left and number $2$
at the bottom right of the diagram. We have $\cat{n}={1}$, and we
denote by $D_\emptyset^{(\epsilon)}$ and $D_{\{1\}}^{(\epsilon)}$ the
two smoothings displayed in~\eqref{eq_crossing}. Then we have
$t_1^{(\emptyset,\epsilon)}=1$, $t_2^{(\emptyset,\epsilon)}=2$, and
$t_1^{(\{1\},\epsilon)}=2$, $t_2^{(\{1\},\epsilon)}=1$.

From the definition of the open-closed cobordisms, see, for example
the cobordisms~\eqref{eq_crossing} to~\eqref{eq_copants}, it is
obvious that the edges $Z(S_{(\alpha,\ell)})$, $\ell\in\cat{n}$,
$\alpha\subseteq\cat{n}\backslash\{\ell\}$, of the cube form
morphisms of $(A^{\otimes r},A^{\otimes r})$-bimodules. Every
`$-$'-boundary point $j\in\{1,\ldots,r\}$ corresponds to one
component of the coloured boundary $\del_1S_{(\alpha,\ell)}$ of the
open-closed cobordism. This component intersects both the source and
the target of $S_{(\alpha,\ell)}$ in one point
(Remark~\ref{rem_intersect}). It intersects the
$s_j^{(\alpha,\epsilon)}$-th component of the source (which is an
arc) and the $s_j^{(\alpha\sqcup\{\ell\},\epsilon)}$-th component of
the target (which is an arc, too). The map $Z(S_{(\alpha,j)})$
therefore commutes with the $A$-left actions defined above due to
the relations that hold in any open-closed TQFT. An analogous result
holds for `$+$'-boundary points if one replaces
$s_j^{(\alpha,\epsilon)}$ and
$s_j^{(\alpha\sqcup\{\ell\},\epsilon)}$ by $t_j^{(\alpha,\epsilon)}$
and $t_j^{(\alpha\sqcup\{\ell\},\epsilon)}$, respectively.

Let us illustrate this result for the example of the crossing
of~\eqref{eq_crossing} and the `$+$'-boundary point number $1$ at the
top left of the diagram, the $A$-right action on the
$t_1^{(\emptyset,\epsilon)}=1$-st tensor factor of
$Z(m^{(\emptyset)})=A\otimes A$ is given by
\begin{equation}
Z\biggl(
\begin{aligned}
\psset{xunit=3mm,yunit=3mm}
\begin{pspicture}(5,5)
  \rput(1.4,0){\multl}
  \rput(4.4,0){\identl}
  \rput(0.4,2.5){\identl}
  \rput(3.4,2.5){\crossl}
\end{pspicture}
\end{aligned}
\biggr)\colon (A\otimes A)\otimes A\to A\otimes A,
\end{equation}
whereas the $A$-right action on the $t_1^{(\{1\},\epsilon)}=2$-nd
tensor factor of $Z(m^{(\{1\})})=A\otimes A$ is given by
\begin{equation}
Z\bigl(
\begin{aligned}
\psset{xunit=3mm,yunit=3mm}
\begin{pspicture}(5,2.5)
  \rput(0.4,0){\identl}
  \rput(3.4,0){\multl}
\end{pspicture}
\end{aligned}
\bigr)\colon (A\otimes A)\otimes A\to A\otimes A.
\end{equation}
The $A$-right action commutes with the differential of the complex
because of the following relation in $\cat{2Cob}^\mathrm{ext}$:
\begin{equation}
\begin{aligned}
\psset{xunit=3mm,yunit=3mm}
\begin{pspicture}(5,10)
  \rput(1,0){\saddler}
  \rput(1,5){\multl}
  \rput(5,5){\curveleftl}
  \rput(0,7.5){\identl}
  \rput(3,7.5){\crossl}
\end{pspicture}
\end{aligned}
\qquad\cong\qquad
\begin{aligned}
\psset{xunit=3mm,yunit=3mm}
\begin{pspicture}(6,7.5)
  \rput(0,0){\identl}
  \rput(5,0){\multl}
  \rput(0,2.5){\saddler}
  \rput(6,2.5){\identl}
  \rput(6,5){\identl}
\end{pspicture}
\end{aligned}
\end{equation}
One can literally slide the $A$-right action down along the
component of $\del_1 M$ that corresponds to the `$+$'-boundary
number $1$. This is the component that connects the right boundary
point of the first component of the source with the right boundary
point of the second component of the target of the open-closed
cobordism of~\eqref{eq_crossing}.
\end{proof}

\begin{rem}
\label{rem_tanglemodule}
For a $(p,q)$-tangle, it is more natural to have $p$ copies of $A$
act from the left and $q$ copies from the right. This can be
achieved by replacing $A$-right actions by $A^\mathrm{op}$-left
actions and vice versa.

We define for $p\in\N_0$,
\begin{equation}
  A^{(p,+1)}:=\underbrace{A\otimes A^\mathrm{op}\otimes A\otimes A^\mathrm{op}\otimes\cdots}_{p},
\end{equation}
and
\begin{equation}
  A^{(p,-1)}:=\underbrace{A^\mathrm{op}\otimes A\otimes A^\mathrm{op}\otimes A\otimes\cdots}_{p}.
\end{equation}
Then the tangle complex ${\tangle{T}}_{Z,\epsilon}$ of a
$(p,q)$-tangle forms a complex of
$(A^{(q,\epsilon)},A^{(p,\epsilon)})$-bimodules.
\end{rem}

\subsection{Obtaining Khovanov's complex from Lee's and Bar-Natan's}
\label{sec_KhfromLee}

In this section, we show in which cases one can compute Khovanov's
chain complex from Lee's or Bar-Natan's, respectively, by exploiting
their filtration and passing to the associated graded complex
(Section~\ref{sect_spectral}). For links, Lee and
Rasmussen~\cite{Lee,Rasmussen} have shown that the filtration on the
chain complex constructed with Lee's Frobenius algebra $C_{0,1}$ gives
rise to a spectral sequence whose $E_0$-page, \ie\ the associated
graded complex, is Khovanov's chain complex, and whose $E_1$-page is
Khovanov's link homology. All pages $E_j$, $j\geq 1$, are invariant
under Reidemeister moves up to isomorphism. Turner~\cite{Turner} has
shown the corresponding result for the chain complex constructed with
Bar-Natan's Frobenius algebra $C_{1,0}$. The following theorem
generalizes these results from links to tangles.

\begin{thm}
\label{thm_spectral} Let $T$ be a plane diagram of an oriented
tangle and $\epsilon\in\{-1,+1\}$ be a checkerboard colouring. Let
$\mathbbm{A}=(A,C,\imath,\imath^\ast)$ be an Euler-filtered
knowledgeable Frobenius algebra and
$\mathbbm{A}^\prime=(A^\prime,C^\prime,\imath^\prime,{\imath^\prime}^\ast)$
be its associated Euler-graded knowledgeable Frobenius algebra (as
in Remark~\ref{rem_spectral}) such that $C^\prime=C_{0,0}$ agrees
with Khovanov's Frobenius algebra. Denote by $Z$ and $Z^\prime$ the
open-closed TQFTs associated with $\mathbbm{A}$ and
$\mathbbm{A}^\prime$, respectively.

The filtration on the tangle complex ${\tangle{T}}_{\epsilon,Z}$ gives
rise to a spectral sequence whose $E_0$-page is the graded complex
${\tangle{T}}_{\epsilon,Z^\prime}$. All pages $E_j$, $j\geq 1$, of
this spectral sequence are invariant under Reidemeister moves up to
isomorphism.
\end{thm}

\begin{rem}
Examples for such pairs $(\mathbbm{A},\mathbbm{A}^\prime)$ of
knowledgeable Frobenius algebras are given by
Examples~\ref{ex_barnatan} and~\ref{ex_khovanov}, by
Examples~\ref{ex_filtered3} and~\ref{ex_modp}, or by any of the
algebras of Theorem~\ref{thm_sskfrob}(2.) and its associated
Euler-graded knowledgeable Frobenius algebra.
\end{rem}

\noindent
For the proof of the theorem, we need one lemma. In order to be able
to state it, we recall the following facts, see, for
example~\cite{McCleary}: Let $C$ and $C^\prime$ be filtered complexes
and $f\colon C\to C^\prime$ be a morphism of filtered complexes, \ie\
a morphism of complexes each of whose components is filtered of degree
$0$. Denote by ${\{E_n\}}_n$ and ${\{E^\prime_n\}}_n$ the spectral
sequences induced by the filtrations of $C$ and $C^\prime$,
respectively. Then the map $f$ induces a morphism of spectral
sequences ${\{f_n\}}_n\colon{\{E_n\}}_n\to{\{E^\prime_n\}}_n$. This
means that for each page $n\geq 0$, there is a morphism of bigraded
differential modules $f_n\colon E_n\to E_n^\prime$, \ie\ a graded
linear map of bidegree $(0,0)$ such that $f_n\circ d_n=d_n^\prime\circ
f_n$, and it means that each $f_{n+1}$ is the map that $f_n$ induces
on homology.

\begin{lem}[Theorem~3.4 of~\cite{McCleary}]
\label{lem_spectral}
Let $C$ and $C^\prime$ be filtered complexes and $f\colon C\to
C^\prime$ be a morphism of filtered complexes. Denote the induced
morphism of spectral sequences by
${\{f_n\}}_n\colon{\{E_n\}}_n\to{\{E^\prime_n\}}_n$. If $f_m\colon
E_m\to E^\prime_m$ is an isomorphism of bigraded differential modules
for some $m$, then so is $f_j$ for any $j\geq m$.
\end{lem}

Notice that it suffices to show that $f_m$ is a linear isomorphism.
It is already linear and graded of bidegree $(0,0)$ and compatible
with the differential since it is part of a morphism of spectral
sequences.

\begin{proof}[Proof of Theorem~\ref{thm_spectral}]
Consider the tangle complex ${\tangle{T}}_{\epsilon,Z}$ with its
filtration. We first show that the associated graded complex is
${\tangle{T}}_{\epsilon,Z^\prime}$. The knowledgeable Frobenius
algebras $\mathbbm{A}$ and $\mathbbm{A}^\prime$ have the same
underlying graded vector spaces, and so both
${\tangle{T}}_{\epsilon,Z}$ and ${\tangle{T}}_{\epsilon,Z^\prime}$
have the same terms. We now compute the differential of the associated
graded complex of the filtered complex ${\tangle{T}}_{\epsilon,Z}$ and
show that it agrees with the differential of
${\tangle{T}}_{\epsilon,Z^\prime}$.

From the definition of the total complex $C_\mathrm{tot}$
(Definition~\ref{def_totalcomplex}), we see that we can do this for
each term of the direct sum and thereby for each smoothing
$\alpha\subseteq\cat{n}$ separately. The differential is in each case
either $+1$ or $-1$ times the linear map associated with one of the
open-closed cobordisms of~\eqref{eq_crossing}
to~\eqref{eq_copants}. Using the knowledgeable Frobenius algebra
$\mathbbm{A}=(A,C,\imath,\imath^\ast)$, these linear maps are (up to
permutations of the tensor factors and up to tensoring with
identities):
\begin{gather}
  (\mu_A\otimes\id_A)\circ(\id_A\otimes\tau_{A,A})\circ(\Delta_A\otimes\id_A),\\
  (\id_A\otimes\mu_A)\circ(\tau_{A,A}\otimes\id_A)\circ(\id_A\otimes\Delta_A),\\
  (\imath^\ast\otimes\id_A)\circ\Delta_A,\\
  \mu_A\circ(\imath\otimes\id_A),\\
  \mu_C,\\
  \Delta_C.
\end{gather}
For each of these compositions, the associated graded map agrees with
the same expression in the knowledgeable Frobenius algebra
$\mathbbm{A}^\prime$. This shows that the associated graded complex
and thereby the $E_0$-page of the spectral sequence is
${\tangle{T}}_{\epsilon,Z^\prime}$.

Now we show that all pages $E_j$, $j\geq 1$, of this spectral sequence
are invariant under Reidemeister moves up to isomorphism.

Let $T$ and $T^\prime$ be plane diagrams of oriented tangles that are
related by a Reidemeister move, let $\epsilon\in\{-1,+1\}$, and let
$C:={\tangle{T}}_{\epsilon,Z}$ and
$C^\prime:={\tangle{T^\prime}}_{\epsilon,Z}$ be the tangle
complexes. These are filtered complexes. We denote by ${\{E_n\}}_n$
and ${\{E^\prime_n\}}_n$ the induced spectral sequences.

The invariance under Reidemeister moves of the homology of $C$ up to
isomorphism is established in the proof of
Theorem~\ref{thm_reidemeister} by constructing chain maps $f\colon
C\to C^\prime$ and $g\colon C^\prime\to C$ which form a homotopy
equivalence. These chain maps are constructed from open-closed
cobordisms by applying $Z$. The cobordisms are displayed, for example,
in~\cite{BN} and in Section~\ref{sec_composition} below.

First, both $f$ and $g$ are homomorphisms of filtered complexes.
This can be seen in a case by case inspection, taking all the
grading shifts in the definition~\eqref{eq_tanglecomplex} into
account. Denote by
${\{f_n\}}_n\colon{\{E_n\}}_n\to{\{E_n^\prime\}}_n$ and by
${\{g_n\}}_n\colon{\{E_n^\prime\}}_n\to{\{E_n\}}_n$ the induced
morphisms of spectral sequences.

Second, on the $E_0$-pages, the $f_0\colon E_0\to E^\prime_0$ and
$g_0\colon E^\prime_0\to E_0$ are the associated graded maps of
$f\colon C\to C^\prime$ and $g\colon C^\prime\to C$, respectively.  By
a similar argument as above, $f_0$ agrees with the chain map that is
obtained from the same open-closed cobordism as $f$, but using the
open-closed TQFT $Z^\prime$ rather than $Z$, and similarly for
$g_0$. The morphisms of differential bigraded modules $f_0\colon
E_0\to E^\prime_0$ and $g_0\colon E^\prime_0\to E_0$ therefore form a
homotopy equivalence between the associated graded complexes. This is
the same homotopy equivalence that appears in the proof of
Reidemeister move invariance for $Z^\prime$ instead of $Z$.  Therefore
both maps induce mutually inverse maps on homology. The homology of
the associated graded complexes, however, is contained in the
$E_1$-pages (\cf~\eqref{eq_gradedhomology}), and therefore the
morphisms of bigraded differential modules $f_1\colon E_1\to
E^\prime_1$ and $g_1\colon E^\prime_1\to E_1$ are mutually inverse.

By Lemma~\ref{lem_spectral}, this implies that each $f_j\colon E_j\to
E^\prime_j$, $j\geq1$, is an isomorphism.
\end{proof}

\begin{rem}\hfill
\label{rem_chartwo}
\begin{enumerate}
\item
  The generalization from links to tangles of Lee's~\cite{Lee} and
  Turner's~\cite{Turner} `endomorphism of the Khovanov invariant' is
  contained in the spectral sequence, namely in the differentials
  $d_j\colon E_j\to E_j^\prime$ of the pages $j\geq 1$.
\item
  Already the simplest example
  $\mathbbm{A}=(A,C,\imath,\imath^\ast)$ of Example~\ref{ex_barnatan}
  for a field of characteristic $2$ is quite powerful. Firstly, by
  the above theorem, it is rich enough to recover Khovanov homology in
  characteristic $2$. Secondly, on links, it agrees with Bar-Natan's
  chain complex in characteristic $2$ which is, by a result of
  Mackaay, Tuner, and Vaz~\cite{MTV}, sufficient to compute
  Rasmussen's $s$-invariant~\cite{Rasmussen}. This invariant is strong
  enough to tell that some knots which are known to be topologically
  slice, are \emph{not} smoothly slice~\cite{Rasmussen2,MTV}. Finally,
  $\mathbbm{A}$ has a strongly separable algebra $A$, and so all of
  the tangle composition technology developed in the following
  section, will be available, too.
\item
 In Section~4 below, we show how one can compute the tangle
 complex for Lee's and Bar-Natan's Frobenius algebras by composing
 arcs and crossings and that one can then in principle compute the
 Khovanov version of the tangle complex from
 Theorem~\ref{thm_spectral}. In order to prove
 Theorem~\ref{thm_spectral}, however, one still needs
 Theorem~\ref{thm_reidemeister}, \ie\ the \emph{global} version of the
 proof of Reidemeister move invariance. So although some computations
 can be done entirely \emph{locally}, the proofs that they are valid
 cannot.
\end{enumerate}
\end{rem}

\subsection{Examples}

In this section, we show that even the tangle complex defined by the
open-closed TQFT $Z_\mathrm{BN}$ associated with the knowledgeable
Frobenius algebra
$\mathbbm{A}_\mathrm{BN}=(A_{1,0},C_{1,0},\imath,\imath^\ast)$ of
Example~\ref{ex_barnatan} in characteristic $2$ for $t=0$ detects
whether the two strands of a $(2,2)$-tangle are non-trivially braided
or not. We therefore compute the tangle homology and the tangle
polynomials of the complexes
\begin{equation}
\llbracket T \rrbracket :=
 \llbracket \;
\begin{pspicture}[.4](2,2.1)
  \rput(-.5,0){ \pspolygon[linecolor=lightgray,
   linewidth=.6pt,fillstyle=solid,fillcolor=lightgray](.5,0)(.5,2.03)(2.5,2.03)(2.5,0)(.5,0)
  \pscustom[linecolor=white,linewidth=0pt,fillstyle=solid,fillcolor=white]{
  \psbezier[linewidth=.8pt](2,0)(2,.5)(1,.5)(1,1)
   \psbezier[linewidth=.8pt](1,1)(1,1.5)(2,1.5)(2,2)
   \psline(1,2)
   \psbezier[linewidth=.8pt](1,2)(1,1.5)(2,1.5)(2,1)
   \psbezier[linewidth=.8pt](2,1)(2,.5)(1,.5)(1,0)
   \psline(2,0)
  }}
 \psbezier[linewidth=.8pt](1.5,0)(1.5,.5)(.5,.5)(.5,1)
 \psbezier[linewidth=.8pt]{<-}(1.5,2)(1.5,1.5)(.5,1.5)(.5,1)
 \pspolygon[linecolor=white,fillstyle=solid,fillcolor=white](.85,0)(.85,2)(1.15,2)(1.15,0)(.85,0)
  \psbezier[linewidth=.8pt]{<-}(.5,2)(.5,1.5)(1.5,1.5)(1.5,1)
  \psbezier[linewidth=.8pt](.5,0)(.5,.5)(1.5,.5)(1.5,1)
\end{pspicture} \;
\rrbracket\quad {\rm and} \quad \llbracket T' \rrbracket :=
 \llbracket \;
\begin{pspicture}[.4](2,2.1)
  \rput(-.5,0){ \pspolygon[linecolor=lightgray,
   linewidth=.6pt,fillstyle=solid,fillcolor=lightgray](.5,0)(.5,2.03)(2.5,2.03)(2.5,0)(.5,0)
  \pscustom[linecolor=white,linewidth=0pt,fillstyle=solid,fillcolor=white]{
  \psbezier[linewidth=.8pt](2,0)(2,.5)(1,.5)(1,1)
   \psbezier[linewidth=.8pt](1,1)(1,1.5)(2,1.5)(2,2)
   \psline(1,2)
   \psbezier[linewidth=.8pt](1,2)(1,1.5)(2,1.5)(2,1)
   \psbezier[linewidth=.8pt](2,1)(2,.5)(1,.5)(1,0)
   \psline(2,0)
  }}
 \psbezier[linewidth=.8pt](1.5,0)(1.5,.5)(.5,.5)(.5,1)
 \psbezier[linewidth=.8pt]{<-}(.5,2)(.5,1.5)(1.5,1.5)(1.5,1)
 \pspolygon[linecolor=white,fillstyle=solid,fillcolor=white](.85,0)(.85,2)(1.15,2)(1.15,0)(.85,0)
  \psbezier[linewidth=.8pt]{<-}(1.5,2)(1.5,1.5)(.5,1.5)(.5,1)
  \psbezier[linewidth=.8pt](.5,0)(.5,.5)(1.5,.5)(1.5,1)
\end{pspicture} \;
\rrbracket.
\end{equation}

By the up to homotopy invariance of the complex $C:=\llbracket T
\rrbracket$ under Reidemeister move two, the first tangle only has
nontrivial homology in cohomological degree zero
 \begin{equation}
\cal{H}^0(C)=(A \otimes A)\{0\},
 \end{equation}
and so
\begin{equation}
 \cal{P}_{Z_\mathrm{BN}}(T) = t^0(A^2 +A^{-2}) = A^2
 +A^{-2}.
\end{equation}

The second complex $C^\prime:=\llbracket T' \rrbracket$ is the filtered complex
\begin{equation}
  \xymatrix{
  0 \ar[r]
  & (A \otimes A)\{4\}
    \ar[r]^-{d^0}
  & (A \otimes A)\{6\} \oplus (A \otimes A)\{6\}
    \ar[r]^-{d^1}
  & (A \otimes C \otimes A)\{8\}
  \ar[r] & 0
  }
\end{equation}
where
\begin{eqnarray}
  d^0 &=& \begin{pmatrix}
            (\id_A\otimes\mu_A)\circ(\tau_{A,A}\otimes\id_A)\circ(\id_A\otimes\Delta_A)\\
            (\id_A\otimes\mu_A)\circ(\tau_{A,A}\otimes\id_A)\circ(\id_A\otimes\Delta_A)
          \end{pmatrix},\\
  d^1 &=& \bigl(((\id_A\otimes\imath^*)\circ\Delta_A)\otimes\id_A,
                \id_A\otimes((\imath^*\otimes\id_A)\circ\Delta_A)\bigr).
\end{eqnarray}
The kernel of $d^0$ is the filtered vector space $\langle1\otimes y +
y \otimes 1, 1 \otimes y + y \otimes y\rangle\{4\}$. A simple
calculation reveals that
\begin{equation}
  \im d^0 = \ker d^1 =
\langle (1 \otimes 1 + 1\otimes y + y \otimes 1) \oplus (1 \otimes 1
+ 1\otimes y + y \otimes 1) , (y \otimes y) \oplus (y \otimes y)
\rangle\{4\}.
\end{equation}
Furthermore,
\begin{eqnarray}
 \im d^1 &=& \langle
 1 \otimes x \otimes y,  y \otimes x \otimes 1,
 y \otimes x \otimes y,
 1 \otimes 1 \otimes 1+1 \otimes x \otimes 1+y \otimes 1 \otimes 1,
\\&&
 1 \otimes 1 \otimes 1+1 \otimes x \otimes 1+1 \otimes 1 \otimes y,
 1 \otimes 1 \otimes y + y \otimes 1 \otimes y
\rangle\{8\}.
\end{eqnarray}
Therefore,
\begin{equation}
  \cal{H}^n\left(C'\right)
  = \left\{
  \begin{array}{ll}
  \langle 1\otimes y + y \otimes 1, 1 \otimes y + y \otimes
  y\rangle \{4\} & {\rm if} \; n=0 , \\
  \langle 1 \otimes 1 \otimes 1, 1\otimes x \otimes 1 \rangle\{8\} & {\rm if} \; n=2 ,\\
  0 &{\rm otherwise}.
  \end{array}
  \right.
\end{equation}
Taking the grading shifts into account, we obtain the following
bigrading for this homology:
\begin{eqnarray}
  \cal{H}^{2,0}(C^\prime) &=& \left<1\otimes y+y\otimes y\right>,\\
  \cal{H}^{4,0}(C^\prime) &=& \left<1\otimes y+y\otimes 1\right>,\\
  \cal{H}^{8,2}(C^\prime) &=& \left<1\otimes x\otimes 1\right>,\\
  \cal{H}^{12,2}(C^\prime)&=& \left<1\otimes 1\otimes 1\right>.
\end{eqnarray}
Hence, the tangle polynomial for $T'$ is
\begin{equation}
\cal{P}_{Z_\mathrm{BN}}(T^\prime) = t^0(A^2+A^4)+t^2(A^{8}+A^{12}) =
A^2+A^4+t^2A^{8}+t^2A^{12}.
\end{equation}

Computing the homology of the associated graded complex for
$\llbracket T' \rrbracket$, i.e. using the open-closed TQFT
$Z_{\mathrm{Kh}}$ associated with the knowledgeable Frobenius algebra
of Example~\ref{ex_khovanov}, we get
\begin{eqnarray}
  \cal{H}^{2,0}(C^\prime) &=& \left<y\otimes y\right>,\\
  \cal{H}^{4,0}(C^\prime) &=& \left<1\otimes y+y\otimes 1\right>,\\
  \cal{H}^{6,1}(C^\prime) &=& \left<
    \begin{pmatrix} 1\otimes y\\ 1\otimes y\end{pmatrix}\right>,\\
  \cal{H}^{8,1}(C^\prime) &=& \left<
    \begin{pmatrix} 1\otimes 1\\ 1\otimes 1\end{pmatrix}\right>,\\
  \cal{H}^{8,2}(C^\prime) &=& \left<y\otimes 1\otimes y\right>,\\
  \cal{H}^{10,2}(C^\prime)&=& \left<1\otimes 1\otimes y,y\otimes 1\otimes 1\right>,\\
  \cal{H}^{12,2}(C^\prime)&=& \left<1\otimes 1\otimes 1\right>,
\end{eqnarray}
so that
\begin{eqnarray}
\cal{P}_{Z_\mathrm{Kh}}(T') &=&
A^4+A^2+t^1(A^8+q^6)+t^2(A^{12}+2A^{10}+A^8).
\end{eqnarray}

%
\section{Composition of tangles}
%
\label{sec_composition}

\subsection{Tensor products of complexes of bimodules}

Let us first introduce an efficient way of describing tensor products
of complexes whose terms are bimodules.

Let $A$ and $B$ be $k$-algebras and $C=(C^i,d^i)$ be a complex of
$(A,B)$-bimodules, \ie\ an object of $\cat{Kom}({}_A\mathcal{M}_B)$.
Denote the $(A,B)$-bimodule structure of each $C^i$ by
$\lambda^i\colon A\otimes C^i\to C^i$ and $\rho^i\colon C^i\otimes
B\to C^i$. If we view $A$ as a $1$-term complex $0\rightarrow
A\rightarrow 0$, then the tensor product of complexes $A\otimes C$ has
the terms ${(A\otimes C)}^i=A\otimes C^i$ and the differential
$d^i_{A\otimes C}=\id_A\otimes d^i$. The ${\{\lambda^i\}}_i$ form a
morphism of complexes $\lambda\colon A\otimes C\to C$, \ie\ a morphism
of $\cat{Kom}(\cat{Vect}_k)$.  Similarly, there is a morphism of
complexes $\rho\colon C\otimes B\to C$.

If both $C$ and $D$ are complexes of $(A,B)$-bimodules with
$\lambda_C\colon A\otimes C\to C$, $\rho_C\colon C\otimes B\to C$,
$\lambda_D\colon A\otimes D\to D$, and $\rho_D\colon D\otimes B\to D$
defined as above, and if $f\colon C\to D$ is a morphism of such
complexes, then we have $f\circ\lambda_C = \lambda_D\circ(\id_A\otimes
f)$ and $f\circ\rho_C=\rho_D\circ(f\otimes\id_B)$.

\begin{lem}
Let $A$, $B$, and $B^\prime$ be $k$-algebras and $C=(C^i,d^i)$ be a
complex of $(A\otimes B^\prime,A\otimes
B^{\prime\prime})$-bimodules, \ie\ an object of
$\cat{Kom}({}_{A\otimes B^\prime}\mathcal{M}_{A\otimes
B^{\prime\prime}})$. Then the following coequalizer in
$\cat{Kom}({}_{B^\prime}\mathcal{M}_{B^{\prime\prime}})$,
\begin{equation}
\label{eq_coeqcomplex}
\xymatrix{
  A\otimes C\ar@<.5ex>[rr]^{\lambda}\ar@<-.5ex>[rr]_{\rho\circ\tau_{A,C}}&&
  C\ar[rr]^{{}_\rho\otimes_\lambda}&&
  {}_\rho C_\lambda
}
\end{equation}
exists. We have denoted by $\lambda\colon A\otimes C\to C$ and
$\rho\colon C\otimes A\to C$ the maps that make $C$ a complex of
$(A,A)$-bimodules. These are morphisms of
$\cat{Kom}({}_{B^\prime}\mathcal{M}_{B^{\prime\prime}})$.

The complex ${}_\rho C_\lambda=({}_\rho C_\lambda^i,d^i_{{}_\rho
C_\lambda})$ is the complex whose terms are the coequalizers
\begin{equation}
\label{eq_coeqleveli}
\xymatrix{
  A\otimes C^i\ar@<.5ex>[rr]^{\lambda^i}\ar@<-.5ex>[rr]_{\rho^i\circ\tau_{A,C}}&&
  C^i\ar[rr]^{{}_\rho\otimes_\lambda^i}&&
  {}_\rho C_\lambda^i
}
\end{equation}
in the category ${}_{B^\prime}\mathcal{M}_{B^{\prime\prime}}$ and
whose differential $d^i_{{}_\rho C_\lambda}$ is uniquely specified as
the morphism in the commutative diagram
\begin{equation}
\xymatrix{
  C^i\ar[rr]^{d^i}\ar[dd]_{{}_\rho\otimes_\lambda^i}&&
    C^{i+1}\ar[dd]^{{}_\rho\otimes_\lambda^{i+1}}\\
  \\
  {}_\rho C_\lambda^i\ar[rr]_{d^i_{{}_\rho C_\lambda}}&&
    {}_\rho C_\lambda^{i+1}
}
\end{equation}
because by the universal property,
${}_\rho\otimes_\lambda^{i+1}\circ d^i$ factors uniquely through the
coequalizer ${}_\rho\otimes_\lambda^i$. The coequalizer map
${}_\rho\otimes_\lambda$ in~\eqref{eq_coeqcomplex} is the morphism
of $\cat{Kom}({}_{B^\prime}\mathcal{M}_{B^{\prime\prime}})$ whose
components are the coequalizers ${}_\rho\otimes_\lambda^i$
of~\eqref{eq_coeqleveli}.
\end{lem}

\subsection{Construction of the tangle complex from constituents}

In this section, we show how one can obtain the tangle complex of
Definition~\ref{def_tanglecomplex} from the iterated tensor product
over $A$ of two basic building blocks: A $1$-term complex
$0\rightarrow A\rightarrow 0$ associated with arcs and a $2$-term
complex $0\rightarrow A\otimes A\rightarrow A\otimes A\rightarrow 0$
associated with crossings which comes in four different variations,
depending on the orientations and on the checkerboard colouring. The
requirement is that we employ a state sum knowledgeable Frobenius
algebra $(A,Z(A),\imath,\imath^\ast)$ that satisfies Bar-Natan's
conditions. Warning: In this section, $C$ denotes a complex whereas
$Z(A)$ is the relevant commutative Frobenius algebra.

\begin{defn}
\label{def_composed}
Let $T$ be a plane diagram of an oriented $(p,q)$-tangle with $n_+$
positive $(\overcrossing)$ and $n_-$ negative crossings
$(\undercrossing)$, $n:=n_++n_-$; let $\epsilon\in\{-1,+1\}$ denote a
checkerboard colouring and let
$\mathbbm{A}=(A,Z(A),\imath,\imath^\ast)$ be an [Euler-graded,
Euler-filtered] state sum knowledgeable Frobenius algebra that
satisfies Bar-Natan's conditions.

\noindent
Choose a finite graph $G\subseteq\R^2$ embedded in the plane such that
\begin{itemize}
\item
  the vertices of $G$ do not intersect $T$, and
\item
  the intersection $G\cap T$ is a finite set, and
\item
  the closure $\overline{R_j}$ of each connected component
  $R_1,\ldots,R_\ell$ of the complement $\R^2\backslash G$ intersects
  $T$ in such a way that $\overline{R_j}\cap T$ is either an arc $(|)$
  or a crossing $(\slashoverback)$.
\end{itemize}
Number the boundary points of the $\overline{R_j}\cap T$ by
$1,\ldots,k$, \ie\ $k=4n+2m$ where $m$ is the number of arcs among the
$\overline{R_j}\cap T$. For each $j\in\{1,\ldots,\ell\}$, we define a
complex $C_j$ as follows.

\begin{enumerate}
\item
  If $\overline{R_j}\cap T$ is an arc $(|)$, we set
\begin{equation}
\label{eq_complexone}
  C_j := \xymatrix{
0\ar[r]&\underline{A}\ar[r]&0
}
\end{equation}
  where we have underlined the term in cohomological degree $0$. Let
  $i_1,i_2\in\{1,\ldots,k\}$ be the numbers of the boundary points of
  the arc. One of them, say $i_1$, is a `$-$'-boundary relative to
  $\epsilon$ (\cf~Proposition~\ref{prop_actions}). The complex $C_j$
  above is a complex of left $A$-modules whose left $A$-action we
  denote by $\lambda^{(i_1)}\colon A\otimes C_j\to C_j$. The other
  point numbered $i_2$ is a `$+$'-boundary, and we denote the right
  $A$-action by $\rho^{(i_2)}\colon C_j\otimes A\to C_j$.
\item
  If $\overline{R_j}\cap T$ is a positive crossing with the
  checkerboard colouring
\begin{equation}
\begin{aligned}
\begin{pspicture}(1.5,1.5)
  \pscustom[fillstyle=solid,fillcolor=lightgray,linecolor=lightgray]{
    \psline(0.0,0.0)(1.5,1.5)(1.5,0.0)(0.0,1.5)(0.0,0.0)
  }
  \psline[linewidth=1pt]{<-}(0.0,1.5)(1.5,0.0)
  \psline[linewidth=3pt,linecolor=white](0.7,0.8)(0.8,0.7)
  \psline[linewidth=1pt]{->}(0.0,0.0)(1.5,1.5)
  \rput(-0.2,-0.2){$i_1$}
  \rput(1.7,-0.2){$i_2$}
  \rput(1.7,1.7){$i_3$}
  \rput(-0.2,1.7){$i_4$}
\end{pspicture}
\end{aligned}
\end{equation}
  whose boundary points are numbered
  $i_1,i_2,i_3,i_4\in\{1,\ldots,k\}$ as indicated,
  we define
\begin{equation}
\label{eq_complextwo}
  C_j:=\xymatrix{
0\ar[r]&\underline{A\otimes A\{2\}}\ar[rrrrr]^{(\mu_A\otimes\id_A)\circ(\id_A\otimes\tau_{A,A})\circ(\Delta_A\otimes\id_A)}&&&&&
  A\otimes A\{4\}\ar[r]&0
}.
\end{equation}
  This is a [graded, filtered] $2$-term complex of $(A\otimes
  A,A\otimes A)$-bimodules. We write
\begin{eqnarray}
  \lambda_1\colon A\otimes (A\otimes A)\to A\otimes A,\quad a\otimes (x\otimes y)&\mapsto&\mu_A(a\otimes x)\otimes y,\\
  \lambda_2\colon A\otimes (A\otimes A)\to A\otimes A,\quad a\otimes (x\otimes y)&\mapsto&x\otimes\mu_A(a\otimes y),\\
  \rho_1\colon (A\otimes A)\otimes A\to A\otimes A,\quad (x\otimes y)\otimes a&\mapsto& \mu_A(x\otimes a)\otimes y,\\
  \rho_2\colon (A\otimes A)\otimes A\to A\otimes A,\quad (x\otimes y)\otimes a&\mapsto& x\otimes\mu_A(y\otimes a)
\end{eqnarray}
  for the two left and the two right actions of $A$ on $A\otimes
  A$. The structure of $C_j$ as a complex of $(A\otimes A,A\otimes
  A)$-bimodules is then given by the following four actions. As above,
  `$-$'-boundaries relative to $\epsilon$ result in left actions and
  `$+$'-boundaries in right actions:
\begin{enumerate}
\item
    $\lambda^{(i_1)}\colon A\otimes C_j\to C_j$ which is defined by
    ${(\lambda^{(i_1)})}^0:=\lambda_1$ and
    ${(\lambda^{(i_1)})}^1:=\lambda_1$ in degree $0$ and $1$,
    respectively.
\item
    $\rho^{(i_2)}\colon C_j\otimes A\to C_j$ with ${(\rho^{(i_2)})}^0:=\rho_2$ and ${(\rho^{(i_2)})}^1:=\rho_1$.
\item
    $\lambda^{(i_3)}\colon A\otimes C_j\to C_j$ with ${(\lambda^{(i_3)})}^0:=\lambda_2$ and ${(\lambda^{(i_3)})}^1:=\lambda_2$.
\item
    $\rho^{(i_4)}\colon C_j\otimes A\to C_j$ with ${(\rho^{(i_4)})}^0:=\rho_1$ and ${(\rho^{(i_4)})}^1:=\rho_2$.
\end{enumerate}
\item
  If $\overline{R_j}\cap T$ is a negative crossing with the
  checkerboard colouring
\begin{equation}
\begin{aligned}
\begin{pspicture}(1.5,1.5)
  \pscustom[fillstyle=solid,fillcolor=lightgray,linecolor=lightgray]{
    \psline(0.0,0.0)(1.5,1.5)(1.5,0.0)(0.0,1.5)(0.0,0.0)
  }
  \psline[linewidth=1pt]{->}(0.0,0.0)(1.5,1.5)
  \psline[linewidth=3pt,linecolor=white](0.7,0.7)(0.8,0.8)
  \psline[linewidth=1pt]{<-}(0.0,1.5)(1.5,0.0)
  \rput(-0.2,-0.2){$i_1$}
  \rput(1.7,-0.2){$i_2$}
  \rput(1.7,1.7){$i_3$}
  \rput(-0.2,1.7){$i_4$}
\end{pspicture}
\end{aligned}
\end{equation}
  we define
\begin{equation}
\label{eq_complexthree}
  C_j:=\xymatrix{
0\ar[r]&A\otimes A\{-4\}\ar[rrrrr]^{(\mu_A\otimes\id_A)\circ(\id_A\otimes\tau_{A,A})\circ(\Delta_A\otimes\id_A)}&&&&&
  \underline{A\otimes A\{-2\}}\ar[r]&0
}.
\end{equation}
  The $A$-actions are:
\begin{enumerate}
\item
    $\lambda^{(i_1)}\colon A\otimes C_j\to C_j$ with ${(\lambda^{(i_1)})}^{-1}:=\lambda_1$ and ${(\lambda^{(i_1)})}^0:=\lambda_1$,
\item
    $\rho^{(i_2)}\colon C_j\otimes A\to C_j$ with ${(\rho^{(i_2)})}^{-1}:=\rho_1$ and ${(\rho^{(i_2)})}^0:=\rho_2$,
\item
    $\lambda^{(i_3)}\colon A\otimes C_j\to C_j$ with ${(\lambda^{(i_3)})}^{-1}:=\lambda_2$ and ${(\lambda^{(i_3)})}^0:=\lambda_2$,
\item
    $\rho^{(i_4)}\colon C_j\otimes A\to C_j$ with ${(\rho^{(i_4)})}^{-1}:=\rho_2$ and ${(\rho^{(i_4)})}^0:=\rho_1$.
\end{enumerate}
\item
  If $\overline{R_j}\cap T$ is a positive crossing with the
  checkerboard colouring
\begin{equation}
\begin{aligned}
\begin{pspicture}(1.5,1.5)
  \pscustom[fillstyle=solid,fillcolor=lightgray,linecolor=lightgray]{
    \psline(0.0,0.0)(1.5,0.0)(0.0,1.5)(1.5,1.5)(0.0,0.0)
  }
  \psline[linewidth=1pt]{<-}(0.0,1.5)(1.5,0.0)
  \psline[linewidth=3pt,linecolor=white](0.7,0.8)(0.8,0.7)
  \psline[linewidth=1pt]{->}(0.0,0.0)(1.5,1.5)
  \rput(-0.2,-0.2){$i_1$}
  \rput(1.7,-0.2){$i_2$}
  \rput(1.7,1.7){$i_3$}
  \rput(-0.2,1.7){$i_4$}
\end{pspicture}
\end{aligned}
\end{equation}
  we define
\begin{equation}
\label{eq_complexfour}
  C_j:=\xymatrix{
0\ar[r]&\underline{A\otimes A\{2\}}\ar[rrrrr]^{(\id_A\otimes \mu_A)\circ(\tau_{A,A}\otimes\id_A)\circ(\id_A\otimes\Delta_A)}&&&&&
  A\otimes A\{4\}\ar[r]&0
}.
\end{equation}
  The $A$-actions are:
\begin{enumerate}
\item
    $\rho^{(i_1)}\colon C_j\otimes A\to C_j$ with ${(\rho^{(i_1)})}^0:=\rho_1$ and ${(\rho^{(i_1)})}^1:=\rho_1$,
\item
    $\lambda^{(i_2)}\colon A\otimes C_j\to C_j$ with ${(\lambda^{(i_2)})}^0:=\lambda_2$ and ${(\lambda^{(i_2)})}^1:=\lambda_1$,
\item
    $\rho^{(i_3)}\colon C_j\otimes A\to C_j$ with ${(\rho^{(i_3)})}^0:=\rho_2$ and ${(\rho^{(i_3)})}^1:=\rho_2$,
\item
    $\lambda^{(i_4)}\colon A\otimes C_j\to C_j$ with ${(\lambda^{(i_4)})}^0:=\lambda_1$ and ${(\lambda^{(i_4)})}^1:=\lambda_2$.
\end{enumerate}
\item
  If $\overline{R_j}\cap T$ is a negative crossing with the
  checkerboard colouring
\begin{equation}
\begin{aligned}
\begin{pspicture}(1.5,1.5)
  \pscustom[fillstyle=solid,fillcolor=lightgray,linecolor=lightgray]{
    \psline(0.0,0.0)(1.5,0.0)(0.0,1.5)(1.5,1.5)(0.0,0.0)
  }
  \psline[linewidth=1pt]{->}(0.0,0.0)(1.5,1.5)
  \psline[linewidth=3pt,linecolor=white](0.7,0.7)(0.8,0.8)
  \psline[linewidth=1pt]{<-}(0.0,1.5)(1.5,0.0)
  \rput(-0.2,-0.2){$i_1$}
  \rput(1.7,-0.2){$i_2$}
  \rput(1.7,1.7){$i_3$}
  \rput(-0.2,1.7){$i_4$}
\end{pspicture}
\end{aligned}
\end{equation}
  we define
\begin{equation}
\label{eq_complexfive}
  C_j:=\xymatrix{
0\ar[r]&A\otimes A\{-4\}\ar[rrrrr]^{(\id_A\otimes\mu_A)\circ(\tau_{A,A}\otimes\id_A)\circ(\id_A\otimes\Delta_A)}&&&&&
  \underline{A\otimes A\{-2\}}\ar[r]&0
}.
\end{equation}
  The $A$-actions are:
\begin{enumerate}
\item
    $\rho^{(i_1)}\colon C_j\otimes A\to C_j$ with ${(\rho^{(i_1)})}^{-1}:=\rho_1$ and ${(\rho^{(i_1)})}^0:=\rho_1$,
\item
    $\lambda^{(i_2)}\colon A\otimes C_j\to C_j$ with ${(\lambda^{(i_2)})}^{-1}:=\lambda_1$ and ${(\lambda^{(i_2)})}^0:=\lambda_2$,
\item
    $\rho^{(i_3)}\colon C_j\otimes A\to C_j$ with ${(\rho^{(i_3)})}^{-1}:=\rho_2$ and ${(\rho^{(i_3)})}^0:=\rho_2$,
\item
    $\lambda^{(i_4)}\colon A\otimes C_j\to C_j$ with ${(\lambda^{(i_4)})}^{-1}:=\lambda_2$ and ${(\lambda^{(i_4)})}^0:=\lambda_1$.
\end{enumerate}
\end{enumerate}
Let $G\cap T=\{x_1,\ldots,x_s\}$ be the set of points at which the
graph $G$ intersects the tangle diagram $T$ and thereby `cuts' it into
its constituents. For each $i\in\{1,\ldots,s\}$, there are two numbers
$p_i,q_i\in\{1,\ldots,k\}$ indicating the numbers of the two boundary
points of some of the $\overline{R_j}\cap T$ that coincide with
$x_i$. One of them, say $p_i$ is a `$+$'-boundary and thus corresponds
to a right $A$-action $\rho^{(p_i)}$ on some $C_j$, whereas the other
one, $q_i$, is a `$-$'-boundary relative to $\epsilon$ and therefore
corresponds to a left $A$-action $\lambda^{(q_i)}$ on some
$C_{j^\prime}$.

We now take the tensor product over $k$ of all complexes $C_j$,
$j\in\{1,\ldots,\ell\}$, and then apply the coequalizer that
`coequalizes' the right $A$-action $\rho^{(p_i)}$ with the left
$A$-action $\lambda^{(q_i)}$ for all $i\in\{1,\ldots,s\}$:
\begin{equation}
\label{eq_compositionresult}
  C(T,\epsilon,\mathbbm{A}) := {}_{\rho^{(p_1)}}\otimes_{\lambda^{(q_1)}}\Biggl(
  {}_{\rho^{(p_2)}}\otimes_{\lambda^{(q_2)}}\biggl(\cdots\bigl(
  {}_{\rho^{(p_s)}}\otimes_{\lambda^{(q_s)}}\bigl(\bigotimes_{j=1}^\ell C_j\bigr)\bigr)\cdots\biggr)\Biggr).
\end{equation}
\end{defn}

\begin{rem}
In each of the five cases, $C_j$ agrees with the tangle complex
$\tangle{\overline{R_j}\cap T}_{\epsilon,Z}$ for the state sum
open-closed TQFT $Z$ associated with $\mathbbm{A}$. Note in
particular that the differential of~\eqref{eq_complextwo} agrees with
$Z(-)$ applied to the open-closed cobordism~\eqref{eq_crossing},
\etc. For each of the complexes, we have in detail specified the
$A$-actions associated with the boundary points
(\cf~Proposition~\ref{prop_actions}). The complex
$C(T,\epsilon,\mathbbm{A})$ is a complex of $(A^{\otimes
r},A^{\otimes r})$-bimodules, $p+q=2r$, associated with the diagram
of the $(p,q)$-tangle $T$.
\end{rem}

\begin{thm}
\label{thm_compositionresult}
Let $T$ be a plane diagram of an oriented $(p,q)$-tangle,
$\epsilon\in\{-1,+1\}$ denote a checkerboard colouring, and let
$Z\colon\cat{2Cob}^\mathrm{ext}\to\cat{Vect}_k$ be a [Euler-graded,
Euler-filtered] state sum open-closed TQFT that satisfies Bar-Natan's
conditions. We denote by $\mathbbm{A}$ its associated knowledgeable
Frobenius algebra. Then the complex of Definition~\ref{def_composed}
is isomorphic to the tangle complex of
Definition~\ref{eq_tanglecomplex},
\begin{equation}
\label{eq_compositionclaim}
  C(T,\epsilon,\mathbbm{A})\cong\tangle{T}_{\epsilon,Z},
\end{equation}
as [graded, filtered] complexes of
$(A^{(q,\epsilon)},A^{(p,\epsilon)})$-bimodules.
\end{thm}

While the proof of this theorem is conceptually simple, making a
precise proof is rather laborious.  The picture that the reader should
have in mind is that the complex constructed via coequalizers is
actually what one would get by using the state sum construction of
\cite{LP2}. Although the proof below is entirely self contained, it
really just amounts to showing that the state sum construction
actually reproduces the knowledgeable Frobenius algebra associated to
the Euler-graded [Euler-filtered] state sum open-closed TQFT $Z$. As
in \cite{LP2}, the key to proving this theorem will be the maps
$P_{ij}$ and $Q_{ij}$ corresponding to triangulated strips and
cylinders whose incoming boundary are triangulated with $i$ arcs and
whose outgoing boundaries are triangulated with $j$ arcs.

\begin{proof}
First, we observe that each of the complexes $C_j$,
$j\in\{1,\ldots,\ell\}$, of~\eqref{eq_complexone},
\eqref{eq_complextwo}, \eqref{eq_complexthree},
\eqref{eq_complexfour}, and~\eqref{eq_complexfive} is the total
complex $C_j\cong C_\mathrm{tot}(\cal{I}_j,{\tilde C}_j)$ of some
commutative cube ${\tilde C}_j$ in $\cat{Vect}_k$ [$\cat{grdVect}_k$,
$\cat{fltVect}_k$]. For~\eqref{eq_complexone}, this is an
$\cal{I}_j$-cube with only one vertex, \ie\ $\cal{I}_j=\emptyset$,
whereas in the other cases, this is an $\cal{I}_j$-cube with two
vertices, \ie\ $\cal{I}=\{z\}$ where $z\in\alpha\subseteq\cat{n}$
indicates the number of the respective crossing. By
Proposition~\ref{prop_total},
\begin{equation}
  \bigotimes_{j=1}^\ell C_j\cong
  C_\mathrm{tot}(\cat{n},\tilde C)\{2n_+-4n_-\}\qquad\mbox{where}\qquad
  \tilde C:=\bigboxtimes_{j=1}^\ell {\tilde C}_j,
\end{equation}
\ie\ the tensor product of complexes that features
in~\eqref{eq_compositionresult}, is the total complex of the external
tensor product of the commutative cubes ${\tilde C}_j$. We have made
some of the grading shifts explicit here.

The vertices of the commutative cube $\tilde C$ are labeled by the
smoothings $\alpha\subseteq\cat{n}$ of the $n$ crossings of $T$. Each
vertex ${\tilde C}_\alpha$ of $\tilde C$ is of the form ${\tilde
C}_\alpha\cong A^{\otimes(2n+m)}\{2|\alpha|\}$ since each of the $n$
crossings contributes a tensor factor of $A\otimes A$ and each of the
$m$ arcs contributes just an $A$. The ordering of the tensor factors
is induced from the numbering $1,\ldots,\ell$ of the components of
$(\R^2\backslash G)\cap T$.

In order to prove the theorem, we permute the tensor factors of each
${\tilde C}_\alpha\cong A^{\otimes(2n+m)}$, $\alpha\subseteq\cat{n}$,
so that an isomorphism of commutative cubes can be explicitly
given. By an argument analogous to that of
Proposition~\ref{prop_indep}(2.), permuting the tensor factors of an
arbitrary vertex of a commutative cube $\tilde C$ results in applying
an isomorphism to the associated total complex
$C_\mathrm{tot}(\cat{n},\tilde C)$.

Consider a vertex $\alpha\subseteq\cat{n}$ of the cube $\tilde C$. The
associated smoothing $D_\alpha\subseteq\R^2$ of the tangle diagram is
a disjoint union of arcs and circles that can be ordered
arbitrarily. Number these components of $D_\alpha$ by
$1,\ldots,k_\alpha$ in the same way as in
Definition~\ref{def_tanglecube} and also define the sequence
$m^{(\alpha)}=(m_1^{(\alpha)},\ldots,m_{k_\alpha}^{(\alpha)})\in{\{0,1\}}^{k_\alpha}$
such that $m_i^{(\alpha)}=0$ if the $i$-th component of $D_\alpha$ is
a circle and $m_i^{(\alpha)}=1$ if it is an arc. The graph $G$ chops
each $D_\alpha\subseteq\R^2$ into several pieces $\overline{R_j}\cap
D_\alpha$, $j\in\{1,\ldots,\ell\}$. If $\overline{R_j}\cap T$ is an
arc, then $\overline{R_j}\cap D_\alpha$ is an arc and ${\tilde
C}_j\cong A$, whereas if $\overline{R_j}\cap T$ is a crossing, then
$\overline{R_j}\cap D_\alpha$ is the disjoint union of two arcs and
${\tilde C}_j\cong A\otimes A$.

We assume that we have permuted the tensor factors of each ${\tilde
C}_\alpha$ in such a way that
\begin{equation}
\label{eq_verticesofcube}
  {\tilde C}_\alpha = {\tilde C}_\alpha^{(1)}\otimes\cdots\otimes{\tilde C}_\alpha^{(k_\alpha)}\{2|\alpha|\},
\end{equation}
where the ${\tilde C}_\alpha^{(i)}$ are defined as follows. If
$m_i^{(\alpha)}=0$, then the $i$-th component $D_\alpha^{(i)}$ of
$D_\alpha$ is a circle which is chopped by the graph complement into
$M_i^{(\alpha)}:=|\{\,j\in\{1,\ldots,\ell\}\mid\,\overline{R_j}\cap
D_\alpha^{(j)}\neq\emptyset\,\}|$ arcs. If $m_i^{(\alpha)}=1$, then
the $i$-th component is an arc which is chopped into $M_i^{(\alpha)}$
arcs. In both cases, we put ${\tilde C}_\alpha^{(i)}:= A^{\otimes
M_i^{(\alpha)}}$.

We also assume that we have permuted the tensor factors of the
${\tilde C}_\alpha^{(i)}$ in such a way that any pair of left and
right action $\lambda^{(q_t)}$ and $\rho^{(p_t)}$ involved in the
coequalizer~\eqref{eq_compositionresult} either
\begin{itemize}
\item
  consists of a right action on the $b$-th and a left action on the
  $(b+1)$-th tensor factor of some ${\tilde C}_\alpha^{(i)}$ if
  $m_i^{(\alpha)}\in\{0,1\}$, or
\item
  consists of a right action on the last and a left action on the
  first tensor factor of some ${\tilde C}_\alpha^{(i)}$ if
  $m_i^{(\alpha)}=0$ (circle only).
\end{itemize}

The complex $C(T,\epsilon,\mathbbm{A})$
of~\eqref{eq_compositionresult} can now be obtained by applying the
composition of coequalizers
\begin{equation}
  \bigotimes_{\rho;\lambda} :=
    {}_{\rho^{(p_1)}}\otimes_{\lambda^{(q_1)}}\circ\cdots\circ
    {}_{\rho^{(p_s)}}\otimes_{\lambda^{(q_s)}}
\end{equation}
of vector spaces to each vertex ${\tilde C}_\alpha$,
$\alpha\subseteq\cat{n}$ of the cube and by then taking the total
complex:
\begin{equation}
  C(T,\epsilon,\mathbbm{A}) \cong C_\mathrm{tot}\bigl(\cat{n},\bigotimes_{\rho;\lambda}(\tilde C)\bigr)\{2n_+-4n_-\}.
\end{equation}
Here we have denoted by $\bigotimes_{\rho;\lambda}(\tilde C)$ the cube
whose vertices are the $\bigotimes_{\rho;\lambda}({\tilde C}_\alpha)$
and whose edges are obtained by repeatedly applying the coequalizer to
maps as in~\eqref{eq_coeqmaps}.

Because of the way we have ordered the tensor factors of the ${\tilde
C}_\alpha^{(i)}$,
\begin{enumerate}
\item
  if $m_i^{(\alpha)}=1$, then the ordering of factors is precisely as
  in~\eqref{eq_standardpkl} and
\begin{equation}
  \bigotimes_{\rho;\lambda}({\tilde C}_\alpha^{(i)}) \cong
    P_{M_i^{(\alpha)}M_i^{(\alpha)}}(A^{\otimes M_i^{(\alpha)}}) \cong
    A,
\end{equation}
\item
  if $m_i^{(\alpha)}=0$, then the ordering of factors is precisely as
  in~\eqref{eq_standardqkl} and
\begin{equation}
  \bigotimes_{\rho;\lambda}({\tilde C}_\alpha^{(i)}) \cong
    Q_{M_i^{(\alpha)}M_i^{(\alpha)}}(A^{\otimes M_i^{(\alpha)}}) \cong
    p(A)\cong Z(A).
\end{equation}
\end{enumerate}
This shows that the vertices~\eqref{eq_verticesofcube} of the cube
$\tilde C$ are isomorphic to the vertices of the cube
$Z(S_{T,\epsilon})$ (Definition~\ref{def_vectcube}). We finally have
to show that the corresponding isomorphisms fit into an isomorphism
of commutative cubes which then implies that the resulting total
complexes in~\eqref{eq_compositionclaim} are isomorphic.

The edges $(\alpha,x)\in r(\cat{n})$ of the cube
$\bigotimes_{\rho;\lambda}(\tilde C)$ are obtained from those of
$\tilde C$ by repeatedly applying coequalizers to maps as
in~\eqref{eq_coeqmaps}. As we know that the coequalizers
$\bigotimes_{\rho;\lambda}\colon{\tilde C}_\alpha^{(i)}\to A$ or
${\tilde C}_\alpha^{(i)}\to Z(A)$ are just $P_{1M_i^{(\alpha)}}$ or
$Q_{1M_i^{(\alpha)}}$ (Corollary~\ref{cor_standardpkl}) that have
inverses $P_{M_i^{(\alpha)}1}$ or $Q_{M_i^{(\alpha)}1}$ if restricted
appropriately (Proposition~\ref{prop_projproperty}), we can write down
the edges of the cube explicitly.

If ${\tilde C}_{(\alpha,x)}\colon {\tilde C}_\alpha\to{\tilde
C}_{\alpha\sqcup\{x\}}$, $x\in\cat{n}$,
$\alpha\subseteq\cat{n}\backslash\{x\}$, is an edge of the cube
$\tilde C$, \ie\ ${\tilde C}_{(\alpha,x)}$ is a tensor product of
either~\eqref{eq_crossing} or~\eqref{eq_crossing2} with identities,
then the induced map
\begin{equation}
  {\hat C}_{(\alpha,x)}\colon
    \bigl(\bigotimes_{\rho;\lambda}(\tilde C)\bigr)_\alpha\to
    \bigl(\bigotimes_{\rho;\lambda}(\tilde C)\bigr)_{\alpha\sqcup\{x\}},
\end{equation}
is given by
\begin{equation}
  {\hat C}_{(\alpha,x)}=
  \bigl(P_{M_1^{(\alpha)}1}^{(m_1^{(\alpha)})}\otimes\cdots\otimes
        P_{M_{k_\alpha}^{(\alpha)}1}^{(m_{k_\alpha}^{(\alpha)})}\bigr)\circ
  {\tilde C}_{(\alpha,x)}\circ
  \bigl(P_{1M_1^{(\alpha)}}^{(m_1^{(\alpha)})}\otimes\cdots\otimes
        P_{1M_{k_\alpha}^{(\alpha)}}^{(m_{k_\alpha}^{(\alpha)})}\bigr),
\end{equation}
where we have used the notation $P_{k\ell}^{(0)}:=Q_{k\ell}$ and
$P_{k\ell}^{(1)}:=P_{k\ell}$. Recall that the composition
in~\eqref{eq_universalmaps} is an example of such a map. We now have
to verify that these edges agree with the edges of the cube
$Z(S_{T,\epsilon})$ of Definition~\ref{def_vectcube}, \ie
\begin{equation}
  {\hat C}_{(\alpha,x)} = Z(S_{(\alpha,x)})
\end{equation}
for all $x\in\cat{n}$, $\alpha\subseteq\cat{n}\backslash\{x\}$. Since
we know that ${\tilde C}_{(\alpha,x)}$ is the tensor product of
either~\eqref{eq_crossing} or~\eqref{eq_crossing2} with identities and
since we have the explicit forms of $P_{k\ell}$ and $Q_{k\ell}$, this
can be done in case by verification.

Up to permutations, duality, $A\leftrightarrow A^\mathrm{op}$ symmetry
and the relations~\eqref{eq_projproperty}, we have to verify four
equations. We write
$d:=(\mu_A\otimes\id_A)\circ(\id_A\otimes\tau_{A,A})\circ(\Delta_A\otimes\id_A)$
for the differential associated with the crossing~\eqref{eq_crossing}.
\begin{enumerate}
\item
  $(P_{12}\otimes\id_A)\circ (\id_A\otimes d)\circ (P_{21}\otimes\id_A) = d$. This
  corresponds to gluing an arc to a crossing in the trivial way and
  can be seen to follow from the relations of the knowledgeable
  Frobenius algebra as follows:
\begin{equation}
  a^{-1}\cdot
\begin{aligned}
\psset{xunit=2mm,yunit=2mm}
\begin{pspicture}(7,10)
  \rput(1,7.5){\comultl}
  \rput(6,7.5){\identl}
  \rput(0,5){\identl}
  \rput(0,2.5){\identl}
  \rput(2,2.5){\saddler}
  \rput(1,0){\multl}
  \rput(6,0){\identl}
\end{pspicture}
\end{aligned}
\qquad\cong\qquad
\begin{aligned}
\psset{xunit=2mm,yunit=2mm}
\begin{pspicture}(5,5)
  \rput(0,0){\saddler}
\end{pspicture}
\end{aligned}
\end{equation}
\item
  $(\id_A\otimes P_{12})\circ
  (d\otimes\id_A)\circ(\id_A\otimes\tau_{A,A})\circ(P_{21}\otimes\id_A) = d$. This
  also corresponds to gluing an arc to a crossing in the trivial way:
\begin{equation}
  a^{-1}\cdot
\begin{aligned}
\psset{xunit=2mm,yunit=2mm}
\begin{pspicture}(7,12.5)
  \rput(6,10){\identl}
  \rput(0,10){\mediumcomultl}
  \rput(0,7.5){\identl}
  \rput(5,7.5){\crossl}
  \rput(0,2.5){\saddler}
  \rput(6,5){\identl}
  \rput(6,2.5){\identl}
  \rput(0,0){\identl}
  \rput(5,0){\multl}
\end{pspicture}
\end{aligned}
\qquad\cong\qquad
\begin{aligned}
\psset{xunit=2mm,yunit=2mm}
\begin{pspicture}(5,5)
  \rput(0,0){\saddler}
\end{pspicture}
\end{aligned}
\end{equation}
\item
  $P_{12}\circ\tau_{A,A}\circ d\circ (Q_{11}\otimes\id_A) =
  \mu_A\circ(\imath\otimes\id_A)$. Note that
  $Q_{11}=p=a^{-1}\cdot\imath$. This corresponds to an arc which is
  glued with both ends to two boundary points of a crossing:
\begin{equation}
\label{eq_compositionwithwindow}
  a^{-1}\cdot
\begin{aligned}
\psset{xunit=2mm,yunit=2mm}
\begin{pspicture}(5,12.5)
  \rput(-0.2,10){\ctl}
  \rput(4,10){\identl}
  \rput(0,5){\saddler}
  \rput(2,2.5){\widecrossl}
  \rput(0,0){\mediummultl}
\end{pspicture}
\end{aligned}
\qquad\cong\qquad
\begin{aligned}
\psset{xunit=2mm,yunit=2mm}
\begin{pspicture}(5,5)
  \rput(-0.2,2.5){\ctl}
  \rput(2,2.5){\identl}
  \rput(1,0){\multl}
\end{pspicture}
\end{aligned}
\end{equation}
\item
  $Q_{11}\circ\mu_A\circ(\imath\otimes Q_{11}) = \mu_{Z(A)}$. Note
  that $Q_{11}=a^{-1}\cdot \imath$ if restricted to $p(A)$ and
  $Q_{11}=\imath^\ast$ on $A$. This corresponds to a crossing which
  already has two boundary points connected by an arc as in case (3.)
  and which gets another arc glued to the remaining two boundary
  points:
\begin{equation}
  a^{-1}\cdot
\begin{aligned}
\psset{xunit=2mm,yunit=2mm}
\begin{pspicture}(5,7.5)
  \rput(-0.2,5){\ctl}
  \rput(1.8,5){\ctl}
  \rput(1,2.5){\multl}
  \rput(1.2,.5){\ltc}
\end{pspicture}
\end{aligned}
\qquad\cong\qquad
\begin{aligned}
\psset{xunit=2mm,yunit=2mm}
\begin{pspicture}(2,2.5)
  \rput(0,0){\multc}
\end{pspicture}
\end{aligned}
\end{equation}
\end{enumerate}
These relations imply that pre- and postcomposing ${\tilde
C}_{(\alpha,x)}$ with the appropriate $P_{k\ell}$ and $Q_{k\ell}$
indeed results in the maps of Remark~\ref{rem_crossings} and thereby
in the edges of the cube $Z(S_{T,\epsilon})$.
\end{proof}

The following corollary is now obvious and settles the question for
the composition rule~\eqref{eq_compositionrule} raised in the
introduction.

\begin{cor}
Let $T$ be a plane diagram of an oriented $(p,q)$-tangle and
$T^\prime$ be one of an oriented $(q,r)$-tangle such that the
orientations match and the composition $T^\prime\circ T$ is well
defined. Let $Z\colon\cat{2Cob}^\mathrm{ext}\to\cat{Vect}_k$ be a
state sum open-closed TQFT that satisfies Bar-Natan's conditions. Then
\begin{equation}
  {\tangle{T^\prime\circ T}}_{\epsilon,Z}\cong
  {\tangle{T^\prime}}_{\epsilon,Z} \otimes_{A^{(q,\epsilon)}} {\tangle{T}}_{\epsilon,Z}
\end{equation}
are isomorphic as complexes of
$(A^{(r,\epsilon)},A^{(p,\epsilon)})$-bimodules,
\cf~Remark~\ref{rem_tanglemodule}.
\end{cor}

\begin{rem}
In order to generalize the composition of tangles to non-strongly
separable algebras $A$, one could try the following. One sticks to the
use of tensor products over $A$ for the composition, but one tries to
avoid the use of anything that involves an $a^{-1}$. This excludes
some of the $P_{jk}$ and $Q_{jk}$, but luckily not those $P_{1k}$ and
$Q_{1k}$ that form the coequalizers. Nevertheless, apart from the
$a^{-1}$, the open-closed cobordism depicted
in~\eqref{eq_compositionwithwindow} contains a window and will always
yield a linear map with non-trivial kernel unless $A$ is strongly
separable.
\end{rem}

\subsection{Invariance under Reidemeister moves}
\label{sec_Reidemeister}

As we now have a construction of the tangle complex in local terms,
\ie\ from a composition of arcs and crossings, we are in a position to
supply a proof of the invariance of its homology under Reidemeister
moves in similarly local terms.

\begin{thm}
Let $T,T^\prime$ be plane diagrams of oriented $(p,q)$-tangles, let
$\epsilon\in\{-1,+1\}$, and let
$Z\colon\cat{2Cob}^\mathrm{ext}\to\cat{Vect}_k$ be an Euler-graded
[Euler-filtered] state sum open-closed TQFT that satisfies Bar-Natan's
conditions. If $T$ and $T^\prime$ are related by a plane isotopy or by
a Reidemeister move, then the complexes $\tangle{T}_{\epsilon,Z}$ and
$\tangle{T^\prime}_{\epsilon,Z}$ are homotopy equivalent as graded
[filtered] complexes of
$(A^{(q,\epsilon)},A^{(p,\epsilon)})$-bimodules.
\end{thm}

The case of plane isotopies is already covered by
Theorem~\ref{thm_compositionresult}. We now provide algebraic at the
same time local proofs of the invariance up to homotopy equivalence
under Reidemeister moves.  These proofs are inspired by Bar-Natan's
picture world, but it is nevertheless worthwhile spelling out the
details in the language of knowledgeable Frobenius algebras. Thanks to
the composition properties of state sum open-closed TQFTs described in
the previous section, we need to check the Reidemeister moves
\begin{equation}
 \Reidemeister
\end{equation}
only in the skein theoretic sense, meaning that any of the above
diagrams can appear inside a possibly larger diagram, employing
tensor products and compositions. Note that because our convention
for orienting the tangle resolutions involves the checkerboard
colouring $\epsilon\in\{0,1\}$ of the tangle, we must take account
for both possible colourings of the above diagrams.

We denote by $(A,C,\imath,\imath^\ast)$ the state sum knowledgeable
Frobenius algebra associated with $Z$. In what follows we will omit
the functor $Z$ around all diagrams of open-closed cobordisms.

\subsubsection{Reidemeister move one}

We show that the complex $\llbracket \Ronei\rrbracket =
\left( 0 \longrightarrow \underline{A\{0\}} \longrightarrow
0\right)$ is homotopy equivalent to the complex
\begin{equation}
\label{eq_reideone}
\llbracket
\psset{xunit=.5cm,yunit=.5cm}
\begin{pspicture}[.4](2,1.5)
 \psbezier[linewidth=.8pt](.5,0)(.5,.5)(1.5,.5)(1.5,1)
 \pspolygon[linecolor=white,fillstyle=solid,fillcolor=white](.9,0)(.9,2)(1.1,2)(1.1,0)(.9,0)
  \psbezier[linewidth=.8pt](1.5,0)(1.5,.5)(.5,.5)(.5,1)
  \psbezier[linewidth=.8pt](.5,1)(.5,1.45)(1.5,1.45)(1.5,1)

   \rput(.3,1){$\scriptstyle 2$}
  \rput(.3,.3){$\scriptstyle 1$}
  \rput(1.65,.3){$\scriptstyle 3$}
\end{pspicture}
\rrbracket=
 \xy
 (0,0)*+{\underline{(A \otimes C)\{2\}}}="1";
 (30,0)*+{A\{4\}}="2";
 {\ar^-{d^0=
     \psset{xunit=.15cm,yunit=.15cm}
        \begin{pspicture}[.4](4,5)
         \rput(2.75,2.5){\ctl}
         \rput(1,2.5){\medidentl}
         \rput(2,0){\multl}
        \end{pspicture}
    } "1";"2"};
 \endxy
\end{equation}
where we have underlined the degree zero term of each complex. Note
that neither of these complexes depends on the choice of
checkerboard colouring $\epsilon\in\{-1,+1\}$. The second complex
always has $(n_+,n_-)=(1,0)$ regardless of the orientation of the
tangle. Hence, we need only consider the one case.

Define the chain map $F \maps \llbracket\Ronei\rrbracket \to
\llbracket\Roneii\rrbracket$ as follows:
\begin{equation}
F^0  :=
    \psset{xunit=.25cm,yunit=.25cm}
    \begin{pspicture}[.4](6,5)
    \rput(1,0){\identl}
    \rput(1,2.5){\medidentl}
    \rput(4,0){\holec}
    \rput(4.2,4){\birthc}
    \end{pspicture}
    \;\;-\;\;
    \begin{pspicture}[.4](6,5)
    \rput(1,0){\medidentl}
    \rput(3.15,0){\ltc}
    \rput(2,2){\comultl}
    \end{pspicture},
\qquad \qquad F^1 :=0.
\end{equation}
The chain map property $d^0\circ F^0=0$ follows from the relations
of the knowledgeable Frobenius algebra $(A,C,\imath,\imath^\ast)$.
In fact, by Corollary~3.24 of~\cite{LP1} it is sufficient to compute
the three topological invariants of the relevant composite
open-closed cobordisms, and then Theorem~3.22 of~\cite{LP1} produces
for us the sequence of relations that constitutes the proof. Next we
define the chain map $G \maps \llbracket\Roneii\rrbracket \to
\llbracket\Ronei\rrbracket$ by
\begin{equation}
G^0 := \psset{xunit=.25cm,yunit=.25cm}
    \begin{pspicture}[.4](4,2.5)
    \rput(1,0){\identl}
    \rput(3.2,1.5){\deathc}
    \end{pspicture},
    \qquad \qquad
G^1:=0.
\end{equation}
There are no conditions that need to be verified to see that this is
a chain map.

The chain maps $F$ and $G$ define a homotopy equivalence between the
complexes $\llbracket \Ronei \rrbracket$ and $\llbracket \Roneii
\rrbracket$. To see this, first note that $G\circ F=\id$
by~\eqref{eq_barnatan1} since
\begin{equation}
G^0\circ F^0 =
    \psset{xunit=.25cm,yunit=.25cm}
    \begin{pspicture}[.4](5,6)
    \rput(1,0){\identl}
    \rput(1,2.5){\identl}
    \rput(4,0){\deathc}
    \rput(4,1){\holec}
    \rput(4.2,5){\birthc}
    \end{pspicture}
   \quad -
    \begin{pspicture}[.4](4,6)
    \rput(1,1){\identl}
    \rput(3,0){\deathc}
    \rput(3.2,1){\ltc}
    \rput(2,3){\comultl}
    \end{pspicture}
    \;\; = \;\; 2
    \psset{xunit=.25cm,yunit=.25cm}
    \begin{pspicture}[.4](2,2.5)
    \rput(1,0){\identl}
    \end{pspicture}
     -
    \psset{xunit=.25cm,yunit=.25cm}
    \begin{pspicture}[.4](2,2.5)
    \rput(1,0){\identl}
    \end{pspicture}
    =
    \psset{xunit=.25cm,yunit=.25cm}
    \begin{pspicture}[.4](2,2.5)
    \rput(1,0){\identl}
    \end{pspicture}
    =
    \id.
\end{equation}
Define the chain homotopy $h\maps\id\to F\circ G$ whose only nonzero
component is given by
\begin{equation}
h^1 := \psset{xunit=.25cm,yunit=.25cm}
   \begin{pspicture}[.4](4,2.5)
    \rput(1,0){\identl}
    \rput(3,0){\birthc}
    \end{pspicture} .
\end{equation}
The equation $\id - F^0\circ G^0=h^1\circ d^0$ is depicted as:
\begin{equation}
    \psset{xunit=.25cm,yunit=.25cm}
    \begin{pspicture}[.4](5,7)
    \rput(1,0){\identl}
    \rput(1,2.5){\identl}
    \rput(1,5){\medidentl}
    \rput(3,0){\identc}
    \rput(3,2.5){\identc}
    \rput(3,5){\medidentc}
    \end{pspicture}
    -
    \begin{pspicture}[.4](6,7)
    \rput(1,0){\identl}
    \rput(1,2.5){\identl}
    \rput(1,5){\medidentl}
    \rput(4,0){\holec}
    \rput(4.2,4){\birthc}
    \rput(4,6){\deathc}
    \end{pspicture}
    +
    \begin{pspicture}[.4](5,7)
    \rput(1,0){\medidentl}
    \rput(3.2,0){\ltc}
    \rput(2,2){\comultl}
    \rput(2,4.5){\identl}
    \rput(4,6){\deathc}
    \end{pspicture}
    =
    \begin{pspicture}[.4](5.5,7)
    \rput(2,0){\medidentl}
    \rput(2,2){\multl}
    \rput(1,4.5){\identl}
    \rput(2.8,4.5){\ctl}
    \rput(5,0){\birthc}
    \end{pspicture} .
\end{equation}
After applying appropriate relations and rearranging the terms, this
equality follows from the relation~\eqref{eq_barnatan3}
\begin{equation}
    \psset{xunit=.25cm,yunit=.25cm}
    \begin{pspicture}[.4](6,12)
    \rput(1,0){\identl}
    \rput(1,2.5){\identl}
    \rput(1,5){\identl}
    \rput(1,7.5){\medidentl}
    \rput(2,9.5){\comultl}
    \rput(4.2,0){\multc}
    \rput(5.2,2.5){\birthc}
    \rput(3.2,2.5){\birthc}
    \rput(4,4){\deathc}
    \rput(4.2,5){\multc}
    \rput(3.2,7.5){\ltc}
     \rput(5.2,7.5){\medidentc}
     \rput(5.2,9.5){\identc}
     \pspolygon[linestyle=dashed](2,2.2)(2,7.5)(6,7.5)(6,2.2)(2,2.2)
    \end{pspicture}
    \quad + \quad
    \begin{pspicture}[.4](6,12)
    \rput(1,0){\identl}
    \rput(1,2.5){\identl}
    \rput(1,5){\identl}
     \rput(1,7.5){\medidentl}
    \rput(2,9.5){\comultl}
    \rput(4.2,0){\multc}
    \rput(4.2,2.5){\comultc}
    \rput(4.2,5){\birthc}
    \rput(3,6.5){\deathc}
    \rput(5,6.5){\deathc}
    \rput(3.2,7.5){\ltc}
     \rput(5.2,7.5){\medidentc}
     \rput(5.2,9.5){\identc}
     \pspolygon[linestyle=dashed](2,2.2)(2,7.5)(6,7.5)(6,2.2)(2,2.2)
    \end{pspicture}
    \quad = \quad
    \begin{pspicture}[.4](6,12)
    \rput(1,0){\identl}
    \rput(1,2.5){\identl}
    \rput(1,5){\identl}
     \rput(1,7.5){\medidentl}
    \rput(2,9.5){\comultl}
    \rput(4.2,0){\multc}
    \rput(5.2,2.5){\birthc}
    \rput(3.2,2.5){\identc}
    \rput(3.2,5){\identc}
    \rput(5,6.5){\deathc}
    \rput(3.2,7.5){\ltc}
    \rput(5.2,7.5){\medidentc}
     \rput(5.2,9.5){\identc}
     \pspolygon[linestyle=dashed](2,2.2)(2,7.5)(6,7.5)(6,2.2)(2,2.2)
    \end{pspicture}
    \quad + \quad
    \begin{pspicture}[.4](6,12)
    \rput(1,0){\identl}
    \rput(1,2.5){\identl}
    \rput(1,5){\identl}
     \rput(1,7.5){\medidentl}
    \rput(2,9.5){\comultl}
    \rput(4.2,0){\multc}
    \rput(3.2,2.5){\birthc}
    \rput(5.2,2.5){\identc}
    \rput(5.2,5){\identc}
    \rput(3,6.5){\deathc}
    \rput(3.2,7.5){\ltc}
    \rput(5.2,7.5){\medidentc}
     \rput(5.2,9.5){\identc}
     \pspolygon[linestyle=dashed](2,2.2)(2,7.5)(6,7.5)(6,2.2)(2,2.2)
    \end{pspicture}
\end{equation}
The equation $\id - F^1\circ G^1 = d^0\circ h^1$ follows from the
knowledgeable Frobenius algebra axiom asserting that the map
$\imath$ is an algebra homomorphism (preserves the unit) together
with the unit axiom for the algebra $A$. Notice that all maps $F^j$,
$G^j$, and $h^j$ are graded [filtered] of degree $0$ if the
appropriate grading shifts $\{-\}$ in~\eqref{eq_reideone} are taken
into account. The other version of the first Reidemeister move can
be treated similarly.

\subsubsection{Reidemeister move two}

We show that the following complexes
\begin{eqnarray} \label{eq_ffvi}
\llbracket \;
\begin{pspicture}[.4](2,2.1)
  \pspolygon[linecolor=lightgray,
   linewidth=.6pt,fillstyle=solid,fillcolor=lightgray](.5,0)(.5,2.03)(0,2.03)(0,0)(.5,0)
   \pspolygon[linecolor=lightgray,
   linewidth=.6pt,fillstyle=solid,fillcolor=lightgray](1.5,0)(1.5,2.03)(2,2.03)(2,0)(1.5,0)
  \psline[linewidth=.8pt](.5,0)(.5,2)
  \psline[linewidth=.8pt](1.5,0)(1.5,2)
  \rput(.35,1){$\scriptstyle 2$}
  \rput(1.65,1){$\scriptstyle 1$}
\end{pspicture}\;
\rrbracket&:=&
 \xy
    (-45,0)*+{0}="1";
    (0,0)*+{\underline{(A \otimes A)\{0\}}}="2";
    (45,0)*+{0}="3";
    {\ar "1";"2"};
    {\ar "2";"3"};
 \endxy
\\
\label{eq_reidetwo}
\llbracket \;
\begin{pspicture}[.4](2,2.1)
  \rput(-.5,0){ \pspolygon[linecolor=lightgray,
   linewidth=.6pt,fillstyle=solid,fillcolor=lightgray](.5,0)(.5,2.03)(2.5,2.03)(2.5,0)(.5,0)
  \pscustom[linecolor=white,linewidth=0pt,fillstyle=solid,fillcolor=white]{
  \psbezier[linewidth=.8pt](2,0)(2,.5)(1,.5)(1,1)
   \psbezier[linewidth=.8pt](1,1)(1,1.5)(2,1.5)(2,2)
   \psline(1,2)
   \psbezier[linewidth=.8pt](1,2)(1,1.5)(2,1.5)(2,1)
   \psbezier[linewidth=.8pt](2,1)(2,.5)(1,.5)(1,0)
   \psline(2,0)
  }}
 \psbezier[linewidth=.8pt](1.5,0)(1.5,.5)(.5,.5)(.5,1)
 \psbezier[linewidth=.8pt](1.5,2)(1.5,1.5)(.5,1.5)(.5,1)
 \pspolygon[linecolor=white,fillstyle=solid,fillcolor=white](.9,0)(.9,2)(1.1,2)(1.1,0)(.9,0)
  \psbezier[linewidth=.8pt](.5,0)(.5,.5)(1.5,.5)(1.5,1)
 \psbezier[linewidth=.8pt](.5,2)(.5,1.5)(1.5,1.5)(1.5,1)
  \rput(.35,1){$\scriptstyle 5$}
  \rput(1.65,1){$\scriptstyle 1$}
  \rput(.35,1.7){$\scriptstyle 2$}
  \rput(1.65,1.7){$\scriptstyle 4$}
  \rput(.35,.3){$\scriptstyle 3$}
  \rput(1.65,.3){$\scriptstyle 6$}
\end{pspicture} \;
\rrbracket&:=&
 \xy
    (-45,0)*+{\scs(A \otimes A)\{-2\}}="1";
    (0,0)*+{\scs \underline{(C \otimes A \otimes A)\{0\}\oplus(A \otimes A)\{0\}}}="2";
    (45,0)*+{\scs (A \otimes A)\{2\}}="3";
    {\ar^-{d^{-1}=
    \left(
      \begin{array}{c}
     \psset{xunit=.15cm,yunit=.15cm}
        \begin{pspicture}(6,5)
         \rput(1,0){\identl}
         \rput(1.25,0){\ltc}
         \rput(4,0){\crossl}
         \rput(2,2.5){\comultl}
         \rput(5,2.5){\identl}
        \end{pspicture} \\
        \psset{xunit=.15cm,yunit=.15cm}
        \begin{pspicture}(4,5)
         \rput(4,0){\saddlel}
        \end{pspicture} \\
      \end{array}
    \right)
    } "1";"2"};
    {\ar^-{d^0=
    \left(
      \begin{array}{c}
     \psset{xunit=.15cm,yunit=.15cm}
        \begin{pspicture}(5,5)
         \rput(2,0){\multl}
         \rput(.75,2.5){\ctl}
         \rput(3,2.5){\identl}
         \rput(5,2.5){\identl}
         \rput(4,0){\curverightl}
        \end{pspicture} \\
        \psset{xunit=.15cm,yunit=.15cm}
        \begin{pspicture}(6,5)
        \rput(-.3,2.5){$-$}
         \rput(2,0){\saddler}
        \end{pspicture} \\
      \end{array}
    \right)^{\mT}
    } "2";"3"};
 \endxy \nn \\ \label{eq_ffvii}
\end{eqnarray}
are homotopy equivalent. Note that in~\eqref{eq_reidetwo} we have
$(n_+,n_-)=(1,1)$ regardless of the orientation of the tangle. We
begin by defining chain maps $F\maps\llbracket
\psset{xunit=.2cm,yunit=.2cm}
\begin{pspicture}[.2](2,2)
  \psline[linewidth=.5pt](.5,0)(.5,2)
  \psline[linewidth=.5pt](1.5,0)(1.5,2)
\end{pspicture}
\rrbracket \to \llbracket \psset{xunit=.2cm,yunit=.2cm}
\begin{pspicture}[.2](2,2)
 \psbezier[linewidth=.5pt](1.5,0)(1.5,.5)(.5,.5)(.5,1)
 \psbezier[linewidth=.5pt](1.5,2)(1.5,1.5)(.5,1.5)(.5,1)
\pspolygon[linecolor=white,fillstyle=solid,fillcolor=white](.85,0)(.85,2)(1.15,2)(1.15,0)(.85,0)
  \psbezier[linewidth=.5pt](.5,0)(.5,.5)(1.5,.5)(1.5,1)
 \psbezier[linewidth=.5pt](.5,2)(.5,1.5)(1.5,1.5)(1.5,1)
\end{pspicture}
\rrbracket$ and $G\maps\llbracket\psset{xunit=.2cm,yunit=.2cm}
\begin{pspicture}[.2](2,2)
 \psbezier[linewidth=.5pt](1.5,0)(1.5,.5)(.5,.5)(.5,1)
 \psbezier[linewidth=.5pt](1.5,2)(1.5,1.5)(.5,1.5)(.5,1)
 \pspolygon[linecolor=white,fillstyle=solid,fillcolor=white](.85,0)(.85,2)(1.15,2)(1.15,0)(.85,0)
  \psbezier[linewidth=.5pt](.5,0)(.5,.5)(1.5,.5)(1.5,1)
 \psbezier[linewidth=.5pt](.5,2)(.5,1.5)(1.5,1.5)(1.5,1)
\end{pspicture}
\rrbracket  \to\llbracket \psset{xunit=.2cm,yunit=.2cm}
\begin{pspicture}[.2](2,2)
  \psline[linewidth=.5pt](.5,0)(.5,2)
  \psline[linewidth=.5pt](1.5,0)(1.5,2)
\end{pspicture}
\rrbracket$ whose only nonzero components are given by the
following diagram:
\begin{equation}
\label{eq_defofF}
 \xy
    (-55,20)*+{0}="1t";
    (0,20)*+{\underline{(A \otimes A)\{0\}}}="2t";
    (60,20)*+{0}="3t";
    (-55,-20)*+{(A \otimes A)\{-2\}}="1";
    (0,-20)*+{\underline{(C \otimes A \otimes A)\{0\}\oplus(A \otimes A)\{0\}}}="2";
    (60,-20)*+{(A \otimes A)\{2\}}="3";
    {\ar_-{d^{-1}=
    \left(
      \begin{array}{c}
     \psset{xunit=.15cm,yunit=.15cm}
        \begin{pspicture}(6,5)
         \rput(1,0){\identl}
         \rput(1.25,0){\ltc}
         \rput(4,0){\crossl}
         \rput(2,2.5){\comultl}
         \rput(5,2.5){\identl}
        \end{pspicture} \\
        \psset{xunit=.15cm,yunit=.15cm}
        \begin{pspicture}(4,5)
         \rput(4,0){\saddlel}
        \end{pspicture} \\
      \end{array}
    \right)
    } "1";"2"};
    {\ar_-{d^0=
    \left(
      \begin{array}{c}
     \psset{xunit=.15cm,yunit=.15cm}
        \begin{pspicture}(5,5)
         \rput(2,0){\multl}
         \rput(.75,2.5){\ctl}
         \rput(3,2.5){\identl}
         \rput(5,2.5){\identl}
         \rput(4,0){\curverightl}
        \end{pspicture} \\
        \psset{xunit=.15cm,yunit=.15cm}
        \begin{pspicture}(6,5)
        \rput(-.3,2.5){$-$}
         \rput(2,0){\saddler}
        \end{pspicture} \\
      \end{array}
    \right)^{\mT}
    } "2";"3"};
    {\ar "1t";"2t"};
    {\ar "2t";"3t"};
   {\ar@{<->} "1";"1t"};
   {\ar@{<->} "3";"3t"};
   {\ar_{
    F^0=\left(
      \begin{array}{c}
     \psset{xunit=.15cm,yunit=.15cm}
        \begin{pspicture}[.5](7,5)
         \rput(2,0){\saddler}
         \rput(0,0){\birthc}
        \end{pspicture}\\
        \psset{xunit=.15cm,yunit=.15cm}
        \begin{pspicture}(6,5)
         \rput(2,0){\identl}
         \rput(4,0){\identl}
        \end{pspicture} \\
      \end{array}
    \right)
    }
   "2t"+(-2,-4);"2"+(-2,4)};
   {\ar_{
    G^0=\left(
      \begin{array}{c}
     \psset{xunit=.15cm,yunit=.15cm}
        \begin{pspicture}[.5](7,5)
        \rput(-.3,2.5){$-$}
         \rput(3,0){\saddler}
         \rput(1,4){\deathc}
        \end{pspicture}\\
        \psset{xunit=.15cm,yunit=.15cm}
        \begin{pspicture}(6,5)
         \rput(2,0){\identl}
         \rput(4,0){\identl}
        \end{pspicture} \\
      \end{array}
    \right)^{\mT}
    }
     "2"+(2,4);"2t"+(2,-4)};
 \endxy
\end{equation}
One can readily check that these maps define chain maps. The equation
$d^0\circ F^0=0$ uses the fact that the zipper $\imath$ is an algebra
homomorphism (preserves the unit) and the left unit axiom for the
algebra $A$. The equation $G^0\circ d^{-1}=0$ follows
from~\eqref{eq_saddleinversion}.

Note that $G\circ F=\id$ by the relation \eqref{eq_barnatan1} so the
homotopy equivalence is established by defining a chain homotopy
$h\maps\id\to F\circ G$. The nonzero components of this chain homotopy
$h$ are given by
\begin{equation}
h^0 := \left(
      \begin{array}{c}
     \psset{xunit=.2cm,yunit=.2cm}
        \begin{pspicture}(7,2.5)
         \rput(4,0){\crossl}
         \rput(1,1.5){\deathc}
        \end{pspicture}\\
        0 \\
      \end{array}
    \right)^{\mT};
    \qquad \qquad
h^1 := \left(
      \begin{array}{c}
     \psset{xunit=.2cm,yunit=.2cm}
        \begin{pspicture}(7,2.5)
         \rput(3,0){\identl}
         \rput(5,0){\identl}
         \rput(1.2,0){\birthc}
        \end{pspicture}\\
        0 \\
      \end{array}
    \right).
\end{equation}
The equation $\id-F^0\circ G^0=h^1\circ d^{0}+d^{-1}\circ h^0$ is
depicted below
\begin{equation}
\left(
      \begin{array}{cc}
     \psset{xunit=.15cm,yunit=.15cm}
        \begin{pspicture}[.5](7,7)
        \rput(2,0){\identc}
        \rput(2,2.5){\identc}
         \rput(4,0){\identl}
         \rput(6,0){\identl}
         \rput(4,2.5){\identl}
         \rput(6,2.5){\identl}
        \end{pspicture} & 0\\
        0 &
        \psset{xunit=.15cm,yunit=.15cm}
        \begin{pspicture}(6,7)
         \rput(2,0){\identl}
         \rput(4,0){\identl}
         \rput(2,2.5){\identl}
         \rput(4,2.5){\identl}
        \end{pspicture} \\
      \end{array}
    \right)
- \left(
      \begin{array}{cc}
     \psset{xunit=.15cm,yunit=.15cm}
        \begin{pspicture}[.5](7,10)
        \rput(-.3,2.5){$-$}
         \rput(3,0){\saddler}
          \rput(3,5){\saddler}
         \rput(.8,9){\deathc}
         \rput(1.2,0){\birthc}
        \end{pspicture}
        &
        \psset{xunit=.15cm,yunit=.15cm}
        \begin{pspicture}[.5](9,5)
         \rput(5,0){\saddler}
         \rput(3.2,0){\birthc}
        \end{pspicture}\\
        \psset{xunit=.15cm,yunit=.15cm}
        \begin{pspicture}[.5](7,7)
        \rput(-.3,2.5){$-$}
         \rput(3,0){\saddler}
         \rput(1.2,4){\deathc}
        \end{pspicture}
       &
        \psset{xunit=.15cm,yunit=.15cm}
        \begin{pspicture}[.5](5,7)
         \rput(3,0){\identl}
         \rput(5,0){\identl}
         \rput(3,2.5){\identl}
         \rput(5,2.5){\identl}
        \end{pspicture} \\
      \end{array}
    \right)
    =
     \left(
      \begin{array}{cc}
     \psset{xunit=.15cm,yunit=.15cm}
        \begin{pspicture}(7,5)
         \rput(3,0){\multl}
         \rput(1.75,2.5){\ctl}
         \rput(4,2.5){\identl}
         \rput(6,2.5){\identl}
         \rput(5,0){\curverightl}
         \rput(1.4,0){\birthc}
        \end{pspicture} &
        \psset{xunit=.15cm,yunit=.15cm}
        \begin{pspicture}(8,7)
        \rput(.7,2.5){$-$}
         \rput(4,0){\saddler}
         \rput(2.2,0){\birthc}
        \end{pspicture}\\ \\
        0&
       0\\
      \end{array}
    \right)+ \left(
      \begin{array}{ccc}
     \psset{xunit=.15cm,yunit=.15cm}
        \begin{pspicture}[.5](6,10)
         \rput(2,0){\identl}
         \rput(2.25,0){\ltc}
         \rput(5,0){\crossl}
         \rput(3,2.5){\comultl}
         \rput(6,2.5){\curveleftl}
         \rput(4,5){\crossl}
         \rput(1,6.5){\deathc}
        \end{pspicture} && 0\\
        \psset{xunit=.15cm,yunit=.15cm}
        \begin{pspicture}[.5](7,9)
         \rput(7,0){\saddlel}
         \rput(1.2,6.5){\deathc}
         \rput(5,5){\widecrossl}
        \end{pspicture}
        &&
        0 \\
      \end{array}
    \right).
\end{equation}
Only the first and third component of this matrix equation are
non-trivial. The third component follows from
\eqref{eq_saddleinversion}. The first component can be written
\begin{equation} \label{eq_fourtuR2}
  \psset{xunit=.15cm,yunit=.15cm}
        \begin{pspicture}[.5](6,10)
         \rput(2,0){\identl}
         \rput(2.25,0){\ltc}
         \rput(5,0){\crossl}
         \rput(3,2.5){\comultl}
         \rput(6,2.5){\curveleftl}
         \rput(4,5){\crossl}
         \rput(1,6.5){\deathc}
        \end{pspicture}
 \quad + \quad
     \psset{xunit=.15cm,yunit=.15cm}
        \begin{pspicture}[.5](5,10)
         \rput(2,2.5){\multl}
         \rput(.75,5){\ctl}
         \rput(3,5){\identl}
         \rput(5,5){\identl}
         \rput(4,2.5){\curverightl}
         \rput(2,0){\identl}
         \rput(4,0){\identl}
         \rput(0.4,0){\birthc}
        \end{pspicture}
 \quad - \quad
  \psset{xunit=.15cm,yunit=.15cm}
        \begin{pspicture}[.5](7,7)
        \rput(2,0){\identc}
        \rput(2,2.5){\identc}
         \rput(4,0){\identl}
         \rput(6,0){\identl}
         \rput(4,2.5){\identl}
         \rput(6,2.5){\identl}
        \end{pspicture}
 \quad -  \quad
     \psset{xunit=.15cm,yunit=.15cm}
        \begin{pspicture}[.5](7,10)
         \rput(3,0){\saddler}
          \rput(3,5){\saddler}
         \rput(.8,9){\deathc}
         \rput(1.2,0){\birthc}
        \end{pspicture}
 \quad = \quad 0.
\end{equation}
Noting that
\begin{equation}
  \psset{xunit=.15cm,yunit=.15cm}
        \begin{pspicture}[.5](5,10)
         \rput(1,0){\saddler}
          \rput(1,5){\saddler}
        \end{pspicture}
  \quad = \quad
  \psset{xunit=.2cm,yunit=.2cm}
        \begin{pspicture}[.5](5,10)
         \rput(2,0){\multl}
          \rput(1,2.5){\medidentl}
          \rput(1,4.5){\medidentl}
          \rput(1,6.5){\identl}
           \rput(5,0){\identl}
          \rput(5,2.5){\medidentl}
          \rput(5,4.5){\medidentl}
          \rput(4,6.5){\comultl}
          \rput(2.8,2.5){\ctl}
          \rput(3.2,4.5){\ltc}
        \end{pspicture}
    \quad = \quad
  \psset{xunit=.2cm,yunit=.2cm}
        \begin{pspicture}[.5](5,12)
         \rput(2,0){\multl}
          \rput(3,2.5){\medidentl}
          \rput(3,4.5){\medidentl}
          \rput(1,9){\identl}
          \rput(2,6.5){\crossl}
          \rput(5,6.5){\identl}
           \rput(5,0){\identl}
          \rput(5,2.5){\medidentl}
          \rput(5,4.5){\medidentl}
          \rput(4,9){\comultl}
          \rput(.8,2.5){\ctl}
          \rput(1.2,4.5){\ltc}
        \end{pspicture}
\end{equation}
then \eqref{eq_fourtuR2} is just \eqref{eq_barnatan3} applied to the
cobordism
\begin{equation}
\psset{xunit=.2cm,yunit=.2cm}
        \begin{pspicture}[.5](11,19)
         \rput(1.2,0){\identc}
         \rput(1.2,2.5){\medidentc}
         \rput(3.9,2.5){\zigc}
         \rput(1.2,4.5){\identc}
         \rput(3.2,4.5){\identc}
         \rput(1.2,7){\medidentc}
         \rput(3.2,7){\medidentc}
         \rput(1.2,9){\identc}
         \rput(4.2,9){\crossc}
         \rput(7.2,9){\identc}
          \rput(2.2,11.5){\zagc}
          \rput(7.1,11.5){\ltc}
          \rput(5.1,11.5){\medidentc}
          \rput(5.1,13.5){\identc}
          \rput(5.1,16){\identc}
         \rput(8,0){\multl}
         \rput(11,0){\identl}
         \rput(11,2.5){\medidentl}
         \rput(11,4.5){\identl}
         \rput(6.8,2.5){\ctl}
         \rput(9,2.5){\medidentl}
         \rput(9,4.5){\identl}
         \rput(9,7){\medidentl}
         \rput(11,7){\medidentl}
         \rput(9,9){\identl}
         \rput(11,9){\identl}
         \rput(9,11.5){\medidentl}
         \rput(11,11.5){\medidentl}
         \rput(8,13.5){\crossl}
         \rput(11,13.5){\identl}
         \rput(10,16){\comultl}
         \rput(7,16){\identl}
         \pspolygon[linestyle=dashed](4,4)(4,9.2)(8,9.2)(8,4)(4,4)
        \end{pspicture}
\end{equation}
The equalities $\id-F^{-1}\circ G^{-1}=h^{0}\circ d^{-1}$ and
$\id-F^{1}\circ G^{1}=h^{0}\circ d^{-1}+d^{0}\circ h^1$ follow from
the identities
\[
  \psset{xunit=.2cm,yunit=.2cm}
        \begin{pspicture}[.5](4,5)
         \rput(2,3){\comultl}
         \rput(1,0){\deathc}
         \rput(1.2,1){\ltc}
         \rput(3,1){\medidentl}
        \end{pspicture}
        \;\; = \;\;
         \psset{xunit=.2cm,yunit=.2cm}
        \begin{pspicture}[.5](4,5)
         \rput(1,0){\identl}
         \rput(1,2.5){\identl}
        \end{pspicture}
        \qquad {\rm and} \qquad
    \psset{xunit=.2cm,yunit=.2cm}
        \begin{pspicture}[.5](4,5)
         \rput(2,0){\multl}
         \rput(.8,2.5){\ctl}
         \rput(1.3,4.5){\birthc}
         \rput(3,2.5){\medidentl}
        \end{pspicture}
        \;\; = \;\;
         \psset{xunit=.2cm,yunit=.2cm}
        \begin{pspicture}[.5](3,5)
         \rput(1,0){\identl}
         \rput(1,2.5){\identl}
        \end{pspicture},
\]
respectively.

For the opposite checkerboard colouring, $\epsilon=-1$, of the
complexes in~\eqref{eq_ffvi} and~\eqref{eq_ffvii} only the complex
\begin{equation}
\llbracket \;
\begin{pspicture}[.4](2,2.1)
  \rput(-.5,0){ \pspolygon[linecolor=white,
   linewidth=.6pt,fillstyle=solid,fillcolor=white](.5,0)(.5,2.03)(2.5,2.03)(2.5,0)(.5,0)
  \pscustom[linecolor=lightgray,linewidth=0pt,fillstyle=solid,fillcolor=lightgray]{
  \psbezier[linewidth=.8pt](2,0)(2,.5)(1,.5)(1,1)
   \psbezier[linewidth=.8pt](1,1)(1,1.5)(2,1.5)(2,2)
   \psline(1,2)
   \psbezier[linewidth=.8pt](1,2)(1,1.5)(2,1.5)(2,1)
   \psbezier[linewidth=.8pt](2,1)(2,.5)(1,.5)(1,0)
   \psline(2,0)
  }}
 \psbezier[linewidth=.8pt](1.5,0)(1.5,.5)(.5,.5)(.5,1)
 \psbezier[linewidth=.8pt](1.5,2)(1.5,1.5)(.5,1.5)(.5,1)
 \pspolygon[linecolor=lightgray,fillstyle=solid,fillcolor=lightgray](.9,0)(.9,2)(1.1,2)(1.1,0)(.9,0)
  \psbezier[linewidth=.8pt](.5,0)(.5,.5)(1.5,.5)(1.5,1)
 \psbezier[linewidth=.8pt](.5,2)(.5,1.5)(1.5,1.5)(1.5,1)
  \rput(.35,1){$\scriptstyle 5$}
  \rput(1.65,1){$\scriptstyle 1$}
  \rput(.35,1.7){$\scriptstyle 2$}
  \rput(1.65,1.7){$\scriptstyle 4$}
  \rput(.35,.3){$\scriptstyle 3$}
  \rput(1.65,.3){$\scriptstyle 6$}
\end{pspicture} \;
\rrbracket
\end{equation}
would change. In particular, one can verify that the complex becomes
\begin{equation}
\llbracket \;
\begin{pspicture}[.4](2,2.1)
  \rput(-.5,0){ \pspolygon[linecolor=white,
   linewidth=.6pt,fillstyle=solid,fillcolor=white](.5,0)(.5,2.03)(2.5,2.03)(2.5,0)(.5,0)
  \pscustom[linecolor=lightgray,linewidth=0pt,fillstyle=solid,fillcolor=lightgray]{
  \psbezier[linewidth=.8pt](2,0)(2,.5)(1,.5)(1,1)
   \psbezier[linewidth=.8pt](1,1)(1,1.5)(2,1.5)(2,2)
   \psline(1,2)
   \psbezier[linewidth=.8pt](1,2)(1,1.5)(2,1.5)(2,1)
   \psbezier[linewidth=.8pt](2,1)(2,.5)(1,.5)(1,0)
   \psline(2,0)
  }}
 \psbezier[linewidth=.8pt](1.5,0)(1.5,.5)(.5,.5)(.5,1)
 \psbezier[linewidth=.8pt](1.5,2)(1.5,1.5)(.5,1.5)(.5,1)
 \pspolygon[linecolor=lightgray,fillstyle=solid,fillcolor=lightgray](.9,0)(.9,2)(1.1,2)(1.1,0)(.9,0)
  \psbezier[linewidth=.8pt](.5,0)(.5,.5)(1.5,.5)(1.5,1)
 \psbezier[linewidth=.8pt](.5,2)(.5,1.5)(1.5,1.5)(1.5,1)
  \rput(.35,1){$\scriptstyle 5$}
  \rput(1.65,1){$\scriptstyle 1$}
  \rput(.35,1.7){$\scriptstyle 2$}
  \rput(1.65,1.7){$\scriptstyle 4$}
  \rput(.35,.3){$\scriptstyle 3$}
  \rput(1.65,.3){$\scriptstyle 6$}
\end{pspicture} \;
\rrbracket :=
 \xy
    (-45,0)*+{\scs(A \otimes A)\{-2\}}="1";
    (0,0)*+{\scs \underline{(C \otimes A \otimes A)\{0\}\oplus(A \otimes A)\{0\}}}="2";
    (45,0)*+{\scs (A \otimes A)\{2\}}="3";
    {\ar^-{d^{-1}=
    \left(
      \begin{array}{c}
     \psset{xunit=.15cm,yunit=.15cm}
        \begin{pspicture}(6,5)
         \rput(1,0){\identl}
         \rput(1.25,0){\ltc}
         \rput(4,0){\crossl}
         \rput(2,2.5){\comultl}
         \rput(5,2.5){\identl}
        \end{pspicture} \\
        \psset{xunit=.15cm,yunit=.15cm}
        \begin{pspicture}(4,5)
         \rput(0,0){\saddler}
        \end{pspicture} \\
      \end{array}
    \right)
    } "1";"2"};
    {\ar^-{d^0=
    \left(
      \begin{array}{c}
     \psset{xunit=.15cm,yunit=.15cm}
        \begin{pspicture}(5,5)
         \rput(2,0){\multl}
         \rput(.75,2.5){\ctl}
         \rput(3,2.5){\identl}
         \rput(5,2.5){\identl}
         \rput(4,0){\curverightl}
        \end{pspicture} \\
        \psset{xunit=.15cm,yunit=.15cm}
        \begin{pspicture}(6,5)
        \rput(-.3,2.5){$-$}
         \rput(5,0){\saddlel}
        \end{pspicture} \\
      \end{array}
    \right)^{\mT}
    } "2";"3"};
 \endxy
\end{equation}
where only the saddle~\eqref{eq_crossing} has been replaced
by~\eqref{eq_crossing2} and vice versa. Performing a similar saddle
switch for the chain maps $F$ and $G$ produces the required homotopy
equivalence. Again, all maps $F^j$, $G^j$, and $h^j$ are graded
[filtered] of degree $0$. Notice that the chain map $G\maps
\llbracket \psset{xunit=.2cm,yunit=.2cm}
\begin{pspicture}[.2](2,2)
 \psbezier[linewidth=.5pt](1.5,0)(1.5,.5)(.5,.5)(.5,1)
 \psbezier[linewidth=.5pt](1.5,2)(1.5,1.5)(.5,1.5)(.5,1)
 \pspolygon[linecolor=white,fillstyle=solid,fillcolor=white](.85,0)(.85,2)(1.15,2)(1.15,0)(.85,0)
  \psbezier[linewidth=.5pt](.5,0)(.5,.5)(1.5,.5)(1.5,1)
 \psbezier[linewidth=.5pt](.5,2)(.5,1.5)(1.5,1.5)(1.5,1)
\end{pspicture}
\rrbracket \to\llbracket\psset{xunit=.2cm,yunit=.2cm}
\begin{pspicture}[.2](2,2)
  \psline[linewidth=.5pt](.5,0)(.5,2)
  \psline[linewidth=.5pt](1.5,0)(1.5,2)
\end{pspicture} \rrbracket $ is actually a \emph{strong
deformation retract} since $h\circ F=0$ by~\eqref{eq_barnatan1} and we
have already shown that $G\circ F=\id$ and $\id-F\circ G=d\circ
h+h\circ d$. Following the standard terminology we say that $F$ is the
\emph{inclusion} in a strong deformation retract.

\subsubsection{Reidemeister move three}

To prove that the complex $\llbracket T \rrbracket$ is invariant up to
homotopy equivalence under the third Reidemeister move, we require the
following results from homological algebra.

\begin{defn}
Let $\Psi \maps (C^r_0,d_0) \to (C^r_1,d_1)$ be a morphism of
complexes. The \emph{cone} $\Gamma(\Psi)$ of $\Psi$ is the complex
with terms $\Gamma^r(\Psi) = C^{r+1}_0 \oplus C_1^r$ and the
differentials $d^r=\left(
\begin{array}{cc}
  -d_0^{r+1} & 0 \\
  \Psi^{r+1} & d_1^r \\
    \end{array}
   \right)$.
\end{defn}

If $\Psi$ is a morphism of graded [filtered] complexes, then the
cone $\Gamma(\Psi)$ is a graded [filtered] complex.

\parpic[r]{
$ \xy
 (0,0)*+{C_{0a}}="t";
 (0,-10)*+{C_{1a}}="b";
 (20,-10)*+{C_{1a}}="b2";
 (20,0)*+{C_{0b}}="t2";
  {\ar_{\Psi} "t";"b"};
  {\ar@<.5ex>^{G_0} "t";"t2"};
  {\ar@<.5ex>^{F_0} "t2";"t"};
  {\ar@<.5ex>^{F_1} "b";"b2"};
  {\ar@<.5ex>^{G_1} "b2";"b"};
 \endxy
$ }
\begin{lem}
\label{lem_natans}
The cone construction is invariant up to homotopy equivalence under
compositions with the inclusions in strong deformation retracts. That
is, consider the diagram of morphisms and complexes on the right.  If
in that diagram $G^0$ is a strong deformation retract with inclusion
$F_0$, then the cones $\Gamma(\Psi)$ and $\Gamma(\Psi\circ F_0)$ are
homotopy equivalent, and if $G_1$ is a strong deformation retract with
inclusion $F_1$, then the cones $\Gamma(\Psi)$ and $\Gamma(F_1\circ
\Psi)$ are homotopy equivalent. Likewise, if $F_{0,1}$ are strong
deformation retracts and $G_{0,1}$ the corresponding inclusions, the
above statements remain true.
\end{lem}

This lemma holds for complexes in any abelian category, \ie\ if
$\Psi$ and $F_0$ are morphisms of graded [filtered] complexes, then
the cones $\Gamma(\Psi)$ and $\Gamma(\Psi\circ F_0)$ are homotopy
equivalent as graded [filtered] complexes, \etc.

We comment here that the proof of the above lemma is constructive,
so that one can explicitly obtain the chain homotopies defining the
above homotopy equivalence. In fact, Bar-Natan explicitly lists the
chain maps and chain homotopy for Reidemeister move three which we
encourage the reader to translate into the language of open-closed
cobordisms.

For each $\xi\in\Z$ define the two chain maps $\cal{S}_1^{\xi}$ and
$\cal{S}_2^{\xi}$ given by the diagram below:
\begin{equation}
\label{eq_S1nS2}
\cal{S}_1^{\xi} := \quad
 \xy
    (-22,10)*+{0}="1t";
    (0,10)*+{\underline{(A \otimes A)\{\xi\}}}="2t";
    (22,10)*+{0}="3t";
    (-22,-10)*+{0}="1";
    (0,-10)*+{\underline{(A \otimes A)\{\xi+2\}}}="2";
    (22,-10)*+{0}="3";
    {\ar "1";"2"};
    {\ar "2";"3"};
    {\ar "1t";"2t"};
    {\ar "2t";"3t"};
   {\ar_{
     \psset{xunit=.15cm,yunit=.15cm}
        \begin{pspicture}[.5](5,5)
         \rput(4,0){\saddlel}
        \end{pspicture}\\
    }
   "2t"+(-2,-4);"2"+(-2,4)};
 \endxy
 \qquad \qquad
 \cal{S}_2^{\xi} := \quad
 \xy
    (-22,10)*+{0}="1t";
    (0,10)*+{\underline{(A \otimes A)\{\xi\}}}="2t";
    (22,10)*+{0}="3t";
    (-22,-10)*+{0}="1";
    (0,-10)*+{\underline{(A \otimes A)\{\xi+2\}}}="2";
    (22,-10)*+{0}="3";
    {\ar "1";"2"};
    {\ar "2";"3"};
    {\ar "1t";"2t"};
    {\ar "2t";"3t"};
   {\ar_{
     \psset{xunit=.15cm,yunit=.15cm}
        \begin{pspicture}[.5](5.5,5)
         \rput(1,0){\saddler}
        \end{pspicture}\\
    }
   "2t"+(-2,-4);"2"+(-2,4)};
 \endxy
\end{equation}
These are morphisms of graded [filtered] complexes.

\begin{lem}
\label{lem_crossings}
With $\cal{S}_1^{\xi}$ and $\cal{S}_2^{\xi}$ as in \eqref{eq_S1nS2} we
have that
\begin{equation}  \psset{xunit=.8cm,yunit=.8cm}
 \llbracket\;\;
 \begin{pspicture}[.3](1,1.2)
 \pspolygon[linecolor=whitegray,fillstyle=solid,fillcolor=whitegray](0,0)(0,1)(.5,.5)(0,0)
 \pspolygon[linecolor=whitegray,fillstyle=solid,fillcolor=whitegray](1,1)(1,0)(.5,.5)(1,1)
 \psline[linewidth=1.0pt](0,0)(.4,.4)
 \psline[linewidth=1.0pt]{->}(.6,.6)(1,1)
 \psline[linewidth=1.0pt]{->}(1,0)(0,1)
  \rput(0,-.2){$\scriptstyle 3$}
 \rput(1,-.2){$\scriptstyle 4$}
 \rput(0,1.2){$\scriptstyle 1$}
 \rput(1,1.2){$\scriptstyle 2$}
 \end{pspicture}\;\;\rrbracket = \Gamma(\cal{S}_1^{(-4)}),
 \quad
  \llbracket\;\;
  \begin{pspicture}[.3](1,1.2)
 \pspolygon[linecolor=whitegray,fillstyle=solid,fillcolor=whitegray](0,0)(1,0)(.5,.5)(0,0)
 \pspolygon[linecolor=whitegray,fillstyle=solid,fillcolor=whitegray](0,1)(1,1)(.5,.5)(0,1)
 \psline[linewidth=1.0pt](0,0)(.4,.4)
 \psline[linewidth=1.0pt]{->}(.6,.6)(1,1)
 \psline[linewidth=1.0pt]{->}(1,0)(0,1)
  \rput(0,-.2){$\scriptstyle 3$}
 \rput(1,-.2){$\scriptstyle 4$}
 \rput(0,1.2){$\scriptstyle 1$}
 \rput(1,1.2){$\scriptstyle 2$}
\end{pspicture} \;\;\rrbracket = \Gamma(\cal{S}_2^{(-4)}),
\end{equation}
\begin{equation}
 \llbracket \;\;
 \begin{pspicture}[.3](1,1.2)
 \pspolygon[linecolor=whitegray,fillstyle=solid,fillcolor=whitegray](0,0)(0,1)(.5,.5)(0,0)
 \pspolygon[linecolor=whitegray,fillstyle=solid,fillcolor=whitegray](1,1)(1,0)(.5,.5)(1,1)
 \psline[linewidth=1.0pt]{->}(0,0)(1,1)
 \psline[linewidth=1.0pt](1,0)(.6,.4)
 \psline[linewidth=1.0pt]{->}(.4,.6)(0,1)
  \rput(0,-.2){$\scriptstyle 3$}
 \rput(1,-.2){$\scriptstyle 4$}
 \rput(0,1.2){$\scriptstyle 1$}
 \rput(1,1.2){$\scriptstyle 2$}
 \end{pspicture}\;\;\rrbracket = \Gamma(\cal{S}_1^{(2)})[1],
\quad
 \llbracket \;\;
\begin{pspicture}[.3](1,1.2)
 \pspolygon[linecolor=whitegray,fillstyle=solid,fillcolor=whitegray](0,0)(1,0)(.5,.5)(0,0)
 \pspolygon[linecolor=whitegray,fillstyle=solid,fillcolor=whitegray](0,1)(1,1)(.5,.5)(0,1)
 \psline[linewidth=1pt]{->}(0,0)(1,1)
 \psline[linewidth=1pt](1,0)(.6,.4)
 \psline[linewidth=1pt]{->}(.4,.6)(0,1)
  \rput(0,-.2){$\scriptstyle 3$}
 \rput(1,-.2){$\scriptstyle 4$}
 \rput(0,1.2){$\scriptstyle 1$}
 \rput(1,1.2){$\scriptstyle 2$}
\end{pspicture} \;\;\rrbracket = \Gamma(\cal{S}_1^{(2)})[1]
\end{equation}
where $[s]$ is the cohomological degree shift.
\end{lem}

By the discussion above the following complex
\begin{equation} \label{eq_R3l}
\llbracket
\begin{pspicture}[.4](3,2)
 \psbezier[linewidth=.8pt](1.5,2)(1.5,1.5)(.5,1.5)(.5,1)
 \psbezier[linewidth=.8pt](1.5,0)(1.5,.5)(.5,.5)(.5,1)
   \pspolygon[linecolor=white,fillstyle=solid,fillcolor=white](.8,0)(.8,2)(1.1,2)(1.1,0)(.8,0)
   \psbezier[linewidth=.8pt](2.5,2)(2.5,1.25)(.5,.5)(.5,0)
   \pspolygon[linecolor=white,fillstyle=solid,fillcolor=white](1.5,0)(1.5,2)(1.8,2)(1.8,0)(1.5,0)
 \psbezier[linewidth=.8pt](.5,2)(.5,1.5)(2.5,.75)(2.5,0)
  \rput(.35,1){$\scriptscriptstyle 5$}
  \rput(1.1,.9){$\scriptscriptstyle 1$}
  \rput(.35,1.7){$\scriptscriptstyle 2$}
  \rput(1.65,1.7){$\scriptscriptstyle 4$}
  \rput(.35,.3){$\scriptscriptstyle 3$}
  \rput(1.65,.3){$\scriptscriptstyle 6$}
  \rput(2.65,1.7){$\scriptscriptstyle 7$}
  \rput(2.65,.3){$\scriptscriptstyle 8$}
  \rput(1.5,1.4){$\scriptscriptstyle 9$}
  \rput(2.1,1){$\scriptstyle \underline{1}$}
  \rput(1,1.85){$\scriptstyle \underline{2}$}
  \rput(1,.15){$\scriptstyle \underline{3}$}
\end{pspicture}
\rrbracket ,
\end{equation}
with the region to the right of the first strand shaded, is homotopy
equivalent to the cone over the chain map
\begin{equation}
\Psi \maps  \llbracket
\begin{pspicture}[.4](3,2)
 \psbezier[linewidth=.8pt](1.5,0)(1.5,.5)(.5,.5)(.5,1)
 \psbezier[linewidth=.8pt](1.5,2)(1.5,1.5)(.5,1.5)(.5,1)
 \pspolygon[linecolor=white,fillstyle=solid,fillcolor=white](1.3,0)(1.3,2)(1.0,2)(1.0,0)(1.3,0)
 \psbezier[linewidth=.8pt](.5,2)(1,2)(1,1.5)(1.5,1.5)
 \psbezier[linewidth=.8pt](2.5,2)(2,2)(2,1.5)(1.5,1.5)
 \psbezier[linewidth=.8pt](.5,0)(1,0)(1,.5)(1.5,.5)
 \psbezier[linewidth=.8pt](2.5,0)(2,0)(2,.5)(1.5,.5)
\end{pspicture}
\rrbracket\quad\to\quad\llbracket
\begin{pspicture}[.4](3,2)
 \psbezier[linewidth=.8pt](1.5,0)(1.5,.5)(.5,.5)(.5,1)
 \psbezier[linewidth=.8pt](1.5,2)(1.5,1.5)(.5,1.5)(.5,1)
 \pspolygon[linecolor=white,fillstyle=solid,fillcolor=white](.9,0)(.9,2)(1.1,2)(1.1,0)(.9,0)
  \psbezier[linewidth=.8pt](.5,0)(.5,.5)(1.5,.5)(1.5,1)
 \psbezier[linewidth=.8pt](.5,2)(.5,1.5)(1.5,1.5)(1.5,1)
 \psbezier[linewidth=.8pt](2.5,0)(2.5,.5)(2,.5)(2,1)
 \psbezier[linewidth=.8pt](2.5,2)(2.5,1.5)(2,1.5)(2,1)
\end{pspicture}
\rrbracket
\end{equation}
 given by
\begin{equation}
\label{eq_psi}
\Psi := \quad
 \xy
    (-55,-22)*+{\scs (A\otimes A \otimes A)\{ -2\}}="1";
    (0,-22)*+{\scs (C\otimes A\otimes A \otimes A)\{ 0\}
                \oplus(A \otimes A \otimes A)     \{ 0\} }="2";
    (50,-22)*+{\scs (A\otimes A \otimes A)\{ 2\}}="3";
    (-55,22)*+{\scs (A\otimes A \otimes A)\{ -4\}}="1'";
    (0,22)*+{\scs (A\otimes A \otimes A)\{ -2\}\oplus (A\otimes A \otimes A)\{ -2\}}="2'";
    (50,22)*+{\scs (A\otimes A \otimes A)\{ 0\}}="3'";
    {\ar_-{
    \left(
      \begin{array}{c}
     \psset{xunit=.15cm,yunit=.15cm}
        \begin{pspicture}(8,5)
         \rput(1,0){\identl}
         \rput(1.25,0){\ltc}
         \rput(4,0){\crossl}
         \rput(2,2.5){\comultl}
         \rput(5,2.5){\identl}
         \rput(7,2.5){\identl}
         \rput(7,0){\identl}
        \end{pspicture} \\
        \psset{xunit=.15cm,yunit=.15cm}
        \begin{pspicture}(6,5)
         \rput(4,0){\saddlel}
         \rput(6,2.5){\identl}
         \rput(6,0){\identl}
        \end{pspicture} \\
      \end{array}
    \right)
    } "1";"2"};
    {\ar_-{
    \left(
      \begin{array}{c}
     \psset{xunit=.15cm,yunit=.15cm}
        \begin{pspicture}(8,5)
         \rput(2,0){\multl}
         \rput(.75,2.5){\ctl}
         \rput(3,2.5){\identl}
         \rput(5,2.5){\identl}
         \rput(4,0){\curverightl}
         \rput(7,2.5){\identl}
         \rput(7,0){\identl}
        \end{pspicture} \\
        \psset{xunit=.15cm,yunit=.15cm}
        \begin{pspicture}(9,5)
        \rput(-.3,2.5){$-$}
         \rput(2,0){\saddler}
         \rput(8,2.5){\identl}
         \rput(8,0){\identl}
        \end{pspicture} \\
      \end{array}
    \right)^{\mT}
    } "2";"3"};
    {\ar^-{
    \left(
      \begin{array}{c}
     \psset{xunit=.15cm,yunit=.15cm}
        \begin{pspicture}(8,10)
         \rput(1,0){\identl}
         \rput(1,7.5){\identl}
         \rput(6,0){\crossl}
         \rput(7,2.5){\identl}
         \rput(7,5){\identl}
         \rput(6,7.5){\crossl}
         \rput(1,2.5){\saddler}
        \end{pspicture} \\
        \psset{xunit=.15cm,yunit=.15cm}
        \begin{pspicture}(8,5)
         \rput(7,0){\saddlel}
         \rput(1,2.5){\identl}
         \rput(1,0){\identl}
        \end{pspicture} \\
      \end{array}
    \right)
    } "1'";"2'"};
    {\ar^-{
    \left(
      \begin{array}{c}
     \psset{xunit=.15cm,yunit=.15cm}
       \begin{pspicture}(8,7.5)
         \rput(1,0){\identl}
         \rput(6,0){\crossl}
         \rput(7,2.5){\identl}
         \rput(7,5){\identl}
         \rput(1,2.5){\saddler}
        \end{pspicture}  \\
        \psset{xunit=.15cm,yunit=.15cm}
        \begin{pspicture}(9,5)
        \rput(-.3,2.5){$-$}
         \rput(2,0){\saddler}
         \rput(8,2.5){\identl}
         \rput(8,0){\identl}
        \end{pspicture} \\
      \end{array}
    \right)^{\mT}
    } "2'";"3'"};
 {\ar_-{
 \psset{xunit=.15cm,yunit=.15cm}
        \begin{pspicture}(9,10)
         \rput(2,0){\crossl}
         \rput(2,7.5){\crossl}
         \rput(7,0){\identl}
         \rput(1,2.5){\identl}
         \rput(1,5){\identl}
         \rput(7,7.5){\identl}
         \rput(7,2.5){\saddlel}
        \end{pspicture}
 } "1'";"1"};
  {\ar^-{
 \psset{xunit=.15cm,yunit=.15cm}
        \begin{pspicture}(9,10)
         \rput(2,0){\crossl}
         \rput(2,7.5){\crossl}
         \rput(7,0){\identl}
         \rput(1,2.5){\identl}
         \rput(1,5){\identl}
         \rput(7,7.5){\identl}
         \rput(7,2.5){\saddlel}
        \end{pspicture}
 } "3'";"3"};
     {\ar_-{
    \left(
      \begin{array}{cc}
     \psset{xunit=.15cm,yunit=.15cm}
       \begin{pspicture}(8,7.5)
         \rput(2,5){\comultl}
         \rput(4,2.5){\crossl}
         \rput(6,0){\crossl}
         \rput(5,5){\identl}
         \rput(3,0){\identl}
         \rput(1,2.5){\identl}
         \rput(1.4,0){\identc}
         \rput(1.3,2.5){\ltc}
         \rput(7,5){\identl}
         \rput(7,2.5){\identl}
        \end{pspicture} & 0 \\
        0 &
         \psset{xunit=.15cm,yunit=.15cm}
        \begin{pspicture}(9,10)
         \rput(2,0){\crossl}
         \rput(2,7.5){\crossl}
         \rput(7,0){\identl}
         \rput(1,2.5){\identl}
         \rput(1,5){\identl}
         \rput(7,7.5){\identl}
         \rput(7,2.5){\saddlel}
        \end{pspicture} \\
      \end{array}
    \right)
    } "2'";"2"};
 \endxy
\end{equation}
but by Lemma~\ref{lem_natans} the
cone over the chain map $\Psi$ is homotopy equivalent to the cone over
the composite of $\Psi$ with any strong deformation retract $G$.
Taking $G$ as in~\eqref{eq_defofF} we have that
$\llbracket\psset{xunit=.2cm,yunit=.2cm}
\begin{pspicture}[.3](3,2)
 \psbezier[linewidth=.6pt](1.5,2)(1.5,1.5)(.5,1.5)(.5,1)
 \psbezier[linewidth=.6pt](1.5,0)(1.5,.5)(.5,.5)(.5,1)
   \pspolygon[linecolor=white,fillstyle=solid,fillcolor=white](.8,0)(.8,2)(1.1,2)(1.1,0)(.8,0)
   \psbezier[linewidth=.6pt](2.5,2)(2.5,1.25)(.5,.5)(.5,0)
   \pspolygon[linecolor=white,fillstyle=solid,fillcolor=white](1.5,0)(1.5,2)(1.8,2)(1.8,0)(1.5,0)
 \psbezier[linewidth=.6pt](.5,2)(.5,1.5)(2.5,.75)(2.5,0)
\end{pspicture} \rrbracket $ is homotopy equivalent to the cone of
the map $\Psi' = G\circ \Psi$ below
\begin{equation} \label{eq_psi'}
\Psi' := \quad
 \xy
    (-65,22)*{\llbracket
    \psset{xunit=.2cm,yunit=.2cm}
\begin{pspicture}[.4](3,2)
 \psbezier[linewidth=.6pt](1.5,0)(1.5,.5)(.5,.5)(.5,1)
 \psbezier[linewidth=.6pt](1.5,2)(1.5,1.5)(.5,1.5)(.5,1)
 \pspolygon[linecolor=white,fillstyle=solid,fillcolor=white](1.3,0)(1.3,2)(1.0,2)(1.0,0)(1.3,0)
 \psbezier[linewidth=.6pt](.5,2)(1,2)(1,1.5)(1.5,1.5)
 \psbezier[linewidth=.6pt](2.5,2)(2,2)(2,1.5)(1.5,1.5)
 \psbezier[linewidth=.6pt](.5,0)(1,0)(1,.5)(1.5,.5)
 \psbezier[linewidth=.6pt](2.5,0)(2,0)(2,.5)(1.5,.5)
\end{pspicture}
\rrbracket =}; (-65,-22)*{\llbracket
    \psset{xunit=.4cm,yunit=.4cm}
\begin{pspicture}[.4](3,2)
 \psline[linewidth=.6pt](.5,0)(.5,2)
 \psline[linewidth=.6pt](1.5,0)(1.5,2)
 \psline[linewidth=.6pt](2.5,0)(2.5,2)
 \rput(.3,1){$\scriptstyle 2$}
 \rput(1.3,1){$\scriptstyle 1$}
 \rput(2.3,1){$\scriptstyle 7$}
\end{pspicture}
\rrbracket =};
    (-45,-22)*+{0}="1";
    (0,-22)*+{\scs (A \otimes A \otimes A)\{ 0\}}="2";
    (45,-22)*+{0}="3";
    (-45,22)*+{\scs (A \otimes A \otimes A)\{ -4\}}="1'";
    (0,22)*+{\scs (A \otimes A \otimes A)\{ -2\}
    \oplus (A \otimes A \otimes A)\{ -2\}}="2'";
    (45,22)*+{\scs (A \otimes A \otimes A)\{ 0\}}="3'";
    {\ar "1";"2"};
    {\ar "2";"3"};
    {\ar^-{
    \left(
      \begin{array}{c}
     \psset{xunit=.15cm,yunit=.15cm}
        \begin{pspicture}(8,10)
         \rput(1,0){\identl}
         \rput(1,7.5){\identl}
         \rput(6,0){\crossl}
         \rput(7,2.5){\identl}
         \rput(7,5){\identl}
         \rput(6,7.5){\crossl}
         \rput(1,2.5){\saddler}
        \end{pspicture} \\
        \psset{xunit=.15cm,yunit=.15cm}
        \begin{pspicture}(8,5)
         \rput(7,0){\saddlel}
         \rput(1,2.5){\identl}
         \rput(1,0){\identl}
        \end{pspicture} \\
      \end{array}
    \right)
    } "1'";"2'"};
    {\ar^-{
    \left(
      \begin{array}{c}
     \psset{xunit=.15cm,yunit=.15cm}
       \begin{pspicture}(8,7.5)
         \rput(1,0){\identl}
         \rput(6,0){\crossl}
         \rput(7,2.5){\identl}
         \rput(7,5){\identl}
         \rput(1,2.5){\saddler}
        \end{pspicture}  \\
        \psset{xunit=.15cm,yunit=.15cm}
        \begin{pspicture}(9,5)
        \rput(-.3,2.5){$-$}
         \rput(2,0){\saddler}
         \rput(8,2.5){\identl}
         \rput(8,0){\identl}
        \end{pspicture} \\
      \end{array}
    \right)^{\mT}
    } "2'";"3'"};
 {\ar "1'";"1"};
  {\ar "3'";"3"};
     {\ar_-{
    \left(
      \begin{array}{c}
     \psset{xunit=.15cm,yunit=.15cm}
       \begin{pspicture}(8,13)
       \rput(-.5,5){$-$}
         \rput(2,10){\crossl}
         \rput(1,5){\identl}
         \rput(1,7.5){\identl}
         \rput(4,7.5){
  \pspolygon[fillcolor=lightgray,fillstyle=solid](1.5,0)(.5,0)(-1.5,2.5)(-.5,2.5)(1.5,0)}
         \rput(1,0){\saddler}
         \rput(6,5){\crossl}
         \rput(7,0){\identl}
         \rput(7,2.5){\identl}
         \rput(7,7.5){\identl}
         \rput(7,10){\identl}
        \end{pspicture}  \\
         \psset{xunit=.15cm,yunit=.15cm}
        \begin{pspicture}(9,10)
         \rput(2,0){\crossl}
         \rput(2,7.5){\crossl}
         \rput(7,0){\identl}
         \rput(1,2.5){\identl}
         \rput(1,5){\identl}
         \rput(7,7.5){\identl}
         \rput(7,2.5){\saddlel}
        \end{pspicture} \\
      \end{array}
    \right)^{\mT}
    } "2'";"2"};
 \endxy
\end{equation}
where we have left the enumeration on the bottom complex for
clarity. Note that depending on the orientation of the tangle
in~\eqref{eq_R3l} and on the number $n_-$ of negative crossings, the
complexes~\eqref{eq_psi} and~\eqref{eq_psi'} can have either of
their three terms in cohomological degree $0$.

Similarly, by Lemma~\ref{lem_crossings} the cone over the chain map
\begin{equation}
 \Phi
  \maps \quad
  \llbracket
\begin{pspicture}[.4](3,2)
\psbezier[linewidth=.8pt](1.5,0)(1.5,.5)(2.5,.5)(2.5,1)
 \psbezier[linewidth=.8pt](1.5,2)(1.5,1.5)(2.5,1.5)(2.5,1)
 \pspolygon[linecolor=white,fillstyle=solid,fillcolor=white](1.65,0)(1.65,2)(1.95,2)(1.95,0)(1.65,0)
 \psbezier[linewidth=.8pt](.5,2)(1,2)(1,1.5)(1.5,1.5)
 \psbezier[linewidth=.8pt](2.5,2)(2,2)(2,1.5)(1.5,1.5)
 \psbezier[linewidth=.8pt](.5,0)(1,0)(1,.5)(1.5,.5)
 \psbezier[linewidth=.8pt](2.5,0)(2,0)(2,.5)(1.5,.5)
\end{pspicture}
\rrbracket \quad\to\quad \llbracket
\begin{pspicture}[.4](3,2)
 \psbezier[linewidth=.8pt](1.5,0)(1.5,.5)(2.5,.5)(2.5,1)
 \psbezier[linewidth=.8pt](1.5,2)(1.5,1.5)(2.5,1.5)(2.5,1)
 \pspolygon[linecolor=white,fillstyle=solid,fillcolor=white](1.9,0)(1.9,2)(2.1,2)(2.1,0)(1.9,0)
\psbezier[linewidth=.8pt](2.5,0)(2.5,.5)(1.5,.5)(1.5,1)
 \psbezier[linewidth=.8pt](2.5,2)(2.5,1.5)(1.5,1.5)(1.5,1)
 \psbezier[linewidth=.8pt](.5,0)(.5,.5)(1,.5)(1,1)
 \psbezier[linewidth=.8pt](.5,2)(.5,1.5)(1,1.5)(1,1)
\end{pspicture}
\rrbracket
\end{equation}
is homotopy equivalent to the complex
\begin{equation}
\label{eq_R3r}\llbracket
\begin{pspicture}[.4](3,2)
   \psbezier[linewidth=.8pt](1.5,0)(1.5,.5)(2.5,.5)(2.5,1)
   \psbezier[linewidth=.8pt](1.5,2)(1.5,1.5)(2.5,1.5)(2.5,1)
   \pspolygon[linecolor=white,fillstyle=solid,fillcolor=white](2.2,0)(2.2,2)(1.9,2)(1.9,0)(2.2,0)
   \psbezier[linewidth=.8pt](2.5,2)(2.5,1.5)(.5,.75)(.5,0)
   \pspolygon[linecolor=white,fillstyle=solid,fillcolor=white](1.5,0)(1.5,2)(1.2,2)(1.2,0)(1.5,0)
 \psbezier[linewidth=.8pt](.5,2)(.5,1.25)(2.5,.5)(2.5,0)
  \rput(1.8,.9){$\scriptscriptstyle 8$}
  \rput(2.65,.9){$\scriptscriptstyle 5$}
  \rput(1.45,1.7){$\scriptscriptstyle 4$}
  \rput(2.65,1.7){$\scriptscriptstyle 7$}
  \rput(1.45,.3){$\scriptscriptstyle 6$}
  \rput(2.65,.3){$\scriptscriptstyle 1$}
  \rput(.35,1.7){$\scriptscriptstyle 2$}
  \rput(.35,.3){$\scriptscriptstyle 3$}
  \rput(1.5,1.4){$\scriptscriptstyle 9$}
  \rput(.9,1){$\scriptstyle \underline{1}$}
  \rput(2,1.85){$\scriptstyle \underline{2}$}
  \rput(2,.15){$\scriptstyle \underline{3}$}
\end{pspicture}
  \rrbracket
\end{equation}
where $\Phi$ is given by
\begin{equation}
\Phi := \quad
 \xy
    (-45,-20)*+{\scs (A \otimes A \otimes A)\{ -2\}}="1";
    (0,-20)*+{\scs (A \otimes A \otimes A)\{ 0\}
    \oplus  (A \otimes A \otimes A \otimes C)\{ 0\}}="2";
    (50,-20)*+{\scs (A \otimes A \otimes A)\{ 2\}}="3";
    (-45,20)*+{\scs (A \otimes A \otimes A)\{-4\}}="1'";
    (0,20)*+{\scs (A \otimes A \otimes A)\{ -2\}
        \oplus (A \otimes A \otimes A)\{ -2\}}="2'";
    (50,20)*+{\scs (A \otimes A \otimes A)\{ 0\}}="3'";
    {\ar_-{
    \left(
      \begin{array}{c}
         \psset{xunit=.15cm,yunit=.15cm}
        \begin{pspicture}(9,12)
         \rput(2,0){\crossl}
         \rput(2,7.5){\crossl}
         \rput(7,0){\identl}
         \rput(1,2.5){\identl}
         \rput(1,5){\identl}
         \rput(7,7.5){\identl}
         \rput(3,2.5){\saddler}
        \end{pspicture} \\
        \psset{xunit=.15cm,yunit=.15cm}
        \begin{pspicture}(8,5)
         \rput(7,0){\identl}
         \rput(7.25,0){\ltc}
         \rput(6,2.5){\comultl}
         \rput(5,0){\identl}
         \rput(1,2.5){\identl}
         \rput(1,0){\identl}
         \rput(3,2.5){\identl}
         \rput(3,0){\identl}
        \end{pspicture}\\
      \end{array}
    \right)
    } "1";"2"};
    {\ar_-{
    \left(
      \begin{array}{c}
    \psset{xunit=.15cm,yunit=.15cm}
        \begin{pspicture}(9,12)
         \rput(2,0){\crossl}
         \rput(2,7.5){\crossl}
         \rput(7,0){\identl}
         \rput(1,2.5){\identl}
         \rput(1,5){\identl}
         \rput(7,7.5){\identl}
         \rput(3,2.5){\saddler}
        \end{pspicture} \\ \psset{xunit=.15cm,yunit=.15cm}
        \begin{pspicture}(9,10)
        \rput(-.3,2.5){$-$}
         \rput(3,0){\multl}
         \rput(3.75,2.5){\ctl}
         \rput(2,2.5){\medidentl}
         \rput(6,2.5){\medidentl}
         \rput(5,0){\curverightl}
         \rput(8,2.5){\medidentl}
         \rput(8,4.5){\identl}
         \rput(2,4.5){\identl}
         \rput(2,7){\identl}
         \rput(4,7){\identl}
         \rput(8,0){\identl}
         \rput(5.4,4.5){\crossmixlc}
         \rput(7.4,7){\crossmixlc}
        \end{pspicture} \\
      \end{array}
    \right)^{\mT}
    } "2";"3"};
    {\ar^-{
    \left(
      \begin{array}{c}
     \psset{xunit=.15cm,yunit=.15cm}
        \begin{pspicture}(8,10)
         \rput(1,0){\identl}
         \rput(1,7.5){\identl}
         \rput(6,0){\crossl}
         \rput(7,2.5){\identl}
         \rput(7,5){\identl}
         \rput(6,7.5){\crossl}
         \rput(1,2.5){\saddler}
        \end{pspicture} \\
        \psset{xunit=.15cm,yunit=.15cm}
        \begin{pspicture}(8,5)
         \rput(7,0){\saddlel}
         \rput(1,2.5){\identl}
         \rput(1,0){\identl}
        \end{pspicture} \\
      \end{array}
    \right)
    } "1'";"2'"};
    {\ar^-{
    \left(
      \begin{array}{c}
     \psset{xunit=.15cm,yunit=.15cm}
       \begin{pspicture}(8,7.5)
         \rput(1,0){\identl}
         \rput(6,0){\crossl}
         \rput(7,2.5){\identl}
         \rput(7,5){\identl}
         \rput(1,2.5){\saddler}
        \end{pspicture}  \\
        \psset{xunit=.15cm,yunit=.15cm}
        \begin{pspicture}(9,5)
        \rput(-.3,2.5){$-$}
         \rput(2,0){\saddler}
         \rput(8,2.5){\identl}
         \rput(8,0){\identl}
        \end{pspicture} \\
      \end{array}
    \right)^{\mT}
    } "2'";"3'"};
 {\ar_-{
 \psset{xunit=.15cm,yunit=.15cm}
        \begin{pspicture}(9,5)
         \rput(1,2.5){\identl}
         \rput(1,0){\identl}
         \rput(7,0){\saddlel}
        \end{pspicture}
 } "1'";"1"};
  {\ar^-{
 \psset{xunit=.15cm,yunit=.15cm}
        \begin{pspicture}(9,5)
         \rput(7,2.5){\identl}
         \rput(7,0){\identl}
         \rput(1,0){\saddler}
        \end{pspicture}
 } "3'";"3"};
     {\ar_-{
    \left(
      \begin{array}{cc}
     \psset{xunit=.15cm,yunit=.15cm}
        \begin{pspicture}(9,5)
         \rput(1,2.5){\identl}
         \rput(1,0){\identl}
         \rput(7,0){\saddlel}
        \end{pspicture} & 0\\
        0 &
              \psset{xunit=.15cm,yunit=.15cm}
       \begin{pspicture}(8,7.5)
         \rput(1,0){\identl}
         \rput(3,0){\identl}
         \rput(6.3,0){\crossmixcl}
         \rput(1,2.5){\medidentl}
         \rput(3,2.5){\medidentl}
         \rput(5.3,2.5){\ltc}
         \rput(7,2.5){\medidentl}
         \rput(1,4.5){\identl}
          \rput(4,4.5){\comultl}
         \rput(7,4.5){\identl}
        \end{pspicture}   \\
      \end{array}
    \right)
    } "2'";"2"};
 \endxy
\end{equation}
but again by Lemma~\ref{lem_natans} the cone $\Gamma(\Phi)$ is equal
to the cone $\Gamma(G\circ \Phi)$ where $G$ is a strong deformation
retract. We take $G$ as in~\eqref{eq_defofF}, but with the arcs
renumbered appropriately. This leads to a map $G$ whose only nonzero
component $G^0$ is given by
\begin{equation}
 G^0 :=      \left(
      \begin{array}{c}
     \psset{xunit=.15cm,yunit=.15cm}
        \begin{pspicture}(7,3)
        \rput(0,1){$-$}
         \rput(2,0){\identl}
         \rput(4,0){\identl}
         \rput(6,0){\identl}
        \end{pspicture}\\
          \psset{xunit=.15cm,yunit=.15cm}
        \begin{pspicture}(9,10)
         \rput(2,0){\crossl}
         \rput(7,0){\identl}
         \rput(2,7.5){\crossl}
         \rput(7,7.5){\identl}
         \rput(1,2.5){\identl}
         \rput(1,5){\identl}
         \rput(3,2.5){\saddler}
         \rput(8.9,6.5){\deathc}
         \rput(9.3,7.5){\identc}
        \end{pspicture}  \\
      \end{array}
    \right)^{\mT}
\end{equation}
In order to make the enumeration of the complex $\llbracket
    \psset{xunit=.2cm,yunit=.2cm}
\begin{pspicture}[.4](3,2)
 \psline[linewidth=.6pt](.5,0)(.5,2)
 \psline[linewidth=.6pt](1.5,0)(1.5,2)
 \psline[linewidth=.6pt](2.5,0)(2.5,2)
\end{pspicture}
\rrbracket$ the same as in \eqref{eq_psi'}, we also postcompose the
composite $G\circ \Phi$ with the chain isomorphism $f$ of
Proposition~\ref{prop_indep}(2) whose only nonzero component is given
as follows:
\begin{equation}
   f^0 :=  \psset{xunit=.2cm,yunit=.2cm}
\begin{pspicture}[.4](6,8)
 \rput(1,0){\identl}
 \rput(4,0){\crossl}
 \rput(5,2.5){\identl}
 \rput(2,2.5){\crossl}
 \rput(1,5){\identl}
 \rput(4,5){\crossl}
\end{pspicture}
\end{equation}
Thus, the complex \eqref{eq_R3r} is homotopy equivalent to the cone of
the composite $\Phi'=f\circ G\circ\Phi$ given by
\begin{equation}
\label{eq_phi'}
\Phi' := \quad
 \xy
 (-60,-46)*{\llbracket
    \psset{xunit=.4cm,yunit=.4cm}
\begin{pspicture}[.4](3,2)
 \psline[linewidth=.6pt](.5,0)(.5,2)
 \psline[linewidth=.6pt](1.5,0)(1.5,2)
 \psline[linewidth=.6pt](2.5,0)(2.5,2)
 \rput(.3,1){$\scriptstyle 2$}
 \rput(1.3,1){$\scriptstyle 1$}
 \rput(2.3,1){$\scriptstyle 7$}
\end{pspicture}
\rrbracket =};
    (-45,-46)*+{0}="1";
    (0,-46)*+{\scs (A \otimes A \otimes A)\{ 0\}}="2";
    (45,-46)*+{0}="3";
    (-45,36)*+{\scs (A \otimes A \otimes A)\{ -4\}}="1'";
    (0,36)*+{\scs (A \otimes A \otimes A)\{ -2\}
    \oplus(A \otimes A \otimes A)\{ -2\}}="2'";
    (45,36)*+{\scs (A \otimes A \otimes A)\{ 0\}}="3'";
    {\ar_-{0} "1";"2"};
    {\ar_-{0} "2";"3"};
    {\ar^-{
    \left(
      \begin{array}{c}
     \psset{xunit=.15cm,yunit=.15cm}
        \begin{pspicture}(8,10)
         \rput(1,0){\identl}
         \rput(1,7.5){\identl}
         \rput(6,0){\crossl}
         \rput(7,2.5){\identl}
         \rput(7,5){\identl}
         \rput(6,7.5){\crossl}
         \rput(1,2.5){\saddler}
        \end{pspicture} \\
        \psset{xunit=.15cm,yunit=.15cm}
        \begin{pspicture}(8,5)
         \rput(7,0){\saddlel}
         \rput(1,2.5){\identl}
         \rput(1,0){\identl}
        \end{pspicture} \\
      \end{array}
    \right)
    } "1'";"2'"};
    {\ar^-{
    \left(
      \begin{array}{c}
     \psset{xunit=.15cm,yunit=.15cm}
       \begin{pspicture}(8,7.5)
         \rput(1,0){\identl}
         \rput(6,0){\crossl}
         \rput(7,2.5){\identl}
         \rput(7,5){\identl}
         \rput(1,2.5){\saddler}
        \end{pspicture}  \\
        \psset{xunit=.15cm,yunit=.15cm}
        \begin{pspicture}(9,5)
        \rput(-.3,2.5){$-$}
         \rput(2,0){\saddler}
         \rput(8,2.5){\identl}
         \rput(8,0){\identl}
        \end{pspicture} \\
      \end{array}
    \right)^{\mT}
    } "2'";"3'"};
 {\ar_-{0} "1'";"1"};
  {\ar^-{0} "3'";"3"};
     {\ar_-{
      \left(
      \begin{array}{c}
        \psset{xunit=.15cm,yunit=.15cm}
         \begin{pspicture}(7,15)
         \rput(-.5,10){$-$}
         \rput(0,0){
            \rput(1,0){\identl}
            \rput(4,0){\crossl}
            \rput(5,2.5){\identl}
            \rput(2,2.5){\crossl}
            \rput(1,5){\identl}
            \rput(4,5){\crossl}
            }
         \rput(3,7.5){\identl}
         \rput(1,7.5){\identl}
         \rput(6,7.5){
  \pspolygon[fillcolor=lightgray,fillstyle=solid](-1.5,0)(-.5,0)(1.5,2.5)(.5,2.5)(-1.5,0)}
         \rput(0,10){\rput(1,2.5){\identl}
            \rput(1,0){\identl}
            \rput(7,0){\saddlel}
         }
        \end{pspicture}\\
              \psset{xunit=.15cm,yunit=.15cm}
       \begin{pspicture}(10,30)
         \rput(0,0){
            \rput(1,0){\identl}
            \rput(4,0){\crossl}
            \rput(5,2.5){\identl}
            \rput(2,2.5){\crossl}
            \rput(1,5){\identl}
            \rput(4,5){\crossl}
            }
         \rput(6,7.5){
  \pspolygon[fillcolor=lightgray,fillstyle=solid](-1.5,0)(-.5,0)(1.5,2.5)(.5,2.5)(-1.5,0)}
         \rput(0,7.5){
            \rput(2,0){\crossl}
             \rput(2,7.5){\crossl}
            \rput(7,7.5){\identl}
            \rput(1,2.5){\identl}
             \rput(1,5){\identl}
            \rput(3,2.5){\saddler}
            \rput(9,6.5){\deathc}
            \rput(9.3,7.5){\identc}
         }
         \rput(1,17.5){\curverightl}
         \rput(3,17.5){\curverightl}
         \rput(7,17.5){\curveleftl}
         \rput(9.3,17.5){\curveleftc}
         \rput(1,20){
            \rput(1,0){\identl}
            \rput(3,0){\identl}
            \rput(6.3,0){\crossmixcl}
            \rput(1,2.5){\medidentl}
            \rput(3,2.5){\medidentl}
            \rput(5.3,2.5){\ltc}
            \rput(7,2.5){\medidentl}
            \rput(1,4.5){\identl}
            \rput(4,4.5){\comultl}
            \rput(7,4.5){\identl}
         }
        \end{pspicture}   \\
      \end{array}
    \right)^{\mT}
    } "2'";"2"};
 \endxy
\end{equation}
As an easy exercise in computing boundary permutations of
open-closed cobordisms, one can see that the chain maps
\eqref{eq_psi'} and \eqref{eq_phi'} are equal. Hence, their cones
$\Gamma(\Psi')=\Gamma(\Phi')$ are equal, making the complexes
\begin{equation}
  \llbracket
\begin{pspicture}[.4](3,2)
 \psbezier[linewidth=.8pt](1.5,2)(1.5,1.5)(.5,1.5)(.5,1)
 \psbezier[linewidth=.8pt](1.5,0)(1.5,.5)(.5,.5)(.5,1)
   \pspolygon[linecolor=white,fillstyle=solid,fillcolor=white](.8,0)(.8,2)(1.1,2)(1.1,0)(.8,0)
   \psbezier[linewidth=.8pt](2.5,2)(2.5,1.25)(.5,.5)(.5,0)
   \pspolygon[linecolor=white,fillstyle=solid,fillcolor=white](1.5,0)(1.5,2)(1.8,2)(1.8,0)(1.5,0)
 \psbezier[linewidth=.8pt](.5,2)(.5,1.5)(2.5,.75)(2.5,0)
\end{pspicture}
\rrbracket \qquad {\rm and} \qquad
 \llbracket
\begin{pspicture}[.4](3,2)
   \psbezier[linewidth=.8pt](1.5,0)(1.5,.5)(2.5,.5)(2.5,1)
   \psbezier[linewidth=.8pt](1.5,2)(1.5,1.5)(2.5,1.5)(2.5,1)
   \pspolygon[linecolor=white,fillstyle=solid,fillcolor=white](2.2,0)(2.2,2)(1.9,2)(1.9,0)(2.2,0)
   \psbezier[linewidth=.8pt](2.5,2)(2.5,1.5)(.5,.75)(.5,0)
   \pspolygon[linecolor=white,fillstyle=solid,fillcolor=white](1.5,0)(1.5,2)(1.2,2)(1.2,0)(1.5,0)
 \psbezier[linewidth=.8pt](.5,2)(.5,1.25)(2.5,.5)(2.5,0)
\end{pspicture}
  \rrbracket
\end{equation}
homotopy equivalent as graded [filtered] complexes by
Lemma~\ref{lem_natans}.

We leave it to the reader to verify the other possible colouring of
the Reidemeister three move as well as the Reidemeister three move for
the other crossing configurations\footnote{Bar-Natan's
paper~\cite{BN2} provides the proof for the other version of
Reidemeister three not proven here.}.

\subsubsection*{Acknowledgements}

We are grateful to Marco Mackaay, Simon Willerton, Christoph
Schweigert, Mikhail Khovanov, Kevin Costello, Ivan Smith, and Paul
Seidel for stimulating discussions. We are grateful to the European
Union Superstring Theory Network for support. A.L.\ would like to
thank the Max Planck Institute in Potsdam, and both of us would like
to thank the Mathematics Department at the University of Chicago for
their hospitality.

The graphics were typeset using \emph{xypic}~3.6 by Kristoffer Rose
and Ross More; using \emph{pstricks}~97 by Timothy van Zandt and
others; and using \emph{dbnsymb} and \emph{makefont} by Dror
Bar-Natan.

%
\end{document}